\documentstyle[amsfonts,graphicx]{report}

\newcommand{\commentout}[1]{}

\newcommand{\ba}{\begin{array}}
		\newcommand{\ea}{\end{array}}
\newcommand{\bdm}{\begin{displaymath}}
		\newcommand{\edm}{\end{displaymath}}
\newcommand{\ben}{\begin{enumerate}}
		\newcommand{\een}{\end{enumerate}}
\newcommand{\beq}{\begin{equation}}
		\newcommand{\eeq}{\end{equation}}
\newcommand{\bfg} {\begin{figure}[p]}
		\newcommand{\efg} {\end{figure}}%Nov 5,99
\newcommand{\bi} {\begin {itemize}}
		\newcommand{\ei} {\end {itemize}}
\newcommand{\bqn}{\begin{eqnarray}}
		\newcommand{\eqn}{\end{eqnarray}}
\newcommand{\bqs}{\begin{eqnarray*}}
		\newcommand{\eqs}{\end{eqnarray*}}
\newcommand{\bsl} {\begin{slide}[8.8in,6.7in]}
		\newcommand{\esl} {\end{slide}}
\newcommand{\bss} {\begin{slide*}[9.3in,6.7in]}
		\newcommand{\ess} {\end{slide*}}
\newcommand{\btb} {\begin {table}}
		\newcommand{\etb} {\end {table}}%Nov 10,99
\newcommand{\bvb}{\begin{verbatim}}
		\newcommand{\evb}{\end{verbatim}}
\newcommand{\bc}{\begin{center}}
		\newcommand{\ec}{\end{center}}
		
\newcommand{\m}{\mbox}
\newcommand {\pd}[2] {{\frac {\partial #1} {\partial #2}}}
\newcommand{\mat}[1]{{{\left[ \ba #1 \ea \right]}}}
\newcommand{\cas}[1]{{{\left \{ \ba #1 \ea \right. }}}
\newcommand {\nm} [1] {{\parallel #1\parallel }}
\newcommand{\reff}[1] {{{\textbf{Figure} \ref {#1}}}}
\newcommand{\refft}[2] {{{\textbf{Figures} \ref {#1} and \ref {#2}}}}
\newcommand{\refe}[1] {{{(\ref {#1})}}}%Nov 5
\newcommand{\reft}[1] {{{\textbf{Table} \ref {#1}}}}
\newcommand{\refet}[2] {{{(\ref {#1},\ref {#2})}}}%Nov 10'99

\def\pmb#1{\setbox0=\hbox{$#1$}%
   \kern-.025em\copy0\kern-\wd0
   \kern.05em\copy0\kern-\wd0
   \kern-.025em\raise.0433em\box0 }

\def\rfp {{{\rho^{\ast j}_{i+1/2} }}}
\def\vfp {{{v^{\ast j}_{i+1/2} }}}
\def\rfm {{{\rho^{\ast j}_{i-1/2} }}}
\def\vfm {{{v^{\ast j}_{i-1/2} }}}
\def\vast {{{v_{\ast}}}}
\def\rl  {{{\rho_l}}}
\def\vl  {{{v_l}}}
\def\rr  {{{\rho_r}}}
\def\vr  {{{v_r}}}
\def\rmm  {{{\rho_m}}}
\def\vm  {{{v_m}}}
\def\mfp {{{m^{\ast j}_{i+1/2} }}}
\def\rfm {{{\rho^{\ast j}_{i-1/2} }}}
\def\mfm {{{m^{\ast j}_{i-1/2} }}}
\def\vfm {{{v^{\ast j}_{i-1/2} }}}
\def\fast {{{f_{\ast}}}}
\def\ml  {{{m_l}}}
\def\mr  {{{m_r}}}
\def\rh {{\rho_h}}

\def\matauZ	{{\left( \begin {array}{cc}
	v &\rho\\
	\rho \vast'^2&v
	\end {array}\right)}}
\def\mattuZ	{{\left( \begin{array} {cc}
	1&1\\
	\vast'(\rho) &-\vast'(\rho)
	\end {array}\right) }}
\def\matluZ	{{\left( \begin{array} {cc}
\lambda_1(u)&0\\
0 &\lambda_2(u)
\end {array}\right) }}

\def\matu 	{{\left (\begin {array}{c}	
	\rho\\
	m
	\end {array} \right)  }}
\def\matf 	{{\left (\begin {array}{c}	
	m\\
	\frac {m^2}{\rho}+c_0^2\rho
	\end {array} \right) }}
\def\mats 	{{\left (\begin {array}{c}	
	0\\
	\frac{f_{\ast}(\rho)-m}{\tau}
	\end {array} \right) }}
\def\matau	{{\left( \begin {array}{cc}
	0 &1\\
	-\frac {m^2} {\rho^2}+c_0^2&\frac {2m}{\rho}
	\end {array}\right)}}
\def\mattu	{{\left( \begin{array} {cc}
	1&1\\
	\frac m {\rho}-c_0 &\frac m {\rho} +c_0
	\end {array}\right) }}
\def\matlu	{{\left( \begin{array} {cc}
\lambda_1(u)&0\\
0 &\lambda_2(u)
\end {array}\right) }}

\def\dx {{\Delta x}}
\def\dt {{\Delta t}}
\def\dxt {{\frac {\dt}{\dx}}}
\def\half {{\frac 12}}

\def\r {\rho}
\def\vs {v_{\ast}}
\def\ra {{\frac {\r}a}}
\def\l {{\lambda}}
\def\vecr{{\textbf{R}}}
\def\vecl{{\textbf{l}}}

\begin{document}
\pagenumbering{roman}
\bc
{\Large \bf Traffic Flow Models and Their Numerical Solutions}\author{{\sc Wenlong Jin}}

\vskip 0.2in
\centerline{By}
\vskip 0.1in
\centerline{\large  Wenlong Jin}
\centerline{\large B.S. (University of Science and Technology of China, Anhui, China) 1998}

\vskip 0.3in
\centerline{\large  THESIS}
\vskip 0.1in
\centerline{\large Submitted in partial satisfaction of the requirements for the
degree of}
\vskip 0.1in
\centerline{\large  MASTER OF SCHIENCE}
\vskip 0.1in
\centerline{\large in}
\vskip 0.1in
\centerline{\large APPLIED MATHEMATICS}
\vskip 0.2in
\centerline{\large in the}
\vskip 0.1in
\centerline{\large  OFFICE OF GRADUATE STUDIES}
\vskip 0.1in
\centerline{\large of the}
\vskip 0.1in
\centerline{\large  UNIVERSITY OF CALIFORNIA}
\vskip 0.1in
\centerline{\large  DAVIS}
\end{center}
\vskip 0.1in
{\large Approved:}

Dr. Michael Zhang

Dr. Elbridge Gerry Puckett

Dr. Albert Fannjiang
\vskip 0.2in
\begin{center}
\centerline{\large Committee in Charge}
\vskip 0.2in
\centerline{\large 2000}
\end{center}
\newpage
\large
\tableofcontents

%%%%%%%%%%%%%%%%%%%%%%%%%%%%%%%%%%%%%%%%%%%%%%%%%%%%%%%%%%%%%%%%%%
%%%%%%%%%%%%%%%%%%%%%%%%%%%%%%%%%%%%%%%%%%%%%%%%%%%%%%%%%%%%%%%%%%

\newpage
\thispagestyle{empty}
\large
{\Large \bf ACKNOWLEDGEMENTS} \\
This thesis is written as a summary of my work in last two years at UC Davis. At the moment when I finish my thesis, I want to express my sincere appreciation to those who have enthusiastically taught and helped me all the time.

First, many thanks go to my advisor, Dr. Michael Zhang. About one year ago, he introduced me to the area of transportation studies. Since then, I've worked on several projects on traffic flow models and ramp metering methods, and this thesis is nothing more than a combination of those projects. During the time when I work with him, he gave me many helpful suggestions on how to carry out researches and how to write down my ideas. 

I must also thank Dr. Elbridge Gerry Puckett, my academic advisor and the committee member of my thesis. He spent a great deal of time in answering all kinds of questions that I concern. He carefully read and revised my thesis. His suggestions on writing and comments on research will continue to influence my research in the future. He also introduced me to his colleagues, including Dr. Randall J. LeVeque and Dr. Phillip Colella, to whom I also owe my thanks for their help.

To Dr. Albert Fannjiang, I gave my thanks for his important comments on my thesis and serving in the committee of this thesis.

I'm grateful to Dr. Zhaojun Bai who gave me many good suggestions. 

I also want to thank my former classmate, Sheng Lu and my former roommate, Dr. John Hong, for sharing their ideas with me.

I should like to take this opportunity to thank Bill Broadley and Zach Johnson, the computer system administrators, whose support is important for so many computations in my projects.

Finally I want to thank my girlfriend, Ling, my parents, Jiasheng and Jiefang, my sister, Xiulian and my brothers, Zhongsheng, Wenhu, and Wenbin for being constant source of encouragement and support.

%%%%%%%%%%%%%%%%%%%%%%%%%%%%%%%%%%%%%%%%%%%%%%%%%%%%%%%%%%%%%%%%%%%%%%%%%%%%%%%%%%%%%%%%%%%%%%%%%%%%%%%%%%%%%%%%%%%%%%%%%%%%%%%%%%%%%%%%%%%%%%%%%%%%%%%%%%

\newpage
\pagenumbering{arabic}
\pagestyle{myheadings} 
\markright{  \rm \normalsize CHAPTER 1. \hspace{0.5cm}
 INTRODUCTION}
\large 
\chapter{Introduction}
\thispagestyle{myheadings}
       
\section{Motivation and Problem Definition}
Traffic networks -- consisting of highways, streets, and other kinds of roadways -- provide convenient and economical conveyance of passengers and goods. The basic activity in transportation is a trip, defined by its origin/destination, departure time/arrival time and travel route. A myriad of trips interact on the network to produce an intricate pattern of traffic flows. Since traffic conditions in many major metropolitan areas are becoming increasingly congested, affecting the operational efficiency of whole networks as well as the travel cost of each trip, traffic flow models are becoming more important in traffic engineering and the transportation policy making process. For example, well-developed traffic models are used in developing advanced ramp metering methods as well as in determining dynamic traffic assignment (JWL \& Zhang, 2000a\nocite{JWL2000a}).

There have been two approaches in mathematical modeling of traffic flow. One approach, from a microscopic view, studies individual movements of vehicles and interactions between vehicle pairs. This approach considers driving behavior and vehicle pair dynamics. But the size of the problem in a microscopic model becomes mathematically intractable when a considerable volume of traffic flow is considered. One example of the microscopic approach is the GM family of car-following models developed in the 1960's (e.g., Gazis et al., 1961\nocite{gazis1961}). The other approach studies the macroscopic features of traffic flows such as flow rate $q$, traffic density $\rho$ and travel speed $v$. The basic relationship between the three variables is: $q=\rho v$. Macroscopic models are more suitable for modeling traffic flow in complex networks since less supporting data and computation are needed. 

In this thesis macroscopic traffic flow models are studied both theoretically and numerically. 
Traffic flows are classified according to traffic conditions, roadway conditions and traffic network structure. Traffic flows are in equilibrium when the travel speed of these flows is uniquely determined as a function of traffic density, otherwise they are in non-equilibrium. Traffic flows are considered inhomogeneous when the roadway has different parameters at different locations. Link flows are flows on road links, and network flows are traffic flows on networks of roadways. Network flows differ from link flows in that vehicles in the former have different characteristics which affect traffic dynamics, such as the origins or destinations. 

Different types of traffic flow are described by different models.
For equilibrium link flow, the celebrated LWR model was developed by Lighthill and Whitham (1955) and Richards (1956). The LWR model has been solved for the homogeneous roadway rigorously. There have been empirical solutions to the inhomogeneous  LWR model. In this work a rigorous procedure to solve the inhomogeneous LWR model is developed. (JWL \& Zhang, 2000b\nocite{JWL2000b}). The LWR model is a first-order model in the sense of PDE system order. In this thesis we also discuss the PW model (Payne, 1971\nocite{payne}; Whitham, 1974\nocite{Whitham}) and Zhang's model (1998\nocite{Zhang1998}, 1999a\nocite{Zhang1999a}) for non-equilibrium link flow. Finally we introduce a multi-commodity model when traffic flow is disaggregated by origins, destinations or departure times.     

All the models we consider are based on conservation of traffic flow, i.e., the increment of vehicles in a section is equal to the difference between upstream in-flux and downstream out-flux in unit time. Mathematically, every model except the multi-commodity model, which is a discrete model, can be written as a continuous system of hyperbolic conservation laws. Solutions to the Riemann problem for all the continuous models are studied analytically, and all the models including the multi-commodity model are solved numerically with Godunov type of methods.

\section{Background and Research Overview}
\subsection{Traffic Flow Models}
In what follows the term ``traffic flow model" means a macroscopic traffic flow model. Many continuum traffic flow models can be described by a system of hyperbolic PDEs. The first of these models was the LWR model. This model relies on the assumption that there exists an equilibrium speed-density relationship $v=\vs (\r)$. Like other dynamic continuum flow models the LWR model is based on the mass conservation, i.e., traffic conservation, and is described by a first-order, nonlinear PDE:
\bqn
\r_t+f(\r)_x&=&0, \label{org:lwr}
\eqn
in which $f(\r)=\r \vs(\r)$ is the traffic flow rate. Equation \refe{org:lwr} is in conservation form. It is also called a kinematic wave model since it shows wave properties analogous to those of gases. There are many numerical methods to solve the LWR model. One approach is to solve the Riemann problem and apply a Godunov method for this model. Both solutions to the Riemann problem and the Godunov method are well-developed for hyperbolic conservation laws (Smoller, 1983\nocite{smoller1983}). Another approach is to use the demand and supply functions (Lebacque, 1996\nocite{lebacque1996}; Daganzo, 1995\nocite{Daganzo95a}), which turns out to be variants of Godunov's method.

The PW model, derived based on microscopic car-following models, discards the equilibrium assumption. It is a second-order system of hyperbolic conservation laws with a source term:
\bqn
\r_t+(\r v)_x&=&0,\label{org:pw1}\\
v_t+vv_x+\frac{c_0^2}\r \r_x &=&\frac {\vs (\r)-v}{\tau},\label{org:pw2}
\eqn
in which the constant $c_0$ is the traffic sound speed and $\tau$ is the relaxation time.  Equation \refe{org:pw1} is the continuity equation, and \refe{org:pw2} is the momentum equation. The PW model relates to driver behavior models better than the LWR model because it accounts for drivers' anticipation and inertia. However it's been shown that the LWR model is an asymptotic approximation of the PW model (Schochet, 1988\nocite{schochet1988}). The PW model better captures non-equilibrium wave phenomena in traffic flow. Another property of the PW model is that it is unstable under certain situations. Since there are no known analytical solutions to the PW model,  we use numerical methods to solve it in this research. All of these methods are developed from Godunov's method. 

Also based on microscopic models, Zhang (1998\nocite{Zhang1998}) developed a new non-equilibrium traffic flow theory :
\bqn
\r_t +(\r v)_x&=&0, \label{org:zhang1}\\
v_t+vv_x+\frac{(\r \vs'(\r))^2}{\r} \r_x &=& \frac {\vs(\r)-v}{\tau}. \label {org:zhang2}
\eqn
In the momentum equation \refe{org:zhang2}, a varying sound speed $c=\r \vs'(\r)$ has been introduced. It has been shown that this new model avoids ``wrong-way travel" which is exhibited in the PW model (Zhang, 1998\nocite{Zhang1998}),  and the model is always stable. Wave solutions to this model were discussed in (Zhang, 1999a\nocite{Zhang1999a}), and finite difference approximations were studied in (Zhang, 2000a\nocite{Zhang2000a}). In this research we perform the numerical simulations of this model.

For a roadway with inhomogeneities such as a change in number of lanes, curvature and slopes, the LWR model can still be used, but the equilibrium speed-density relationship varies with location. By introducing an inhomogeneity function $a(x)$ which is a profile of the roadway at the location $x$, we can write the inhomogeneous LWR model as
\bqn
\r_t+f(a,\r)_x&=&0.
\eqn
Here the traffic flow rate $f(a,\r)$ is a function of the inhomogeneity $a(x)$. By writing 
\bqn
a_t&=&0,
\eqn
the inhomogeneous LWR model is a non-strictly hyperbolic system. There have been empirical methods to solve the inhomogeneous LWR model (Lebacque, 1995\nocite{lebacque1995}; Daganzo, 1995a\nocite{Daganzo95a}). In this research, we develop a rigorous procedure to solve the LWR model based on the work by Isaacson et al. (1992\nocite{isaacson_t1992}) and Lin et al. (1995\nocite{lin_t_w1995}). We find the solutions that are consistent with those given by Lebacque. However,  our method can be extended to solve higher-order inhomogeneous models while those of Lebacque and Daganzo cannot.

Multi-commodity models are discussed in (Daganzo, 1994\nocite{Daganzo94} and 1995\nocite{Daganzo95}; Jayakrishnan, 1991\nocite{jayakrishnan91}; Vaughan et al., 1984\nocite{vaughan84}). In these models traffic flow is disaggregated by origins, destinations or departure times. All of these models are based on traffic conservation.  A First-In-First-Out (FIFO) discipline is assumed in all of these multi-commodity models. The model studied by Vaughan et al. is a continuous model. The models studied by Jayakrishnan and Daganzo are discrete models. In these two discrete models vehicles that are close to each other (in the sense of location or time) and have the same origin or destination or some other common characteristics are combined as a macroparticle. The macroparticles in a zone are ordered by time or location. Macroparticles are moved according to traffic conditions, which are solved with a link flow model. In the thesis, we introduce a new multi-commodity model that is more efficient in moving macroparticles. 

\subsection{Hyperbolic Systems of Conservation Laws and Godunov Methods}
A system of hyperbolic conservation laws (Smoller, 1983\nocite{smoller1983}) takes the following form:
\bqn
u_t+f(u)_x&=&0,\label{conse}
\eqn
where $u=(u_1,\cdots, u_n)\in {\mathbb R}^n, n\geq 1$, and $(x,t)\in {\mathbb R}\times {\mathbb R}_+$. The $n$ eigenvalues of the differential of $f(u)$, $D f(u)$, are denoted as $\lambda_1,\cdots,\lambda_n$. The solutions related to $i$-th eigenvalue are called $i$-family wave solutions. If the eigenvalues are distinct, the system \refe{conse} is a strictly hyperbolic system. To solve \refe{conse}, initial and boundary conditions are needed. The Riemann problem is to solve \refe{conse} with the following jump initial condition:
\bqn
u(x,t=0)&=&\cas{{ll}u_l, & x<0\\u_r, &x> 0},
\eqn
where $u_l$ and $u-r$ are constants.

It is well-known that the weak solutions to the Riemann problem exist and are unique for the system \refe{conse} under the so-called ``Lax's entropy condition" (Lax, 1972\nocite{lax1972}). The system admits discontinuous solutions, i.e., shock waves, and the wave speed $s$ is determined by Rankine-Hugoniot condition:
\bqn
s[u]&=[f(u)],
\eqn
where $[u]=u_l-u_r$, and similarly, $[f(u)]=f(u_l)-f(u_r)$. For valid $i$-family shock wave solutions, the entropy inequalities hold, i.e., $\lambda_i(u_l)>s>\lambda_i(u_r)$. There are continuous solutions $u=u(\xi)$, $\xi=x/t$, which satisfy the ordinary differential equations
\bqn
-\xi u_{\xi}+f(u)_{\xi}&=&0.
\eqn

For most general initial and boundary conditions, there are no analytical solutions to \refe{conse}, hence, one must use numerical methods to solve it. The Godunov method is one of the most efficient numerical methods, which combines solutions to the Riemann problem and conservation laws. In a Godunov method, the space region $[a,b]$ is divided into $N$ grids $x_0=a, x_1,\cdots, x_{N-1}, x_N=b$; the time scale $[t_0, t_1]$ are partitioned into $M$ time steps $t^0=t_0, t^1, \cdots, t^{M-1}, t^M=t_1$. The hyperbolic system of conservation laws \refe{conse} can be approximated by the finite difference equations:
\bqn
\frac{ U_i^{j+1}-U_i^j}{\dt}+\frac {f(U^{\ast}_{i-1/2})-f(U^{\ast}_{i+1/2})}{\dx}&=&0, \label{diff}
\eqn
where $U_i^j$ is the average of $u(x,t)$ in grid $[x_{i-1/2},x_{i+1/2}]$ at time step $t_j$, similarly $U_i^{j+1}$ is the average at time step $t_{j+1}$; $U^{\ast}_{i-1/2}$ is the average of $u(x,t)$ during time interval $[t_j, t_{j+1}]$ at grid boundary $x_{i-1/2}$, similarly $U^{\ast}_{i+1/2}$ is the average at boundary $x_{i+1/2}$. The boundary flux $f(U^{\ast}_{i-1/2})$ is calculated by solving a Riemann problem at each cell edge with the following initial conditions:
\bqs
u(x,t^j)&=&\cas{{ll}U^j_{i-1}, &x<x_{i-1/2}\\U^j_i, &x> x_{i+1/2}}.
\eqs
Equations \refe{diff} says that the increment of $u$ is equal to the difference between both boundary fluxes, which is the general idea of conservation.

In this thesis our second-order models are not exact conservation laws since they have a source term, therefore our methods  have been extended to address the effect of source terms.

\section{Layout of the Thesis}
In chapter 2, we study the homogeneous LWR model through theoretical discussions and numerical simulations. This is the first step for us to understand traffic flow models and associated numerical methods. In chapter 3, Godunov-type methods are developed for Zhang's model. The Riemann problem is solved in detail, and computational results are discussed. In chapter 4, the PW model is studied with several different Godunov-type finite difference methods, and the stability of the model is tested. In chapter 5, the Riemann problem for the inhomogeneous LWR model is solved rigorously. In chapter 6, a multi-commodity model based on the LWR model is studied. In chapter 7, possible extensions of this research are discussed.

\newpage
\pagestyle{myheadings} 
\markright{  \rm \normalsize CHAPTER 2. \hspace{0.5cm}
 The LWR Model and Its Numerical Solutions}
\large 
\chapter{The LWR Model and Its Numerical Solutions}
\thispagestyle{myheadings}

The landmark paper by Lighthill and Whitham (1955) set the tone for many researchers' investigations into the theory of traffic flow. Introduced for traffic flows on a single, long, and rather idealized road, the LWR model proposes that their dynamics is described by the following PDE:
\bqn
\r_t+f(\r)_x&=&0, \label{lwr:eqn1}
\eqn
where subscript $t$ means the partial derivative with respect to time $t$, and subscript $x$ means the partial derivative with respect to location $x$.
In \refe{lwr:eqn1}, the function $f(\r)=\r \vs(\r)$ is called the fundamental diagram of traffic flow, in which $\vs(\r)$ reflects the equilibrium relationship between travel speed and traffic density. It's generally assumed that $f(\r)$ is a concave function, i.e., $f_{\r\r}(\r)<0$. The characteristic wave speed $\lambda(\r)=f_{\r}(\r)=\vs(\r)+\r \vs'(\r)$. $\lambda(\r)$ can be positive or negative.

Weak solutions (Smoller, 1983\nocite{smoller1983}) for \refe{lwr:eqn1} satisfy the integral form of the conservation law:
\bqn
\frac{\partial}{\partial t}\int_{x_1}^{x_2} \r(x,t)\m{dx}&=&f(x_1,t)-f(x_2,t).
\eqn
First, we discuss the Riemann problem for \refe{lwr:eqn1} theoretically and then give the numerical solutions.

\section{The Riemann problem}
We consider the Riemann problem for the LWR model with the following jump initial condition:
\bqn
\r(x,t=0)&=&\cas{{ll}\rl, & x<0\\\rr, &x> 0}. \label{lwr:ini}
\eqn

There are two types of wave solutions to the Riemann problem. The discontinuous solution is a shock wave
\bqn
\r(x,t>0)&=&\cas{{ll}\rl, & x/t<s\\\rr, &x/t> s},
\eqn
where $s$ is the shock wave speed. The shock wave speed is determined by the Rankine-Hugoniot jump condition,
\bqn
s&=&\frac {[f(\r)]}{[\r]}.
\eqn

The wave speed of a valid shock wave solution has to satisfy the entropy condition:
\bqn
\lambda(\rl)>s>\lambda(\rr).
\eqn
Specifically, for a concave fundamental diagram, a shock wave is a solution to the Riemann problem for \refe {lwr:eqn1} with initial conditions \refe{lwr:ini}  when $\rl<\rr$; i.e., the upstream traffic density is lower.

When the upstream traffic density is higher, solution to the LWR model with initial data \refe{lwr:ini} is a continuous rarefaction wave.  The rarefaction wave is given by $\r=\r(\xi), \xi=x/t$, where $\r$ satisfies the ordinary differential equation
\bqn
-\xi \r_{\xi}+f(\r)_{\xi}&=&0.
\eqn
If $\r(\xi)\neq 0$, we obtain
\bqn
\lambda(\r(\xi))=\xi,
\eqn
from which we can find $\r(\xi)$ uniquely. On any characteristic $x/t=\xi$, $\r$ is constant.

\section {Computation of boundary fluxes}
Given a wave solution to the Riemann problem, we can compute the average $\r^{\ast}$ at $x=0$ and therefore the flux $f(\r^{\ast})$ through the boundary. There are the following five cases:
\bi
\item[Case 1] When the solution to the Riemann problem is a shock wave with wave speed $s>0$, $\r^{\ast}=\rl$.
\item[Case 2] When the solution to the Riemann problem is a shock wave with wave speed $s\leq0$, $\r^{\ast}=\rr$.
\item[Case 3] When the solution to the Riemann problem is a rarefaction wave, and $\lambda(\rl)>0$, $\r^{\ast}=\rl$.
\item[Case 4] When the solution to the Riemann problem is a rarefaction wave, and $\lambda(\rr)<0$, $\r^{\ast}=\rr$.
\item[Case 5] When the solution to the Riemann problem is a rarefaction wave, $\lambda(\rl)<0$ and $\lambda(\rr)>0$, $\r^{\ast}$ is the solution of $\lambda(\r^{\ast})=0$.
\ei

Given initial and boundary conditions, we can use a first-order Godunov method to calculate the traffic conditions. The numerical solutions are given in next section.

%Tue Dec 21 14:53:58 PST 1999
\section {Numerical solutions to the LWR model}
In this section, we solve the Riemann problem for \refe{conse} numerically. We use Newell's equilibrium model,
\bqn
v_{\ast}(\rho)&=&v_f\left(1-\exp\{\frac {|c_j|}{v_f}(1-\rho_j/\rho\}\right)
\eqn
Without loss of generality, we set $v_f=1, c_j=1, \rho_j=1$ to obtain
\bqn
v_{\ast}(\rho)&=&1-\exp{(1-\frac 1{\rho})}
\eqn
The valid range for $\r$ and $v$ is $0<\r, v \leq 1$.  The corresponding flow rate, $f_{\ast}=\rho v_{\ast}$, is a normalized fundamental diagram.

\ben
\item Given the initial conditions:
\bqn
\rho(x,0)&=&\cas {{ll} 0.65& x\in [0l, 200l]\\0.4& \m{otherwise}}\\
v(x,0)&=&v_{\ast}(\rho(x,0)), \qquad \forall x\in [0l, 800l]
\eqn
we obtain the solutions shown in  \reff{LWRRiemann3_3d},\reff{LWRRiemann3_2d}.

\item Given the initial conditions:
\bqn
\rho(x,0)&=&\cas {{ll} 0.65& x\in [0l, 200l]\\0.9& \m{otherwise}}\\
v(x,0)&=&v_{\ast}(\rho(x,0)), \qquad \forall x\in [0l, 800l]
\eqn
we obtain the solutions shown in \reff{LWRRiemann4_3d},\reff{LWRRiemann4_2d}.

\een

From these numerical solutions we can see the formation of shock waves and rarefaction waves. These solutions are only an approximation of real solutions. For example, the solutions after $t=0$ in \reff{LWRRiemann4_2d} are not exact jumps, while the theoretical solution to the LWR model with initial conditions \refe{lwr:ini} is still a jump at any time. However as $\dx\to0$ and $\dt\to 0$, the solutions given by Godunov's method converge to the exact solutions.

\bfg
\bc\includegraphics [height=8cm] {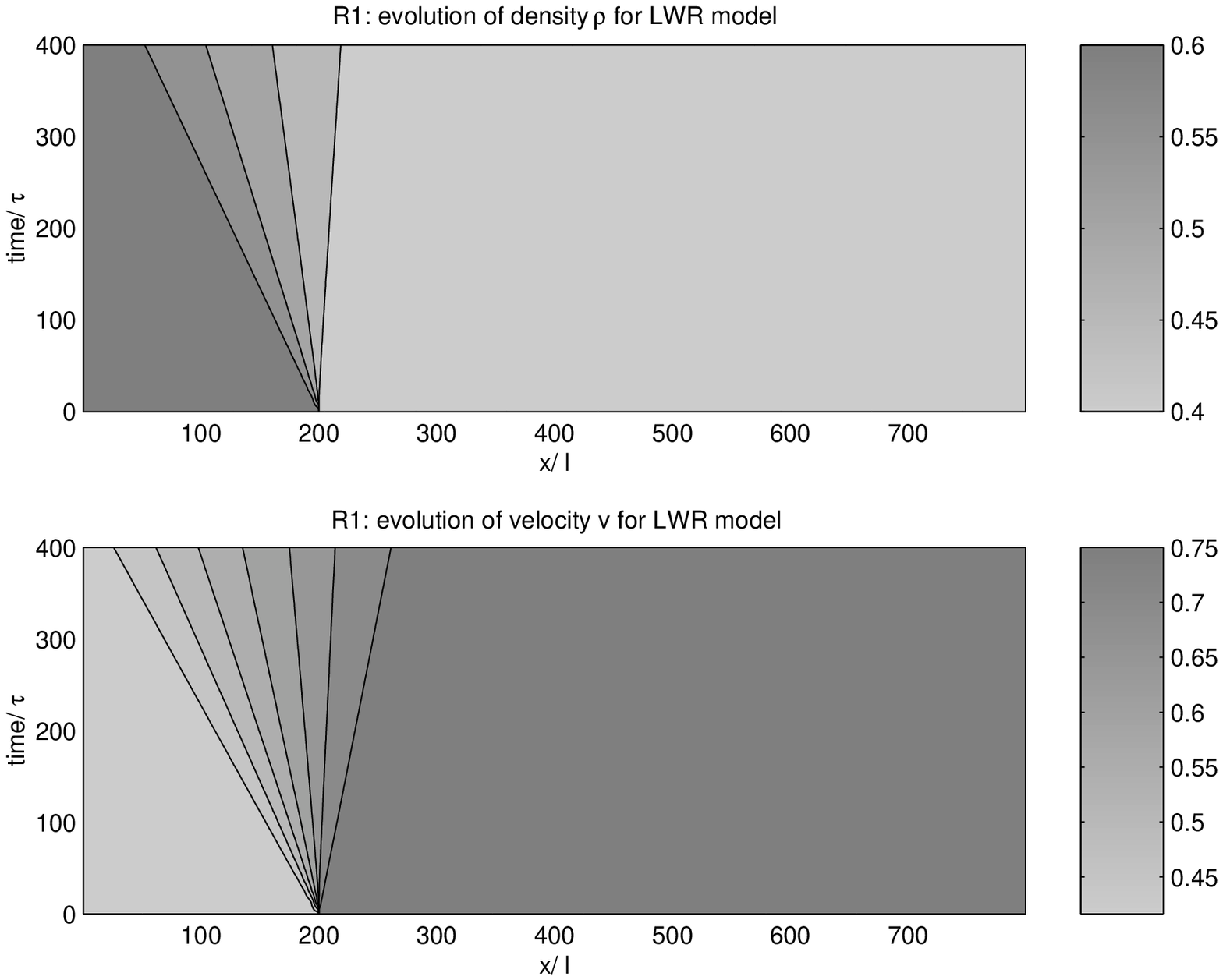}\ec
\caption {Rarefaction wave solution of the LWR model with initial data \refe{lwr:ini}} \label {LWRRiemann3_3d}
\efg
\bfg
\bc\includegraphics [height=8cm] {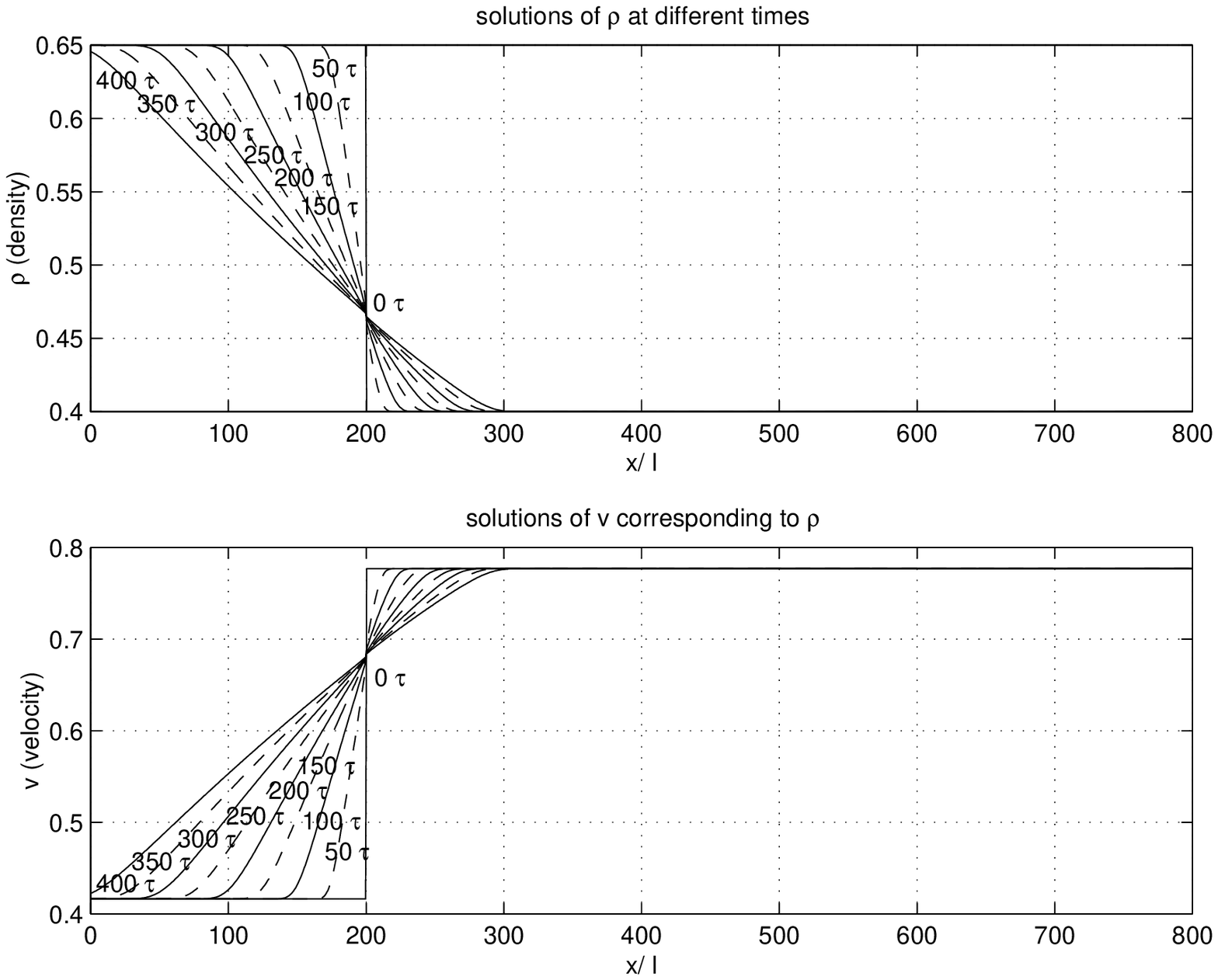}\ec
\caption {Solutions from \reff{LWRRiemann3_3d} at selected times} \label {LWRRiemann3_2d}
\efg

\bfg
\bc\includegraphics [height=8cm]{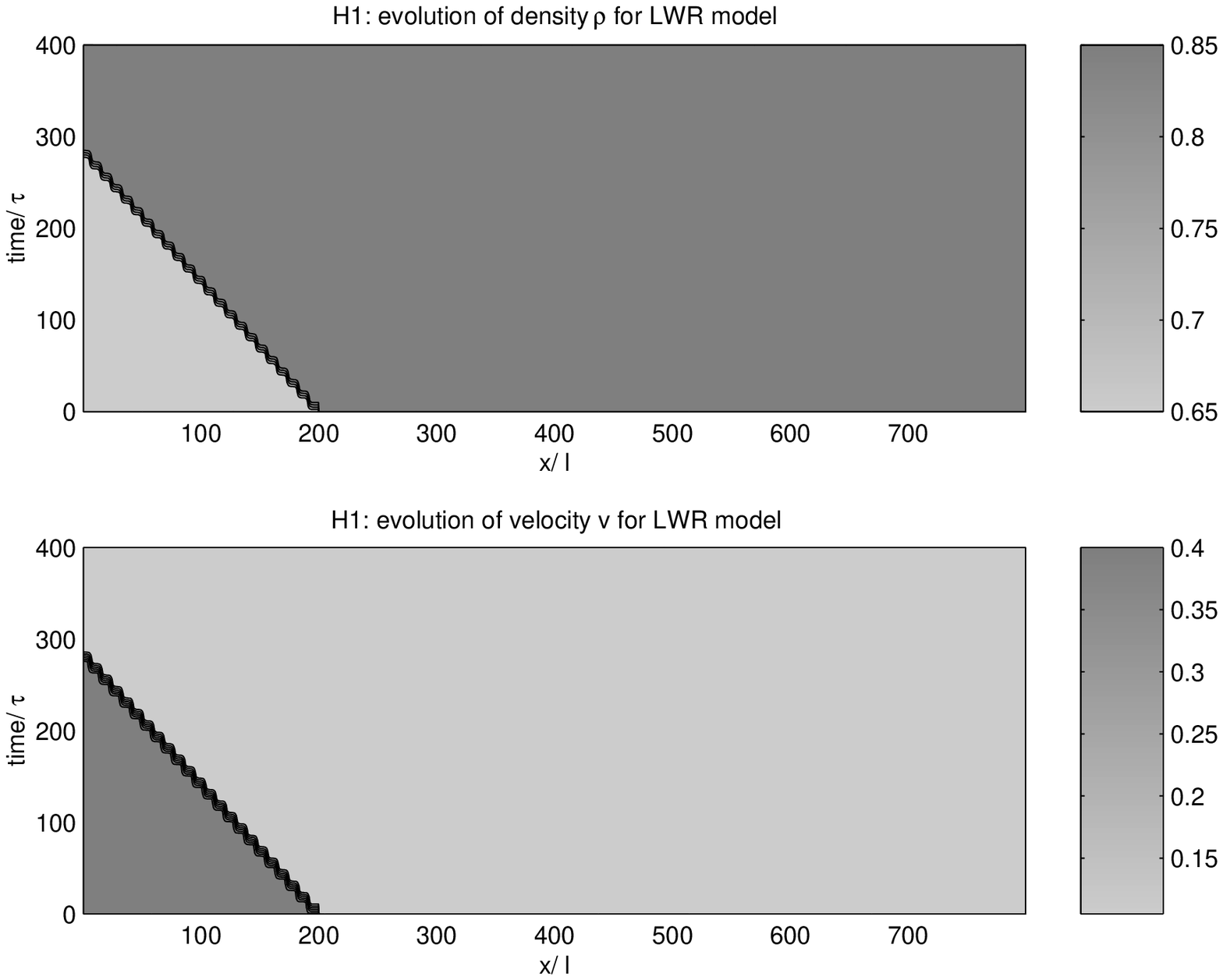}\ec
\caption {Shock wave solution of the LWR model with initial data \refe{lwr:ini}} \label {LWRRiemann4_3d}
\efg
\bfg
\bc\includegraphics [height=8cm] {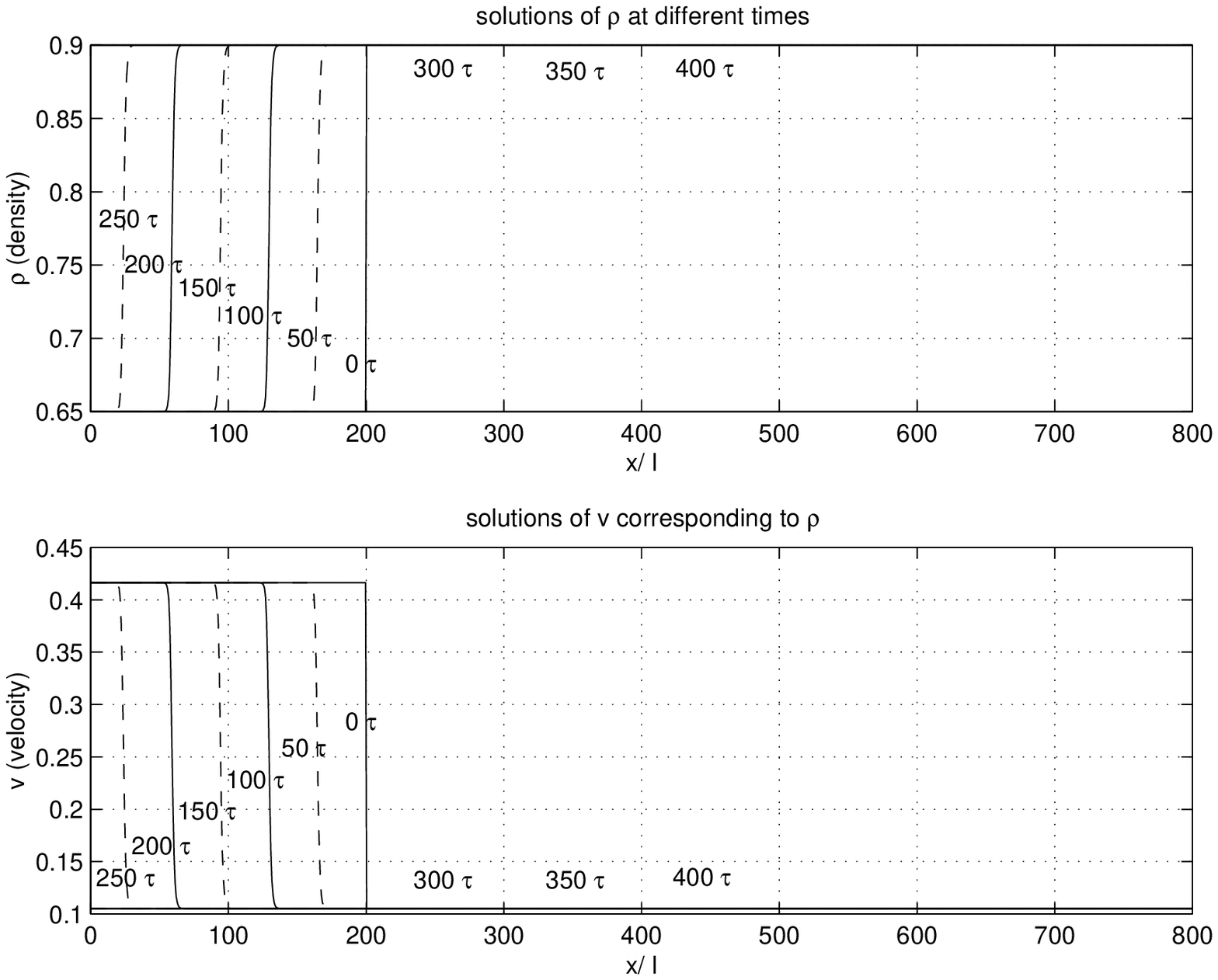}\ec
\caption {Solutions from \reff{LWRRiemann4_3d} at selected times} \label {LWRRiemann4_2d}
\efg

\newpage
\pagestyle{myheadings} 
\markright{  \rm \normalsize CHAPTER 3. \hspace{0.5cm}
 Zhang's Second-Order Traffic Flow Model and Its Numerical Solutions}
\large
\chapter{Zhang's Second-Order Traffic Flow Model and Its Numerical Solutions}
\section {Introduction}
A theory of non-equilibrium traffic flow has been developed by Zhang (1998\nocite{Zhang1998}), and in Zhang (2000a\nocite{Zhang2000a}) he developed the Godunov-type finite difference equations (FDE) for this model. Zhang's model is a second-order model, and can be written in the conservation form
\begin {eqnarray}
\left ( \begin {array} {c}
\rho\\
v
\end {array} \right )_t
+
\left ( \begin {array} {c}
\rho v\\
\frac {v^2} 2+\phi(\rho)
\end {array} \right )_x
&=&
\left ( \begin {array} {c}
0\\
\frac {v_{\ast}(\rho)-v}{\tau}
\end {array} \right ), \label{zhang.1}
\end {eqnarray}
where $\phi(\rho)$ is a velocity flux function and defined as
\begin {eqnarray}
\phi ' (\rho) &=& \frac {c^2(\rho)} {\rho} =\rho(v_{\ast}'(\rho))^2.
\end {eqnarray}
Here $c(\r)=-\r \vs'(\r)$ is the traffic sound speed.

In \refe{zhang.1}, $\vs(\r)$ is the equilibrium speed. Some often-used equilibrium travel speed functions are listed along with the corresponding velocity flux functions $\phi(\rho)$ as follows:
\begin {center}
\begin {tabular} {||c|c|c||}\hline
Functions&$\vast(\rho)$ &$\phi(\rho)$ \\\hline
Greenshields & $v_f(1-\rho/\rho_j)$ & $\frac {{v_f}^2}{2{\rho_j}^2} \rho^2$\\\hline
Polynomial & $v_f(1-(\rho/\rho_j)^n), n>1$ & $\frac {v_f^2} {2} (\rho/\rho_j)^{2n} $\\\hline
Greenberg & $v_0 \ln(\rho_j/\rho)$ & $v_0^2 \ln (\rho)$ \\\hline
Underwood & $v_f \exp (-\rho/\rho_j)$ & $-v_f^2 (1+\rho/\rho_j) \exp(-\rho/\rho_j)$\\\hline
Newell& $v_f[1-\exp(\frac {|c_j|}{v_f}(1-\rho_j/\rho)) ]$ & $\frac {v_f^2}{2}(\frac {\rho_j|c_j|}{v_f} \frac 1 {\rho}+\frac 12)\exp(2\frac {|c_j|}{v_f}(1-\rho_j/\rho))  $\\\hline
\end {tabular}
\end {center}

The equilibrium travel speed $\vs$ is decreasing with respect to traffic density; i.e., $\vast'(\rho)<0$. The fundamental diagram $f_{\ast}(\rho)\equiv \rho\vast(\rho)$ is concave; i.e., $f_{\ast}''(\rho)<0$. In \refe{zhang.1}, $\tau$ is the relaxation time and the relaxation term $\frac {v_{\ast}(\rho)-v}{\tau}$ constrains the difference between the real travel speed $v$ and the equilibrium travel speed $v_{\ast}$. When the relaxation term is 0, Zhang's model reduces to the LWR model:
\bqn
\rho_t+(\rho v_{\ast})_x&=&0.
\eqn

 Zhang's model has three different wave velocities (relative to the road): A first-order wave velocity  and two second-order wave velocities. The first-order wave velocity is the wave velocity of the corresponding LWR model:
\begin {eqnarray}
\lambda_{\ast} (\rho)&=&\vs(\r)-c(\r)=\vast(\rho)+\rho \vast'(\rho).
\end {eqnarray}
The two second-order wave velocities are
\begin {eqnarray}
\lambda_{1,2}(\rho,v)&=&v\mp c(\rho)=v\pm\rho\vast'(\rho).
\end {eqnarray}
The relationship between these three wave speeds along $(\rho,\vast)$ phase curves is that
\begin {eqnarray}
\lambda_1=\lambda_{\ast}<\lambda_2.
\end {eqnarray}
 The waves with wave speed $\lambda_1$ are called 1-waves; similarly the waves with wave speed $\lambda_2$ is called 2-waves. Since $\r,v\geq0$ and $v'_{\ast}<0$, the 2-wave speed $\lambda_2>0$ for any $\rho,v$.

Since Zhang's model is a hyperbolic system of conservation laws with a relaxation term, the system is stable when (Liu, 1979\nocite {Liu79}, 1987\nocite{Liu87}; Chen et al., 1994\nocite{Chen94})
\bdm
\lambda_1\leq \lambda_{\ast}\leq \lambda_2.
\edm
This condition is satisfied by Zhang's model. Thus Zhang's model is always stable.

In the following sections we study Godunov-type methods and use them to solve \refe{zhang.1} numerically. In section 2 we discuss Godunov's method and properties of Zhang's model. In section 3 we present a second-order Godunov method. In Section 4 we solve the Riemann problems numerically and discuss the order of accuracy for different methods.

\section {Godunov's method}
A Godunov-type finite difference method for Zhang's model was first presented in (Zhang, 2000a\nocite {Zhang2000a}). In this section we review this Godunov method and solve the associated Riemann problem.

The Godunov-type FDEs for Zhang's model are
\begin {eqnarray}
\frac {\rho^{j+1}_i-\rho^j_i }{k}+\frac {\rfp \vfp-\rfm \vfm } {h} &=&0 ,
\end {eqnarray}
\begin {eqnarray}
\frac {v^{j+1}_i-v^j_i }{k}
+\frac {\frac{(\vfp)^2 }2 +\phi(\rfp) -\frac{(\vfm)^2 }2 -\phi(\rfm)}{h}\nonumber\\
&= \frac {\vast(\rho^{j+1}_i)-v^{j+1}_i} {\tau}.
\end {eqnarray}
In these FDEs, $\rho^{j}_i$ is the average of $\r$ in cell $i$ at time step $j$; i.e.,
\begin {eqnarray}
\rho^j_i&=&\frac 1 h \int_{x_{i-1/2}}^{x_{i+1/2}} \rho(x,t_j) dx.
\end {eqnarray}
Similarly $v^j_i$ is the average of $v$. We use $\rfp$ as the average of $\r$ through the cell boundary $x_{i+1/2}$ over the time interval $(t_j,t_{j+1})$, i.e.,
\begin {eqnarray}
\rfp &=& \frac 1 k \int _{t_j}^{t_{j+1}} \rho(x_{i+1/2},t) dt.
\end {eqnarray}
Similarly we define $\vfp,\rfm,\vfm$ as boundary averages.

By measuring the source term with values at time $t_{j+1}$, we write the evolution equations for Zhang's model as
\begin {eqnarray}
&\rho^{j+1}_i& = \rho^j_i-\frac kh (\rfp\vfp-\rfm\vfm )\label{ev.1}\\
&v_i^{j+1}&= \frac 1 {(1+\frac k {\tau})  } \{v_i^j -\frac kh[\frac{(\vfp)^2 }2 +\phi(\rfp) -\frac{(\vfm)^2 }2 -\phi(\rfm) ] \nonumber\\&&\qquad+\frac k {\tau} \vast (\rho^{j+1}_i)\} \label{ev.2}
\end {eqnarray}

Provided traffic conditions $(\r,v)$ at time $t_j$, traffic conditions at time $t_{j+1}$ can be calculated if we know the boundary averages $\rfp,\vfp,\rfm,\vfm$. The computation of $\rfp,\vfp$ at the cell boundary $x_{i+1/2}$ during the time interval $(t_j,t_{j+1})$ depends on  a Riemann problem for \refe{zhang.1} with the following initial conditions
\begin {eqnarray}
u_{i+1/2}(x,t_j)=\left\{\begin {array} {l c}
U_l, & \mbox{ if } x-x_{i+1/2}<0 \\
U_r, & \mbox{ if } x-x_{i+1/2}> 0
\end {array}\right. ,
\end {eqnarray}
where we define the state variable $u(x,t)=(\rho,v)$, $U_i^j=(\rho_i^j,v_i^j)=u(x_i,t_j)$ and left and right states $U_l=(\rl,\vl),U_r=(\rr,\vr)$.
In a first-order Godunov method, we use the cell averages $\rl=\rho_i^j,\vl=v_i^j$ as the left side (upstream) traffic conditions and $\rr=\rho_{i+1}^j,\vr=v_{i+1}^j$ as the right side (downstream) conditions. In a second-order Godunov method, we use higher-order approximations to the left and right states. (For details for a second-order Godunov method, refer to Section \ref{secondGodunov}.)

Here we have neglected the relaxation term in \refe{zhang.1} when solving the Riemann problem. This Riemann problem has been discussed by Zhang (1999a\nocite{Zhang1999a}), and the solutions to the boundary averages are provided there. The solutions are self-similar and can be expressed in the form of
\begin {displaymath}
\psi(\frac {x-x_{i+1/2}} {t}; U_r,U_l).
\end {displaymath}

\subsection {Solutions of the boundary averages}\label{sub3.2.1}
There are 8 types of wave solutions to the Riemann problem, which are combinations of two 1-waves and two 2-waves. The calculation of the boundary averages depend on the type of solutions. The formula for calculating the boundary averages are listed as follows.

\begin {enumerate}

\item The wave solution is a 1-shock when the initial conditions satisfy
\begin {eqnarray}
\mbox{ H1: } \vr-\vl&=& -\sqrt {\frac {2(\rl-\rr)(\phi(\rl)-\phi(\rr))}{\rl+\rr}}, \quad \rr>\rl,\: \vr<\vl.
\end {eqnarray}
The wave speed is
\begin {eqnarray}
s&=&\frac {\rr\vr-\rl\vl}{\rr-\rl}
\end {eqnarray}
The boundary averages $(\rfp,\vfp)$ are given in the following table:
\\\begin {center}
\begin {tabular} {||c||c|c|c||}\hline
&$s=\frac {\rr\vr-\rl\vl}{\rr-\rl}$ &$\rfp$&$\vfp$ \\\cline{2-4}
&$s>0$ & $\rl$ &$\vl$ \\\cline{2-4}
H1&$s<0$ & $\rr$ & $\vr$ \\\cline{2-4}
&$s=0$ & $\frac {\rl+\rr} 2$ & $\frac {\vl+\vr}2$ \\\hline
\end {tabular}
\end {center}

\item The wave solution is a 2-shock when the initial states satisfy
\begin {eqnarray}
\mbox{ H2: } \vr-\vl&=& -\sqrt {\frac {2(\rl-\rr)(\phi(\rl)-\phi(\rr))}{\rl+\rr}}, \quad \rr<\rl,\: \vr<\vl
\end {eqnarray}
The wave speed is
\begin {eqnarray}
s&=&\frac {\rr\vr-\rl\vl}{\rr-\rl} >0 .
\end {eqnarray}
The boundary averages $(\rfp,\vfp)$ are given in the following table:
\\\begin {center}
\begin {tabular} {||c||c|c|c||}\hline
&$s=\frac {\rr\vr-\rl\vl}{\rr-\rl}$ &$\rfp$&$\vfp$ \\\cline{2-4}
H2&$s>0$ & $\rl$ &$\vl$ \\\hline
\end {tabular}
\end {center}

\item The wave solution is a 1-rarefaction when the initial states satisfy
\begin {eqnarray}
\mbox{ R1: } \vr-\vl&=&\vast(\rr)-\vast(\rl), \quad \rr<\rl,\: \vr>\vl
\end {eqnarray}
The characteristic speed of a 1-rarefaction wave is
\begin {eqnarray}
\lambda_1(\rho,v)&=&v+\rho \vast'(\rho)
\end {eqnarray}
The boundary averages are the left state when $\lambda_1(\rl,\vl)>0$, similarly they are the right state when $\lambda_1(\rr,\vr)<0$. Otherwise, $(\rfp,\vfp)$  are the solutions of the equations
\begin {eqnarray}
\lambda_1(\rfp,\vfp)&=&\rfp \vast'(\rfp)+\vfp=0\label{a1}\\
\vfp-\vl &=& \vast(\rfp)-\vast(\rl). \label{a2}
\end {eqnarray}
We simplify equations \refet{a1}{a2} as
\begin {eqnarray}
\lambda_{\ast} (\rfp) &=&\vast(\rl)-\vl \equiv \Delta v \label {9c}\\
\vfp &=& \vast(\rfp) -\Delta v \label {9d}.
\end {eqnarray}

The boundary averages $(\rfp,\vfp)$ are given in the following table:
\\\begin {center}
\begin {tabular} {||c||c|c|c||}\hline
&$\lambda_1$ &$\rfp$&$\vfp$ \\\cline{2-4}
&$\lambda_1(\rl,\vl)>0$ & $\rl$ &$\vl$ \\\cline{2-4}
R1&$\lambda_1(\rr,\vr)<0$ & $\rr$ & $\vr$ \\\cline{2-4}
&o.w. & \multicolumn {2} {c||} {solution to (\ref{9c},\ref{9d}) } \\\hline
\end {tabular}
\end {center}

\item The wave solution is a 2-rarefaction when the initial states satisfy
\begin {eqnarray}
\mbox{ R2: } \vr-\vl&=&\vast(\rl)-\vast(\rr), \quad \rr>\rl,\: \vr>\vl
\end {eqnarray}
The characteristic speed of the 2-rarefaction wave is
\begin {eqnarray}
\lambda_2(\rho,v)&=&v-\rho \vast'(\rho)>0.
\end {eqnarray}

The solutions of $(\rfp,\vfp)$ are given in the following table:
\\\begin {center}
\begin {tabular} {||c||c|c|c||}\hline
&$\lambda_2$ &$\rfp$&$\vfp$ \\\cline{2-4}
R2&$\lambda_2>0$ & $\rl$ &$\vl$  \\\hline
\end {tabular}
\end {center}

\item The wave solution is a 1-rarefaction + 2-rarefaction when there exists an intermediate state $(\rmm,\vm)$ satisfying
\begin {eqnarray}
\mbox{ R1: } \vm-\vl&=&\vast(\rmm)-\vast(\rl), \quad \rmm<\rl,\: \vm>\vl \\
\mbox{ R2: } \vr-\vm&=&\vast(\rmm)-\vast(\rr), \quad \rr>\rmm,\: \vr>\vm.
\end {eqnarray}
That is to say, $\rmm$ satisfies
\begin {eqnarray}
2*\vast(\rmm)-\vast(\rl)-\vast(\rr)-(\vr-\vl)&=&0 \label{m1}
\end {eqnarray}
in which $\rmm<\rl,\rmm<\rr$. We can write $\vm$ as
\begin {eqnarray}
\vm=\vast(\rmm)+\vl-\vast(\rl).
\end {eqnarray}

The boundary averages $(\rfp,\vfp)$ are given in the following table:
\\\begin {center}
\begin {tabular} {||c||c|c|c||}\hline
&$\lambda_1$ &$\rfp$&$\vfp$ \\\cline{2-4}
&$\lambda_1(\rl,\vl)>0$ & $\rl$ &$\vl$ \\\cline{2-4}
R1-R2&$\lambda_1(\rmm,\vm)<0$ & $\rmm$ & $\vm$ \\\cline{2-4}
&o.w. & \multicolumn {2} {c||} {solution to (\ref{9c},\ref{9d}) } \\\hline
\end {tabular}
\end {center}

\item The wave solution is a 1-rarefaction + 2-shock when there exists an intermediate state $(\rmm,\vm)$ satisfying
\begin {eqnarray}
\mbox{ R1: } \vm-\vl=&\vast(\rmm)-\vast(\rl) &, \rmm<\rl,\: \vm>\vl \\
\mbox{ H2: } \vr-\vm=&-\sqrt {\frac {2(\rmm-\rr)(\phi(\rmm)-\phi(\rr))}{\rmm+\rr}} &, \rr<\rmm,\: \vr<\vm.
\end {eqnarray}
That is to say, $\rmm$ satisfies
\begin {eqnarray}
&\vast(\rmm)-\vast(\rl)-\sqrt {\frac {2(\rmm-\rr)(\phi(\rmm)-\phi(\rr))}{\rmm+\rr}}-(\vr-\vl)&
      \lefteqn{=0} \label{m2}
\end {eqnarray}
in which $\rr<\rmm<\rl$. We can write $\vm$ as
\begin {eqnarray}
\vm=\vast(\rmm)+\vl-\vast(\rl).
\end {eqnarray}

The boundary averages $(\rfp,\vfp)$ are given as the following
\\\begin {center}
\begin {tabular} {||c||c|c|c||}\hline
&$\lambda_1$ &$\rfp$&$\vfp$ \\\cline{2-4}
&$\lambda_1(\rl,\vl)>0$ & $\rl$ &$\vl$ \\\cline{2-4}
R1-H2&$\lambda_1(\rmm,\vm)<0$ & $\rmm$ & $\vm$ \\\cline{2-4}
&o.w. & \multicolumn {2} {c||} {solution to (\ref{9c},\ref{9d}) } \\\hline
\end {tabular}
\end {center}

\item The wave solution is a 1-shock + 2-shock when there exists an intermediate state $(\rmm,\vm)$ satisfying
\begin {eqnarray}
\mbox{ H1: } \vm-\vl=&-\sqrt {\frac {2(\rl-\rmm)(\phi(\rl)-\phi(\rmm))}{\rl+\rmm}} &, \rmm>\rl,\: \vm<\vl\\
\mbox{ H2: } \vr-\vm=&-\sqrt {\frac {2(\rmm-\rr)(\phi(\rmm)-\phi(\rr))}{\rmm+\rr}} &, \rr<\rmm,\: \vr<\vm.
\end {eqnarray}
That is to say, $\rmm$ satisfies
\begin {eqnarray}
&-\sqrt {\frac {2(\rl-\rmm)(\phi(\rl)-\phi(\rmm))}{\rl+\rmm}}-\sqrt {\frac {2(\rmm-\rr)(\phi(\rmm)-\phi(\rr))}{\rmm+\rr}}-(\vr-\vl)&
      \lefteqn{=0} \label{m3}
\end {eqnarray}
in which $\rmm>\rl,\rmm>\rr$. We can compute $\vm$ as
\begin {eqnarray}
\vm=-\sqrt {\frac {2(\rl-\rmm)(\phi(\rl)-\phi(\rmm))}{\rl+\rmm}}+\vl.
\end {eqnarray}

The boundary averages $(\rfp,\vfp)$ are given in the following table:
\\\begin {center}
\begin {tabular} {||c||c|c|c||}\hline
&$s=\frac {\rmm\vm-\rl\vl}{\rr-\rl}$ &$\rfp$&$\vfp$ \\\cline{2-4}
&$s>0$ & $\rl$ &$\vl$ \\\cline{2-4}
H1-H2&$s<0$ & $\rmm$ & $\vm$ \\\cline{2-4}
&$s=0$ & $\frac {\rl+\rmm} 2$ & $\frac {\vl+\vm}2$ \\\hline
\end {tabular}
\end {center}

\item The wave solution is a 1-shock + 2-rarefaction when there exists an intermediate state $(\rmm,\vm)$ satisfying
\begin {eqnarray}
\mbox{ H1: } \vm-\vl=&-\sqrt {\frac {2(\rl-\rmm)(\phi(\rl)-\phi(\rmm))}{\rl+\rmm}} &, \rmm>\rl,\: \vm<\vl\\
\mbox{ R2: } \vr-\vm=&\vast(\rmm)-\vast(\rr) &, \rr>\rmm,\: \vr>\vm.
\end {eqnarray}
That is to say, $\rmm$ satisfies
\begin {eqnarray}
&-\sqrt {\frac {2(\rl-\rmm)(\phi(\rl)-\phi(\rmm))}{\rl+\rmm}}+\vast(\rmm)-\vast(\rr)-(\vr-\vl)&
      \lefteqn{=0} \label{m4}
\end {eqnarray}
in which $\rmm>\rl,\rmm>\rr$. We compute $\vm$ as
\begin {eqnarray}
\vm=-\sqrt {\frac {2(\rl-\rmm)(\phi(\rl)-\phi(\rmm))}{\rl+\rmm}}+\vl.
\end {eqnarray}

The boundary averages $(\rfp,\vfp)$  are given in the following table:
\\\begin {center}
\begin {tabular} {||c||c|c|c||}\hline
&$s=\frac {\rmm\vm-\rl\vl}{\rr-\rl}$ &$\rfp$&$\vfp$ \\\cline{2-4}
&$s>0$ & $\rl$ &$\vl$ \\\cline{2-4}
H1-R2&$s<0$ & $\rmm$ & $\vm$ \\\cline{2-4}
&$s=0$ & $\frac {\rl+\rmm} 2$ & $\frac {\vl+\vm}2$ \\\hline
\end {tabular}
\end {center}

\end {enumerate}

\subsection {Some points concerning the implementation of the Godunov method}
In subsection \ref{sub3.2.1}  above, we studied the solutions of the boundary averages. We now discuss how to get a stable, convergent and efficient numerical method.

In a linear hyperbolic system, Godunov's method is stable and convergent if the Courant-Friedrichs-Lewy (CFL) number is less than unity. Similarly we require the CFL number be less than unity for Zhang's model; i.e.,
\begin {eqnarray}
\max \left |\frac k h \lambda_2(\rho,v) \right | &\leq & 1
\end {eqnarray}
where $\lambda_2(\rho,v)=v-\rho \vast'(\rho)$ has the bigger magnitude of two wave velocities.

There are two types of boundary conditions (BC) : Dirichlet BC and Neumann (or natural) BC. If the interval for $x$ is $[a,b]$, we will impose boundary condition on $u(a-\frac h2,t_j)$ and $u(b+\frac h2)$ instead of the real boundary $x=a$ and $x=b$. This is to say we have to solve a Riemann problem to get the fluxes on both the end boundaries instead of imposing boundary conditions directly on those fluxes. \footnote{For the detailed discussions on treatment of boundary conditions, refer to (Zhang, 2000a\nocite{Zhang2000a}).}

The source term, $s(u)=(0,\frac {\vast(\rho)-v}{\tau})$, doesn't involve spatial gradients and therefore remains bounded as terms on the left hand-side of \refe{zhang.1} go to $\infty$. We approximate the source term implicitly  with $s(U_i^{j+1})$ or $s(\frac {U_i^{j+1}+U_i^j}{2})$ in order to improve the stability property of the Godunov method. However, when the source term is stiff; i.e., when $\tau$ is small, the problem of numerical instability may arise. In this case we use a much smaller time increment $k\ll\tau$.

The solution to the Riemann problem is important  both theoretically and computationally. One can compute the numerical solutions when  all of the Riemann problems are well-posed and solvable at each step of the iteration. However, when the left and right states for a Riemann problem are far from each other, the intermediate state $(\r_m,v_m)$ for the wave solutions may be out of domain of validity, e.g., $\r_m<0$. In this case, we have a ``vacuum problem" and hence the numerical solutions can not be uniquely determined.

The cost of solving the Riemann problem determines the computational efficiency of the numerical method. Here we propose improvements to the computational efficiency for Zhang's model.
In Godunov's methods for Zhang's model, most of the calculations are basic arithmetic computations except the calculation of the intermediate state $(\rmm,\vm)$ in equations \refe{m1},\refe{m2},\refe{m3} and \refe{m4} and the solutions of equations \refe{9c} and \refe{9d}).  Given $(\r_l,v_l)$ and $(\r_r,v_r)$, the nonlinear algebraic equations \refe{m1},\refe{m2},\refe{m3} and \refe{m4} can all be written in the form of $g(\r_m)=0$. The functions $g(\r_m)$ for these equations are all monotonically decreasing in the interval of validity of $\r_m$.
To find $\rfp$ from equation \refe{9c}, we define a function
\begin {eqnarray}
g(\rho)&=&\lambda_{\ast} (\rho) -\Delta v, \quad  \rr \: (\rmm \mbox{ for R1-R2,R1-H2})<\rho<\rl .
\end {eqnarray}
Since $f_{\ast}(\rho)$ is concave and $\lambda_{\ast}(\rho)=f_{\ast}'(\rho)$, we find $\lambda_{\ast}(\r)$ and $g(\r)$ are decreasing. Here we chose secant method to solve these equations since it is very efficient when the associated functions are monotonically decreasing.

\section {A Second-order Godunov Method}\label{secondGodunov}
In this section we introduce a second-order finite difference method for Zhang's model. This method is a two-stage predictor/corrector method. In this method, \refe{zhang.1} is decoupled into two nonlinear scalar equations in each interval $(x_{i-1/2},x_{i+1/2})$ at time $t_j$ and the predictor/corrector procedures are applied to those scalar functions.

We can write Zhang's model as:
\begin {eqnarray}
u_t+A(u) u_x &=& s(u)
\end {eqnarray}
where
\begin {eqnarray}
u&=&\left ( \begin {array} {c}
\rho(x,t)\\
v(x,t)
\end {array} \right )
\end {eqnarray}
and
\begin {eqnarray}
A(u)&=&\matauZ
\end {eqnarray}

The two eigenvalues and corresponding eigenvectors are
\begin {eqnarray}
\begin {array}{ll}
\lambda_1(u)=v+\rho\vast(\rho),&r_1(u)=[1,\vast'(\rho)]^t,\\
\lambda_2(u)=v-\rho\vast'(\rho),&r_2(u)=[1,-\vast'(\rho)]^t.
\end {array}
\end {eqnarray}

We diagonalize $A(u)$ by
\begin {eqnarray}
T^{-1}(u)A(u)T(u)&=&\matluZ\equiv\Lambda(u),
\end {eqnarray}
where the transformation matrix $T(u)$ is
\begin {eqnarray}
T(u)&=&\mattuZ.
\end {eqnarray}
Letting $W=T^{-1}(u)u$, Zhang's model \refe{zhang.1} under the transformation becomes
\begin {eqnarray}
W_t+\Lambda(u)W_x&=& T^{-1}(u)s(u). \label{second.1}
\end {eqnarray}

For the solution $w(x,t)$ to a scalar equation $w_t+\lambda(w) w_x= 0$, the first-order Godunov method uses a step function $w_I(x,t_j)$ to interpolate the solution, i.e.,
\begin {eqnarray}
w_I(x,t_j) &=& w_i^j, \mbox {  if } x_{i-1/2}<x\leq x_{i+1/2}.
\end {eqnarray}
In a first-order Godunov method, the Riemann problem has the following jump initial conditions:
\begin {eqnarray}
\begin {array} {lcl}
w_{i+1/2}^{j,L}& =&w_i^j \\
w_{i-1/2}^{j,R} &=& w_i^j.
\end {array}
\end {eqnarray}
For a second-order Godunov method, we interpolate the initial condition with a piecewise linear function
\begin {eqnarray}
w_I(x,t_j) &=& w_i^j+\frac {(x-ih)}{h} \Delta^{VL} w_i^j, \mbox {  if } x_{i-1/2}<x\leq x_{i+1/2}.
\end {eqnarray}
Then we do a half-step prediction
\begin {eqnarray}
\begin {array} {lcl}
w_{i+1/2}^{j+1/2,L}& =& w_i^j+\frac 12 (1-\lambda(w_i^j)\frac k h) \Delta ^{VL} w_i^j \\
w_{i-1/2}^{j+1/2,R} &=& w_i^j-\frac 12 (1+\lambda(w_i^j)\frac k h) \Delta ^{VL} w_i^j,
\end {array}
\end {eqnarray}
where $\Delta ^{VL} w_i^j$ is the van Leer slope defined as (all subscripts $j$ have been suppressed)
\begin {eqnarray*}
&\Delta ^{VL}w_i& = \left \{\begin {array} {l}
S_i\cdot \min(2|w_{i+1}-w_i|,2|w_i-w_{i-1}|,\frac 12|w_{i+1}-w_{i-1}|) ,\: \lefteqn{\xi>0} \\
0, \qquad\qquad\qquad \mbox{ otherwise }
\end {array}\right.
\end {eqnarray*}
\begin {eqnarray}
S_i&=& \mbox{sign} (w_{i+1}-w_{i-1})\\
\xi&=&(w_{i+1}-w_i)\cdot(w_i-w_{i-1})
\end {eqnarray}
The van Leer slope limiter ensures that the method remains second order when the solution $w(x,t)$ is smooth  and eliminates Gibb's phenomenon at discontinuities.

We apply the procedure above to the two scalar equations of the related homogeneous $2\times 2$ system in \refe{second.1} to obtain $W_{i+1/2}^{j+1/2,L}$ and $W_{i-1/2}^{j+1/2,R}$. Given the half-step values of $W_{i+1/2}^{j+1/2,L}$ and $W_{i-1/2}^{j+1/2,R}$,  $U_{i+1/2}^{j+1/2,L}$ and $U_{i-1/2}^{j+1/2,R}$ can be calculated by an inverse transformation:
\begin {eqnarray}
\begin {array} {ccc}
U_{i+1/2}^{j+1/2,L}&=&T(U_i^j) W_{i+1/2}^{j+1/2,L} \\
U_{i-1/2}^{j+1/2,R}&=&T(U_i^j) W_{i-1/2}^{j+1/2,R}
\end {array}
\end {eqnarray}
We then solve the Riemann problem to find the boundary averages.

\section {Numerical Solutions of Zhang's model}
%Tue Nov  2 14:10:13 PST 1999
Based on the discussions in the former sections, we carry out some numerical computations to test the validity of the Godunov method and the properties of Zhang's model.

Here we use Newell's equilibrium model,
\bqn
v_{\ast}(\rho)&=&v_f\left(1-\exp\{\frac {|c_j|}{v_f}(1-\rho_j/\rho)\}\right).
\eqn
 and set $v_f=1, c_j=1, \rho_j=1$ to get the standardized equilibrium relationship:
\bqn
v_{\ast}(\rho)&=&1-\exp{(1-\frac 1{\rho})}
\eqn
The domain for $\rho$ is $0\leq\rho\leq 1$\footnote{We use the limit of $\vs$ when $\r\to 0$ as its value at $\r=0$.}, the range for $v_{\ast}$ is the same. The traffic flow rate is $f_{\ast}=\rho v_{\ast}$. The equilibrium travel speed $\vs(\r)$ and flow rate $f_{\ast}(\r)$ are shown in \reff{zfunc}.

The subcharacteristic; i.e., the first-order characteristic velocity, (i.e., the wave velocity of the corresponding LWR model) is
\bdm
\lambda_{\ast}=1-(1+\frac 1{\rho})\exp { (1-\frac 1{\rho})},
\edm
and the eigenvalues are
\bdm
\lambda_{1,2}=v\pm \frac 1{\rho} \exp {(1-\frac 1{\rho})}.
\edm

When $0\leq\rho\leq 1$, we have
\bdm
|\lambda_2|<|v|+1.
\edm
The CFL condition number is defined as:
\bqn
\max|\frac kh\lambda_2(\rho,v)|\leq \frac kh(\max v+1).
\eqn
Since the CFL number is no larger than 1, we find
\bqn
k\leq \frac h{\max v+1}.
\eqn

In all of the computations that follow, we let $x\in[0 l,800 l]$, where $l$ is the unit of length. Here the number of grid points is denoted by $N$, and $h=\frac {800 l} N$ is the space step. We let $T_0=K \tau$ denote the final time ($\tau$ is the unit of time), and $m$ the number of time steps and $k=\frac {K \tau}m$. From the CFL condition, we have
\bqn
\frac{(\max v+1)K}{800} \frac N m \frac {\tau}l\leq 1.
\eqn
Setting $\tau=l=10.0$, and $m=N$, CFL condition is not violated when $K\leq 400$ since $\max v<1$.

\bfg
\bc\includegraphics[height=12cm] {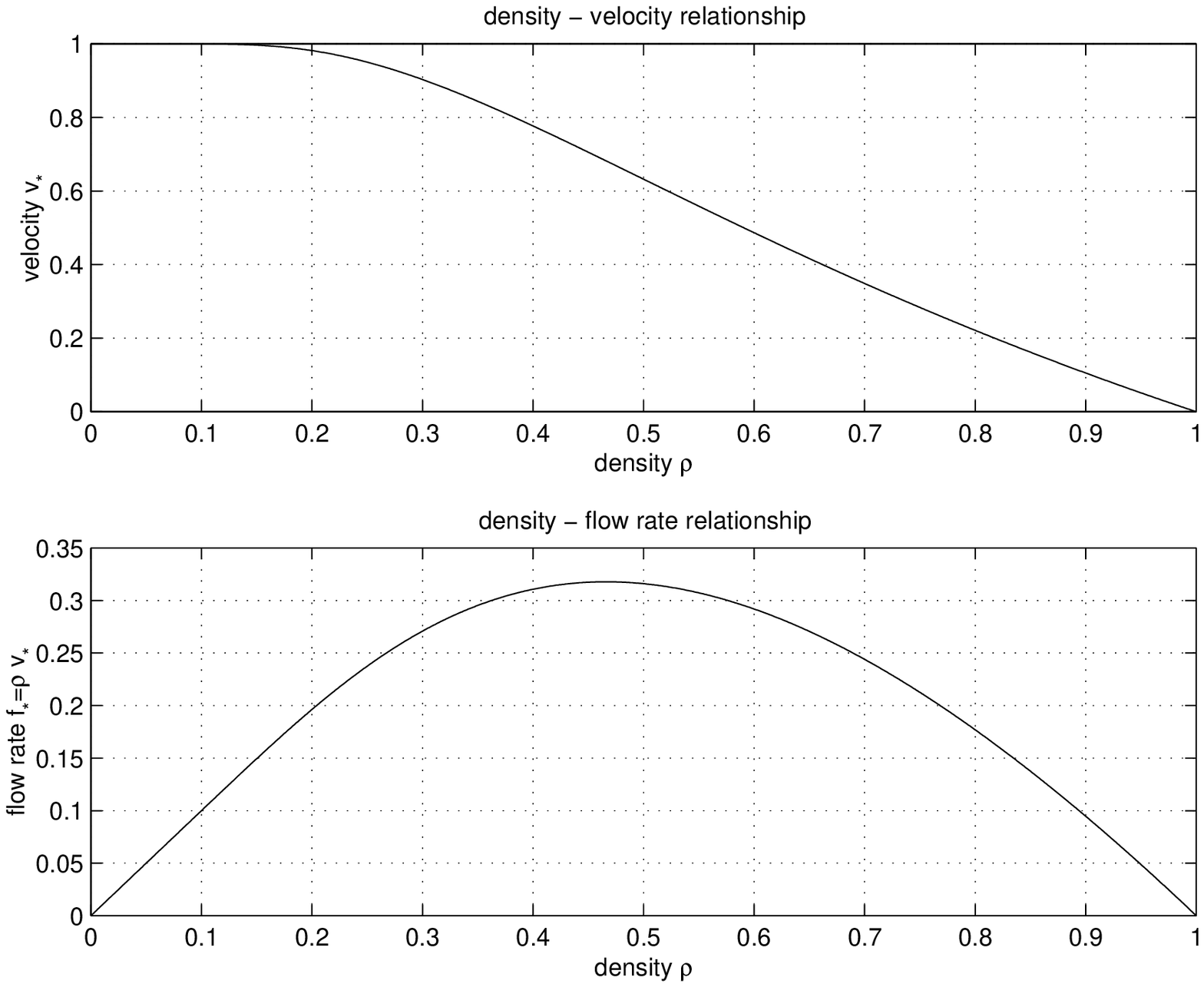}\ec
\caption {Newell's Equilibrium $\r$--$\vs$/ $\r$--$f_{\ast}$ Relationship} \label {zfunc}
\efg

\subsection {Riemann solutions}
%Thu Dec  9 12:03:39 PST 1999
In the following four computations we use the first-order Godunov method to examine different types of waves solutions. With four well-chosen jump initial conditions, we can observe four different type of waves H1-H2, R1-R2, R1-H2 and H1-R2. To prevent the 2-waves from relaxing to 1-waves in a short time, here we rescale the relaxation time $\tau\to1000\tau$. For each computation we present  a contour plot  of both $\r$ and $v$ until $T_0=400\tau$, and several 2-D curves at different time $t$ have been selected from the contour plot.

\bi
\item[Computation 1] We use the following  initial conditions:
\bqn
\rho(x,0)&=&0.65, \qquad \forall x\in [0l, 800l]\\
v(x,0)&=&\cas {{ll} v_{\ast}(0.65)& x\in [0l, 200l]\\v_{\ast}(0.65)-0.2& \m{otherwise}}
\eqn
The solutions are of H1-H2 type, shown in \refft{Riemann1_3d}{Riemann1_2d}. These figures show that the downstream travel speed keeps increasing along with the propagation of the 2-shock wave. This is due to the effect of the relaxation term. We can predict that when the downstream travel speed  reaches $\vs(0.65)$, which is the equilibrium travel speed with respect to the downstream traffic density $\r=0.65$, the 2-shock disappears and a 1-rarefaction wave forms.

\item[Computation 2] We use the following initial conditions:
\bqn
\rho(x,0)&=&0.65, \qquad \forall x\in [0l, 800l]\\
v(x,0)&=&\cas {{ll} v_{\ast}(0.65)& x\in [0l, 200l]\\v_{\ast}(0.65)+0.2& \m{otherwise}}
\eqn
The solutions are of R1-R2 type, shown in \refft{Riemann2_3d}{Riemann2_2d}. These figures show that the downstream travel speed keeps decreasing along with the propagation of the 2-rarefaction wave. This is also due to the effect of the relaxation term, which will relax the 2-rarefaction wave to a 1-shock wave.

\item[Computation 3] We use the following  initial conditions:
\bqn
\rho(x,0)&=&\cas {{ll} 0.65& x\in [0l, 200l]\\0.4& \m{otherwise}}\\
v(x,0)&=&v_{\ast}(0.65), \qquad \forall x\in [0l, 800l]
\eqn
The solutions are of R1-H2 type, shown in \refft{Riemann3_3d}{Riemann3_2d}. Here the effect of the relaxation term will relax the 2-shock wave to a 1-rarefaction wave.

\item[Computation 4] We use the following initial conditions:
\bqn
\rho(x,0)&=&\cas {{ll} 0.65& x\in [0l, 200l]\\0.9& \m{otherwise}}\\
v(x,0)&=&v_{\ast}(0.65), \qquad \forall x\in [0l, 800l]
\eqn
The the solutions are of H1-R2 type, shown in \refft{Riemann4_3d}{Riemann4_2d}. Here the effect of the relaxation term will relax the 2-rarefaction wave to a 1-shock wave.
\ei

\bfg
\bc\includegraphics[height=8cm] {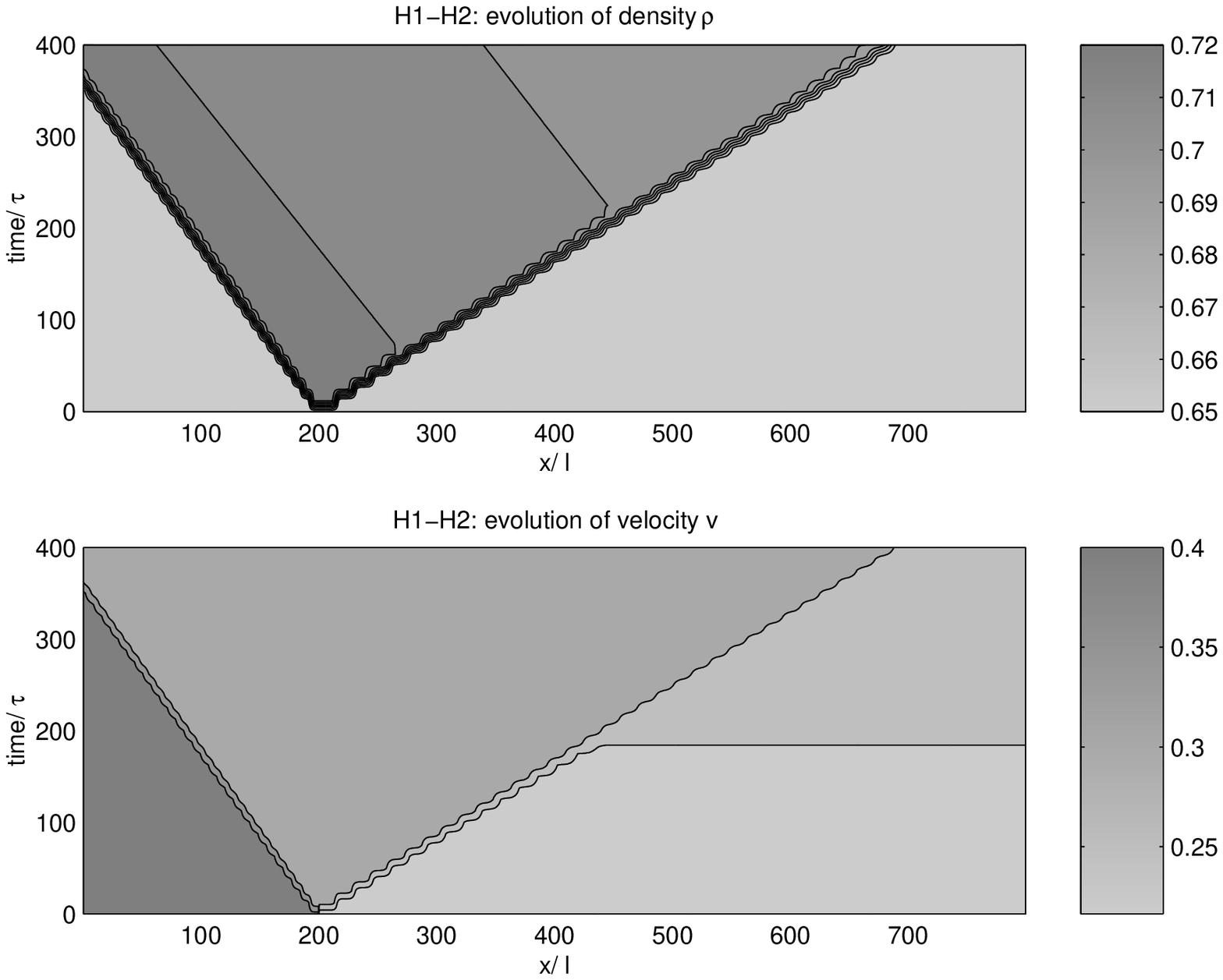}\ec
\caption {Left-Shock-Right-Shock Wave Solution} \label {Riemann1_3d}
\efg
\bfg
\bc\includegraphics[height=8cm] {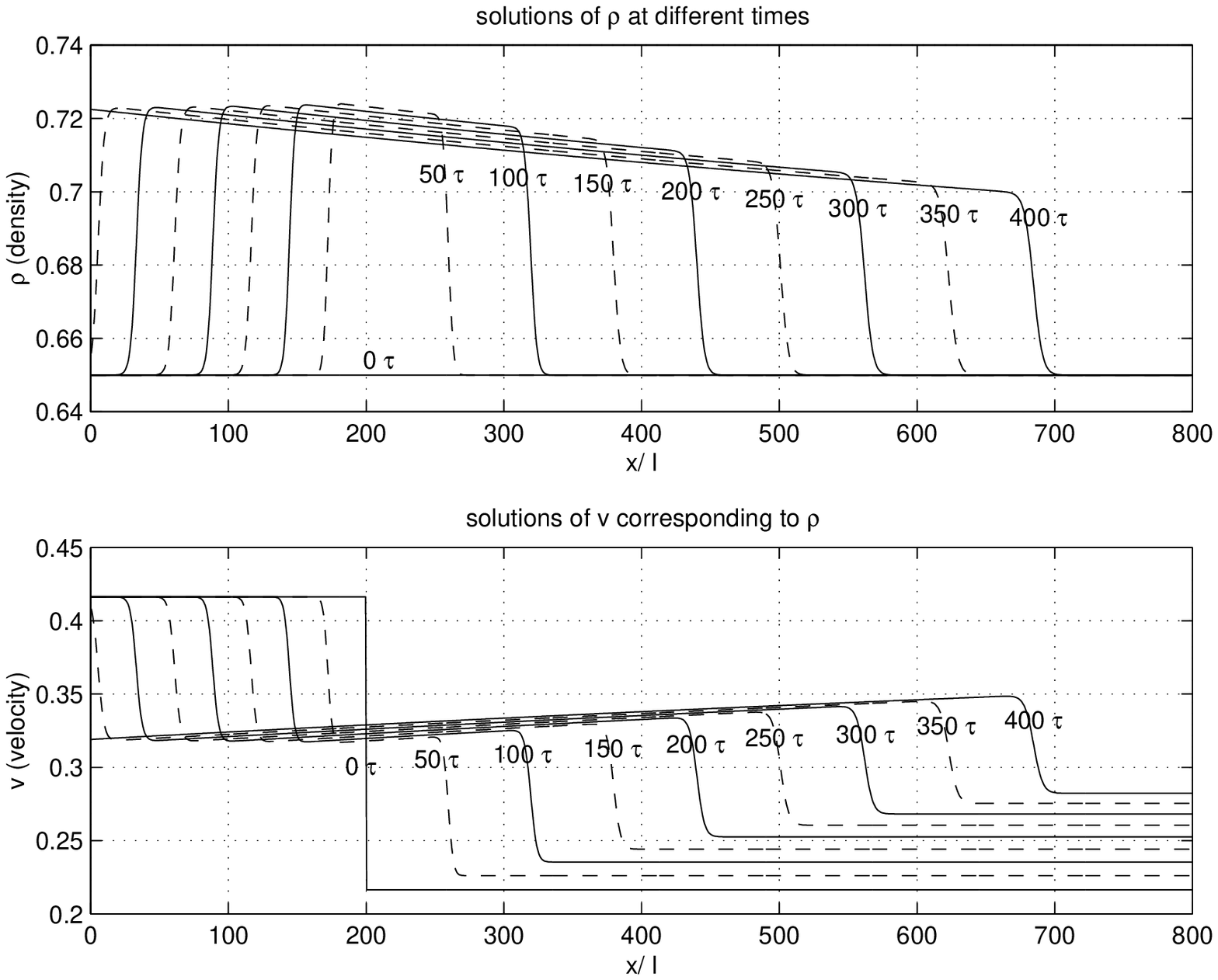}\ec
\caption {H1-H2 Solutions from  at Selected Times} \label {Riemann1_2d}
\efg

\bfg
\bc\includegraphics[height=8cm] {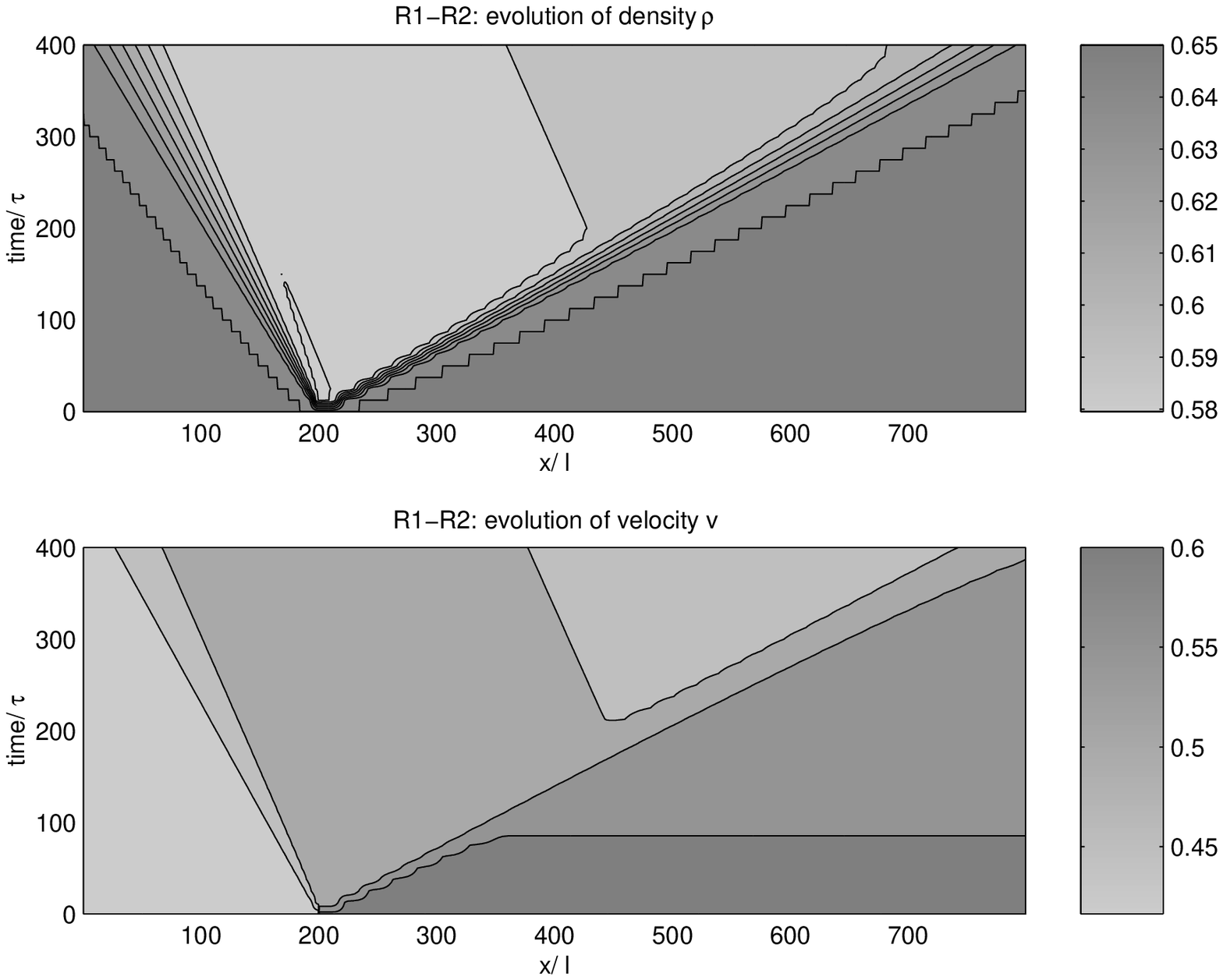}\ec
\caption {Left-Rarefaction-Right-Rarefaction Wave Solution} \label {Riemann2_3d}
\efg
\bfg
\bc\includegraphics[height=8cm] {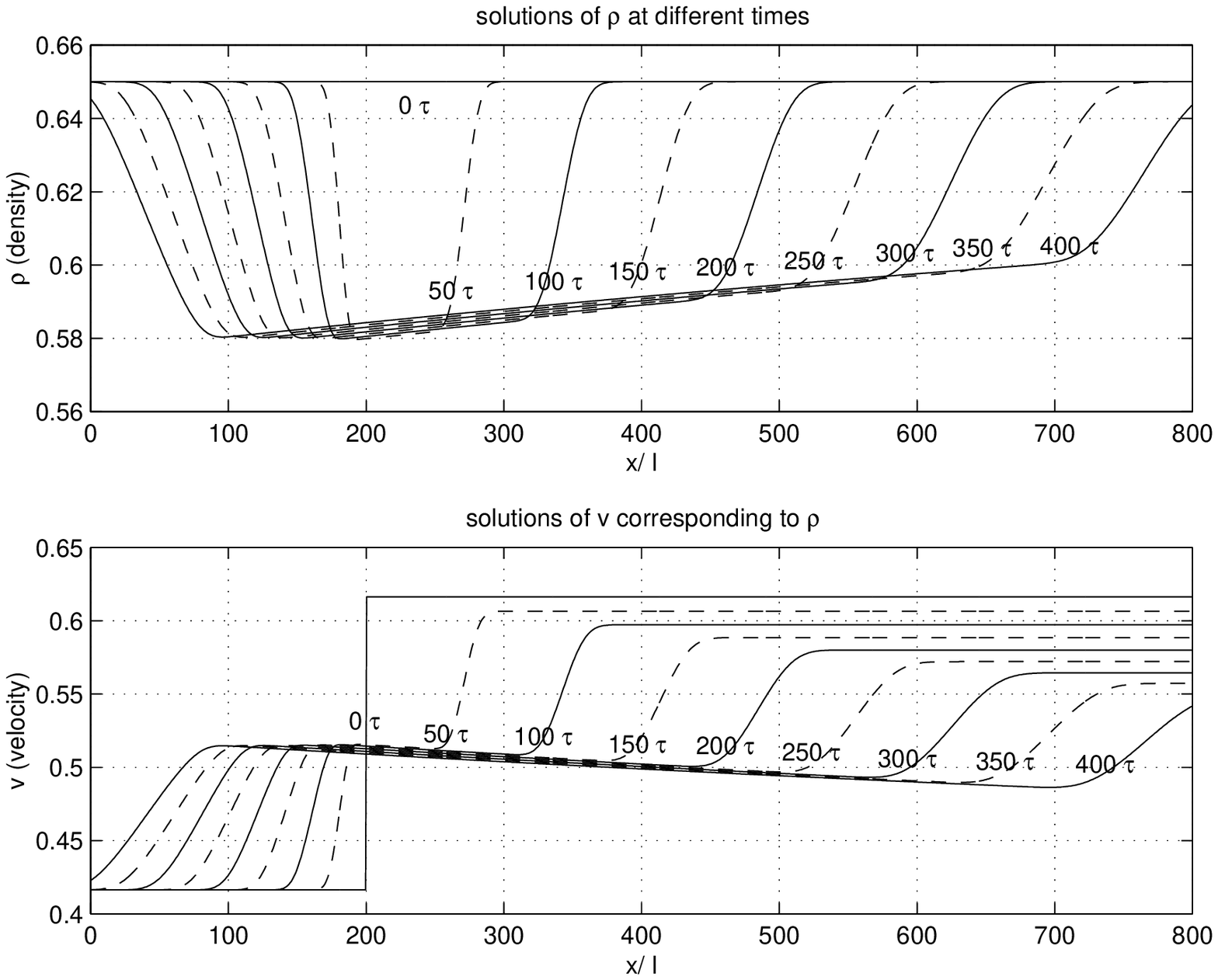}\ec
\caption {R1-R2 Solutions at Selected Times} \label {Riemann2_2d}
\efg

\bfg
\bc\includegraphics[height=8cm] {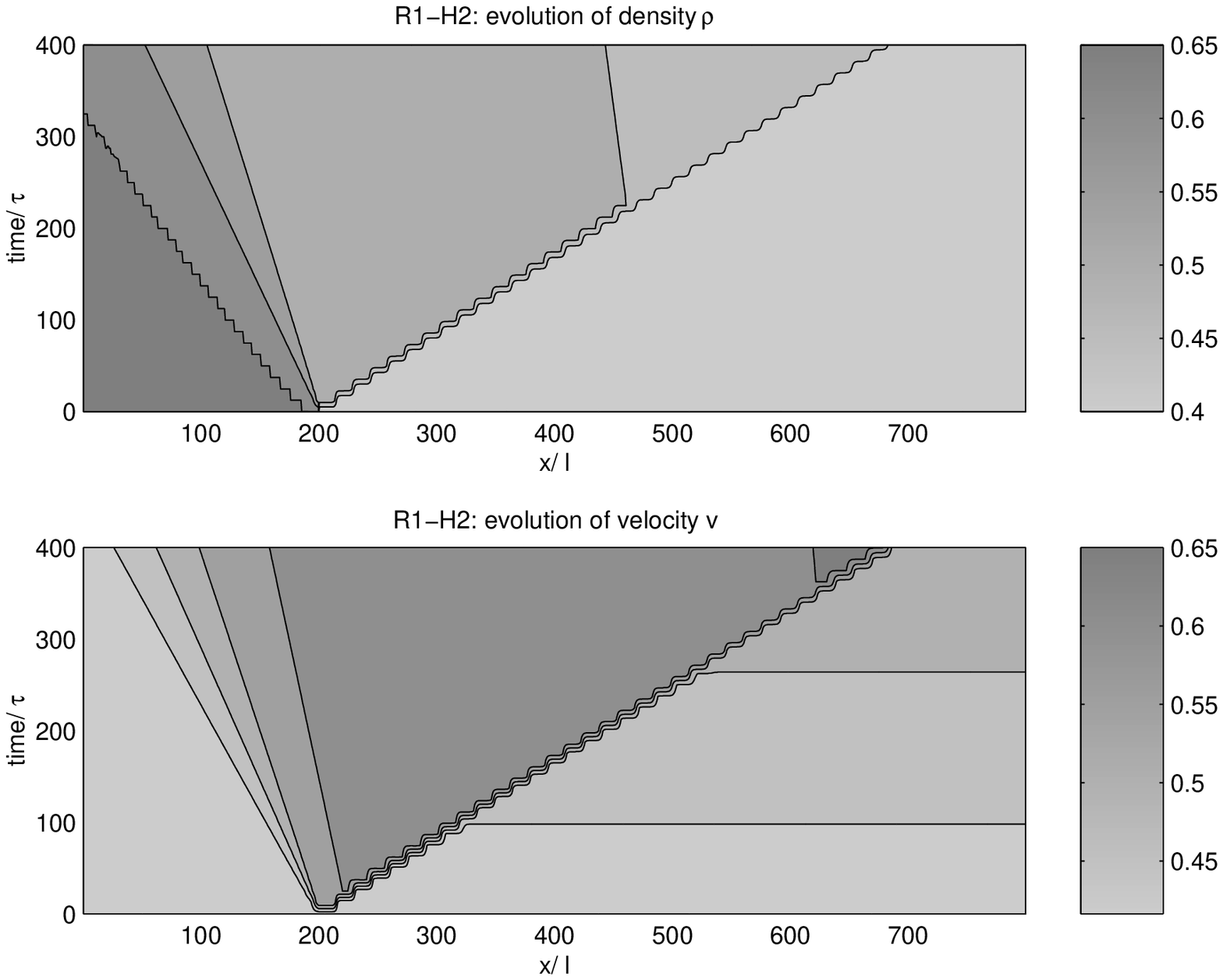}\ec
\caption {Left-Rarefaction-Right-Shock Wave Solution} \label {Riemann3_3d}
\efg
\bfg
\bc\includegraphics[height=8cm] {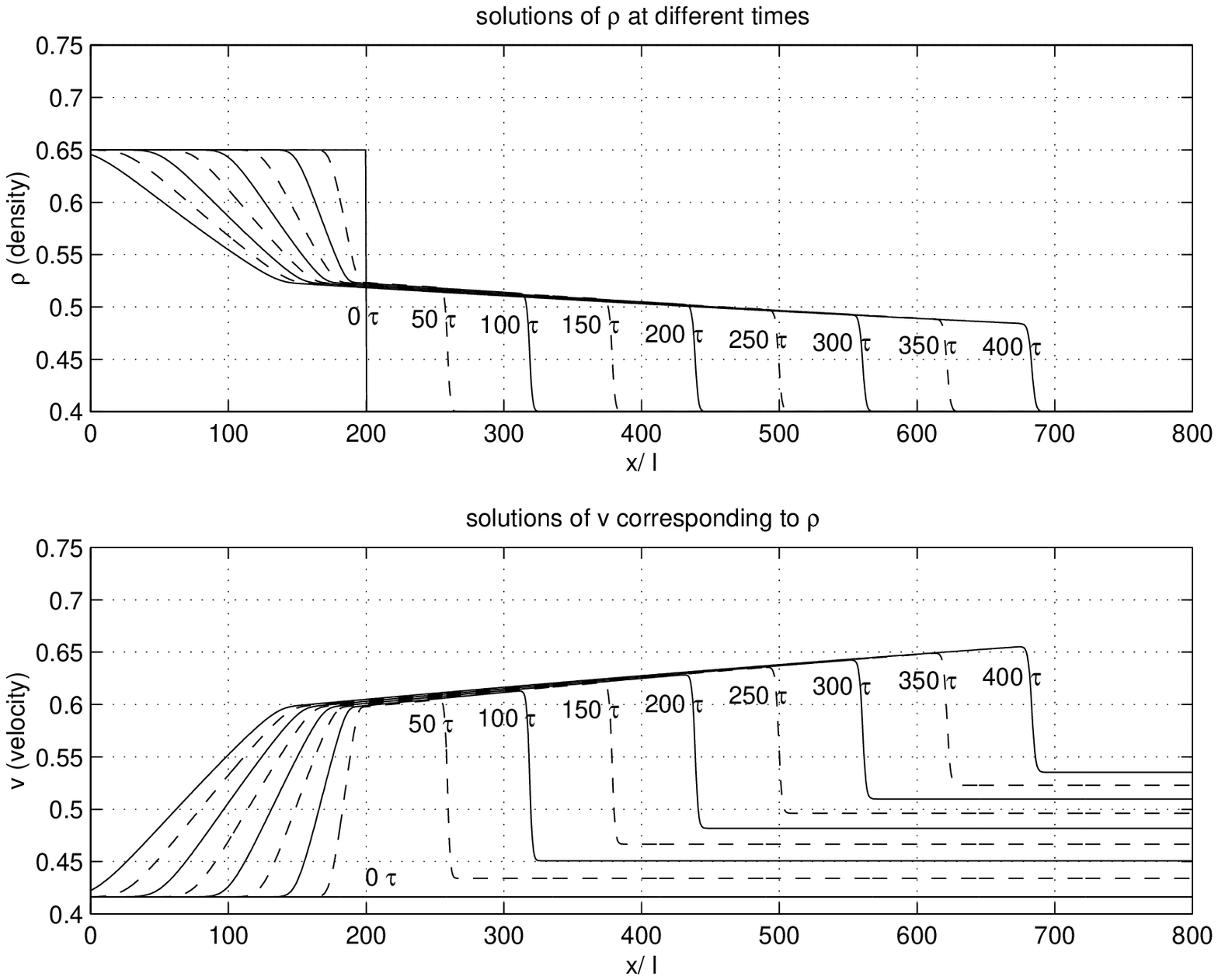}\ec
\caption {R1-H2 Solutions at Selected Times} \label {Riemann3_2d}
\efg

\bfg
\bc\includegraphics[height=8cm] {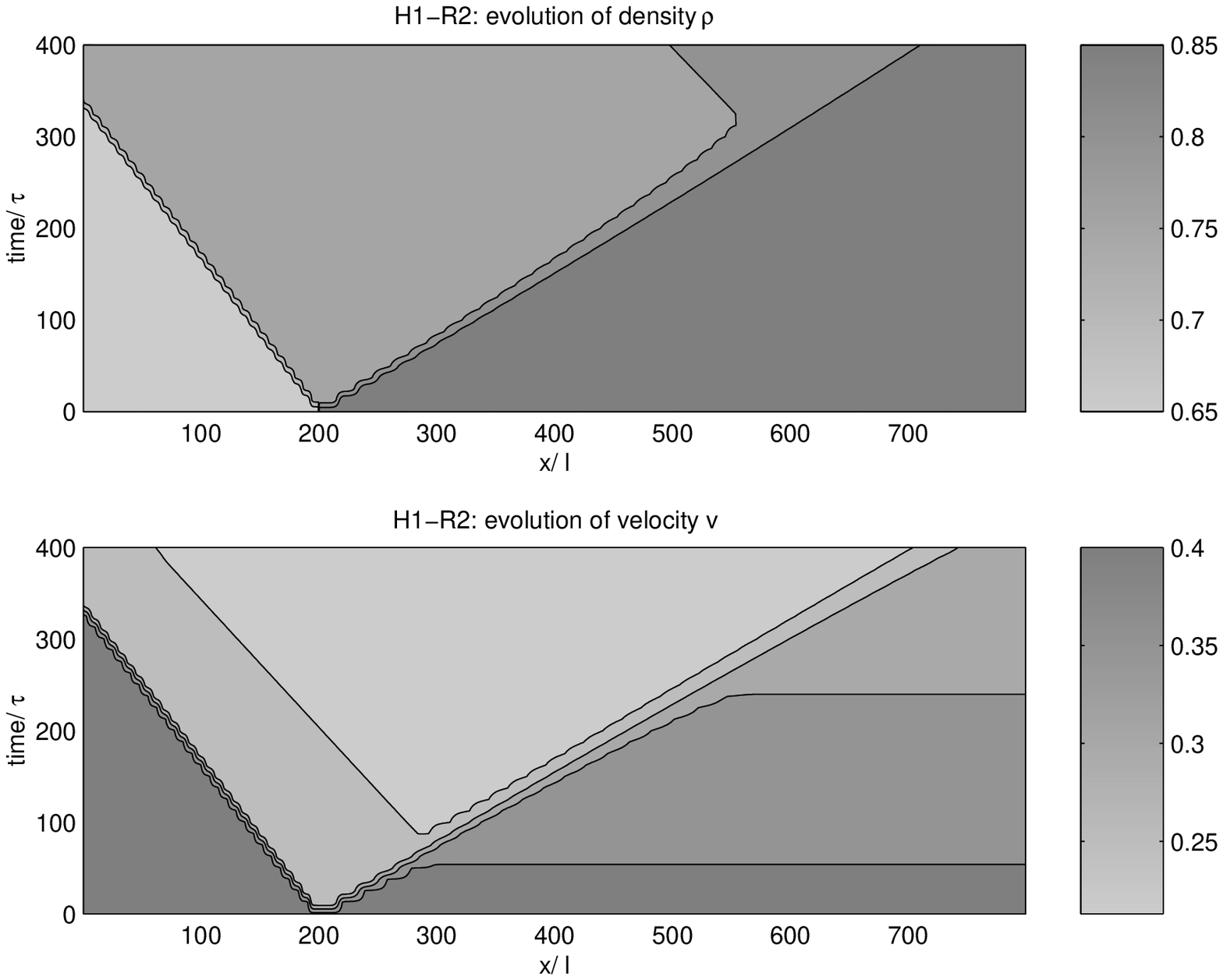}\ec
\caption {Left-Shock-Right-Rarefaction Wave Solution} \label {Riemann4_3d}
\efg
\bfg
\bc\includegraphics[height=8cm] {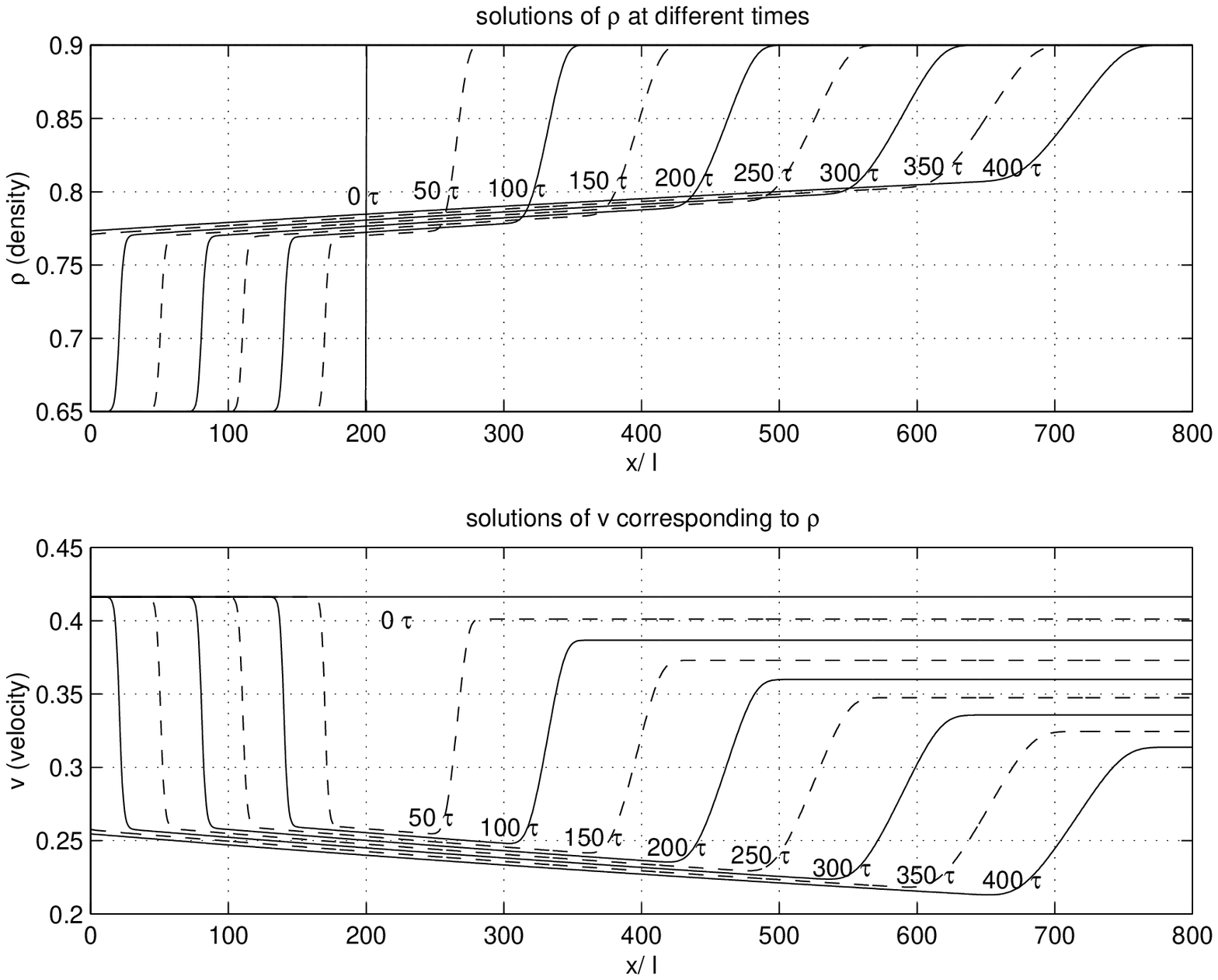}\ec
\caption {H1-R2 Solutions at Selected Times} \label {Riemann4_2d}
\efg

\subsection {A General Solution and Convergence Rates}
In this subsection the relaxation time is no longer rescaled since we wish to observe the solutions to \refe{zhang.1} for a longer time. Here we are interested in the solution for $0\leq t\leq T_0=400\tau$.

Using the general initial condition
\bqn
\rho(x,0)&=&0.65+\sin {(\frac {2\pi x} {800l})/4},\label{ini.1}\\
v(x,0)&=&v_{\ast}(\rho(x,0))+0.1, \qquad \forall x\in [0l, 800l],\label{ini.2}
\eqn
we use the first-order Godunov method to obtain the solutions shown in \refft{general1_3d}{general1_2d}. From these figures we see a shock wave forms at the downstream section, and a complicated combination of rarefaction waves forms in the upstream section. The solutions also show that $\r$-$v$ are in equilibrium by $t=100\tau$, i.e., $v=\vs(\r)$, due to the effect of the relaxation term, although the initial condition is not in equilibrium.

Next, we compute the convergence rate for the first- and second-order methods with Neumann boundary conditions with initial conditions \refet{ini.1}{ini.2}. The convergence rate is calculated from the comparison between the solution on different grids.  The grid numbers $2N$ and $N$ generate different grid sizes: $\frac h 2$ and $h$.  We denote the solutions at time $T_0$ as $(U^{2N}_i)_{i=1}^{2N}$ and $(U^N_i)_{i=1}^N$ respectively, and define the difference vector $(e^{2N-N})_{i=1}^N$ between these two solutions as
\bqn
\textbf e_i^{2N-N}&=&\frac12 (U^{2N}_{2i-1}+U^{2N}_{2i})-U^{N}_i, i=1,\cdots,N \label{def:error}.
\eqn
Then the relative error between the two solutions is defined as the norm of the difference vector:
\bqn
\epsilon^{2N-N}&=&\nm{{\textbf e}^{2N-N}}. \label{epsilon}
\eqn
The convergence rate is defined as
\bqn
r&=&\log_2(\frac{\epsilon^{2N-N}}{\epsilon^{4N-2N}}).\label{rate}
\eqn
In \refe{epsilon}, the norm can be $L^1$-, $L^2$- or $L^{\infty}$-norm.

For $N$ equal to 64, 128, 256, 512 and 1024, the relative errors and convergence rates for the first-order method are given in \reft{first}, and those for the second-order Godunov's method are computed and given in \reft {second}. For the first-order Godunov's method, the convergence rates related to $L^1$-norm errors are around 1, which is larger than that related to $L^2$-norm errors and even larger than that related $L^{\infty}$-norm errors. From \reft{first} we also see that the rates for $\rho$ and $v$ are consistent since $\rho$ and $v$ are in equilibrium at $400\tau$.
In \reft{second}, the convergence rates for the second-order method, although better slightly than those for first-order method, are not second order due to the effect of the relaxation term.

\bfg
\bc\includegraphics[height=8cm] {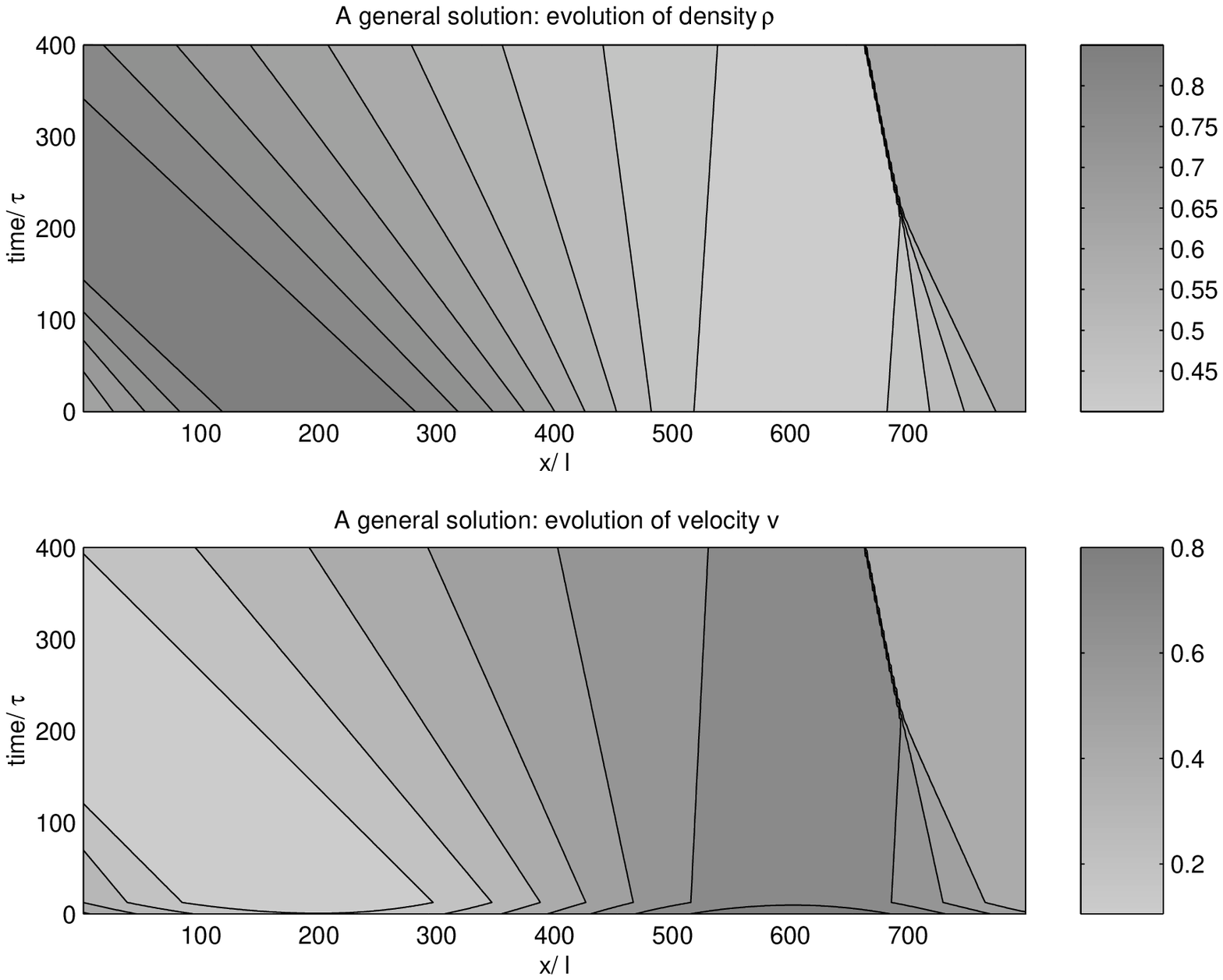}\ec
\caption {A General Solution of Zhang's Model with Neumann BC} \label {general1_3d}
\efg
\bfg
\bc\includegraphics[height=8cm] {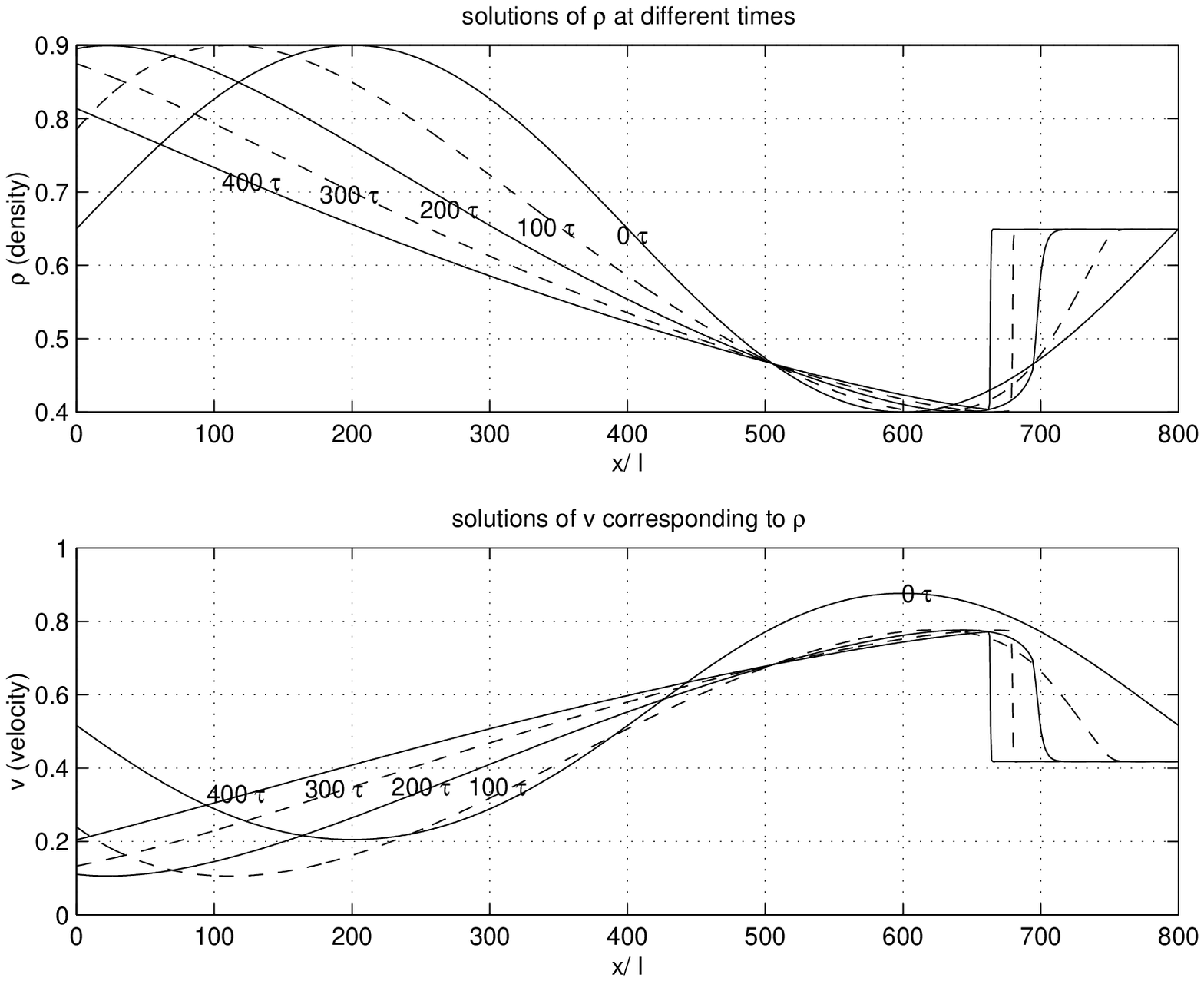}\ec
\caption {The Solutions from \reff{general1_3d} at Selected Times} \label {general1_2d}
\efg

%Wed Dec 22 15:26:26 PST 1999
\begin {table} [p]
\begin{tabular}{lccccccc}\\\hline
$\rho$&128-64&Rate&256-128&Rate&512-256&Rate&1024-512\\\hline
$L^1$&3.43e-03&0.969&1.75e-03&0.985&8.84e-04&0.994&4.44e-04\\\hline
$L^2$&6.94e-03&0.632&4.48e-03&0.635&2.88e-03&0.718&1.75e-03\\\hline
$L^{\infty}$&4.56e-02&0.0219&4.49e-02&0.0973&4.20e-02&0.355&3.28e-02\\\hline\hline
$v$&128-64&Rate&256-128&Rate&512-256&Rate&1024-512\\\hline
$L^1$&4.83e-03&0.973&2.46e-03&0.985&1.24e-03&0.991&6.25e-04\\\hline
$L^2$&9.86e-03&0.631&6.37e-03&0.628&4.12e-03&0.695&2.54e-03\\\hline
$L^{\infty}$&6.53e-02&0.0243&6.42e-02&0.0903&6.03e-02&0.313&4.86e-02\\\hline
\end {tabular}
\caption {Convergence rates for the first-order method with initial conditions \refet{ini.1}{ini.2}} \label {first}
\etb

\begin {table} [p]
\begin{tabular}{lccccccc}\\\hline
$\rho$&128-64&Rate&256-128&Rate&512-256&Rate&1024-512\\\hline
$L^1$&5.81e-03&0.966&2.98e-03&1.03&1.46e-03&1.05&7.06e-04\\\hline
$L^2$&9.85e-03&0.823&5.57e-03&0.836&3.12e-03&0.885&1.69e-03\\\hline
$L^{\infty}$&3.73e-02&0.215&3.22e-02&0.232&2.74e-02&0.616&1.79e-02\\\hline\hline
$v$&128-64&Rate&256-128&Rate&512-256&Rate&1024-512\\\hline
$L^1$&8.24e-03&0.973&4.20e-03&1.03&2.05e-03&1.04&9.93e-04\\\hline
$L^2$&1.40e-02&0.842&7.80e-03&0.865&4.28e-03&0.833&2.40e-03\\\hline
$L^{\infty}$&5.40e-02&0.250&4.54e-02&0.354&3.56e-02&0.527&2.47e-02\\\hline
\end {tabular}
\caption {Convergence rates for the second-order method with initial conditions \refe{ini.1}{ini.2}} \label {second}
\etb

\newpage
\pagestyle{myheadings} 
\markright{  \rm \normalsize CHAPTER 4. \hspace{0.5cm}
 The PW Model and Its Numerical Solutions}
\large
\chapter{The PW Model and Its Numerical Solutions}
\section {Introduction}
The Payne-Whitham (PW) model, suggested by Payne (1971\nocite {payne}) and Whitham (1974\nocite {Whitham}), is one of the first non-equilibrium traffic models. It can be written as 
\bqn
\rho_t+\rho v_x+v \rho_x&=&0 \label {eqn:1a} \\
v_t+vv_x+\frac {c_0^2}{\rho} \rho_x&=&\frac{v_{\ast}(\rho)-v}{\tau}\label {eqn:1b},
\eqn
where the traffic sound speed $c_0>0$ is constant, and $v_{\ast}(\rho)$ is the relationship between velocity $v$ and density $\rho$ for equilibrium states. The fundamental diagram $f_{\ast}(\rho)=\rho v_{\ast}(\rho)$, reflecting basic features of a roadway, is assumed to be known. Setting $m=\rho v$, \refe{eqn:1a} and \refe{eqn:1b} can be written as
\bqn
\matu_t+\matf_x=\mats \label {eqn:cons},
\eqn
which is in conservative form
\bqn
u_t+f(u)_x&=&s(u) \label {eqn:gene},
\eqn
with state $u=(\rho, m)$ of traffic flow.  

For $c_0=\rho v'_{\ast}(\rho)$, system \refe{eqn:cons} becomes Zhang's model, which is discussed in Chapter 3. Both PW model and Zhang's model are ``systems of hyperbolic conservation laws with relaxation'' in the sense of Whitham (1959\nocite {Whitham59}; 1974\nocite {Whitham}) and Liu (1987\nocite {Liu87}). 

Schochet (1988\nocite{schochet1988}) has shown that, as $\tau\to 0$, the system \refe{eqn:cons} admits the limit 
\bqn
\rho_t+(\rho \vast)_x&=&\nu\frac{\partial^2 \r}{\partial x^2}. \label{lim}
\eqn
Note that \refe{lim} is the LWR model with viscous right handside.

Generally, the equilibrium velocity is assumed to be decreasing with respect to density; i.e., $\vast'(\rho)<0$; the equilibrium flow rate is concave; i.e., $f_{\ast}''(\rho)<0$. Equation \refe{eqn:cons} has three different wave velocities (relative to the road): $\lambda_{\ast}$, $\lambda_1$ and $\lambda_2$. The first-order wave speed $\lambda_{\ast}$ is 
\begin {eqnarray}
\lambda_{\ast} (\rho)&=& \vast(\rho)+\rho \vast'(\rho).
\end {eqnarray}
$\lambda_{\ast}$ can be positive or negative. Since the wave speed of the degenerate system (i.e., the LWR model), it is called a sub-characteristic.
The second-order wave speeds, or ``frozen characteristic speeds'' (Pember, 1993\nocite {Pember93a}), are 
\begin {eqnarray}
\lambda_{1}(\rho,v)&=&v- c_0,\\
\lambda_{2}(\rho,v)&=&v+ c_0. 
\end {eqnarray}
Since $v,c_0\geq 0$, $\lambda_1<\lambda_2$ and $v+c_0>0$. 

Whitham (1959\nocite{Whitham59}; 1974\nocite{Whitham}) showed that the stability condition for the linearized system with a relaxation term is 
\bqn
\lambda_1< \lambda_{\ast} < \lambda_2  \qquad \m { when } t=0 \label{stable:con}.
\eqn
Liu (1987\nocite{Liu87}) showed that if condition \refe{stable:con} is always satisfied, then the corresponding LWR model is stable under small perturbations and the time-asymptotic solutions of the system \refe{eqn:cons} are completed determined by the equilibrium LWR model. Chen, Levermore, and Liu (1994\nocite{Chen94}) showed that if condition \refe{stable:con} is satisfied, then solutions of the system tend to solutions of the equilibrium equations as the relaxation time tends to zero. Besides \refe{stable:con}, that $\l_1\neq 0$ for all $x$ and $t$ can serve as another stability condition, since, otherwise, the standing wave $f(u)_x= s(u)$ may be singular and can't always be solved.

This chapter is organized as follows. After discussing the boundary averages in section 2, we study several different Godunov methods for the PW model in section 3. In section 4, we present numerical solutions to the PW model.

\section {Computation of the boundary averages of $\r$ and $v$}
To develop Godunov's methods for the PW model, we first partition a piece of roadway, e.g., an interval of $[a,b]$, into $N$ zones, and then approximate current traffic conditions $\r$ and $v$ with certain types of functions. At each zone boundary, the averages of $\r$ and $v$ over time have to be computed at every time step for computing $\r$ and $v$ at next time step. 

For a system of conservation laws, we can compute the boundary averages by solving a Riemann problem, which has been discussed by Smoller (1983\nocite{smoller1983}). For the homogeneous version of the PW model, we can adopt the solutions by Zhang (2000a\nocite {Zhang2000a}). Zhang developed solutions to the Riemann problem for the homogeneous version of Zhang's model, which is similar to the PW model, and obtained 8 types of wave solutions and the formula for the boundary averages of $\r$ and $v$ related to each type of solutions. However, the PW model has one relaxation term. The Riemann problem for \refe{eqn:cons} with a source term is still open. In this section, we still calculate the boundary averages based on the solutions to a Riemann problem for the homogeneous version of the PW model.

Liu (1979\nocite{Liu79}) discussed the Cauchy problem for a hyperbolic system of conservation laws with source terms. Inspired by his study, we present another approach of computing the boundary averages of $\r$ and $v$ by solving a Cauchy problem for \refe{eqn:cons}. This method is presented in the second part of this section.

\subsection{Computing the boundary averages from the Riemann problem}
In this subsection, we study the Riemann problem for the homogeneous version of \refe{eqn:cons} with the following jump initial conditions:
\bqn
u_{i+1/2}(x,t=0)&=&\cas{{ll}U_l,&\m{ if } x<x_{i+1/2}\\U_r,&\m{ if } x\geq x_{i+1/2}}, \label{ini.ri}
\eqn
where the left state $U_l=(\r_l, v_l)$ and the right state $U_r=(\r_r,v_r)$ are constant. For computational purpose, we are interested in the averages of $\r$ and $v$ at the boundary $x=x_{i+1/2}$ over a time interval $\dt$,  which are denoted by $\rfp$ and $\vfp$; i.e., we want to find
\bqn
\rfp=\frac 1 {\dt} \int_0^{\dt} \r(x=0,t) \m{dt} &\m{ and }& \vfp=\frac 1 {\dt} \int_0^{\dt} v(x=0,t) \m{dt}.
\eqn

For the homogeneous PW model, the equations of the characteristic curves are written in terms of $(\rho,v)$ instead of $(\rho,m)$, and the left and right initial values for $v$ are given as $\vl=\frac {\ml} {\rl}$ and $\vr=\frac {\mr}{\rr}$. Thus the boundary average of $m$, $\mfp=\rfp\vfp$, can be easily obtained once we computed $\rfp$ and $\vfp$.

Determined by the relationship between left state $(\rl,\vl)$ and right state $(\rr,\vr)$, there are 8 types of wave solutions to the Riemann problem, including 4 first-order waves and 4 second-order waves. The four first-order wave solutions are a 1-shock, a 2-shock, a 1-rarefaction and a 2-rarefaction. The four second-order wave solutions are a H1-H2 (Left-Shock-Right-Shock) wave, a R1-R2 (Left-Rarefaction-Right-Rarefaction) wave, a R1-H2 (Left-Rarefaction-Right-Shock) wave and a H1-R2 (Left-Shock-Right-Rarefaction) wave.

In the remaining part of this subsection, we define the velocity flux function $\phi(\r)$ as $\phi(\r)=c_0^2 \r$, discuss the eight types of wave solutions to the Riemann problem one by one, and present a table containing the boundary averages $\rfp$ and $\vfp$ for each case.

\begin {enumerate}
\item The wave solution is a 1-shock if the left and right states satisfy 
\begin {eqnarray}
\mbox{ H1: } \vr-\vl&=& -\sqrt {\frac {2(\rl-\rr)(\phi(\rl)-\phi(\rr))}{\rl\rr}}, \quad \rr>\rl,\: \vr<\vl.
\end {eqnarray}
The wave speed is 
\begin {eqnarray}
s&=&\frac {\rr\vr-\rl\vl}{\rr-\rl}
\end {eqnarray}
The solutions $(\rfp,\vfp)$ for case 1 are summarized in the following table: 
\\\begin {center}
\begin {tabular} {||c||c|c|c||}\hline
&$s=\frac {\rr\vr-\rl\vl}{\rr-\rl}$ &$\rfp$&$\vfp$ \\\cline{2-4}
&$s>0$ & $\rl$ &$\vl$ \\\cline{2-4}
H1&$s<0$ & $\rr$ & $\vr$ \\\cline{2-4}
&$s=0$ & $\frac {\rl+\rr} 2$ & $\frac {\vl+\vr}2$ \\\hline
\end {tabular}
\end {center}

\item The wave solution is a 2-shock if the left and right states satisfy 
\begin {eqnarray}
\mbox{ H2: } \vr-\vl&=& -\sqrt {\frac {2(\rl-\rr)(\phi(\rl)-\phi(\rr))}{\rl\rr}} \quad \rr<\rl,\: \vr<\vl
\end {eqnarray}
The wave speed is 
\begin {eqnarray}
s&=&\frac {\rr\vr-\rl\vl}{\rr-\rl} >0 
\end {eqnarray}
The solutions $(\rfp,\vfp)$ for case 2 are summarized in the following table: 
\\\begin {center}
\begin {tabular} {||c||c|c|c||}\hline
&$s=\frac {\rr\vr-\rl\vl}{\rr-\rl}$ &$\rfp$&$\vfp$ \\\cline{2-4}
H2&$s>0$ & $\rl$ &$\vl$ \\\hline
\end {tabular}
\end {center}

\item The wave solution is a 1-rarefaction if the left and right states satisfy 
\begin {eqnarray}
\mbox{ R1: } \vr-\vl&=&\vast(\rr)-\vast(\rl) \quad \rr<\rl,\: \vr>\vl
\end {eqnarray}
The characteristic velocity is determined by the first eigenvalue of the system: \begin {eqnarray}
\lambda_1(\rho,v)&=&v-c_0.
\end {eqnarray}
If $\lambda_1(\rl,\vl)>0$, the boundary averages $\rfp$ and $\vfp$ are the left initial values for $\r$ and $v$ , similarly, if  $\lambda_1(\rr,\vr)<0$, they are the right initial values. Otherwise, $(\rfp,\vfp)$ are solutions to the equations:
\begin {eqnarray}
\lambda_1(\rfp,\vfp)&=&\vfp-c_0=0,\\
\vfp-\vl &=& \vast(\rfp)-\vast(\rl),
\end {eqnarray}
 which can be simplified as follows,
\begin {eqnarray}
\vfp&=&c_0 \label {eqn:star1},\\
\vast(\rfp)&=&c_0-v_l+\vast(\rl)\label {eqn:star2}.
\end {eqnarray}

The solutions $(\rfp,\vfp)$ to \refet{eqn:star1}{eqn:star2} for case 3 are summarized in the following table: 
\\\begin {center}
\begin {tabular} {||c||c|c|c||}\hline
&$\lambda_1$ &$\rfp$&$\vfp$ \\\cline{2-4}
&$\lambda_1(\rl,\vl)>0$ & $\rl$ &$\vl$ \\\cline{2-4}
R1&$\lambda_1(\rr,\vr)<0$ & $\rr$ & $\vr$ \\\cline{2-4}
&o.w. & \multicolumn {2} {c||} {solution to equations \ref {eqn:star1}, \ref {eqn:star2} } \\\hline
\end {tabular}
\end {center}

\item The wave solution is a 2-rarefaction  if the left and right states satisfy 
\begin {eqnarray}
\mbox{ R2: } \vr-\vl&=&\vast(\rl)-\vast(\rr) \quad \rr>\rl,\: \vr>\vl
\end {eqnarray}
The characteristic velocity is the second eigenvalue of the system: 
\begin {eqnarray}
\lambda_2(\rho,v)&=&v+c_0>0.
\end {eqnarray}

The solutions $(\rfp,\vfp)$ for case 4 are summarized in the following table:
\\\begin {center}
\begin {tabular} {||c||c|c|c||}\hline
&$\lambda_2$ &$\rfp$&$\vfp$ \\\cline{2-4}
R2&$\lambda_2>0$ & $\rl$ &$\vl$  \\\hline
\end {tabular}
\end {center}

\item The wave solution is a 1-rarefaction + 2-rarefaction with an intermediate state $(\rmm,\vm)$ if the left, right and intermediate states satisfy
\begin {eqnarray}
\mbox{ R1: } \vm-\vl&=&\vast(\rmm)-\vast(\rl) \quad \rmm<\rl,\: \vm>\vl; \label{l.1} \\
\mbox{ R2: } \vr-\vm&=&\vast(\rmm)-\vast(\rr) \quad \rr>\rmm,\: \vr>\vm. \label{l.2}
\end {eqnarray}
Adding \refe{l.1} to \refe{l.2} we find
\begin {eqnarray}
2*\vast(\rmm)-\vast(\rl)-\vast(\rr)-(\vr-\vl)&=&0,
\end {eqnarray}
for $\rmm<\rl,\rmm<\rr$. Thus
\begin {eqnarray}
\vm=\vast(\rmm)+\vl-\vast(\rl).
\end {eqnarray}

The solutions $(\rfp,\vfp)$ for case 5 are summarized in the following table:
\\\begin {center}
\begin {tabular} {||c||c|c|c||}\hline
&$\lambda_1$ &$\rfp$&$\vfp$ \\\cline{2-4}
&$\lambda_1(\rl,\vl)>0$ & $\rl$ &$\vl$ \\\cline{2-4}
R1-R2&$\lambda_1(\rmm,\vm)<0$ & $\rmm$ & $\vm$ \\\cline{2-4}
&o.w. & \multicolumn {2} {c||} {solution to equations \ref {eqn:star1},
\ref {eqn:star2} } \\\hline
\end {tabular}
\end {center}

\item The wave solution is a 1-rarefaction + 2-shock with an intermediate state $(\rmm,\vm)$ if the left, right and intermediate states satisfy
\begin {eqnarray}
\mbox{ R1: } \vm-\vl=&\vast(\rmm)-\vast(\rl), & \rmm<\rl,\: \vm>\vl \\
\mbox{ H2: } \vr-\vm=&-\sqrt {\frac {2(\rmm-\rr)(\phi(\rmm)-\phi(\rr))}{\rmm\rr}} ,& \rr<\rmm,\: \vr<\vm.
\end {eqnarray}
These two equations yield
\begin {eqnarray}
&\vast(\rmm)-\vast(\rl)-\sqrt {\frac {2(\rmm-\rr)(\phi(\rmm)-\phi(\rr))}{\rmm\rr}}-(\vr-\vl)&
      =0,
\end {eqnarray}
for $\rr<\rmm<\rl$. Thus
\begin {eqnarray}
\vm=\vast(\rmm)+\vl-\vast(\rl).
\end {eqnarray}

The solutions $(\rfp,\vfp)$ for case 6 are summarized in the following table: 
\\\begin {center}
\begin {tabular} {||c||c|c|c||}\hline
&$\lambda_1$ &$\rfp$&$\vfp$ \\\cline{2-4}
&$\lambda_1(\rl,\vl)>0$ & $\rl$ &$\vl$ \\\cline{2-4}
R1-H2&$\lambda_1(\rmm,\vm)<0$ & $\rmm$ & $\vm$ \\\cline{2-4}
&o.w. & \multicolumn {2} {c||} {solution to equations \ref {eqn:star1},
\ref {eqn:star2} } \\\hline
\end {tabular}
\end {center}

\item The wave solution is a 1-shock + 2-shock with an intermediate state $(\rmm,\vm)$ if the left, right and intermediate states satisfy
\begin {eqnarray}
\mbox{ H1: } \vm-\vl=&-\sqrt {\frac {2(\rl-\rmm)(\phi(\rl)-\phi(\rmm))}{\rl\rmm}}, & \rmm>\rl,\: \vm<\vl;\\
\mbox{ H2: } \vr-\vm=&-\sqrt {\frac {2(\rmm-\rr)(\phi(\rmm)-\phi(\rr))}{\rmm\rr}}, & \rr<\rmm,\: \vr<\vm.
\end {eqnarray}
These two equations imply
\begin {eqnarray}
\sqrt {\frac {2(\rl-\rmm)(\phi(\rl)-\phi(\rmm))}{\rl\rmm}}+\sqrt {\frac {2(\rmm-\rr)(\phi(\rmm)-\phi(\rr))}{\rmm\rr}}&&\nonumber\\{}+(\vr-\vl)&=&0,
\end {eqnarray}
for $\rmm>\rl,\rmm>\rr$. Thus
\begin {eqnarray}
\vm=-\sqrt {\frac {2(\rl-\rmm)(\phi(\rl)-\phi(\rmm))}{\rl\rmm}}+\vl.
\end {eqnarray}

The solutions $(\rfp,\vfp)$ for case 7 are summarized in the following table: 
\\\begin {center}
\begin {tabular} {||c||c|c|c||}\hline
&$s=\frac {\rmm\vm-\rl\vl}{\rr-\rl}$ &$\rfp$&$\vfp$ \\\cline{2-4}
&$s>0$ & $\rl$ &$\vl$ \\\cline{2-4}
H1-H2&$s<0$ & $\rmm$ & $\vm$ \\\cline{2-4}
&$s=0$ & $\frac {\rl+\rmm} 2$ & $\frac {\vl+\vm}2$ \\\hline
\end {tabular}
\end {center}

\item The wave solution is a 1-shock + 2-rarefaction with an intermediate state $(\rmm,\vm)$ if the left, right and intermediate states satisfy
\begin {eqnarray}
\mbox{ H1: } \vm-\vl=&-\sqrt {\frac {2(\rl-\rmm)(\phi(\rl)-\phi(\rmm))}{\rl\rmm}}, & \rmm>\rl,\: \vm<\vl;\\
\mbox{ R2: } \vr-\vm=&\vast(\rmm)-\vast(\rr), & \rr>\rmm,\: \vr>\vm.
\end {eqnarray}
These two equations imply
\begin {eqnarray}
&-\sqrt {\frac {2(\rl-\rmm)(\phi(\rl)-\phi(\rmm))}{\rl\rmm}}+\vast(\rmm)-\vast(\rr)-(\vr-\vl)&
      \lefteqn{=0},
\end {eqnarray}
for $\rmm>\rl,\rmm>\rr$. Thus
\begin {eqnarray}
\vm=-\sqrt {\frac {2(\rl-\rmm)(\phi(\rl)-\phi(\rmm))}{\rl\rmm}}+\vl.
\end {eqnarray}

The solutions $(\rfp,\vfp)$ for case 8 are summarized in the following table:
\\\begin {center}
\begin {tabular} {||c||c|c|c||}\hline
&$s=\frac {\rmm\vm-\rl\vl}{\rr-\rl}$ &$\rfp$&$\vfp$ \\\cline{2-4}
&$s>0$ & $\rl$ &$\vl$ \\\cline{2-4}
H1-R2&$s<0$ & $\rmm$ & $\vm$ \\\cline{2-4}
&$s=0$ & $\frac {\rl+\rmm} 2$ & $\frac {\vl+\vm}2$ \\\hline
\end {tabular}
\end {center}

\end {enumerate}

\subsection {Computing the boundary averages from the Cauchy problem}
Liu (1979\nocite {Liu79}) studied the Cauchy problem for a hyperbolic system of conservation laws with source terms.  In this subsection we apply his theory to find the boundary averages by solving the Cauchy problem for the PW model. 

The homogeneous version of the PW model can be written as follows, 
\bqn
u_t+f(u)_x&=&0,
\eqn
where $u=(\r,v)$. The Riemann problem for this system at $x=0$ has the following initial conditions:  
\bqn
u(x,t_j)&=&\left \{ \begin {array} {l l}
U_l, &\mbox{ if } x < 0\\
U_r, &\mbox{ if } x > 0
\end {array}\right.\label {cauchy:1}.
\eqn
The wave solutions to the Riemann problem consist of two basic waves with an intermediate state $U_1$.  By denoting $U_0\equiv U_l$ and $U_2\equiv U_r$, we define $(U_{i-1},U_i)$ as the $i^{th} (i=1,2)$ propagation wave. Each basic wave $(U_{i-1},U_i)$ may be a shock or a rarefaction. These wave solutions, which have been discussed in previous subsection, serve as the basis for solutions to a related Cauchy problem.

The Cauchy problem for the PW model has the following initial conditions: 
\bqn
u(x,t_j)&=&\left \{ \begin {array} {l l}
u_l(x), &\mbox{ if } x < 0\\
u_r(x), &\mbox{ if } x > 0
\end {array}\right.\label {cauchy:2} ,
\eqn
in which, 
\bqn
f(u_l(x))_x&=& s(u_l(x)), \qquad u_l(0)=U_l \label {standingwave:1}\\
f(u_r(x))_x&=& s(u_r(x)), \qquad u_r(0)=U_r \label {standingwave:2}.
\eqn
Here the intermediate state $u_1(x)$ is determined from 
\bqn
f(u_1(x))_x&=& s(u_1(x)), \qquad u_1(0)=U_1,
\eqn
in which $U_1$ is the intermediate state solution to the corresponding Riemann problem.
We denote $u_0(x) \equiv u_l(x)$ and $u_2(x)\equiv u_r(x)$, and define $(u_{i-1}(x),u_i(x))$ as the $i^{th} (i=1,2 )$ propagation wave of the Cauchy problem.

For computational purposes, we still want to compute the boundary average of $u(x,t)$ at $x=0$; i.e., $U^{\ast}\equiv \int_0^{\dt} u(0,t) dt$. The boundary average varies when $x=0$ is covered by different states or waves. Since solutions to the Cauchy problem for the PW model consist of three states and four types of waves, the boundary $x=0$ may be covered by the left, right or intermediate state, and may be crossed by 1-shock, 2-shock, 1-rarefaction or 2-rarefaction wave. 

If the boundary $x=0$ is covered by the intermediate state, according to the definition of the intermediate state, the solutions of the boundary average $U^{\ast}$ is equal to $U_1$. Similarly, $U^{\ast}=U_l$ or $U_r$ when the boundary is covered by the left or right state.

Of the 4 types of basic waves, the 2-H and 2-R waves never cross $x=0$ since $\lambda_2(U)>0$, and 1-shock still propagates along a line described by $\frac x t=\sigma(U_{i-1}, U_i)$. Therefore, these three waves don't affect the boundary averages. In the remaining part, we consider the adjustment to the 1-rarefaction wave. 

The 1-rarefaction wave was shown by Liu to be a perturbation on the solutions to the corresponding Riemann problem. 
Assume that $(u_{i-1}(x),u_i(x))$ is a 1-R wave, $u_{i-1}(x)$ and $u_i(x)$ are separated in a region $x_{i-1}(t)<x<x_i(t)$. The wave solution of the Cauchy problem (\ref{cauchy:2}) approaches the 1-rarefaction wave $(U_{i-1},U_i)$ as $t\rightarrow 0$,
\bqn
\left \{\begin {array} {l} \lim _{t\rightarrow 0} \frac {x_{i-1}(t)}t =\lambda_1(u_{i-1}(x))\equiv \xi_0 \\
 \lim _{t\rightarrow 0} \frac {x_{i}(t)}t =\lambda_1(u_{i}(x))\equiv \xi_1 \end {array}\right . 
\eqn
\bqn
\lim _{t\rightarrow 0;\frac xt=\eta} u(x,t)=v(\eta), \quad x_{i-1}(t)<x<x_i(t) ,
\eqn 
where $v(\eta)\in R_1(U_{i-1})$ with $\lambda_1(v(\eta))=\eta$. 

Using the coordinate of time $t$, we define the initial $i$-characterized speed $\xi$ as follows:
\bqn
\frac {\partial x(\xi,t) }{\partial t} &=& \lambda_1(w(\xi,t))\label{cauchy:t1}\\ 
w(\xi,t)&=& u(x(\xi,t),t)\label {cauchy:t2}
\eqn
and 
\bqn
w(\xi,0)&=&v(\xi)\equiv \phi(\xi) \\
x_{i-1}&=&x(\xi_0,t)\\
u(x(\xi_0,t),t)&=&u_{i-1}(x(\xi_0,t),t)=u_{i-1}(x(\xi_0,t))\\
u(x(\xi_1,t),t)&=&u_i(x(\xi_1,t),t)=u_i(x(\xi_1,t))
\eqn
For a homogeneous system $\xi$ will be a constant slope for characteristics. However, $\xi$ is not constant for the PW model with a source term. 

With the transformation (\ref {cauchy:t1},\ref {cauchy:t2}), the original system
\bqn
\pd {u(x(\xi,t),t)}{t} +A(u) \pd ux &=& s(u) \qquad \m {in which } A(u)=\nabla f(u)
\eqn
becomes 
\bqn
\pd {x(\xi,t)}{t} &=&\lambda_1(w(\xi,t))\\
\pd x {\xi} \pd w t+(A-\lambda_1) \pd w {\xi}&=& \pd x {\xi} s.
\eqn

Since 
\bqn
x(\xi,0)=0, 
\eqn
we obtain
\bqn
x(\xi,t)&=&\pd x t\big |_{t=0} t+\frac 12 \pd {{}^2x}{t^2}\big |_{t=0} t^2+\textbf O(t^3) \nonumber
\\
&=&\xi t+\frac 12 \pd {{}^2x}{t^2}\big |_{t=0} t^2+\textbf O(t^3)
\eqn
We calculate $\pd {{}^2x}{t^2}$ as follows:
\bqn
\pd {{}^2x}{t^2}&=&\pd {\lambda_1(w)}{t} = \nabla \lambda_1 \pd wt \nonumber\\
&=&\nabla \lambda_1 (\pd ux \pd xt+\pd ut)=\nabla \lambda_1(s+(\xi I-A)\pd ux)\big |_{t=0} \nonumber\\
&=&\nabla \lambda_1 s +\textbf {O} (t),
\eqn
where we assume the difference between $\xi \pd ux$ and $A\pd ux$ is small since $-\xi u_{\xi}+f(u)_{\xi}=0$ for a rarefaction wave.
Then, the rarefaction wave path is   
\bqn
x(t)&=& \lambda_1 t+  \frac 12\nabla \lambda_1(u) s(u) t^2 +\textbf {O} (t^3)
 \nonumber\\
&=&(\frac m {\rho} -c_0) t + \frac {f_{\ast} -m} {2\tau\rho} t^2 +\textbf {O} (t^3) .\label {cauchy:para}  
\eqn
which is a parabola. The linear term of the rarefaction wave is the 1-R wave for the Riemann problem. The second-order term is determined by the source term. For small time scales; i.e., $t$ small, this rarefaction wave is a perturbation of the rarefaction for the corresponding Riemann problem. The 2-rarefaction wave functions can be derived similarly, although they are not necessary for computing the boundary averages.

Recall that, by assumption, the characteristic curve of 1-rarefaction wave is
\bqn
v_l&=&v_{\ast}(\rho)-v_{\ast}(\rho_l).
\eqn
Let $x(t)=0$. From (\ref {cauchy:para}) we get
\bqn
v&=&c_0-\frac {f_{\ast} -m} {2\tau\rho} t \\
\rho&=&v_{\ast}^{-1}(v-v_l+v_{\ast}(\rl))=v_{\ast}^{-1}(\rho_0)-\frac 1 {v'_{\ast}(\rho_0)} \frac {f_{\ast} -m} {2\tau\rho} t, 
\eqn
from which we can calculate $U^{\ast}\equiv \int_0^{\dt} u(0,t) dt$.

The above analysis shows that the 1-rarefaction wave fans of the Cauchy problem consists of parabolic curves instead of lines. When we use these parabolic rarefaction waves to compute $U^{\ast}$, the numerical solution of the PW model improves. However, the adjustment for a 1-rarefaction wave is needed only when it crosses the boundary $x=0$, and 1-rarefaction waves are adjusted with a lower order perturbation in a short time step $\dt$. Thus this improvement doesn't appear to be significant. Besides, we have $\l_1(U)=0$ for some state $U$ when 1-rarefaction wave cross the boundary. Therefore, the state $U$ is in the unstable region for the PW model. This may be another reason that this adjustment doesn't yield significantly better solutions.

\section {Godunov methods}
In this section, we study Godunov methods for solving the PW model (\ref {eqn:cons}). For a general system (\ref {eqn:gene}), the finite difference equations are
\bqn
U_i^{j+1}=U_i^j-\dxt (f(U_{i+1/2}^{j+1/2})-f(U_{i-1/2}^{j+1/2}))+\dt \tilde s(U)
,\label {eqn:disc}\eqn 
in which $U_i^{j+1}$ $U_i^j$ are both averages of $u(x,t)$ over $i^{th}$ cell at time $(j+1)\dt$ and $j\dt$, $U_{i\pm 1/2}^{j+1/2}$ are the boundary averages calculated as shown in the preceding section, and $\tilde s(U)$ is the source average over $((i-1/2)\dx,(i+1/2)\dx)\times (j\dt,(j+1)\dt)$. 

When we treat the source term implicitly, the system is discretized as 
\begin {eqnarray}
\frac {\rho^{j+1}_i-\rho^j_i }{k}+\frac {\mfp-\mfm} {h} &=&0 \label{ll.1}
\end {eqnarray}
\begin {eqnarray}
\frac {m^{j+1}_i-m^j_i }{k}
+\frac {\frac{\mfp^2 }{\rfp} +c_0^2\rfp -\frac{\mfm^2 }{\rfm} -c_0^2\rfm }{h} \nonumber\\
&=& \frac {\fast(\rho^{j+1}_i)-m^{j+1}_i} {\tau}  \label{ll.2}
\end {eqnarray}
in which, $\rho^{j}_i$ and $m_i^j$ are the cell average of $\r$ and $m$ respectively over the $i$th cell, and 
 $\rfp$ and $\mfp$ are the averages of $\rho$ and $m$ respectively on the cell boundary $x_{i+1/2}$ in the time interval $(t_j,t_{j+1})$.
In \refe{ll.2}, the source term $\tilde s(U)$ is treated implicitly. 

From \refet{ll.1}{ll.2}, we can write the evolution equations for the PW model as
\begin {eqnarray}
&\rho^{j+1}_i& = \rho^j_i-\frac kh (\mfp-\mfm )\\ 
&m_i^{j+1}&= \frac 1 {(1+\frac k {\tau})  } \{m_i^j -\frac kh[\frac{\mfp^2 }{\rfp} +c_0^2\rfp -\frac{(\mfm)^2 }{\rfp} -c_0^2\rfm ] \\&&\qquad+\frac k {\tau} \fast (\rho^{j+1}_i)\}\nonumber
\end {eqnarray}

In a first-order Godunov method we use the cell averages $\rl=\rho_i^j,\ml=m_i^j$ and $\rr=\rho_{i+1}^j,\mr=m_{i+1}^j$ as left/right states as the initial condition for the Riemann problem 
\begin {eqnarray}
u_{i+1/2}(x,t_j)=\left\{\begin {array} {l c}
U_l, & \mbox{ if } x-x_{i+1/2}<0 \\
U_r, & \mbox{ if } x-x_{i+1/2}\geq 0
\end {array}\right. .
\end {eqnarray}

\subsection {The Second-order Godunov Method}
In this subsection we introduce a second-order Godunov method for the PW model in the process similar to that in section \ref{secondGodunov}. 

We begin by writing the PW model in the linearized form
\begin {eqnarray}
u_t+A(u) u_x &=& s(u),
\end {eqnarray}
where 
\begin {eqnarray}
u&=&\left ( \begin {array} {c}
\rho(x,t)\\
m(x,t)
\end {array} \right )
\end {eqnarray}
and 
\begin {eqnarray}
A(u)&=&\matau.
\end {eqnarray}

The eigenvalues and eigenvectors of $A(u)$ are
\begin {eqnarray}
\begin {array}{ll}
\lambda_1(u)=\frac m {\rho}-c_0,&r_1(u)=[1,\lambda_1]^t=[1,\frac m {\rho} -c_0]^t\\
\lambda_2(u)=\frac m {\rho}+c_0,&r_1(u)=[1,\lambda_2]^t=[1,\frac m {\rho} +c_0]^t
\end {array}
\end {eqnarray}
We diagonalize $A(u)$ by 
\begin {eqnarray}
T^{-1}(u)A(u)T(u)&=&\matlu\equiv\Lambda(u),
\end {eqnarray}
where the transformation matrix $T(u)$ is 
\begin {eqnarray}
T(u)&=&\mattu .
\end {eqnarray}

Letting $W=T^{-1}(u)u$, the PW model under the transformation becomes 
\begin {eqnarray}
W_t+\Lambda(u)W_x&=& T^{-1}(u)s(u).
\end {eqnarray}
Therefore, the PW model is transformed into two separated scalar equations in $W=(w_1,w_2)$. For the scalar equation in $w_i$ ($i=1,2$), we introduce an interpolation for $w_i(x,t)$ which yields a second-order Godunov method for solving this equation. In the remaining part of this subsection, we first introduce an interpolation for $w_i(x,t)$ and then apply inverse transformation on them in order to develop a second-order method for the whole system.

For a scalar equation $w_t+\lambda(w) w_x= 0$, in a first-order Godunov method we use a step function $u_I(x,t_j)$ to interpolate the data at time $t_j$, 
\begin {eqnarray}
w_I(x,t_j) &=& w_i^j, \mbox {  if } x_{i-1/2}<x\leq x_{i+1/2} ,
\end {eqnarray}
and solve the Riemann problem with the following initial conditions: 
\begin {eqnarray}
\begin {array} {lcl}
w_{i+1/2}^{j,L}& =&w_i^j \\
w_{i-1/2}^{j,R} &=& w_i^j
\end {array}
\end {eqnarray}

In a second-order Godunov method, we interpolate the data at time $t_j$ with a piecewise linear function,
\begin {eqnarray}
w_I(x,t_j) &=& w_i^j+\frac {(x-ih)}{h} \Delta^{VL} w_i^j, \mbox {  if } x_{i-1/2}<x\leq x_{i+1/2} .
\end {eqnarray}

With a half-step prediction , the Riemann problem has the initial conditions of
\begin {eqnarray}
\begin {array} {lcl}
w_{i+1/2}^{j+1/2,L}& =& w_i^j+\frac 12 (1-\lambda(w_i^j)\frac k h) \Delta ^{VL} w_i^j \\
w_{i-1/2}^{j+1/2,R} &=& w_i^j-\frac 12 (1+\lambda(w_i^j)\frac k h) \Delta ^{VL} w_i^j
\end {array},
\end {eqnarray}
where $\Delta ^{VL} w_i^j$ is the van Leer slope (all superscripts $j$ are suppressed)
\begin {eqnarray*}
&\Delta ^{VL}w_i& = \left \{\begin {array} {l}
S_i\cdot \min(2|w_{i+1}-w_i|,2|w_i-w_{i-1}|,\frac 12|w_{i+1}-w_{i-1}|) ,\: \lefteqn{\varphi>0} \\
0, \qquad\qquad\qquad \mbox{ otherwise }
\end {array}\right. 
\end {eqnarray*}
\begin {eqnarray}
S_i&=& \mbox{sign} (w_{i+1}-w_{i-1})\\
\varphi&=&(w_{i+1}-w_i)\cdot(w_i-w_{i-1})
\end {eqnarray}

We apply the above procedure twice to compute the half-step values $W_{i+1/2}^{j+1/2,L}$ and $W_{i-1/2}^{j+1/2,R}$. Then we apply the inverse transformation on these half-step values to obtain $U_{i+1/2}^{j+1/2,L}$ and $U_{i-1/2}^{j+1/2,R}$:
\begin {eqnarray}
\begin {array} {ccc}
U_{i+1/2}^{j+1/2,L}&=&T(U_i^j) W_{i+1/2}^{j+1/2,L} \\
U_{i-1/2}^{j+1/2,R}&=&T(U_i^j) W_{i-1/2}^{j+1/2,R}
\end {array}
\end {eqnarray}

With this new interpolation of $\r$ and $m$, we can solve the Riemann problem or the Cauchy problem on a cell boundary. For the inhomogeneous system, this second-order method has a convergence rate of 2. However, it is different for the PW model with a relaxation term.

\subsection {Some other Godunov-type variant methods}
In this subsection, we review other variants of Godunov method for \refe{eqn:cons} with a source term.

The first variant was suggested by Pember (1993a\nocite {Pember93a},1993b\nocite {Pember93b}). He treated the source term as  $\tilde s(U)=\half(s(U_{i-1/2}^{j+1/2})+s(U_{i+1/2}^{j+1/2}))$, where both $U_{i-1/2}^{j+1/2}$ and $U_{i+1/2}^{j+1/2}$ are the boundary averages solved in the Riemann problem. 

The second variant is the fractional step splitting method, in which each time step $\dt$ is split in-to three steps. In the first and third fractional steps, a first-order implicit method is used to integrate
\bqn
\matu_t=\mats 
\eqn
for time steps of $\dt/2$. In the second step,
 we solve the corresponding homogeneous system of (\ref {eqn:cons}), i.e.,    
\bqn
\matu_t+\matf_x= 0
\eqn
for a time step of $\dt$.

A third variant is the quasi-steady wave-propagation algorithm suggested by LeVeque (1998a\nocite {LeVeque98a},1998b\nocite {LeVeque98b}). This method introduces a new discontinuity in the center of each grid, i.e.,
\bqn
\frac 12 (U_i^{-}+U_i^{+})&=& U_i,
\eqn
and
\bqn
f(U_i^{+})-f(U_i^{-})&=&s(U_i)\dx,
\eqn
where $U_i^{\pm}=\left (\begin {array}{c} \rho_i^{\pm}\\m_i^{\pm}\end {array}\right )$.
Then we get
\bqn
\left (\begin {array} {c}
m_i^{+}-m_i^{-}\\
{m_i^{+}}^2/\rho_i^{+}+c_0^2\rho_i^{+}-{m_i^{-}}^2/\rho_i^{-}-c_0^2\rho_i^{-}
\end {array}\right )
&=&
\left (\begin {array}{c}
0\\
\frac {f_{\ast}(\rho_i)-m_i}{\tau}
\end {array} \right ) \dx
\eqn
so $m_i^{+}=m_i^{-}=m_i$. Setting 
\bqn
\rho_i^{+}=\rho_i+\delta ,\qquad \rho_i^{-}=\rho_i -\delta
\eqn
we find
\bqn
2c_0^2\delta^3-K\delta ^2-(2m_i^2+2c_0^2\rho_i^2)\delta+K\rho_i^2=0, \label {var:e1}
\eqn
in which $K=\frac {f_{\ast}(\rho_i)-m_i}{\tau}\dx$. 

Given the solution to (\ref {var:e1}), LeVeque solved the Riemann problem of the homogeneous system with the initial conditions
\bqn
u_{i+1/2}(x,t_j)&=&\left \{ \begin {array} {l l} 
U_i^{+} & \m { if } x-x_{i+1/2}<0 \\
U_i^{-} & \m { if } x-x_{i+1/2}\geq0
\end {array}\right.
\eqn
Then in the evolution equations, the source term is not considered. 

In LeVeque's method, the data at time $t_j$ are interpolated with a solution to the standing wave for the PW model in each cell. This method may be used together with the Cauchy problem we discussed before.

\section{Numerical Solutions to the PW Model}
In this section, we use the model parameters given in (Kerner, 1994\nocite {Kerner}), i.e., $c_0=2.48445 l/\tau$, $v_{\ast}(\rho)\equiv V(\rho)=5.0461[(1+\exp\{[\rho-0.25]/0.06\})^{-1}-3.72\times 10^{-6}] l/\tau$ and $L=800l$. Here $l$ is the unit of length, $\tau$ is the relaxation coefficient, and the section of the roadway is from $0 l$ to $800 l$.  We set $l=10$ (m) $\tau=10$ (sec) and $\rho_h=0.172$. The equilibrium functions $\vs(\r)$ and $f_{\ast}(\r)$ are given in \reff {fig:func1}. In the figure, $\rho_{c1},\rho_{c2}$ are two critical densities and the region $\rho_{c1}<\r<\r_{c2}$ is the unstable region. According to (\ref {stable:con}), $\rho_{c1}$ and $\rho_{c2}$ satisfy the equation:
\bqn
\rho v'_{\ast}(\rho)+c_0&=&0.
\eqn
From this equation we get $\rho_{c1}=0.173$,$\rho_{c2}=0.396$.

For numerical computation purpose, we use two initial conditions:
\bqn
\rho(x,0)&=&\rh+0.02\sin(2\pi x/800/l) \label {ini1:1}\\
v(x,0)&=&v_{\ast}(\rh)-0.02\cos(2\pi x/800/l)\label {ini1:2},
\eqn
which is a global perturbation, and
\bqn
\rho(x,0)&=&\cas {{ll} \rh+\delta &\m { when } x\in [37.5 l,48.4 l]\\\rh-\delta/3 &\m { when } x\in [50.0 l,82.8 l]\\
\rh&\m { otherwise }} \label {ini2:1}\\
v(x,0)&=&v_{\ast}(\rho(x,0)) \label {ini2:2},
\eqn
which is a local perturbation. Setting $\delta=0.02$ and $\rh=0.16$, the two initial conditions are given in \reff {fig:ini1}. The global perturbation is not in equilibrium, however the local perturbation is in equilibrium.
\bfg
\bc\includegraphics[height=8cm] {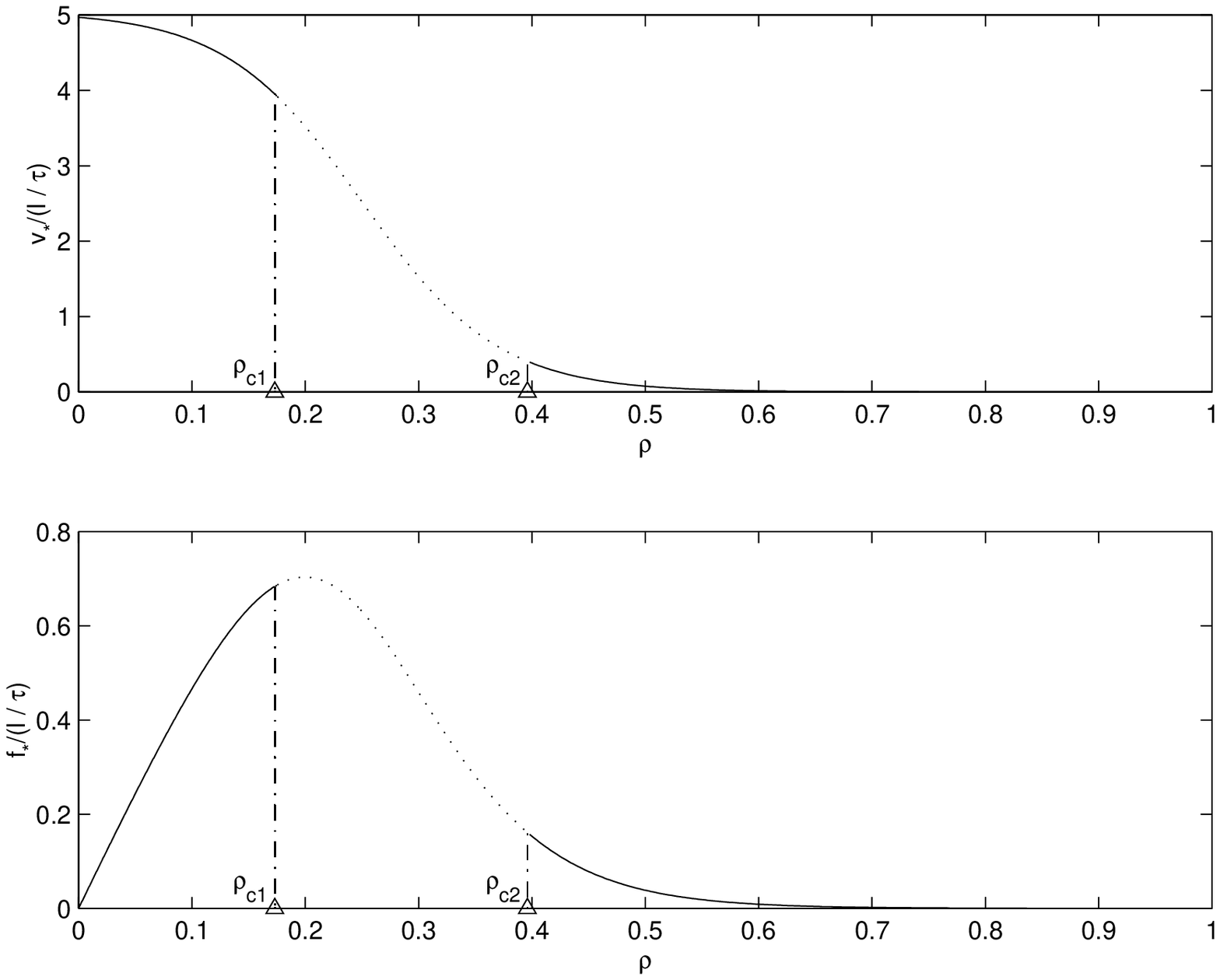}\ec
\caption {One selection of the equilibrium velocity and flow rate} \label {fig:func1}
\efg

\bfg
\bc\includegraphics[height=8cm] {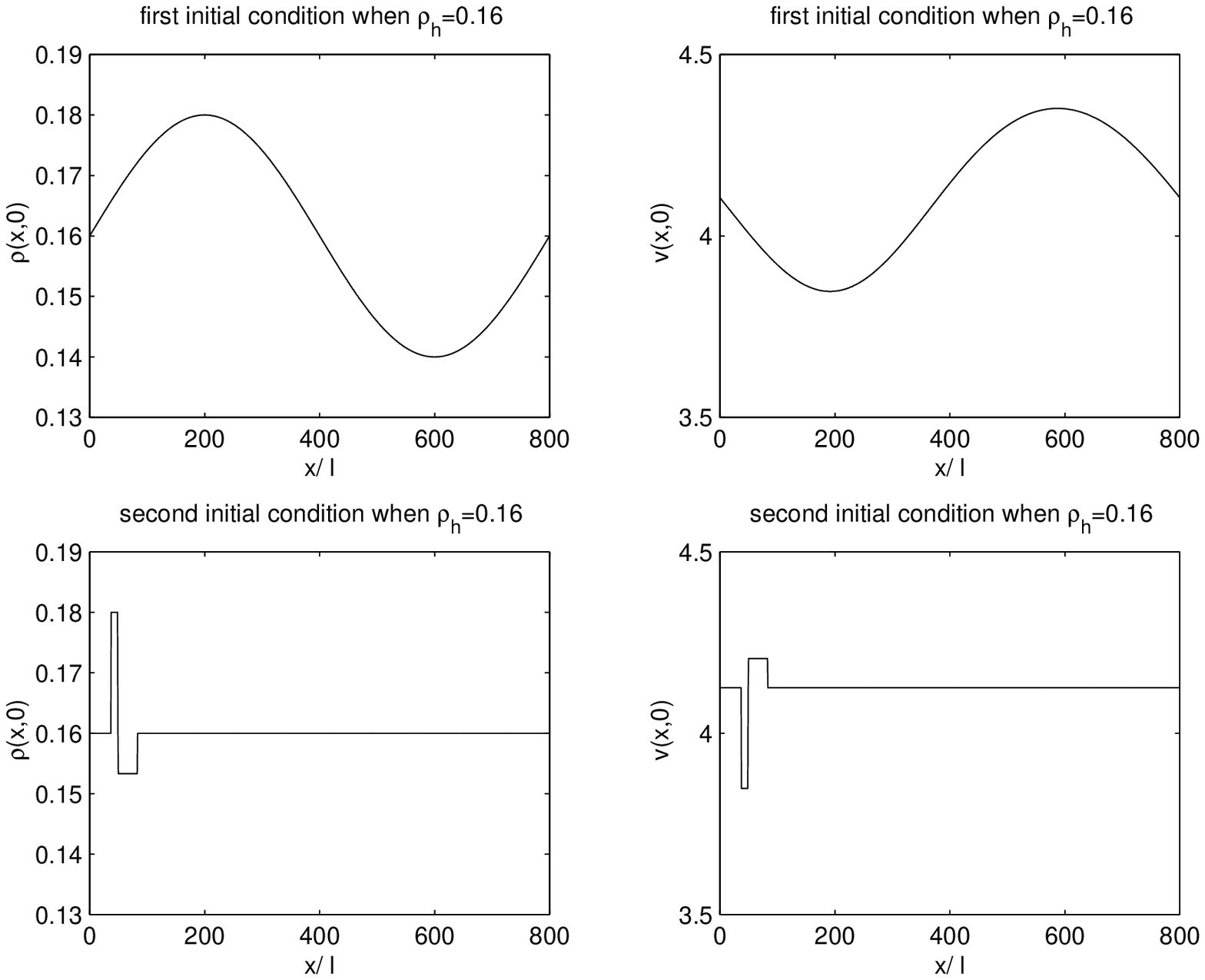}\ec
\caption {Two initial conditions: global non-equilibrium perturbation v.s. local equilibrium perturbation} \label {fig:ini1}
\efg

\bfg
\bc\includegraphics[height=8cm] {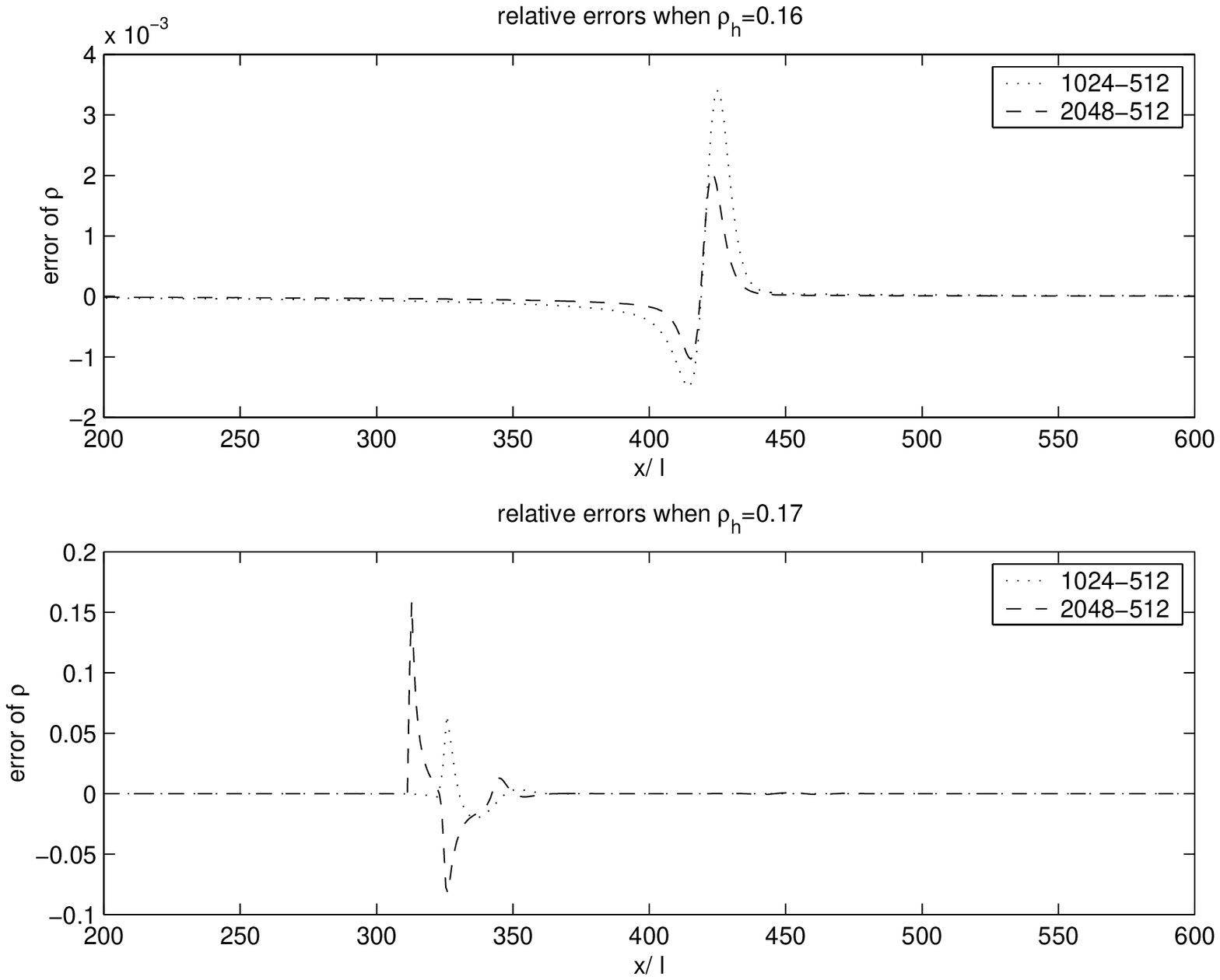}\ec
\caption {Stability and instability for (\ref {ini1:1},\ref {ini1:2}) with Periodic boundary condition}\label {fig:stable1}
\efg
\bfg
\bc\includegraphics[height=8cm] {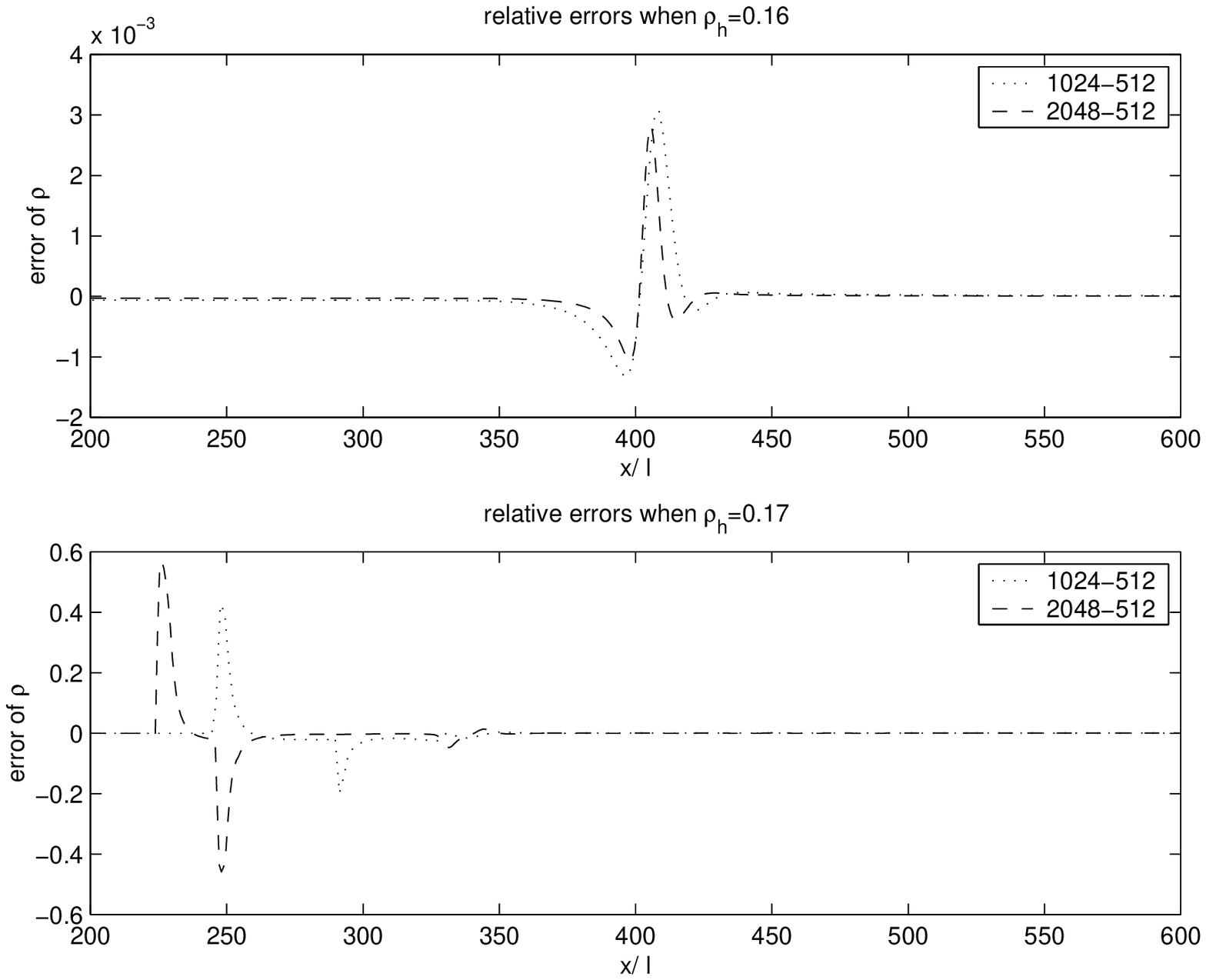}\ec
\caption {Stability and instability for (\ref {ini1:1},\ref {ini1:2}) with Neumann boundary condition} \label {fig:stable2}
\efg

\bfg
\bc\includegraphics[height=10cm] {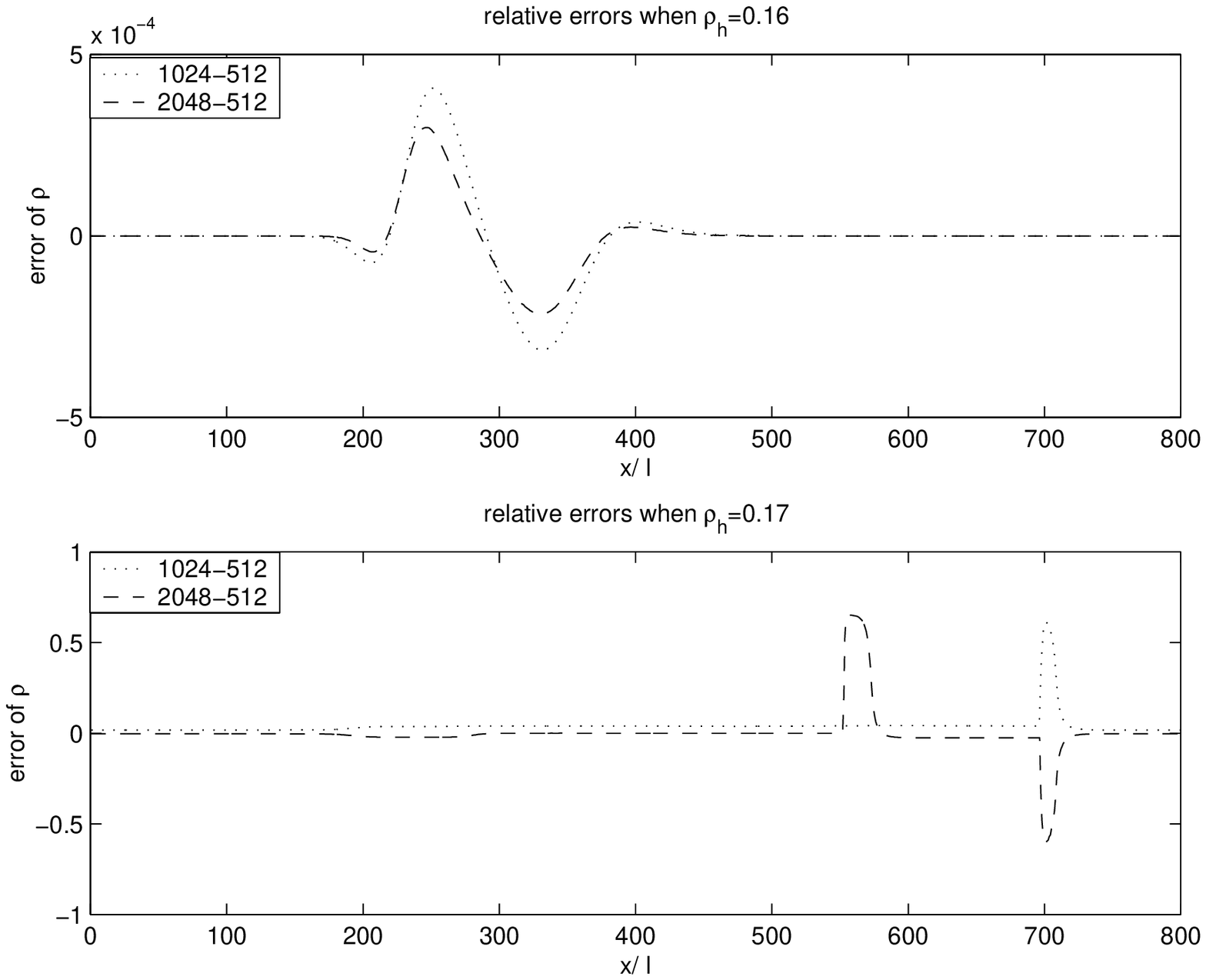}\ec
\caption {Stability and instability for (\ref {ini2:1},\ref {ini2:2}) with Periodic boundary condition}\label {fig:stable3}
\efg

\subsection {Stability test}
The PW model is unstable in some regions, which is different from Zhang's model. In this subsection we test the stability property of the PW model at time $T_0=200.0\tau$. We use a first-order Godunov's method to compute the relative differences, defined in \refe{def:error}, between solutions with $N=512$ grids and those with 1024 grids, and the difference between solutions with 1024 grids and those with 2048 grids. These differences are drawn as curves labeled as  $1024-512$ and $2048-1024$ respectively in the following figures. Since the solutions $(\r,v)$ are close to the equilibrium state, only the differences of $\r$ are given.   Given the convergent Godunov's method, the difference caused by different grid numbers decreases if the PW model is stable. We test the PW model with different initial conditions and different boundary conditions, which are in the stable or unstable region for the PW model, and the results are the same as predicted.

Setting $\rh=0.16$ or $\rh=0.17$ for the initial conditions (\ref {ini1:1},\ref {ini1:2}), we solve the PW model with 512, 1024 and 2048 grids with periodic boundary conditions. The differences are shown in \reff {fig:stable1}. In the figure, the difference decreases when we increase the number of grids when $\rh=0.16$; but the difference increases when $\rh=0.17$. The figure proves that the PW model is stable for $\rh=0.16$, and unstable for $\rh=0.17$.

Using the same initial conditions as in last example, we solve the PW model with Neumann boundary conditions. The relative differences are given in \reff{fig:stable2}. In this case we get the same results as in last example.

Next we use the initial conditions (\ref {ini2:1},\ref {ini2:2}) with $\delta=0.02$. We solve the PW model with periodic boundary conditions. According to the differences shown in \reff{fig:stable3}, we conclude that the PW model also stable when $\rh=0.16$ and unstable when $\rh=0.17$. We don't test the PW model for the local perturbation with Neumann boundary conditions, since we expect the same stability property.

According to the numerical tests, the system is stable when $\max \rho(x,0)\leq0.18$ and unstable when $\max \rho(x,0)\geq 0.19$, which is consistent with the theoretical prediction.

\subsection {The Riemann problem and steady-state solutions}
In this subsections, we solve the PW model with four well-selected initial conditions so that we observe second-order waves. The Neumann boundary conditions are used and the number of grids is set as $N=1024$. We also change the relaxation term to $\frac {f_{\ast}-m}{1000\tau}$. The relaxation time $1000\tau$ is  long enough for us to watch the second-order waves at $T_0=100\tau$, before all 2-waves relax to 1-waves.

\begin {figure}
\bc\includegraphics[height=8cm] {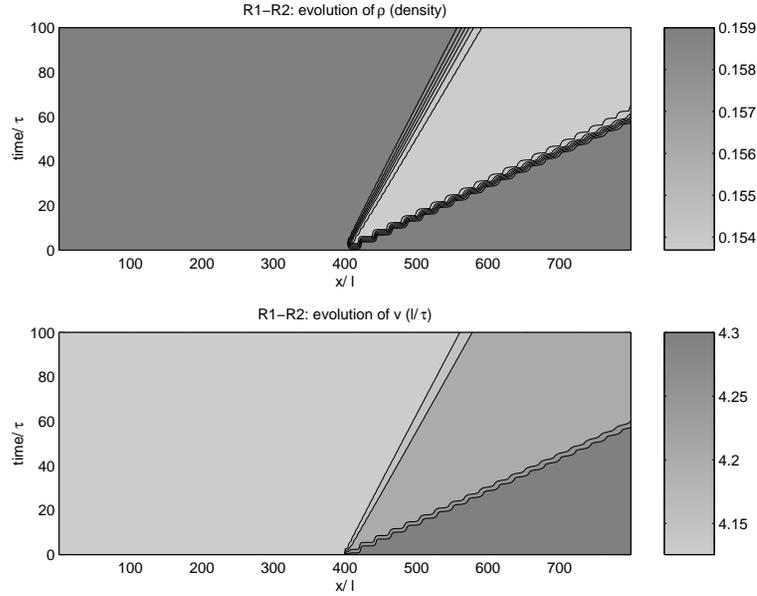}\ec
\caption{ The R1-R2 wave solutions to the Riemann problem}\label {riemann1}
\end {figure}

\begin {figure}
\bc\includegraphics[height=8cm] {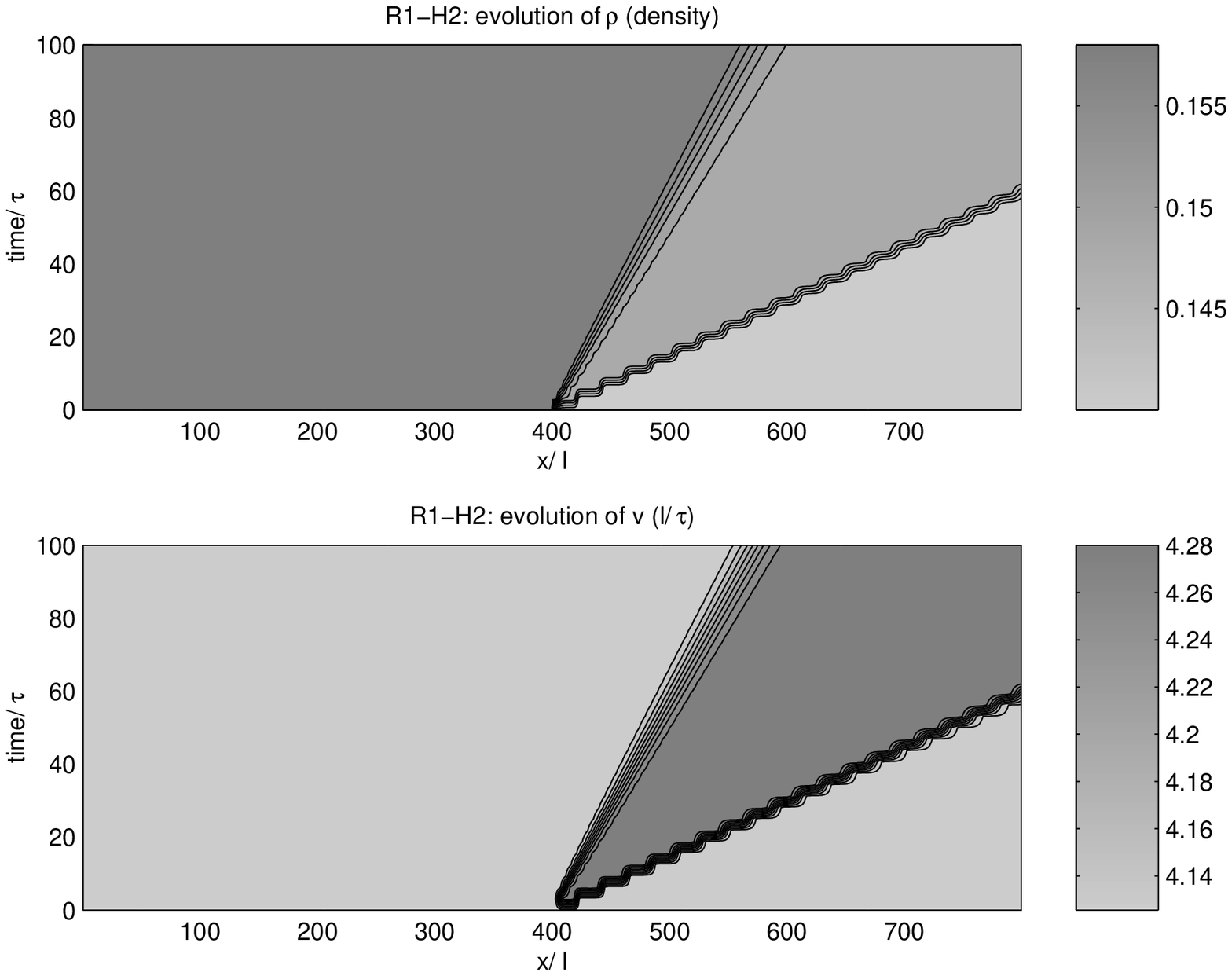}\ec
\caption{The R1-H2 wave solutions to the Riemann problem}\label {riemann2}
\end {figure}

\begin {figure}
\bc\includegraphics[height=8cm] {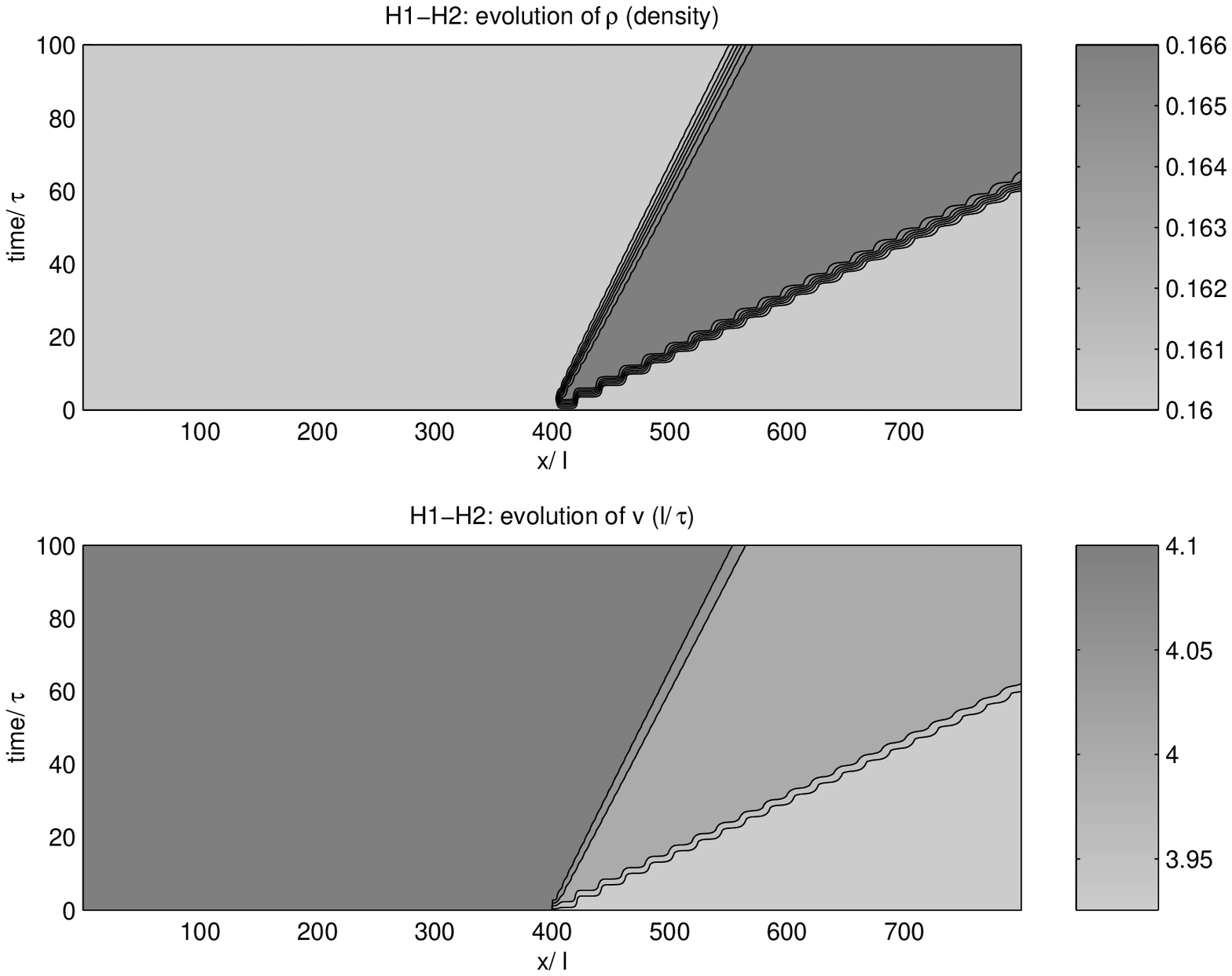}\ec
\caption{The H1-H2 wave solutions to the Riemann problem}\label {riemann3}
\end {figure}

\begin {figure}
\bc\includegraphics[height=8cm] {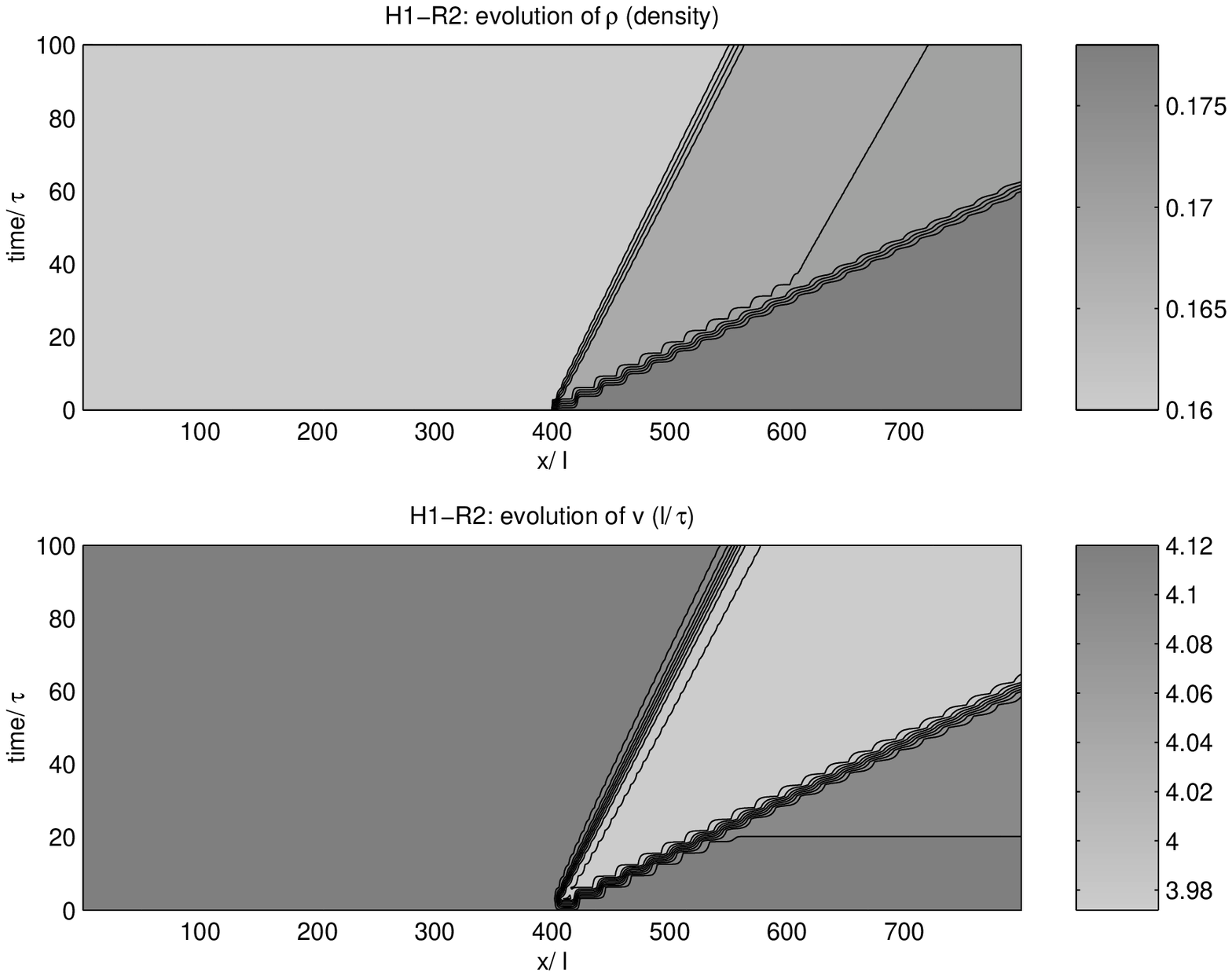}\ec
\caption{The H1-R2 wave solutions to the Riemann problem}\label {riemann4}
\end {figure}

\ben
\item We use the following jump initial conditions
\bqn
\rho(x,0)&=& 0.16 \label{riemann1:1}\\
v(x,0)&=&\cas{{ll} v_l=v_{\ast}(0.16)& \m { when } x<=400 l\\v_r=v_{\ast}(0.16)+0.2 l/\tau & \m { when } x>400 l}.\label{riemann1:2}
\eqn
 The solutions to the Riemann problem are a R1-R2 wave, given in \reff{riemann1}.

\item We use the following jump initial conditions
\bqn
\rho(x,0)&=&\cas{{ll} \rl=0.16& \m { when } x<=400 l\\\rr=0.16-0.02 & \m { when } x>400 l}\label{riemann2:1} \\
v(x,0)&=&v_{\ast}(0.16).\label{riemann2:2}
\eqn
The wave solutions to the Riemann problem of the PW model are a R1-H2 wave, given in \reff{riemann2}.

\item We use the following jump initial conditions
\bqn
\rho(x,0)&=& 0.16 \label{riemann3:1}\\
v(x,0)&=&\cas{{ll} v_l=v_{\ast}(0.16)& \m { when } x<=400 l\\v_r=v_{\ast}(0.16)-0.2 l/\tau & \m { when } x>400 l}.\label{riemann3:2}
\eqn
The wave solutions are a H1-H2 wave, given in \reff{riemann3}.

\item We use the following jump initial conditions
\bqn
\rho(x,0)&=&\cas{{ll} \rl=0.16& \m { when } x<=400 l\\\rr=0.16+0.02 & \m { when } x>400 l} \label{riemann4:1}\\
v(x,0)&=&v_{\ast}(0.16).\label{riemann4:2}
\eqn
The wave solutions are a H1-R2 wave, given in \reff{riemann4}.
\een

The solutions above show four different type of second-order waves, which consist of two basic waves. However, 2-waves relax to 1-waves if the relaxation time is shorter or the observing time is longer. A 2-shock wave relaxes to a 1-rarefaction wave; and a 2-rarefaction wave relaxes to a 1-shock wave, due to the effect of the source term. In next part, we show how the 2-waves relax to 1-waves and how the free regions and cluster regions form. We set relaxation time as $10\tau$, and observe the solutions until $T_0=100\tau$.

%Fri Nov  5 19:33:00 PST 1999

\ben
\item
With the initial conditions (\ref{riemann1:1}, \ref{riemann1:2}), we get the solutions  shown in \reff {cluster1}. At around $5\tau$ a downstream 1-shock forms when the traffic conditions relax to the equilibrium state, i.e., $v=v_{\ast}(0.16)$. After that traffic flow forms a free region with lower density and higher travel speed. The free region travels in the speed of $\lambda_{\ast}$. However the free region will disappear as the R1-H1 wave propagates, and finally the traffic flow will become uniform.
\item With the initial conditions (\ref{riemann2:1}, \ref{riemann2:2}), we get the solutions shown in \reff{cluster2}. A new 1-rarefaction wave forms when the H2-wave disappears. As these two rarefaction waves propagate, the traffic conditions will become uniform.

\item With the initial conditions (\ref{riemann3:1}, \ref{riemann3:2}), we get the solutions shown in \reff{cluster3}. At around $5\tau$, the 2-shock wave disappears and a 1-rarefaction wave forms. After that a cluster, with higher density and lower travel speed, travels with at the speed of $\lambda_{\ast}$. However, the traffic conditions will become uniform as the H1-R1 wave propagates.
\item With the initial conditions (\ref{riemann4:1}, \ref{riemann4:2}), we get the solutions shown in \reff{cluster4}. As long as the 2-rarefaction disappears, a new 1-shock is formed. After that, both shock waves travel in the speed of $\lambda_{\ast}(0.16)=2.12 l/ \tau$.
\een

In the remaining part of this subsection we consider the steady-state solutions of the system after a long time $T_0=1600 \tau$ with the relaxation time $\tau$ .
\ben
\item Using the initial conditions (\ref {ini1:1},\ref {ini1:2}) with $\rh=0.16$ and periodic boundary conditions, we get the solutions shown in \reff{fig:steady1}. The figure shows that an upstream shock wave and a downstream rarefaction wave form. After that, the traffic conditions become more and more uniform.

\item Using the initial conditions (\ref {ini2:1},\ref {ini2:2}) with $\rh=0.16$ and periodic boundary conditions, we have the solutions shown in \reff {fig:steady2}. The solutions show that the traffic conditions get more and more uniform, quicker than the case shown in \reff {fig:steady1}.

\een

\bfg
\bc\includegraphics[height=8cm] {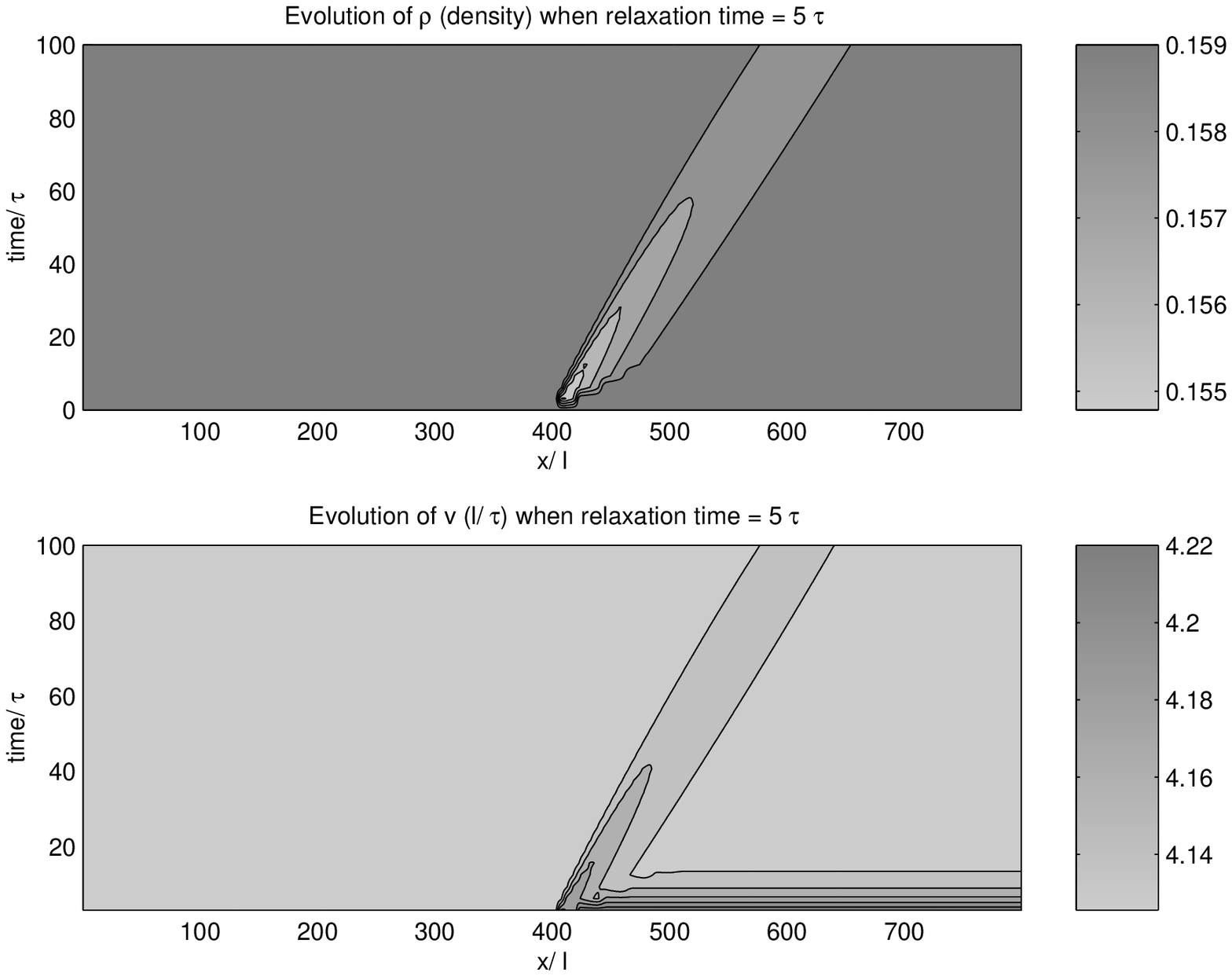}\ec
\caption {Formation of a free region}\label {cluster1}
\efg

\bfg
\bc\includegraphics[height=8cm] {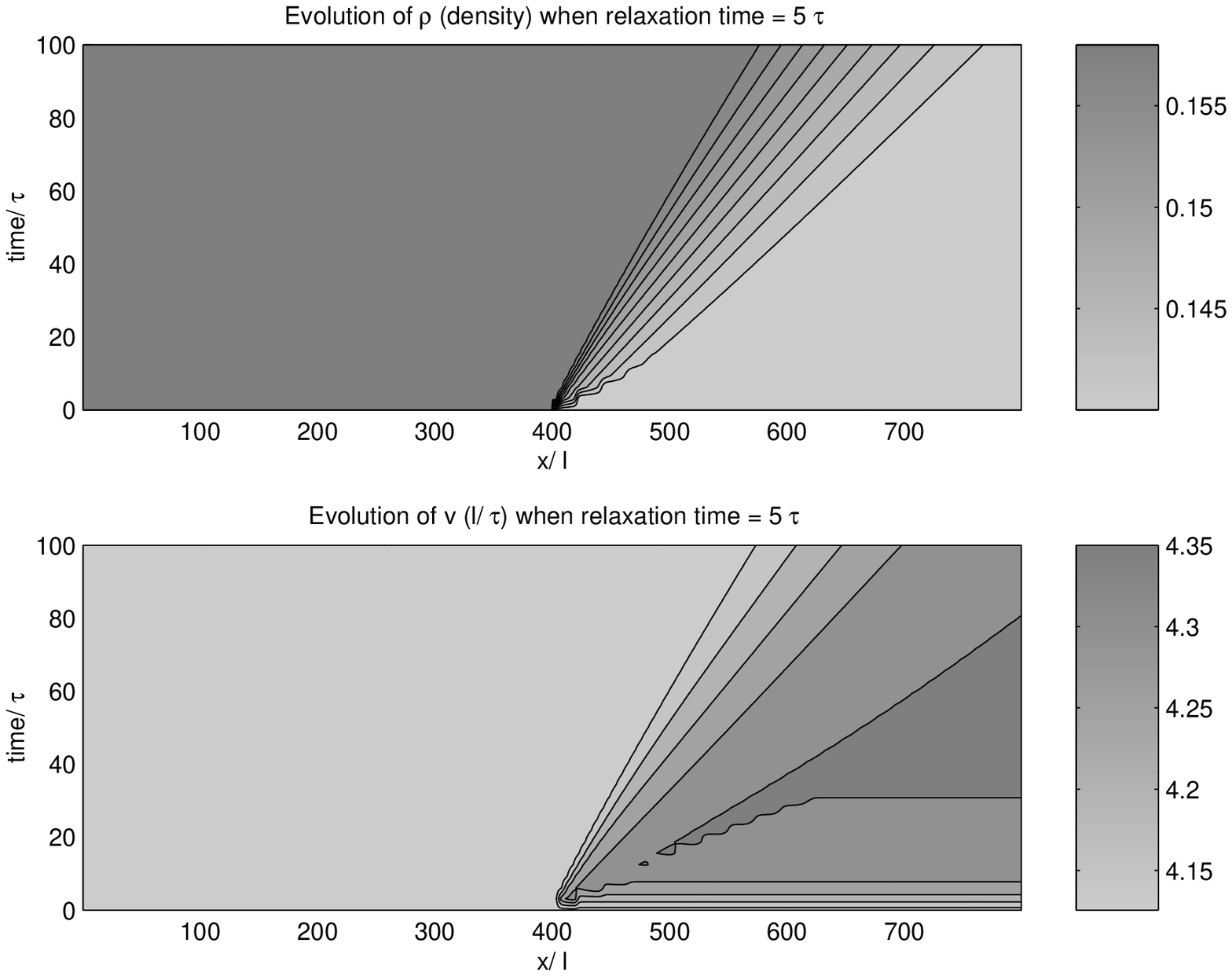}\ec
\caption {Formation of two 1-rarefaction waves}\label {cluster2}
\efg

\bfg
\bc\includegraphics[height=8cm] {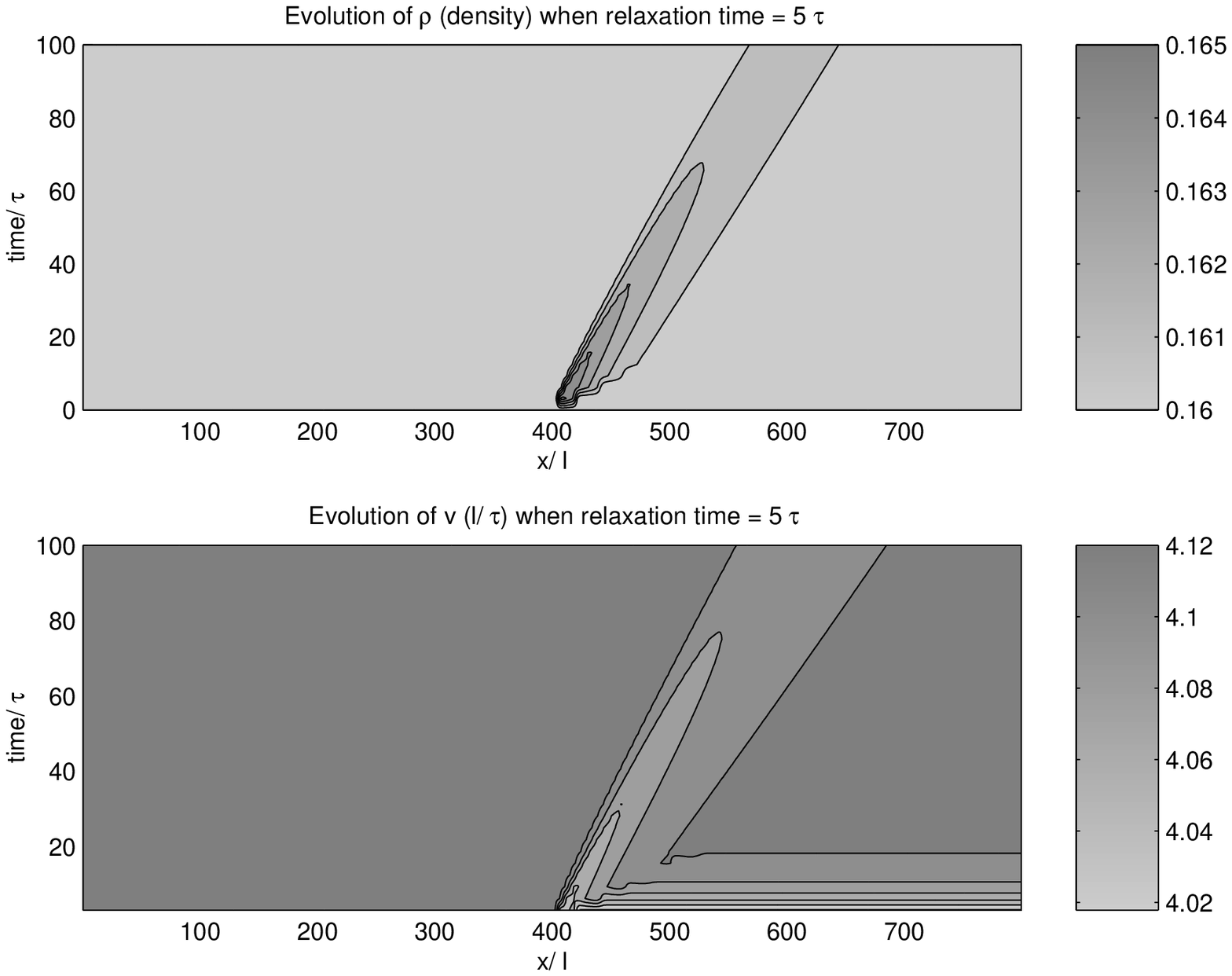}\ec
\caption {Formation of a cluster}\label {cluster3}
\efg

\bfg
\bc\includegraphics[height=8cm] {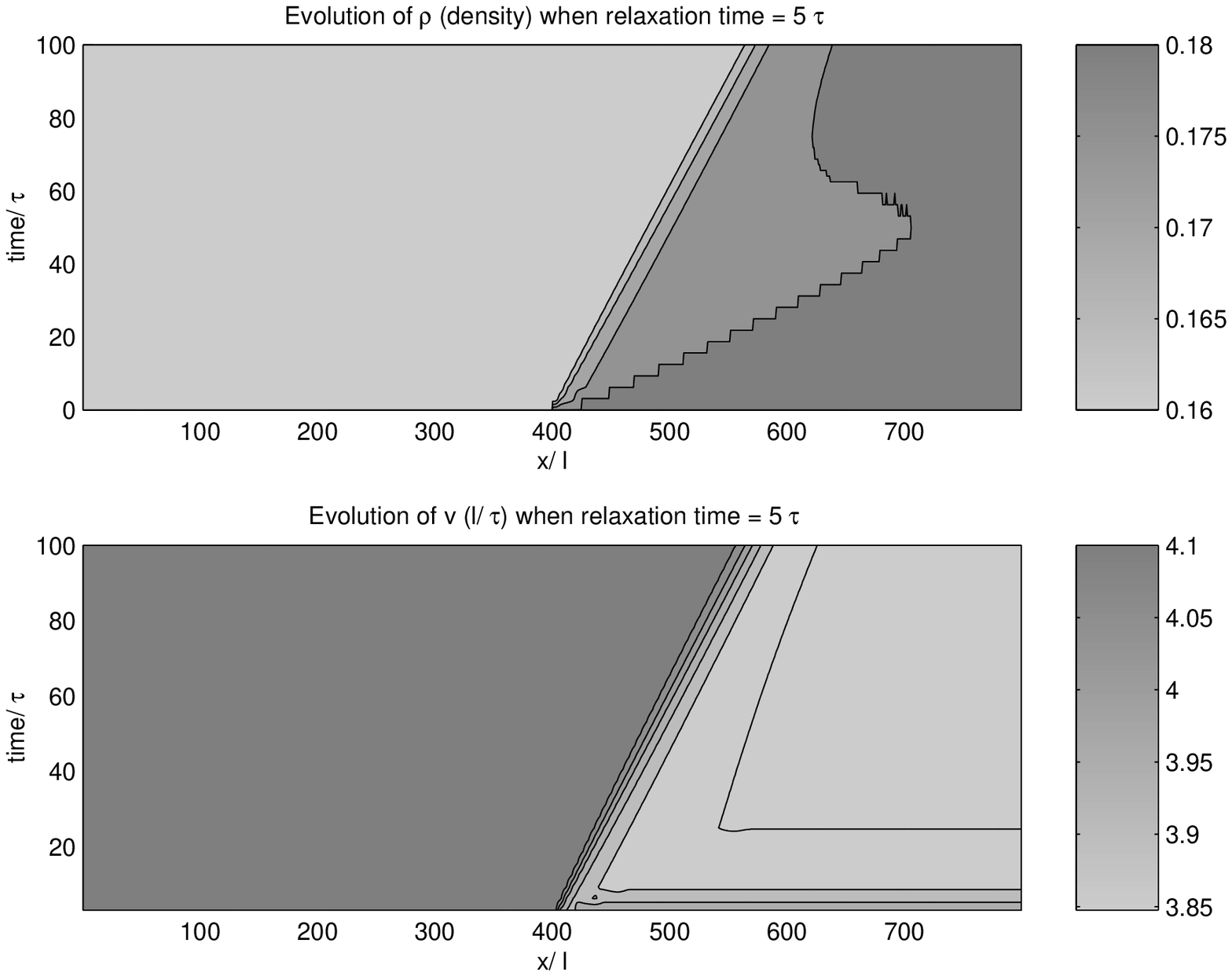}\ec
\caption {Formation of two 1-shock waves}\label {cluster4}
\efg

%Wed Nov 10 20:54:53 PST 1999
\subsection {General solutions and convergence rates}
In this subsection, different numerical methods for the PW model are discussed with the initial conditions \refet {ini1:1}{ini1:2} ($\rh=0.16$). Using Neumann boundary conditions, we solve the PW model until $T_0=400\tau$. With the number of grids as 64, 128, 256, 512 or 1024, we carry out 5 separate computations. For 1024 grids, we plot the contour of $\r$ and $v$ on the $x-t$ phase plane as well as their solutions at selected times. We also compute the convergence rates based on solutions with different number of grids. The convergence rate is defined in (\ref{def:error} -- \ref{epsilon}).

\bfg
\bc\includegraphics[height=8cm] {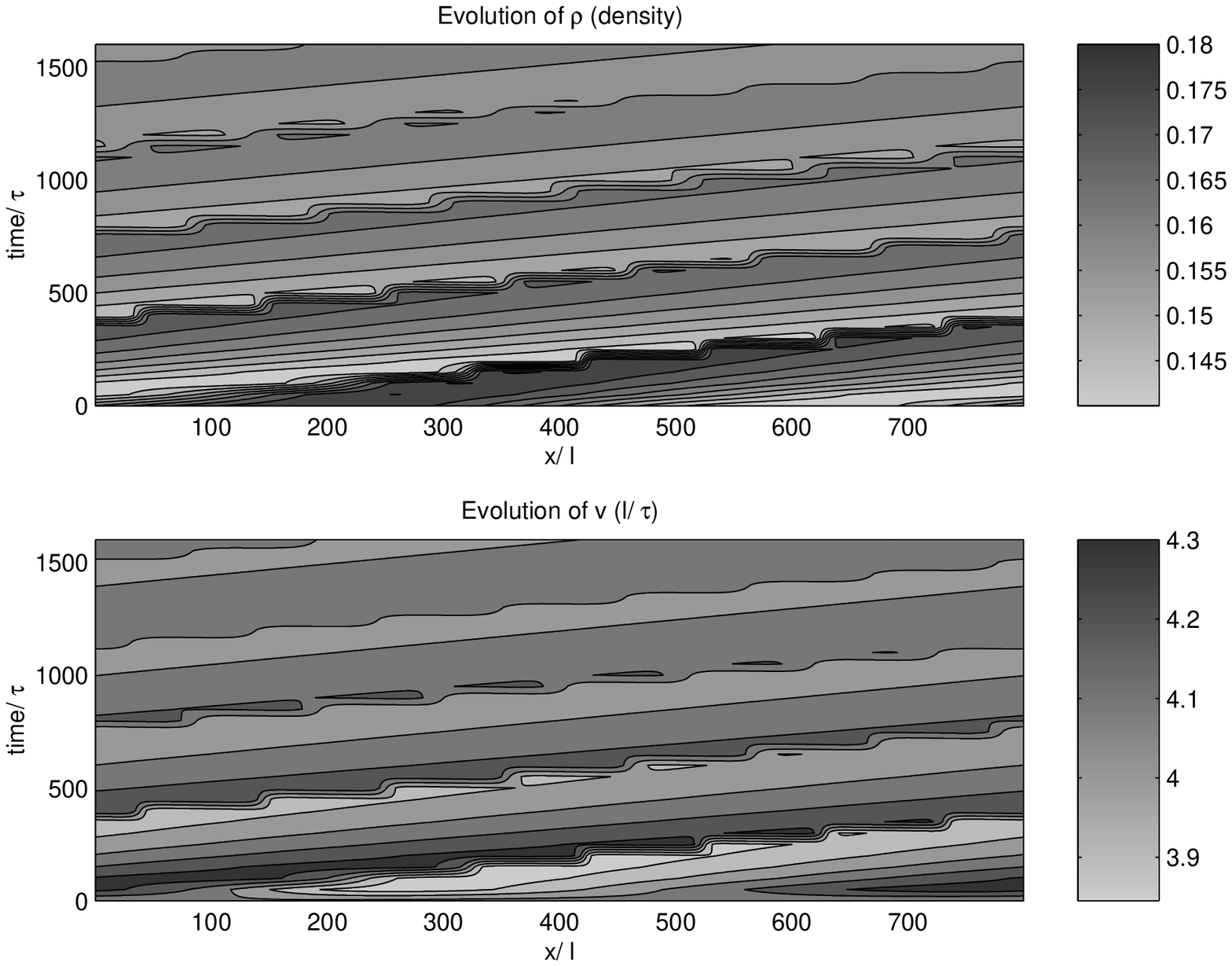}\ec
\caption {Solution for (\ref {ini1:1},\ref {ini1:2}) till 1600$\tau$}\label {fig:steady1}
\efg

\bfg
\bc\includegraphics[height=8cm] {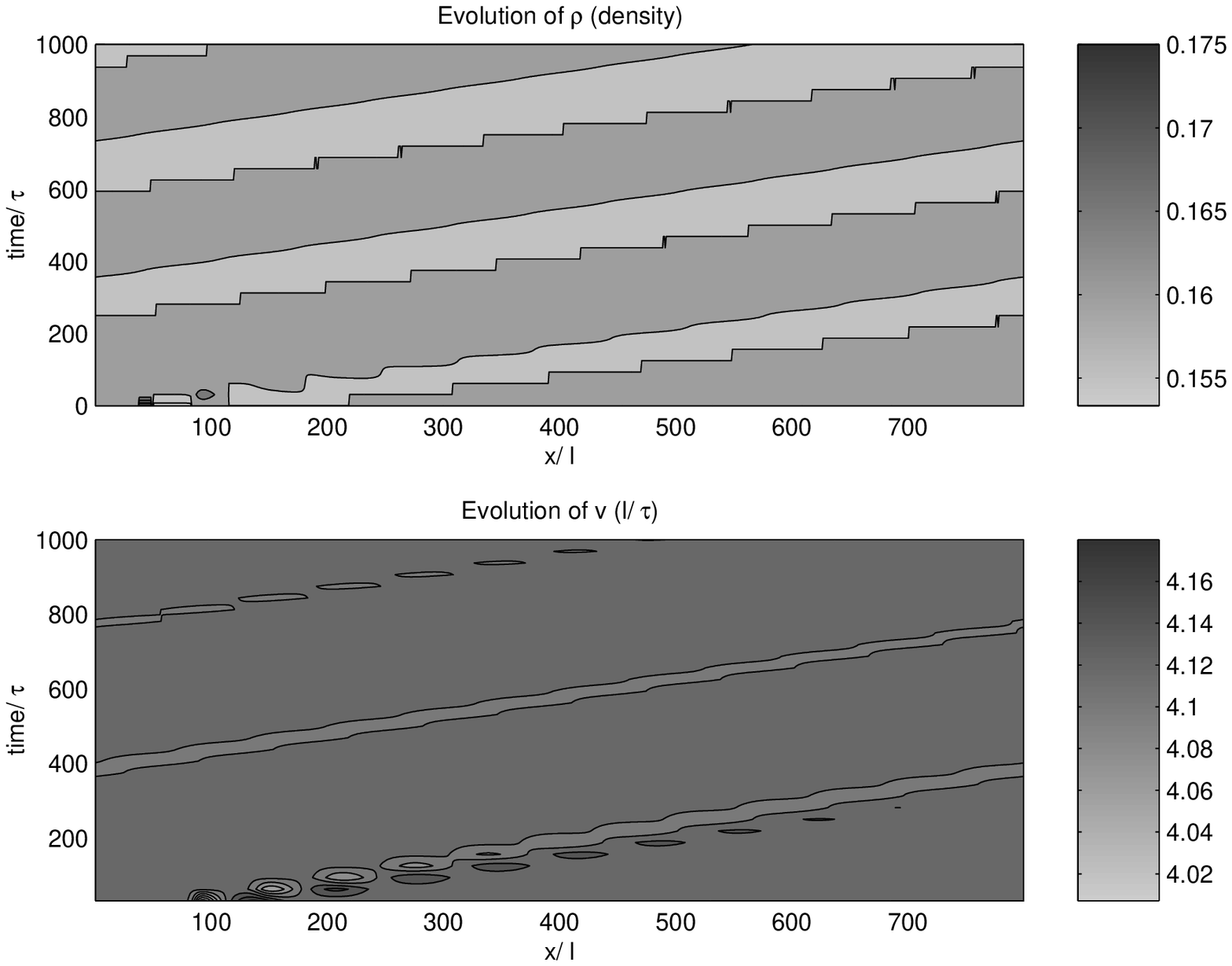}\ec
\caption {Solution for (\ref {ini2:1},\ref {ini2:2}) till 1000$\tau$}\label {fig:steady2}
\efg

\bfg
\bc\includegraphics[height=8cm] {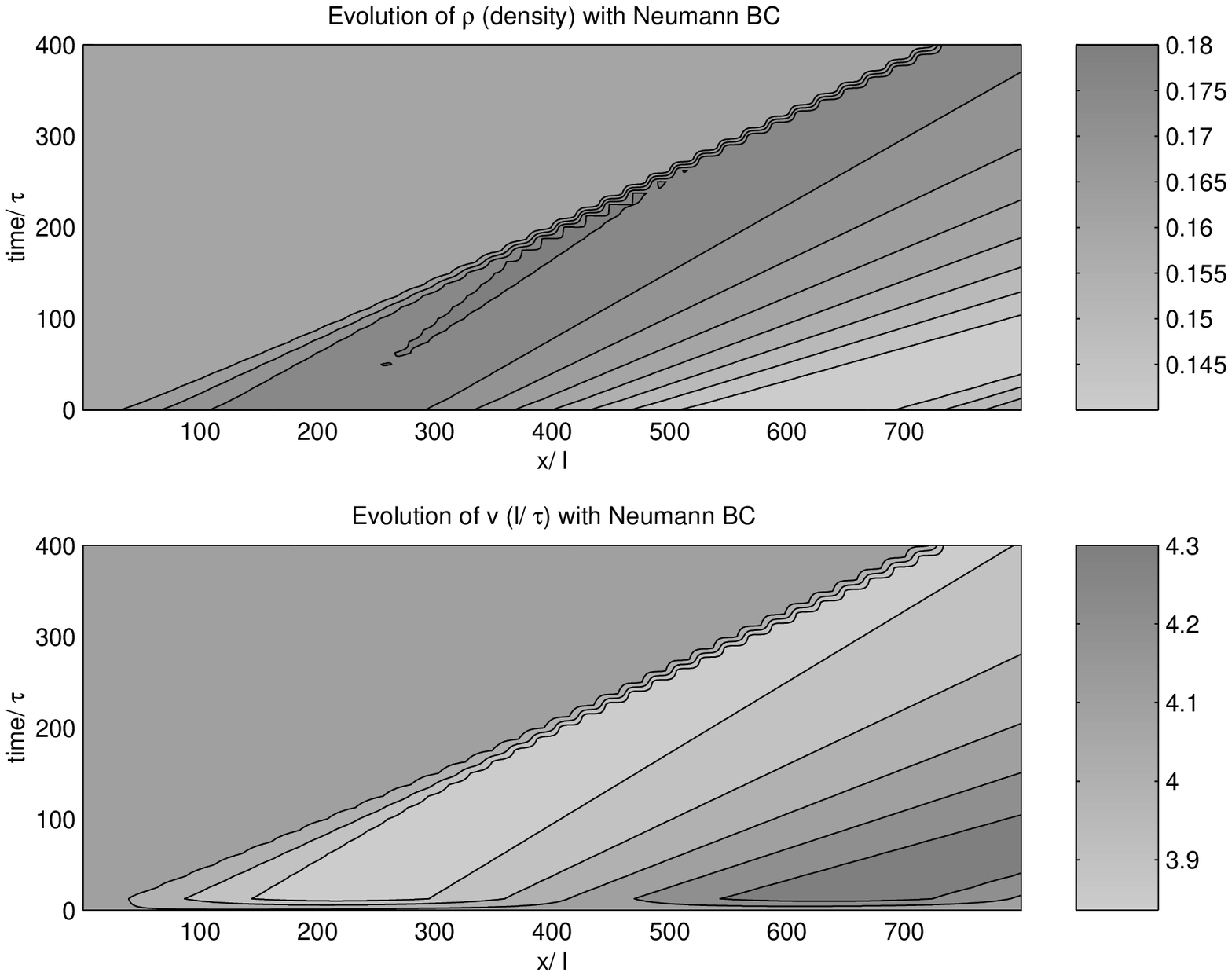}\ec
\caption {Solutions by a first-order Godunov's method for 1024 grids}\label {fig:general13d}
\efg

\bfg
\bc\includegraphics[height=8cm] {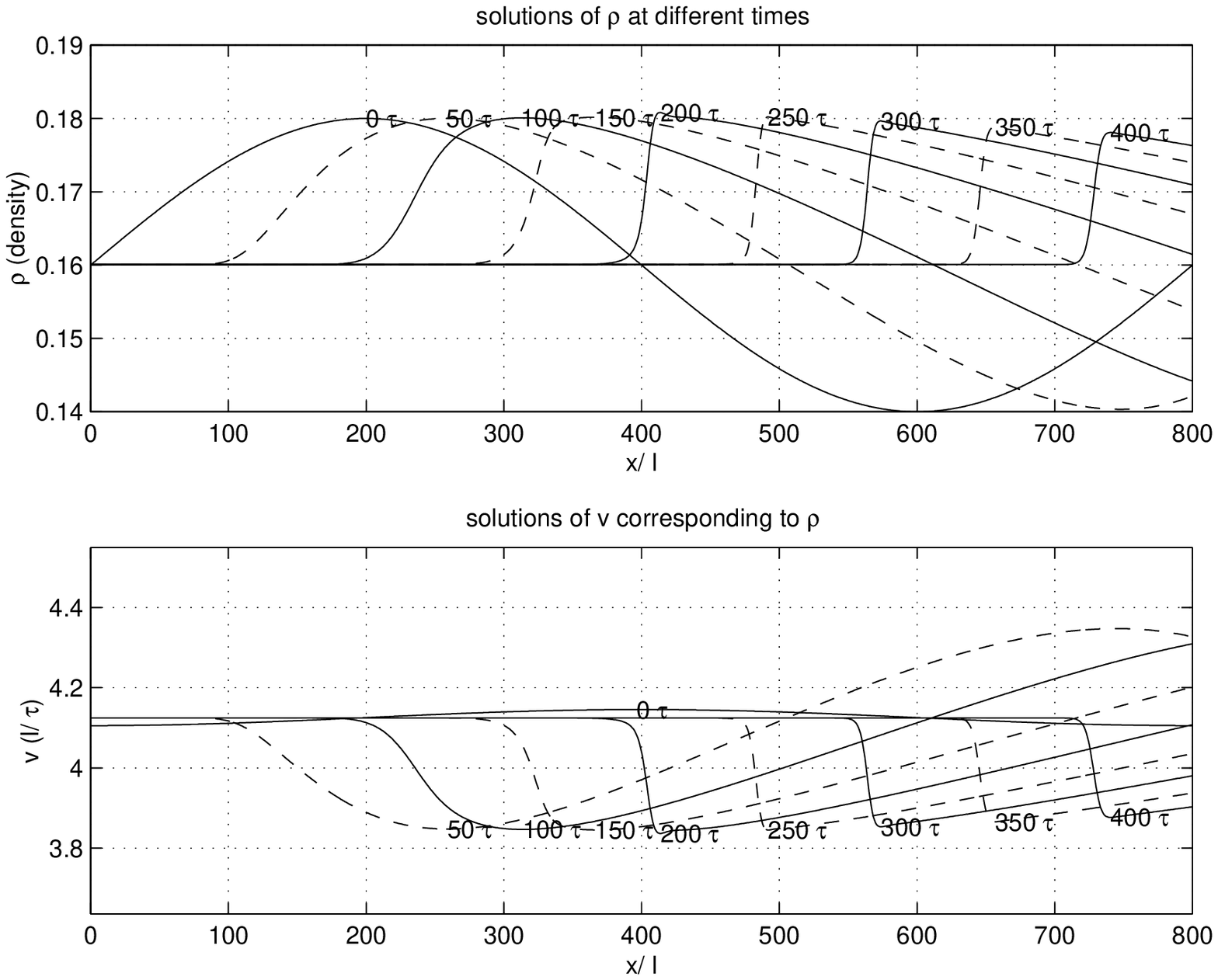}\ec
\caption {Solutions from \reff{fig:general13d} at selected times}\label {fig:general12d}
\efg
The first-order Godunov's method gives the solutions shown in \reff {fig:general13d} and \reff {fig:general12d}.
 The relative errors and convergence rates are given by \reft {tb:general1}. From the table, we can see that for $L^1$ norm the method is almost of first order, but for $L^2$ or $L^{\infty}$ norms, the rate of convergence is lower.

\btb\bc
\begin{tabular}{lccccccc}\\\hline
$\rho$&128-64&Rate&256-128&Rate&512-256&Rate&1024-512\\\hline
$L^1$&7.36e-04&9.72e-01&3.75e-04&9.95e-01&1.88e-04&9.69e-01&9.62e-05\\\hline
$L^2$&9.38e-04&6.65e-01&5.92e-04&6.81e-01&3.69e-04&6.74e-01&2.31e-04\\\hline
$L^{\infty}$&3.07e-03&4.02e-02&2.99e-03&1.79e-01&2.64e-03&3.27e-01&2.10e-03\\\hline\hline
$v$&128-64&Rate&256-128&Rate&512-256&Rate&1024-512\\\hline
$L^1$&9.68e-03&9.67e-01&4.95e-03&9.99e-01&2.48e-03&9.59e-01&1.27e-03\\\hline
$L^2$&1.27e-02&6.43e-01&8.14e-03&6.63e-01&5.14e-03&6.56e-01&3.27e-03\\\hline
$L^{\infty}$&4.41e-02&4.81e-02&4.27e-02&1.76e-01&3.78e-02&3.15e-01&3.04e-02\\\hline
\end {tabular}
\caption {Convergence rates for  the first-order Godunov's method} \label {tb:general1}\ec
\etb

\bfg
\bc\includegraphics[height=8cm] {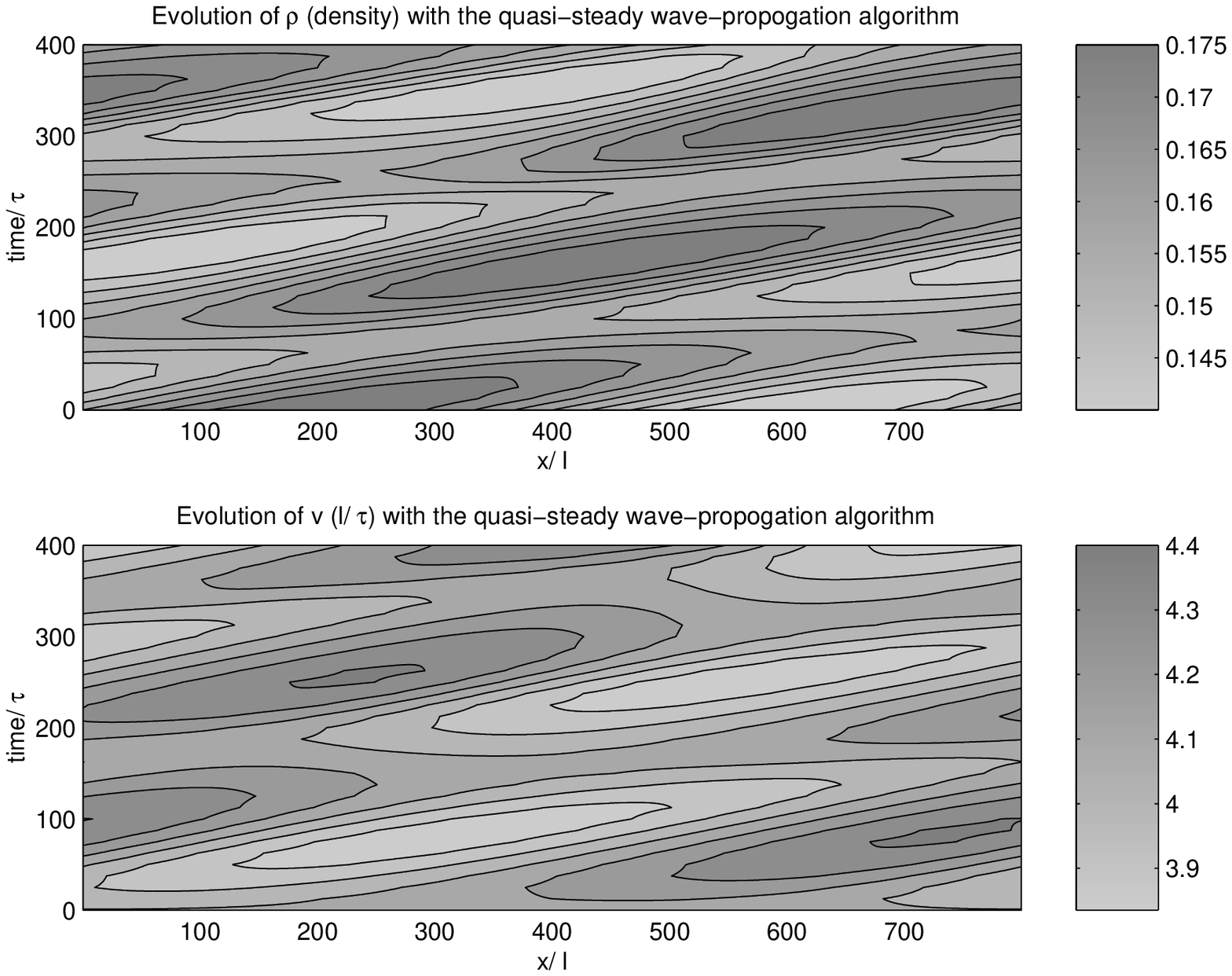}\ec
\caption {Solutions by LeVeque's method for 1024 grids} \label {fig:LeVeque13d}
\efg
\bfg
\bc\includegraphics[height=8cm] {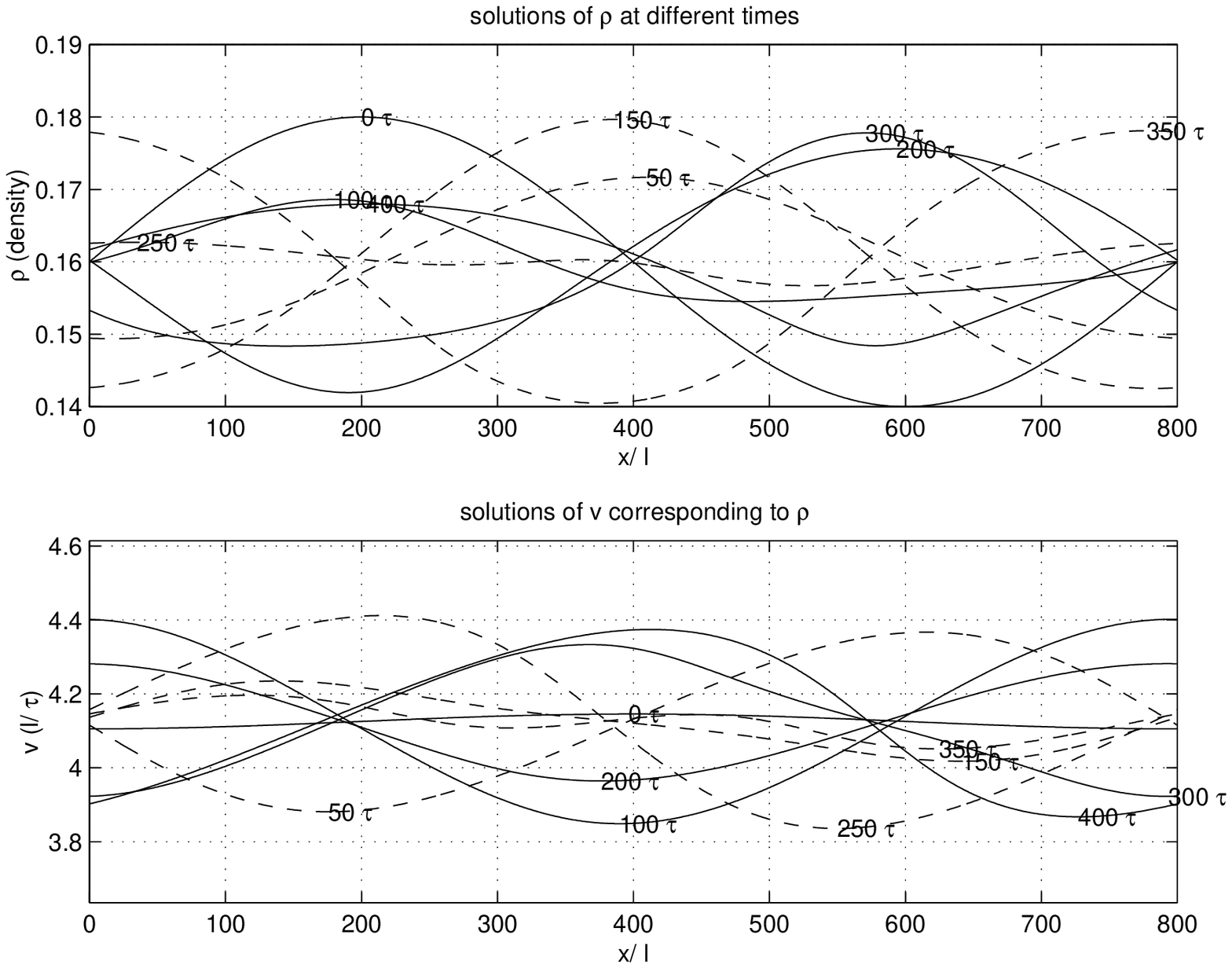}\ec
\caption {Solutions from \reff {fig:LeVeque13d} at selected times} \label {fig:LeVeque12d}
\efg

The quasi-steady wave-propagation algorithm by LeVeque (1998a\nocite {LeVeque98a}, 1998b\nocite {LeVeque98b}) gives the solutions shown in \reff {fig:LeVeque13d} and \reff {fig:LeVeque12d} for 1024 grids. The relative errors and convergence rates are given in \reft {tb:LeVeque1}. This scheme de-estimate the effects of the source term, since the convergence rates of $\rho$ and $v$ are totally different.

\btb\bc
\begin{tabular}{lccccccc}\\\hline
$\rho$&256-128&Rate&512-256&Rate&1024-512\\\hline
$L^1$&2.09e-02&5.10e-01&1.47e-02&1.16e+00&6.54e-03\\\hline
$L^2$&2.28e-02&4.88e-01&1.63e-02&1.14e+00&7.40e-03\\\hline
$L^{\infty}$&3.49e-02&5.81e-01&2.33e-02&8.86e-01&1.26e-02\\\hline\hline
$v$&256-128&Rate&512-256&Rate&1024-512\\\hline
$L^1$&3.42e-02&-1.61e+00&1.05e-01&6.18e-02&1.00e-01\\\hline
$L^2$&3.91e-02&-1.58e+00&1.17e-01&5.78e-02&1.12e-01\\\hline
$L^{\infty}$&6.96e-02&-1.30e+00&1.71e-01&-1.51e-01&1.90e-01\\\hline
\end {tabular}\ec

\caption {Convergence rates for the quasi-steady wave-propagation algorithm} \label {tb:LeVeque1}
\etb

\bfg
\bc\includegraphics[height=8cm] {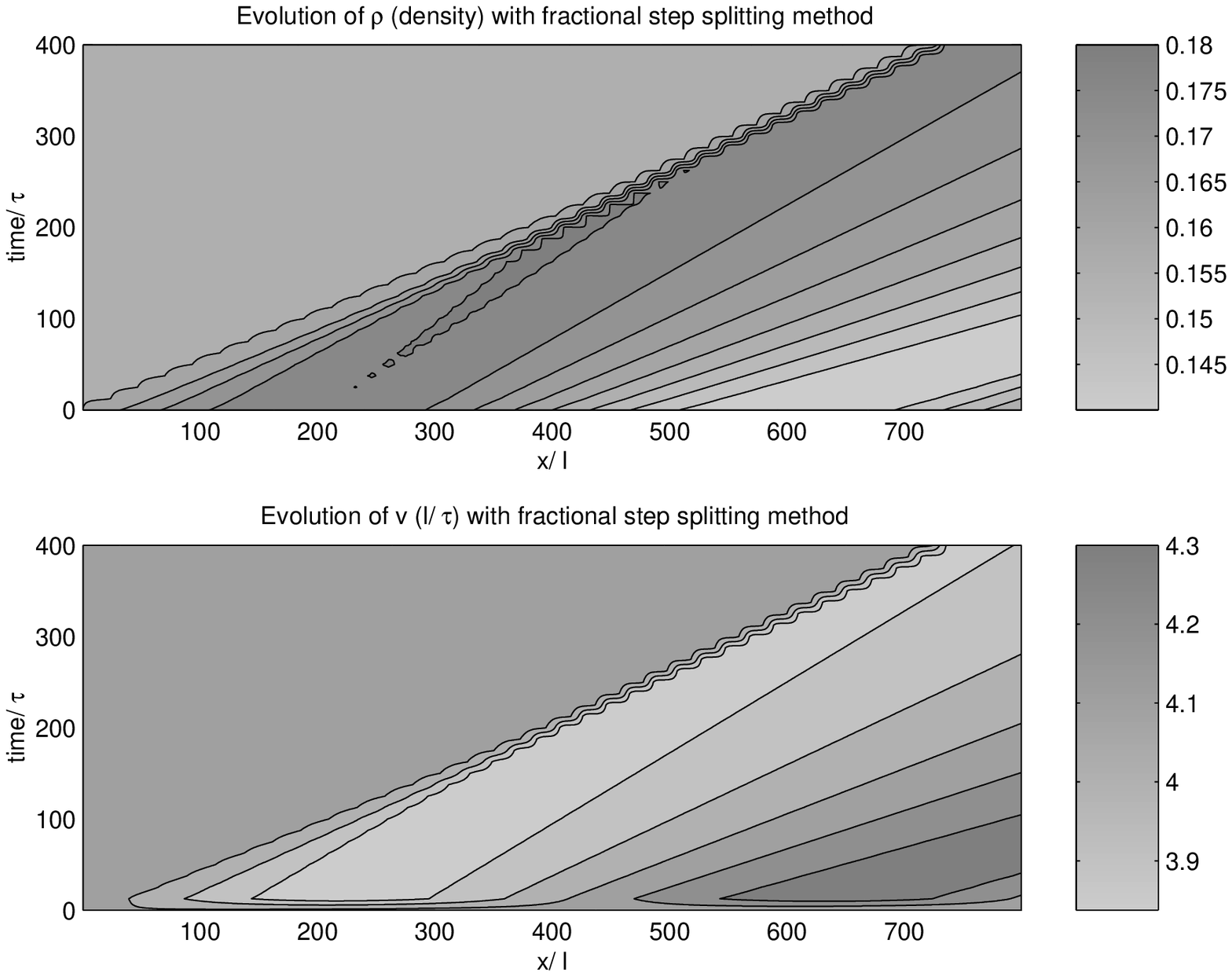}\ec
\caption {Solutions by fractional splitting method for 1024 grids} \label {fig:Fractional13d}
\efg
\bfg
\bc\includegraphics[height=8cm] {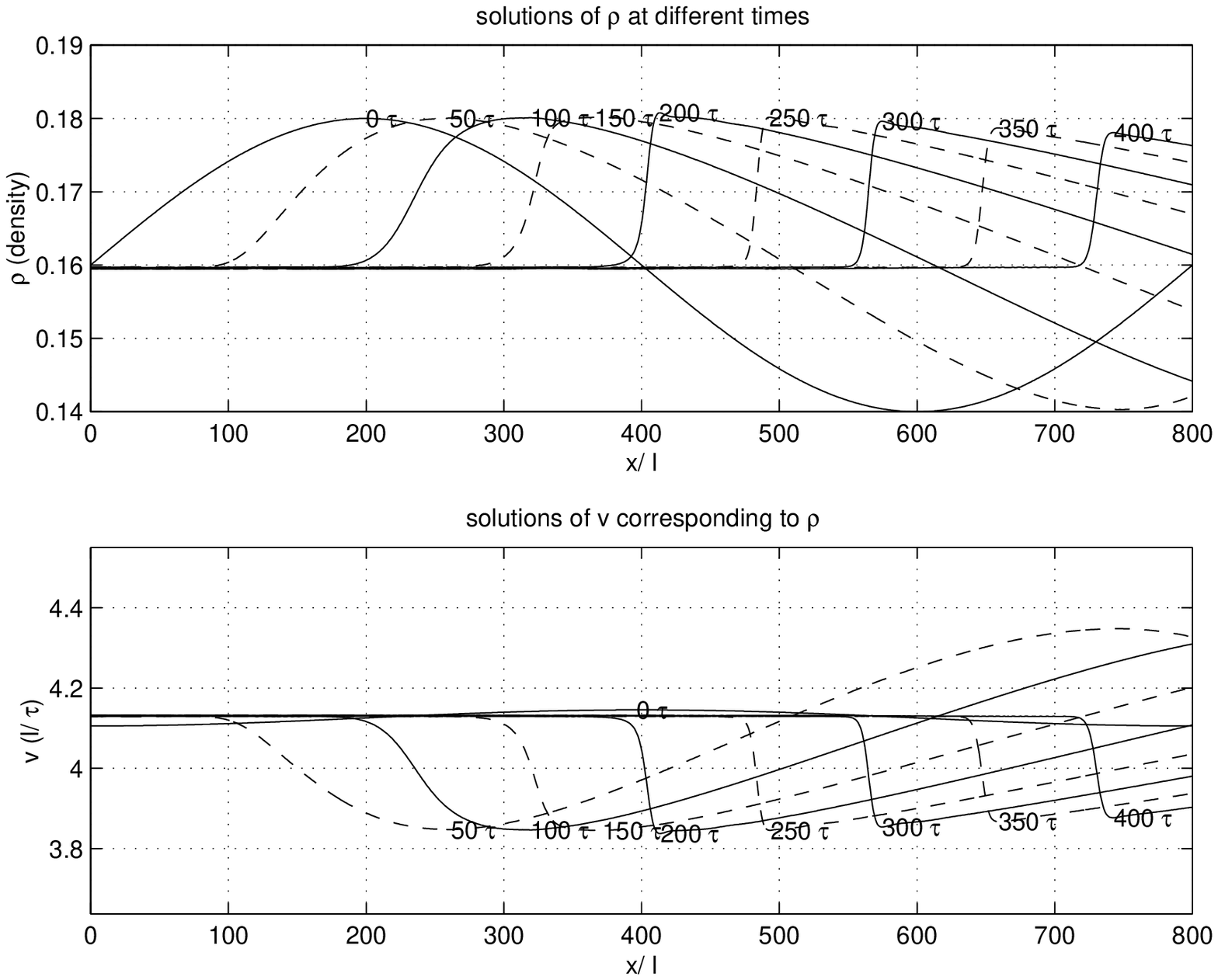}\ec
\caption {Solutions from \reff {fig:Fractional13d} at different times} \label {fig:Fractional12d}
\efg

The fractional step splitting method gives solutions of the PW model shown in \reff {fig:Fractional13d} and \reff {fig:Fractional12d} for 1024 grids. The relative errors and convergence rates are given in \reft{tb:Fractional1}. We find that the fractional step splitting method is not so stable as the first-order method, since convergence rates have big oscillations.
\btb\bc
\begin{tabular}{lccccccc}\\\hline
$\rho$&128-64&Rate&256-128&Rate&512-256&Rate&1024-512\\\hline
$L^1$&8.27e-04&9.49e-01&4.28e-04&-5.45e-01&6.25e-04&1.38e+00&2.41e-04\\\hline
$L^2$&1.00e-03&7.35e-01&6.01e-04&-2.72e-01&7.26e-04&5.81e-01&4.85e-04\\\hline
$L^{\infty}$&2.83e-03&1.95e-01&2.47e-03&-4.07e-01&3.28e-03&-4.48e-01&4.47e-03\\\hline\hline
$v$&128-64&Rate&256-128&Rate&512-256&Rate&1024-512\\\hline
$L^1$&1.08e-02&9.46e-01&5.58e-03&-5.01e-01&7.90e-03&1.35e+00&3.11e-03\\\hline
$L^2$&1.33e-02&7.14e-01&8.13e-03&-1.92e-01&9.28e-03&4.76e-01&6.67e-03\\\hline
$L^{\infty}$&3.86e-02&1.19e-01&3.55e-02&-2.77e-01&4.30e-02&-5.32e-01&6.22e-02\\\hline
\end {tabular}\ec

\caption {Convergence rates for the fractional step splitting method} \label {tb:Fractional1}
\etb

\bfg
\bc\includegraphics[height=8cm] {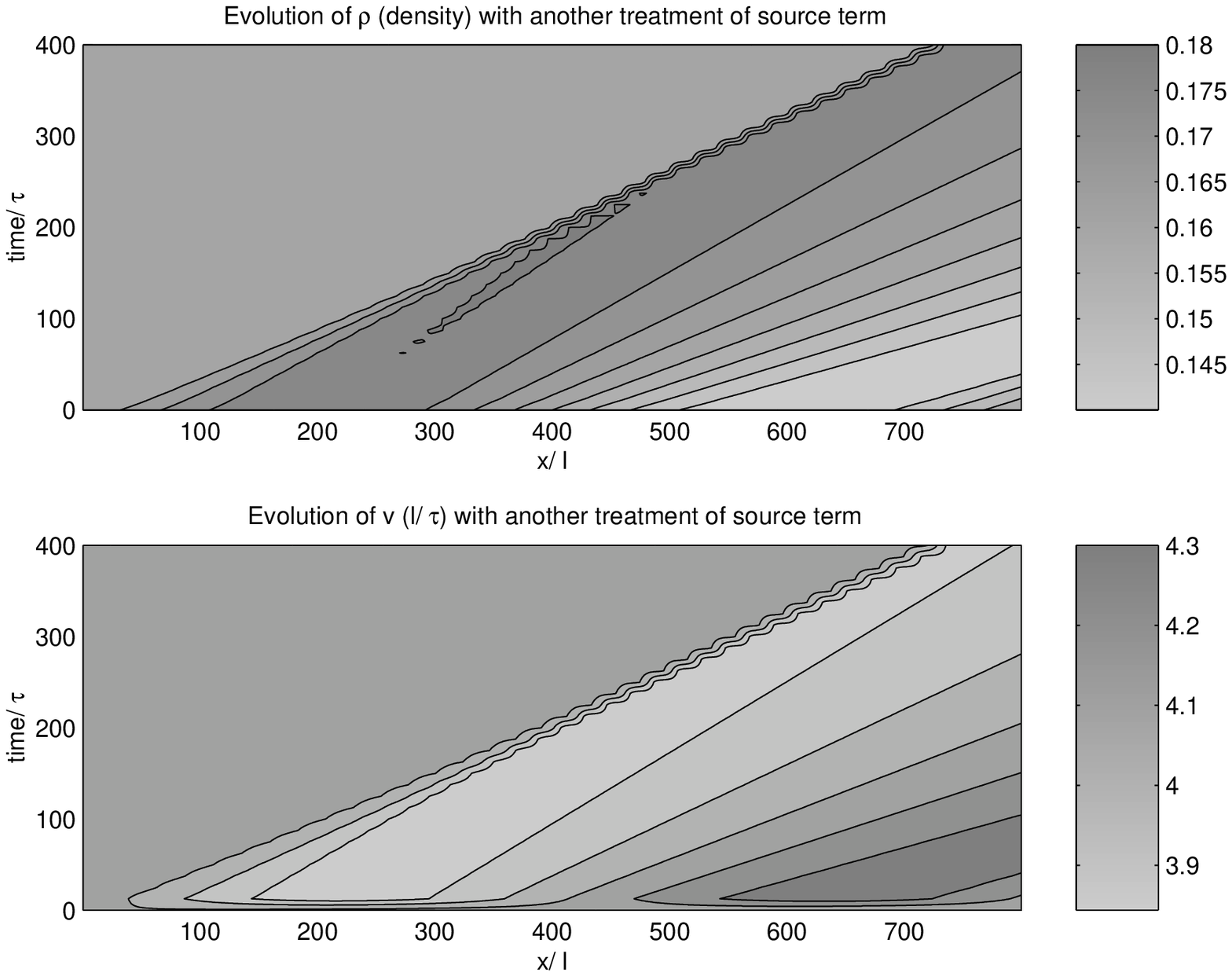}\ec
\caption {Solutions by Pember's method for 1024 grids} \label {fig:Pember13d}
\efg
\bfg
\bc\includegraphics[height=8cm] {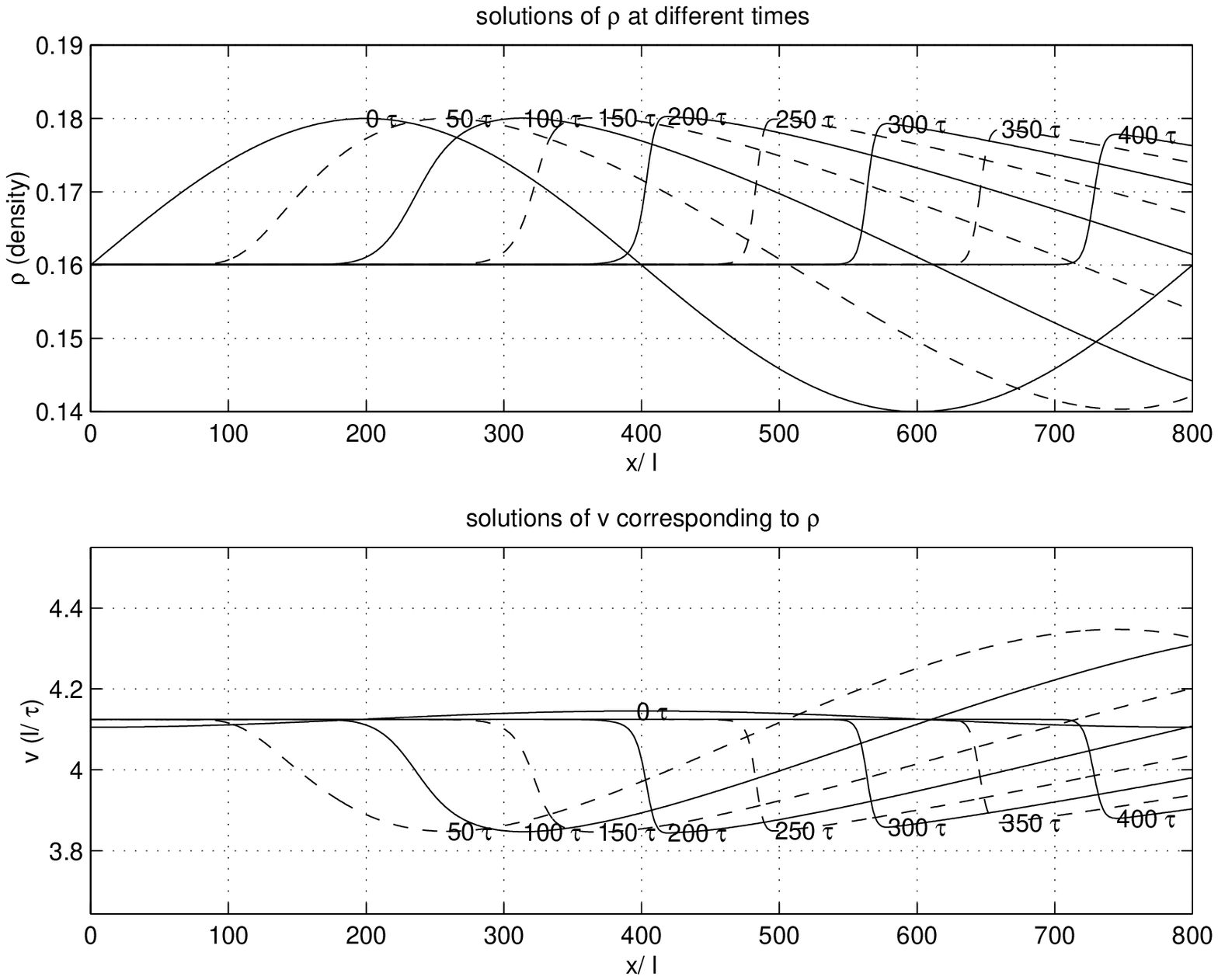}\ec
\caption {Solutions from \reff {fig:Pember13d} at different times} \label {fig:Pember12d}
\efg

%Fri Nov 12 22:13:18 PST 1999
Pember's method (1993a\nocite {Pember93a}, 1993b\nocite {Pember93b}) give the solutions shown in \reff {fig:Pember13d} and \reff {fig:Pember12d} for 1024 grids. The relative errors and convergence rates are given in \reft{tb:Pember1}.

\btb\bc
\begin{tabular}{lccccccc}\\\hline
$\rho$&128-64&Rate&256-128&Rate&512-256&Rate&1024-512\\\hline
$L^1$&8.31e-04&7.69e-01&4.88e-04&9.73e-01&2.49e-04&1.02e+00&1.22e-04\\\hline
$L^2$&1.07e-03&4.67e-01&7.76e-04&6.39e-01&4.98e-04&6.82e-01&3.10e-04\\\hline
$L^{\infty}$&2.94e-03&-4.29e-02&3.03e-03&9.29e-02&2.84e-03&2.85e-01&2.33e-03\\\hline\hline
$v$&128-64&Rate&256-128&Rate&512-256&Rate&1024-512\\\hline
$L^1$&1.10e-02&7.60e-01&6.50e-03&9.73e-01&3.31e-03&1.02e+00&1.64e-03\\\hline
$L^2$&1.46e-02&4.52e-01&1.06e-02&6.33e-01&6.87e-03&6.74e-01&4.30e-03\\\hline
$L^{\infty}$&4.09e-02&-4.25e-02&4.21e-02&9.82e-02&3.94e-02&2.78e-01&3.24e-02\\\hline
\end {tabular}\ec

\caption {Convergence rates for Pember's method} \label {tb:Pember1}
\etb

\bfg
\bc\includegraphics[height=8cm] {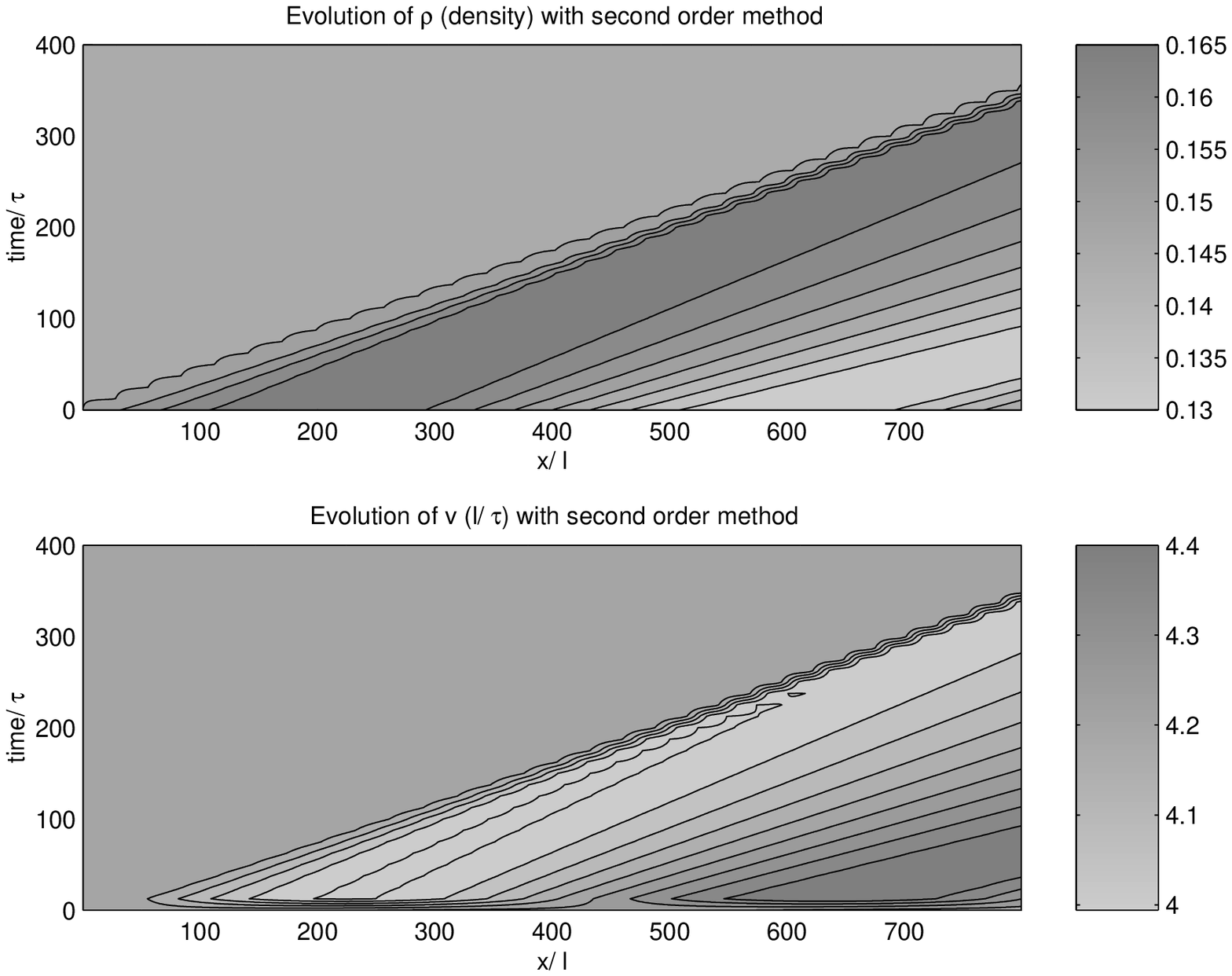}\ec
\caption {Solutions by the second-order Godunov's method for 1024 grids} \label {fig:second13d}
\efg
\bfg
\bc\includegraphics[height=8cm] {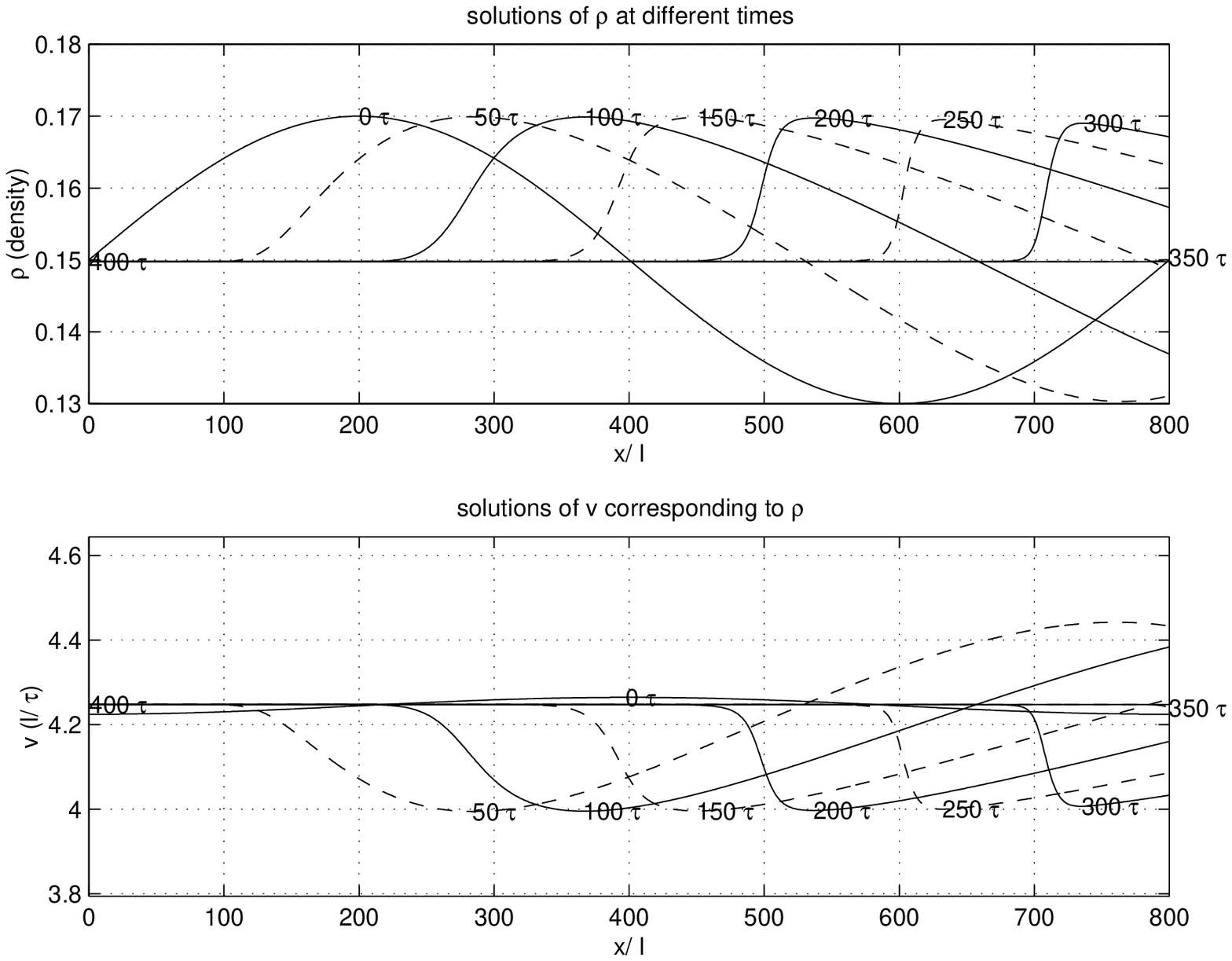}\ec
\caption {Solutions from \reff {fig:second13d} at selected times} \label {fig:second12d}
\efg

In the remaining part of this subsection we consider the second-order Godunov's method. We use $\rh=0.15$ for initial conditions \refet {ini1:1}{ini1:2}, since the second-order method is not stable for $\rh=0.16$. For 1024 grids, we get the solutions shown in \reff {fig:second13d} and \reff {fig:second12d}. The relative errors and convergence rates are given in \reft {tb:second1}.
\btb\bc
\begin{tabular}{lccccccc}\\\hline
$\rho$&128-64&Rate&256-128&Rate&512-256&Rate&1024-512\\\hline
$L^1$&5.34e-04&1.29e+00&2.18e-04&1.19e+00&9.55e-05&1.03e+00&4.69e-05\\\hline
$L^2$&5.34e-04&1.29e+00&2.18e-04&1.19e+00&9.55e-05&1.03e+00&4.69e-05\\\hline
$L^{\infty}$&5.34e-04&1.29e+00&2.18e-04&1.17e+00&9.65e-05&1.01e+00&4.78e-05\\\hline\hline
$v$&128-64&Rate&256-128&Rate&512-256&Rate&1024-512\\\hline
$L^1$&6.02e-03&1.30e+00&2.45e-03&1.19e+00&1.07e-03&1.03e+00&5.26e-04\\\hline
$L^2$&6.02e-03&1.30e+00&2.45e-03&1.19e+00&1.07e-03&1.03e+00&5.26e-04\\\hline
$L^{\infty}$&6.02e-03&1.30e+00&2.45e-03&1.18e+00&1.08e-03&1.01e+00&5.36e-04\\\hline
\end {tabular}\ec

\caption {Convergence rate for second-order method} \label {tb:second1}
\etb

Comparing the convergence rates with those for first-order Godunov's method, we see no significant improvement. This is different from the case for Zhang's model. This is a special property of the PW model.

%\input{figurePW.tex}

%Tue Nov 30 14:58:01 PST 1999
\subsection {Unstable solutions of the PW model}
In this subsection we check the unstable solutions of the PW model. We use the first-order Godunov's method with periodic boundary conditions, and we observe the solutions at  $T_0=200\tau$ for different number of grids.

We use the initial conditions \refet{ini1:1}{ini1:2} with $\rh=0.175$, which is in the unstable region of the PW model. Solutions of the PW model are listed as the following  for different number of grids:
\ben
\item For 512 grids, solutions are given in \reff {fig:unst1_512_3d} and \reff {fig:unst1_512_2d}.

\item For 1024 grids, solutions are given in \reff {fig:unst1_1024_3d} and \reff {fig:unst1_1024_2d}.

\item For 2048 grids, solutions are given in \reff {fig:unst1_2048_3d} and \reff {fig:unst1_2048_2d}.
\een

The figures show that the number and position of spikes are different for different number of grids. This difference is caused by different approximations used by Godunov's method for different number of grids, since the PW model is unstable in the region where the initial conditions are.    

\bfg
\bc\includegraphics[height=8cm] {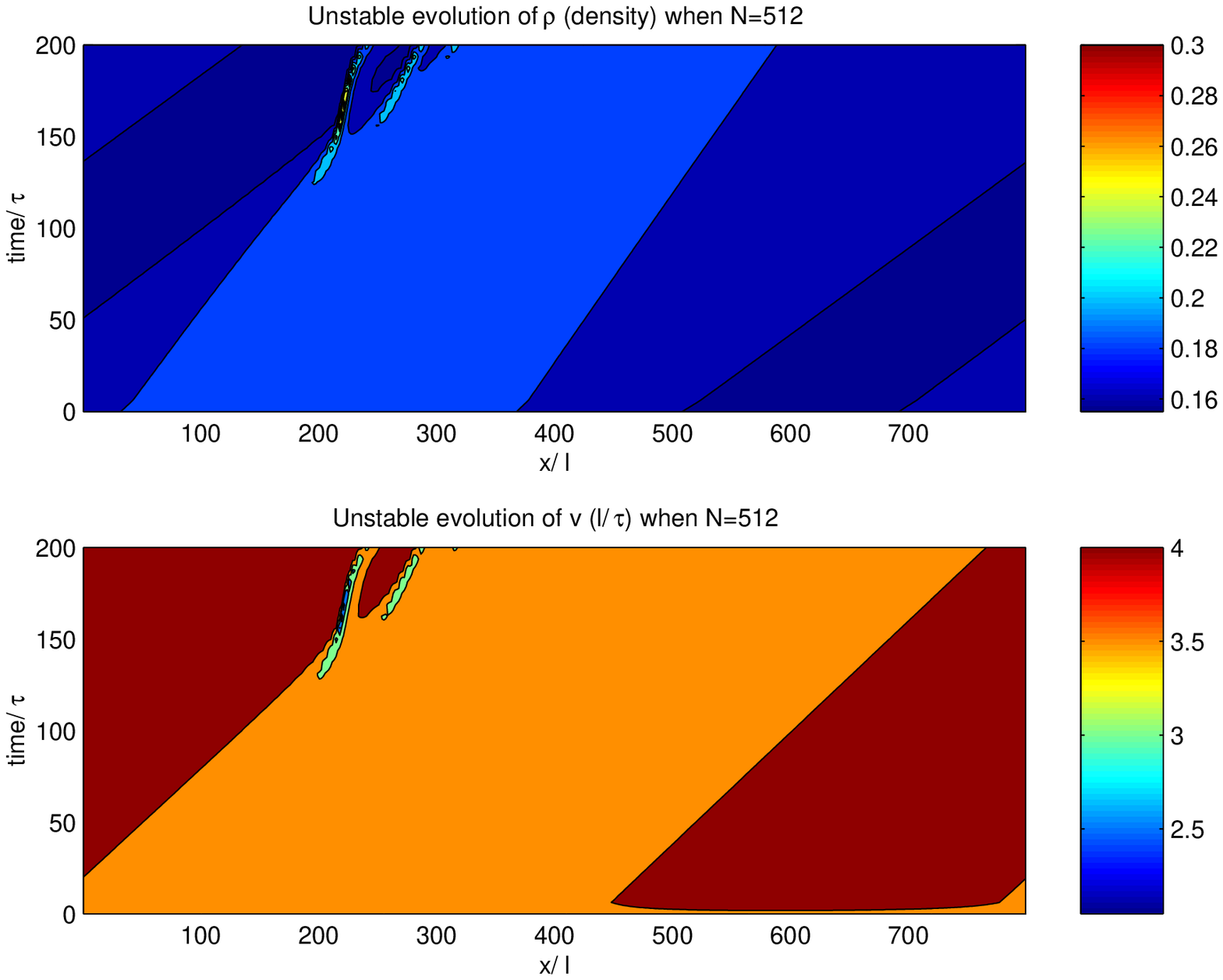}\ec
\caption {Solutions for 512 grids with initial conditions \refet {ini1:1}{ini1:2}} \label {fig:unst1_512_3d}
\efg
\bfg
\bc\includegraphics[height=8cm] {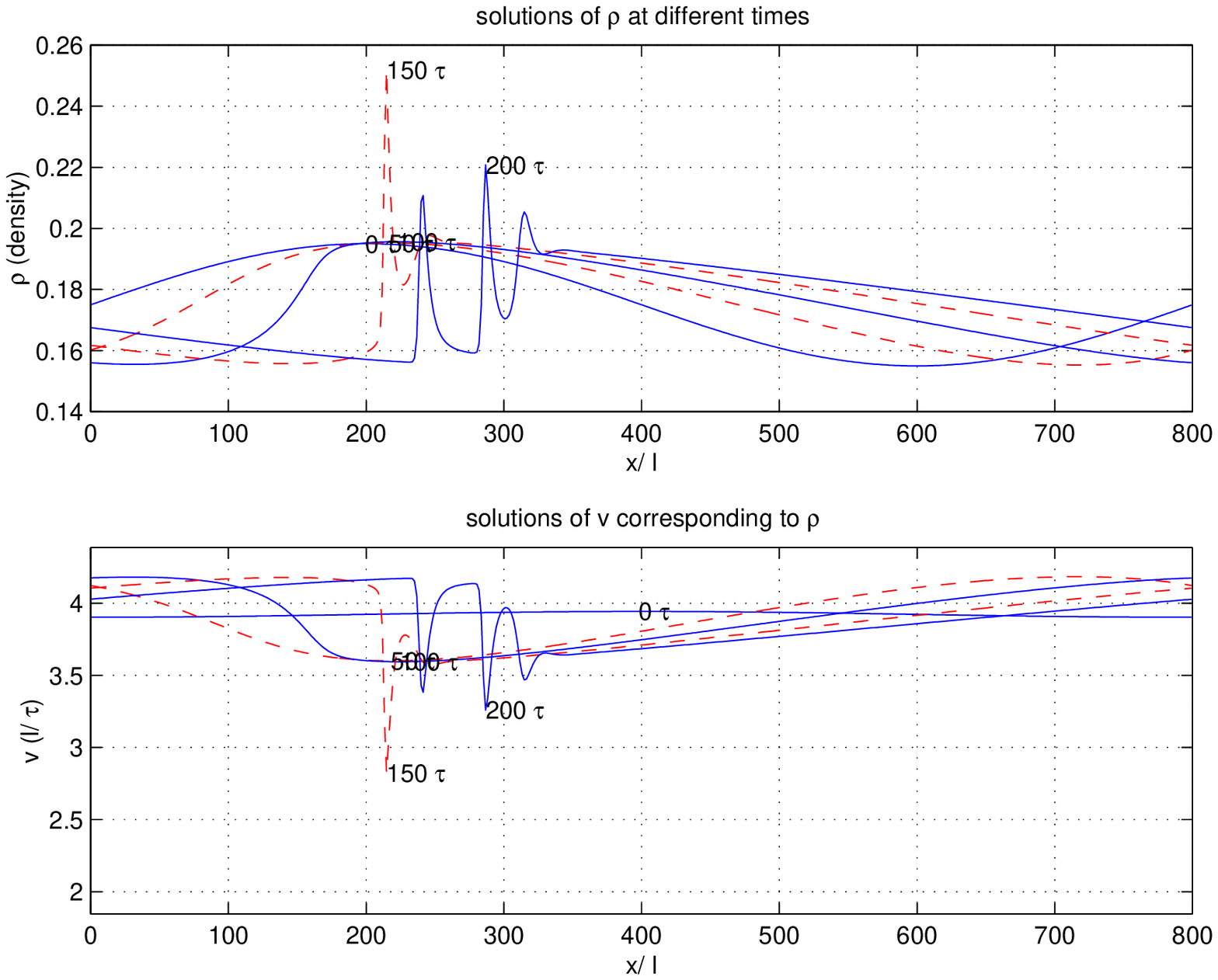}\ec
\caption {Solutions from \reff {fig:unst1_512_3d} at selected times} \label {fig:unst1_512_2d}
\efg

\bfg
\bc\includegraphics[height=8cm] {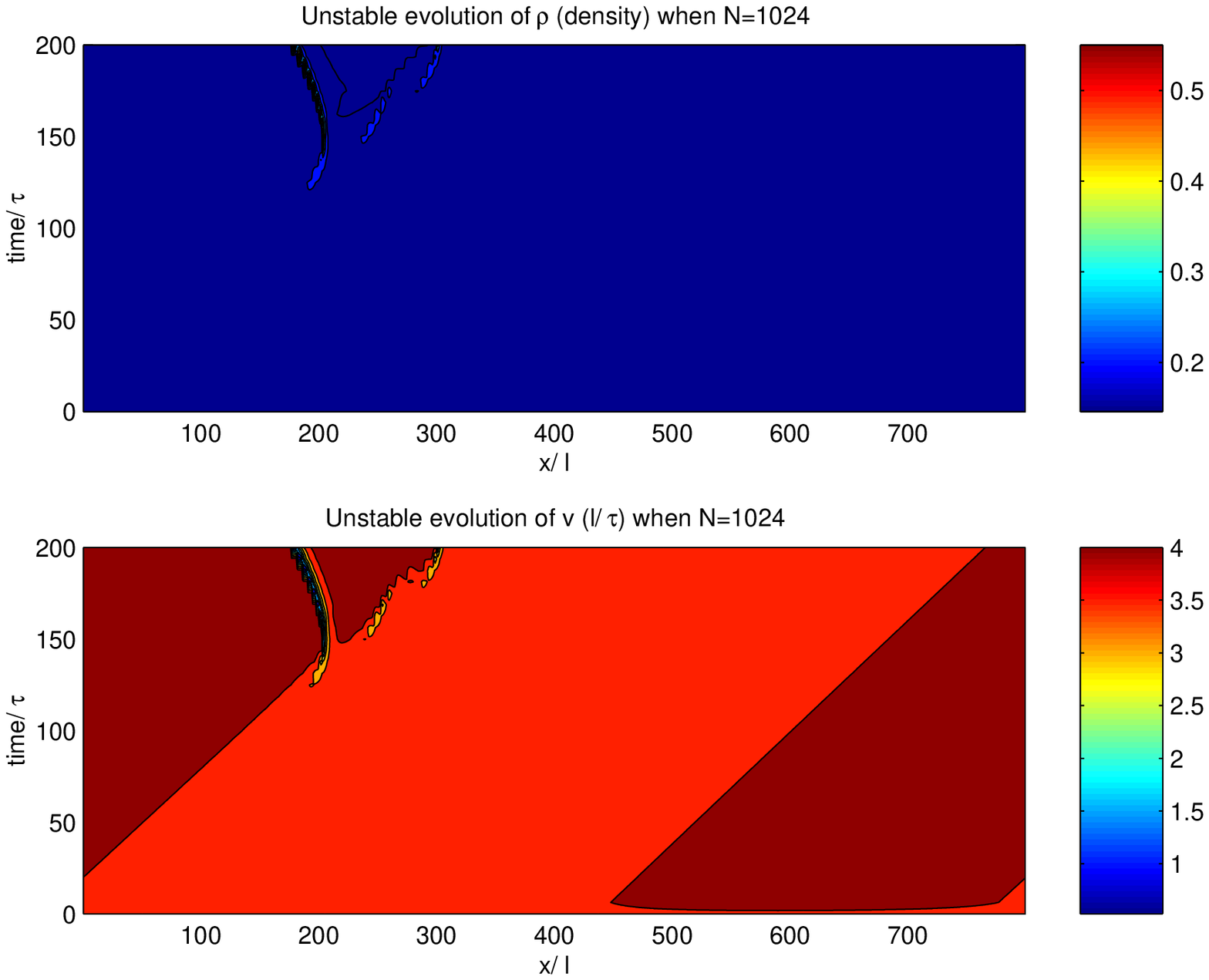}\ec
\caption {Solutions for 1024 grids with initial conditions \refet {ini1:1}{ini1:2}} \label {fig:unst1_1024_3d}
\efg
\bfg
\bc\includegraphics[height=8cm] {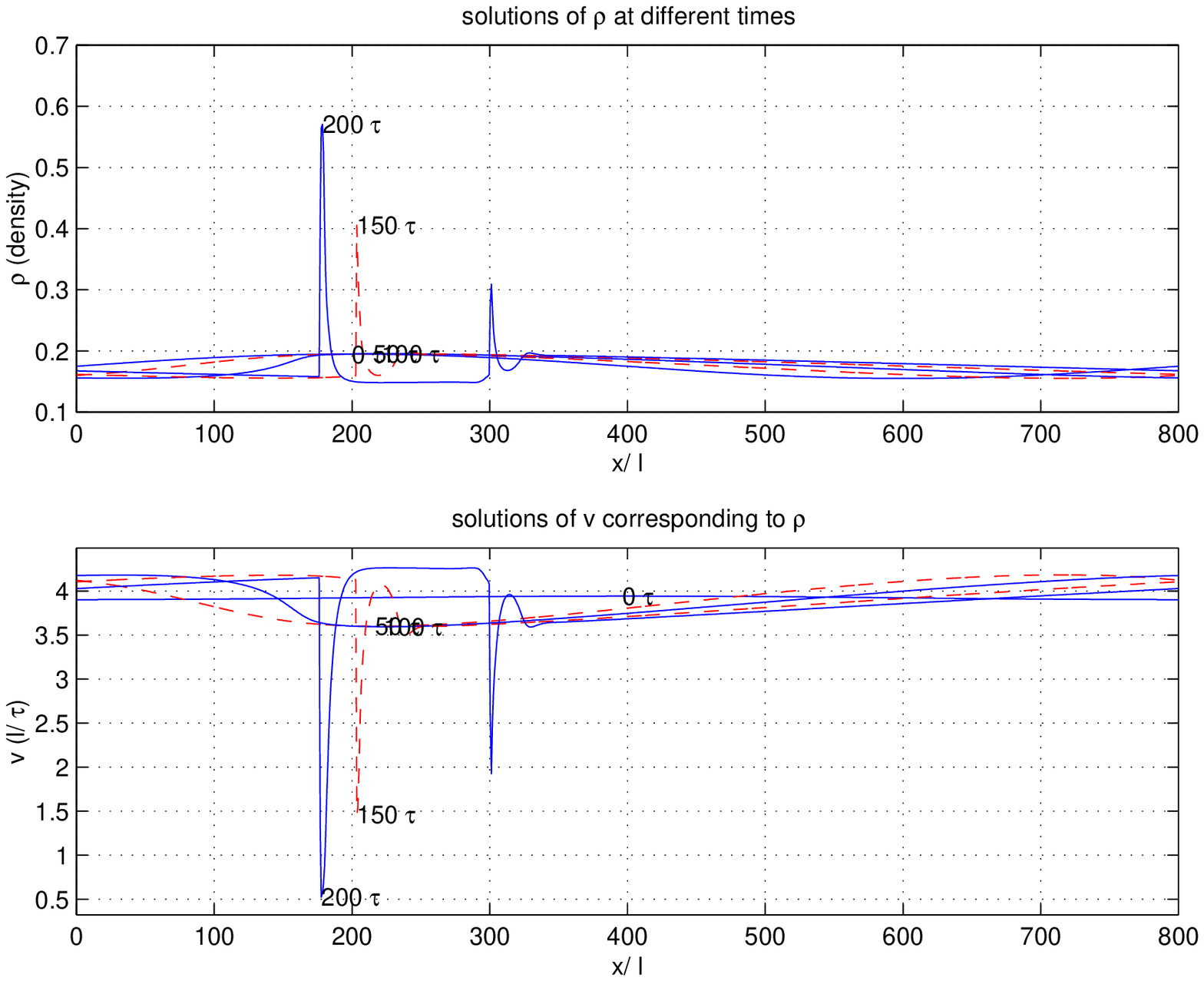}\ec
\caption {Solutions from \reff {fig:unst1_1024_3d} at selected times} \label {fig:unst1_1024_2d}
\efg

\bfg
\bc\includegraphics[height=8cm] {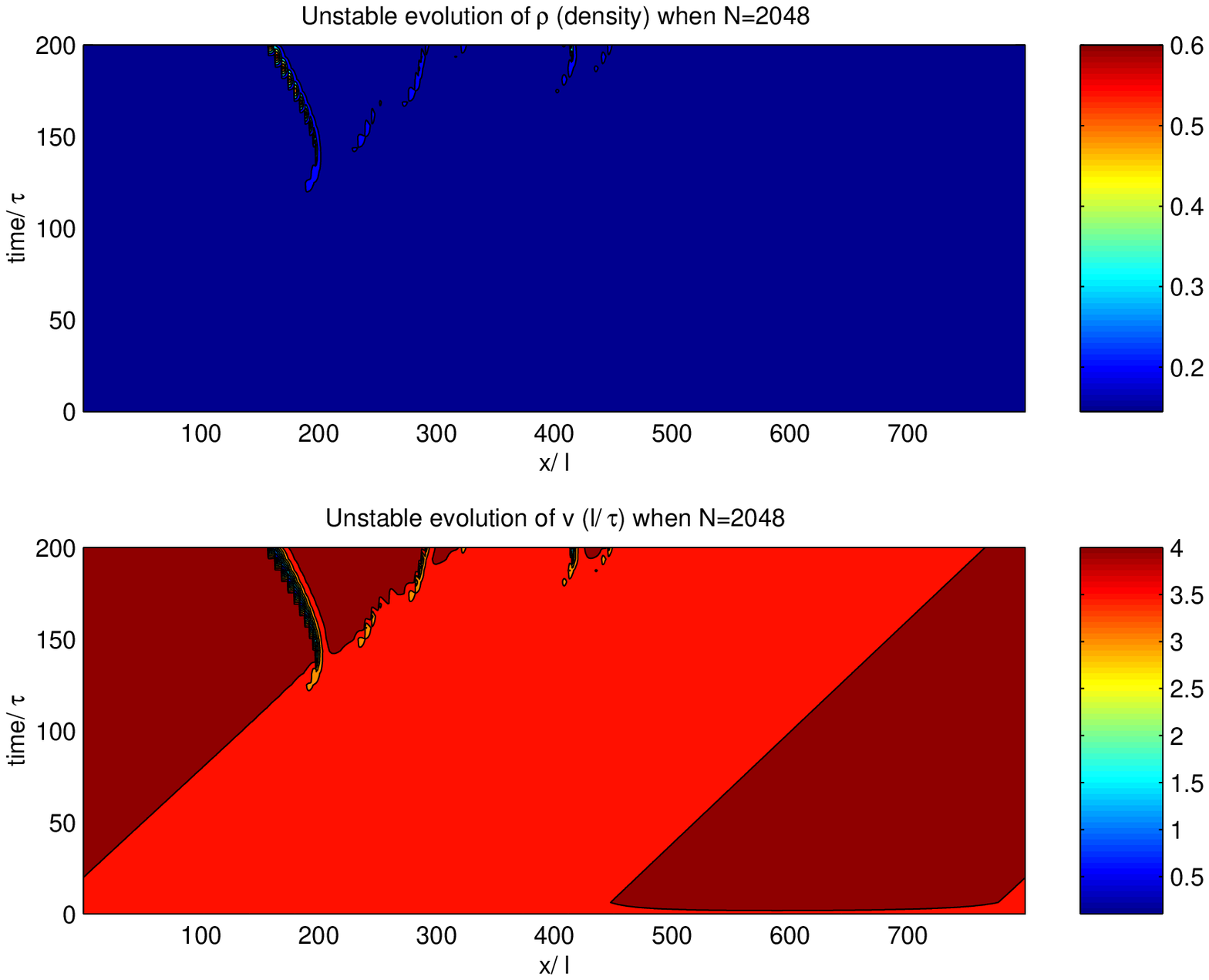}\ec
\caption {Solutions for 2048 grids with initial conditions \refet {ini1:1}{ini1:2}} \label {fig:unst1_2048_3d}
\efg

\bfg
\bc\includegraphics[height=8cm] {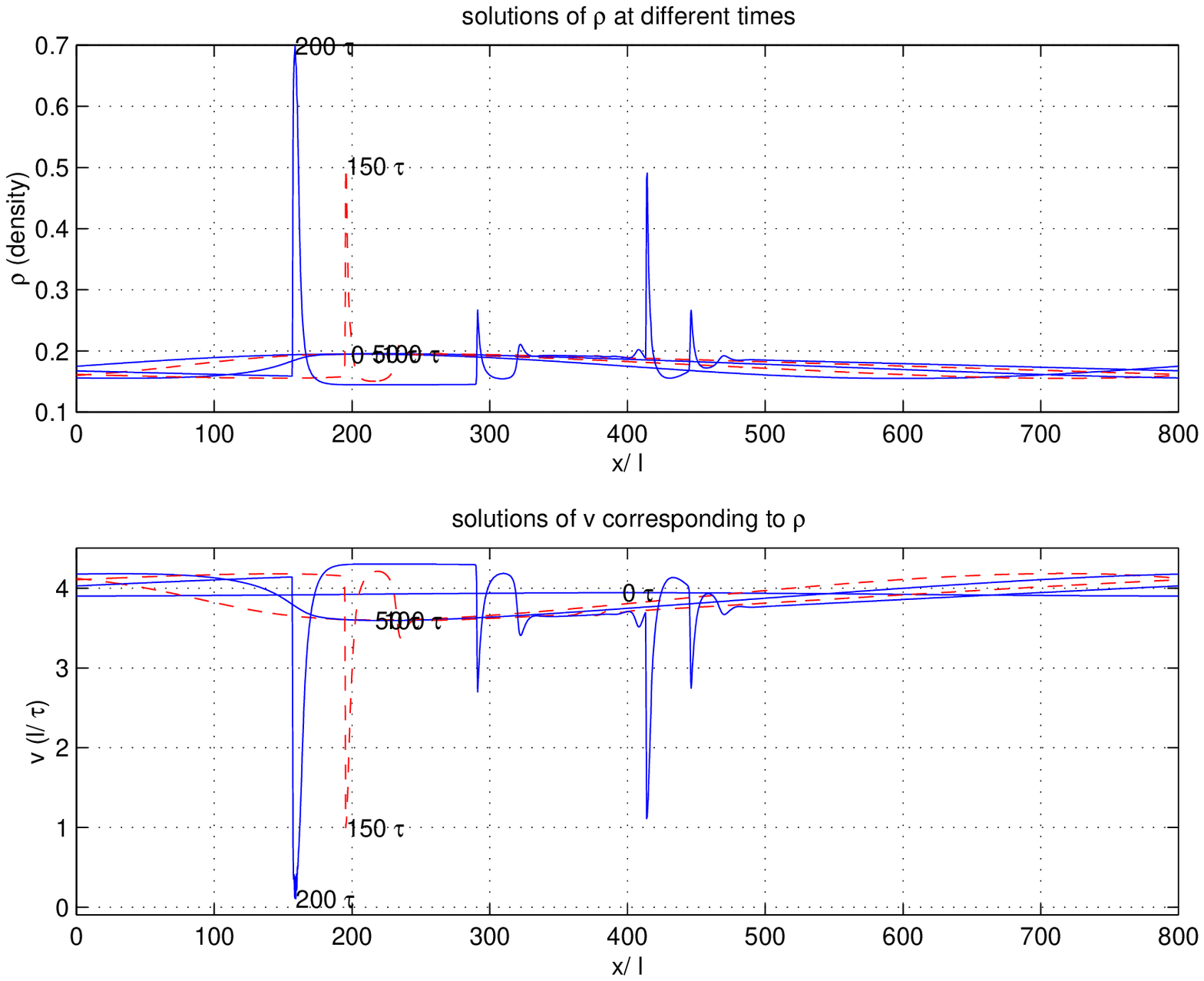}\ec
\caption {Solutions from \reff {fig:unst1_2048_3d} at selected times} \label {fig:unst1_2048_2d}
\efg

\newpage
\pagestyle{myheadings} 
\markright{  \rm \normalsize CHAPTER 5. \hspace{0.5cm}
 The Inhomogeneous LWR Model and Its Numerical Solutions}
\large
\chapter{The Inhomogeneous LWR Model and Its Numerical Solutions}

\section{Introduction}
The LWR model (Lighthill and Whitham, 1955\nocite{lighthill};Richards, 1956\nocite{richards}) was introduced based on the conservation of traffic flow, and is written as:
\bqn
\pd{}{t} \r+\pd{}{x} f&=&0.
\eqn

The LWR model assumes that traffic flow is in equilibrium, or equivalently that travel speed $v$ (unit: mph) is defined as a function of traffic density $\r$ (unit: vpm) at any location $x$:
\bqn
v=\vs(x,\r).
\eqn
The function of flow rate $f=\r \vs(x,\r)$ (unit: vph) is called a fundamental diagram. For typical equilibrium traffic flows, travel speed is a decreasing function; i.e., $v_{\ast\r}<0$, and traffic flow rate is a concave function; i.e., $f_{\r\r}<0$.

The homogeneous LWR model is for modeling a homogeneous road, where travel speed $\vs$ is uniform with respect to location $x$. The homogeneous LWR model can thus be written as
\bqn
\r_t+f(\r)_x&=&0. \label{hom_LWR}
\eqn

The homogeneous LWR model is a scalar conservation law, and there have been many methods to compute its entropy solutions (Lebacque, 1995\nocite{lebacque1995}). It's well-known that the entropy solutions exist and are unique under the so-called ``Lax entropy condition". For computation of these entropy solutions, the Godunov method is most often used.

For a road with inhomogeneity, such as variable number of lanes, curvature and slope, we can formulate the inhomogeneous LWR model as
\bqn
\r_t+f(a,\r)_x&=&0, \label{inh_1st}
\eqn
where $a=a(x)$ is a time-invariant variable. The inhomogeneity factor $a(x)$ gives a profile of a piece of roadway, for example $a(x)$ can be the number of lanes at location $x$.

The inhomogeneous LWR model has been studied by Daganzo (1994\nocite{Daganzo94}) and Lebacque (1995\nocite{lebacque1995}) and Daganzo (1994\nocite{Daganzo94}). Both of these authors suggested solutions to the inhomogeneous LWR model, and their solutions are consitent. However, these studies only presented empirical solutions without rigorous proof. 

The difficulty of dealing with the inhomogeneous LWR model is due to the extra variable $a(x)$. Here, by introducing $a(x)$ as an additional conservation law, we consider the inhomogeneous LWR model as a resonant nonlinear system, which has been discussed in (Isaacson \& Temple, 1992\nocite{isaacson_t1992}; Lin et al., 1995\nocite{lin_t_w1995}). In this chapter, we follow the procedures provided in those researches to study the inhomogeneous LWR model. 

This chapter is organized as follows. In section 2, we formulate the inhomogeneous traffic model as a resonant nonlinear system, and its properties are discussed. In section 3, we solve Riemann problem for the inhomogeneous LWR model. In section 4 we present the numerical methods for this model. We conclude the discussions of this chapter in section 5.

\section{The Properties of the inhomogeneous LWR model}
To apply the results for hyperbolic systems of conservation laws, we express the inhomogeneity factor $a(x)$ as an additional conservation law, i.e., $a_t=0$. Hence, we can write the inhomogeneous LWR model as
\bqn
U_t+F(U)_x&=&0,\label{system}
\eqn
where $U=(a,\r), F(U)=(0,f(a,\r)),x\in R,t\geq 0$. In this chapter we consider one type of inhomogeneity -- variable number of lanes. The fundamental diagram is thus written as $f(a,\r)=\r\vs(\ra)$, given that all the lanes are of the same condition.

The inhomogeneous LWR model \refe{system} can be linearized as
\bqn
U_t+DF(U) U_x&=&0, \label{linear}
\eqn
where the differential $DF(U)$ of the flux vector $F(U)$ is  
\bqn
DF&=&\mat{{cc}0&0\\-\frac {\r^2}{a^2} \vs'(\ra)&\vs(\ra)+\ra\vs'(\ra)}.
\eqn
The two eigenvalues for $DF$ are
\bqn
\l_0=0\qquad \l_1=\vs(\ra)+\ra\vs'(\ra).
\eqn
The corresponding right eigenvectors are
\bqs
\vecr_0=\mat{{c}\vs(\ra)+\ra\vs'(\ra)\\(\ra)^2\vs'(\ra)}
\qquad \vecr_1=\mat{{c}0\\1},
\eqs
and the left eigenvector of $\partial f/\partial \r$ as $\vecl_1=1$.
The system \refe{system} is a non-strictly hyperbolic system, since $\l_1$ may be equal to $\l_0=0$.

We consider a traffic state $U_{\ast}=(a_{\ast},\r_{\ast})$ as critical if
\bqn
\lambda_1(U_{\ast})=0\label{assum_1};
\eqn
i.e., at critical states, the two wave speeds are the same and system \refe{system} is singular. 
For a critical traffic state $U_{\ast}$ we also have
\bqn
\frac {\partial}{\partial \r} \l_1(U_{\ast})&=&f_{\r\r}<0 \label{assum_2},
\eqn
and
\bqn
\frac{\partial}{\partial a} f(U_{\ast})&=&-(\ra)^2\vs'(\ra)\big|_{U_{\ast}}=\ra \vs(\ra)\big |_{U_{\ast}}>0. \label{assum_3}
\eqn

A consequence of properties \refe{assum_2} and \refe{assum_3} is that the linearized system \refe{linear} at $U_{\ast}$ has the following normal form
\bqn
\mat{{c} \delta a\\\delta \r}_t+\mat{{cc}0&0\\1&0}\mat{{c} \delta a\\\delta \r}_x&=&0. \label{nlinear}
\eqn
The system \refe{nlinear} has the solution $\delta \r(x,t)=\delta a'(x)t+c$, and the solution goes to infinity as $t$ goes to infinity. Therefore \refe{nlinear} is a linear resonant system, and the original inhomogeneous LWR model \refe{system} is a nonlinear resonant system.

For \refe{system}, the smooth curve $\Gamma$ in $U$-space formed by all critical states $U_{\ast}$ are  named a transition curve. Therefore $\Gamma$ is defined as 
\bqs
\Gamma&=&\left\{ U \big | \l_1(U)=0\right\}. 
\eqs
Since $\l_1(U)=\vs(\ra)+\ra\vs'(\ra)$, we obtain
\bqn
\Gamma &=&\left\{(a,\r)\big|\ra=\alpha, \m{ where } \alpha \m{ uniquely solves }\vs(\alpha)+\alpha \vs'(\alpha)=0\right\}; \label{trancur}
\eqn
i.e., the transition curve for \refe{system} is a straight line passing through the origin in $U$-space. In \refe{trancur}, $\alpha$ is unique since $f(a,\r)$ is concave in $\r$.

The entropy solutions to a nonlinear resonant system are different from those to a strict hyperbolic system of conservation laws. Isaacson \& Temple (1992\nocite{isaacson_t1992}) proved that solutions to the Riemann problem  for system \refe{system} exist and are unique with the conditions \refe{assum_1}-\refe{assum_3}. Lin et al. (1995\nocite{lin_t_w1995}) also presented solutions to a scalar nonlinear resonant system, which is similar to our system \refe{system} except that $f$ is convex. In the next section we apply those results to solve the Riemann problem for the inhomogeneous LWR model.

\section{Solutions to the Riemann problem}\label{ri}
\begin{figure}[t]\bc\includegraphics [height=10cm]{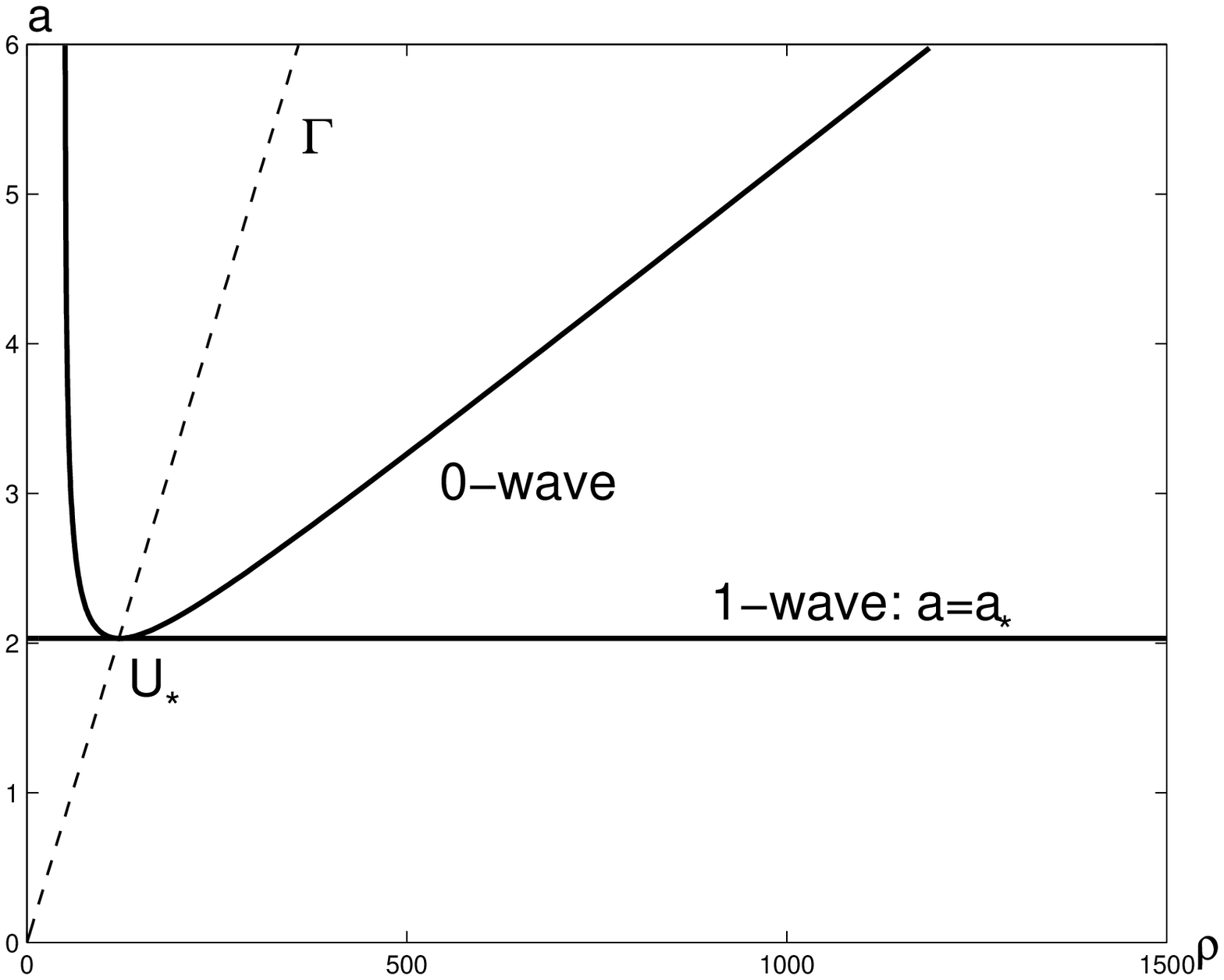}\ec\caption{Integral curves}\label{F_waves}\efg
In this section we study the wave solutions to the Riemann problem for \refe{system} with  the following jump initial conditions
\bqn
U(x,t=0)&=&\cas{{cc} U_L&\m{ if }x<0\\U_R &\m{ if }x>0}, \label{inh.ini}
\eqn
where the initial values of $U_L, U_R$ are constant. For computational purpose, we are interested in  the average flux at the boundary $x=0$ over a time interval $\dt$, which is denoted by $f^{\ast}_0$; i.e., we want to find
\bqn
f^{\ast}_0&=&\frac 1 {\dt} \int_0^{\dt} f(U(x=0,t)) \m{dt}.
\eqn

The inhomogeneous LWR model \refe{system} has two families of basic wave solutions associated to the two eigenvalues. The solutions whose wave speed is $\l_0$ are in the 0-family, and the waves are called 0-waves. Similarly the solutions whose wave speed is $\l_1$ are in the 1-family, and the waves are called 1-waves. In $U$-space, the wave curves for \refe{system} are the integral curves of the right eigenvectors $\vecr_0$ and $\vecr_1$. Hence the 0-wave curves are given by $f(U)$=const, and the 1-wave curves are given by $a=\bar a$, where $\bar a$ is constant. The 0-wave is also called a standing wave since its wave speed is always 0. The 1-wave solutions are determined by the solutions of the scalar conservation law $\r_t+f(\bar a,\r)_x=0$.
A 0-wave curve, a 1-wave curve passing a critical state $U_{\ast}$ and the transition curve $\Gamma$ are shown in \reff{F_waves}, where $a$ is set as the vertical axis and $\r$ is set as the horizontal axis.

As shown in \reff{F_waves}, the 0-wave curve is convex, and the 1-wave curve is tangent to the 0-wave curve at the critical state $U_{\ast}$. The transition curve $\Gamma$ intersects the 0-wave and 1-wave curves  transversely at $U_{\ast}$, and there is only one critical state on one 0-wave or 1-wave curve. For any point $U$, there is only one 0-wave curve and only one 1-wave curve passing it. In \reff{F_waves}, the states left to the transition curve are undercritical since $\r/a<\alpha$; and the states right to the transition curve are overcritical since $\r/a>\alpha$.

\begin{figure}[t]\bc\includegraphics [height=8cm] {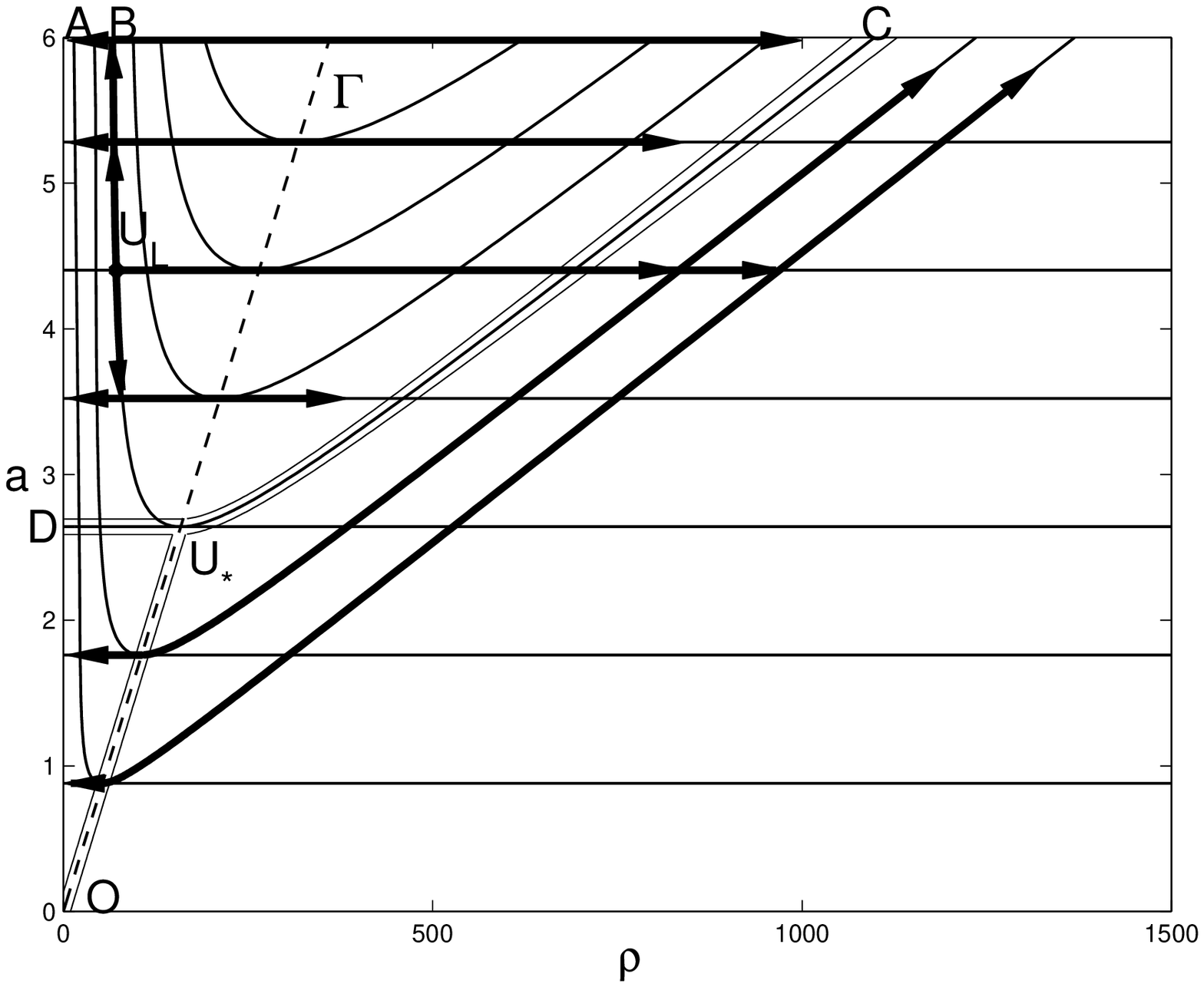}\ec\caption{The Riemann problem for $U_L$ left of $\Gamma$}\label{F_riemann1}\efg
\begin{figure}[t]\bc\includegraphics [height=8cm] {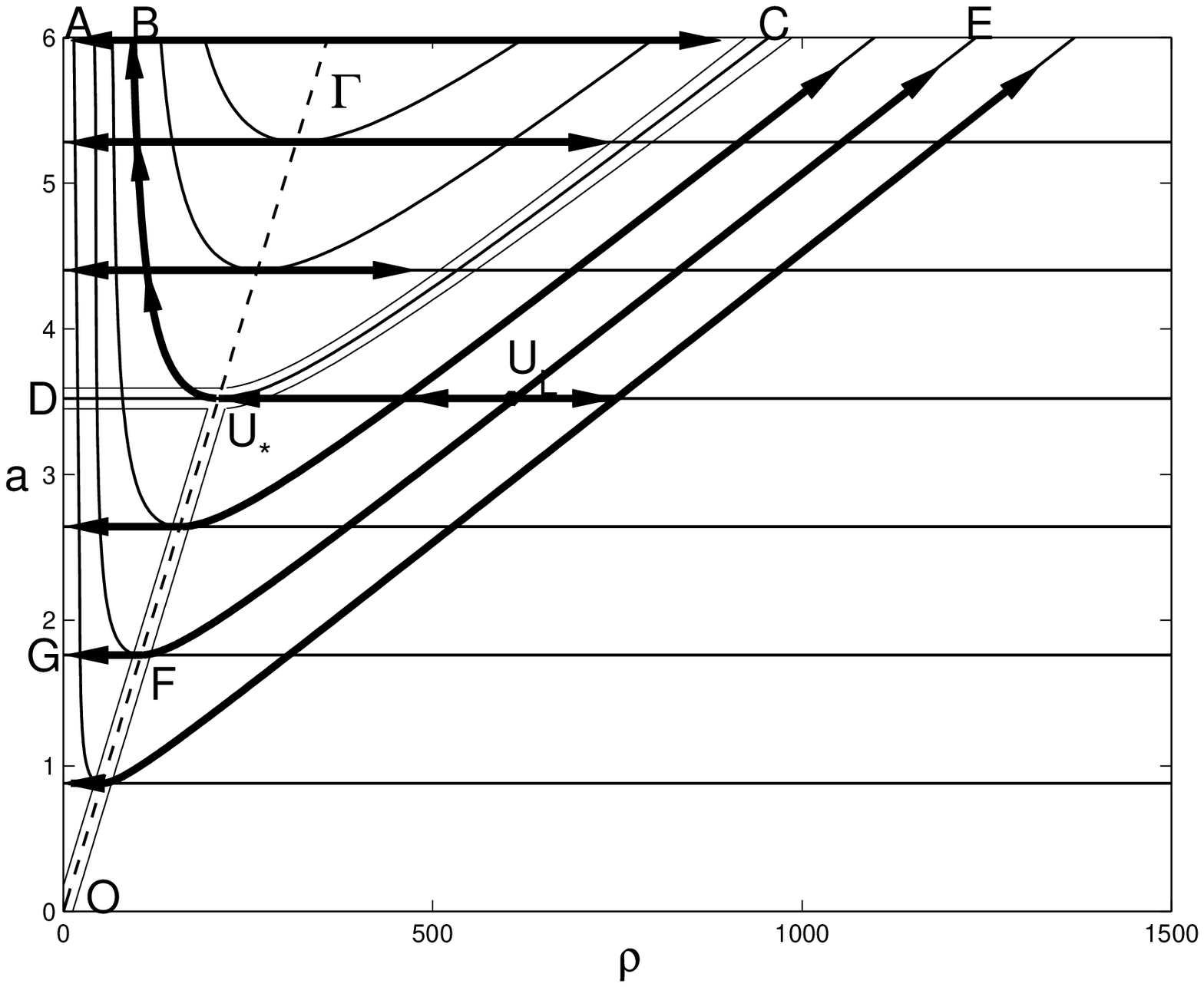}\ec\caption{The Riemann problem for $U_L$ right of $\Gamma$}\label{F_riemann2}\efg

 The wave solutions to the Riemann problem for \refe{system} are combinations of basic 0-waves and 1-waves. If we order the waves with respect to space $x$ at any time $t$, all the waves must satisfy Lax's entropy condition; i.e., the waves from left (upstream) to right (downstream) should increase their wave speeds so that they don't cross each other. This condition is imposed on all hyperbolic systems of conservation laws.  For a nonlinear resonant system \refe{system}, an additional condition has to be imposed; i.e., as long as the standing wave is not interrupted by a shock-wave, its associated density must vary continuously. To guarantee this, the following entropy condition is required:
\bqn
\m{ The standing wave can NOT cross the transition curve } \Gamma.
\eqn

With the two entropy conditions, the solutions to the inhomogeneous LWR model exist and are unique. The wave solutions for undercritical left state $U_L$ is  shown in \reff{F_riemann1}, and those for overcritical left state $U_L$ is shown in \reff{F_riemann2}.

In the remaining part of this section, we discuss wave solutions to the Riemann problem for \refe{system}, present the formula for the boundary flux $f^{\ast}_0$ related to each type of solutions, summarize our results and compare them with those existing in literature.

\subsection{Solutions of the boundary fluxes} \label{bf}
\bfg\bc\includegraphics[width=9cm] {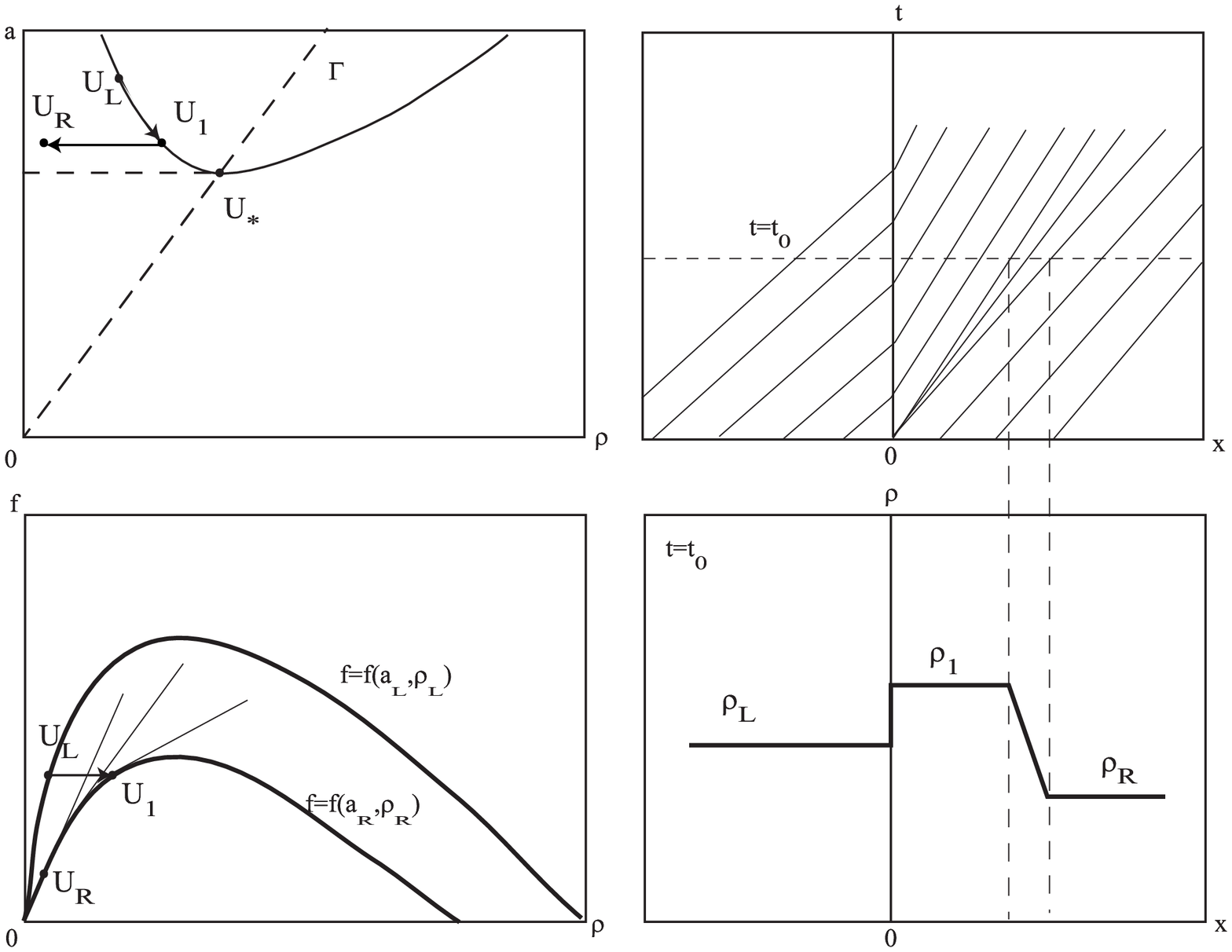}\ec\caption{One example for wave solutions of type 1 for \refe{system} with initial conditions \refe{inh.ini}}\label{F_case1}\efg
\bfg\bc\includegraphics[width=9cm] {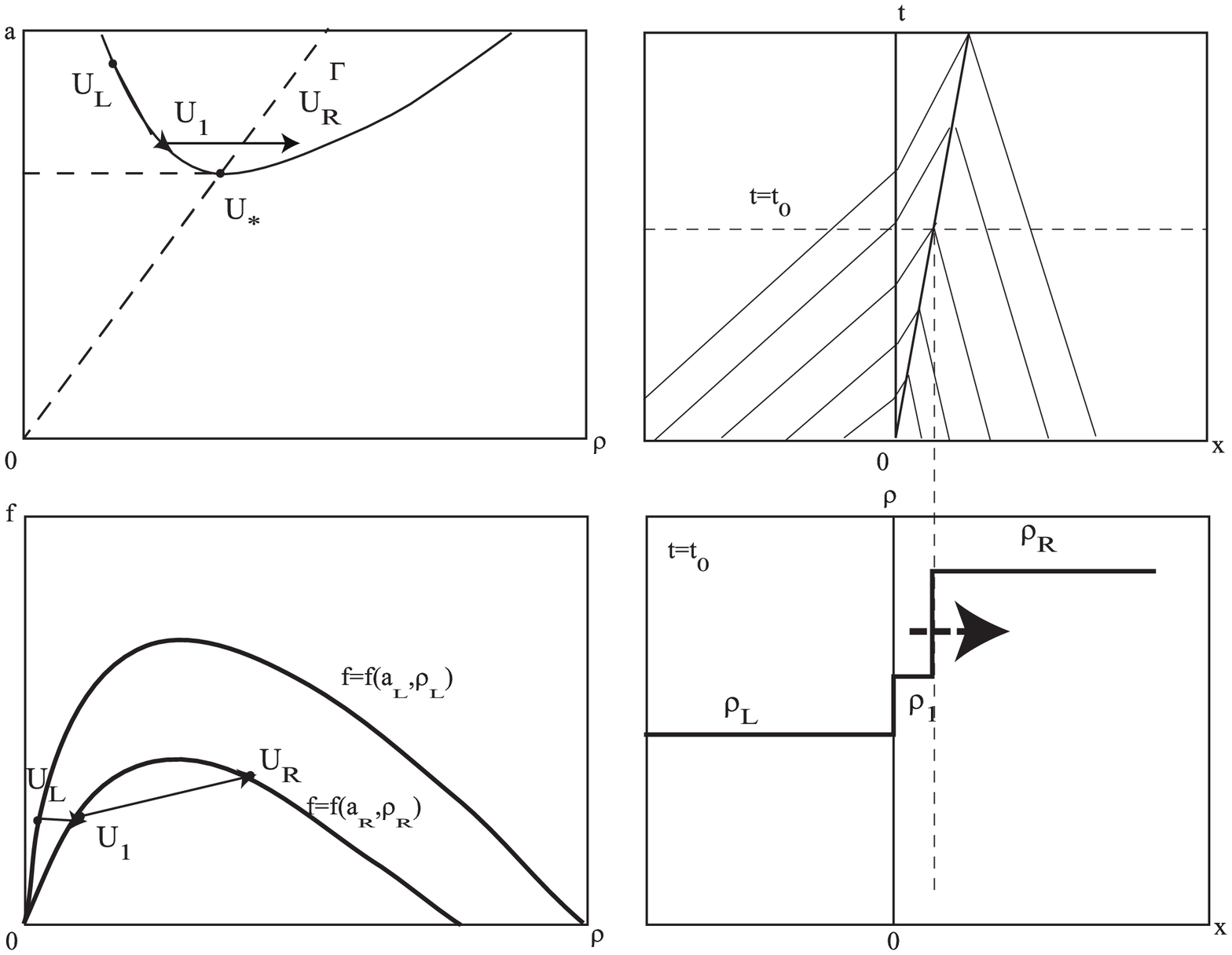}\ec\caption{One example for wave solutions of type 2 for \refe{system} with initial conditions \refe{inh.ini}}\label{F_case2}\efg
\bfg\bc\includegraphics[width=9cm] {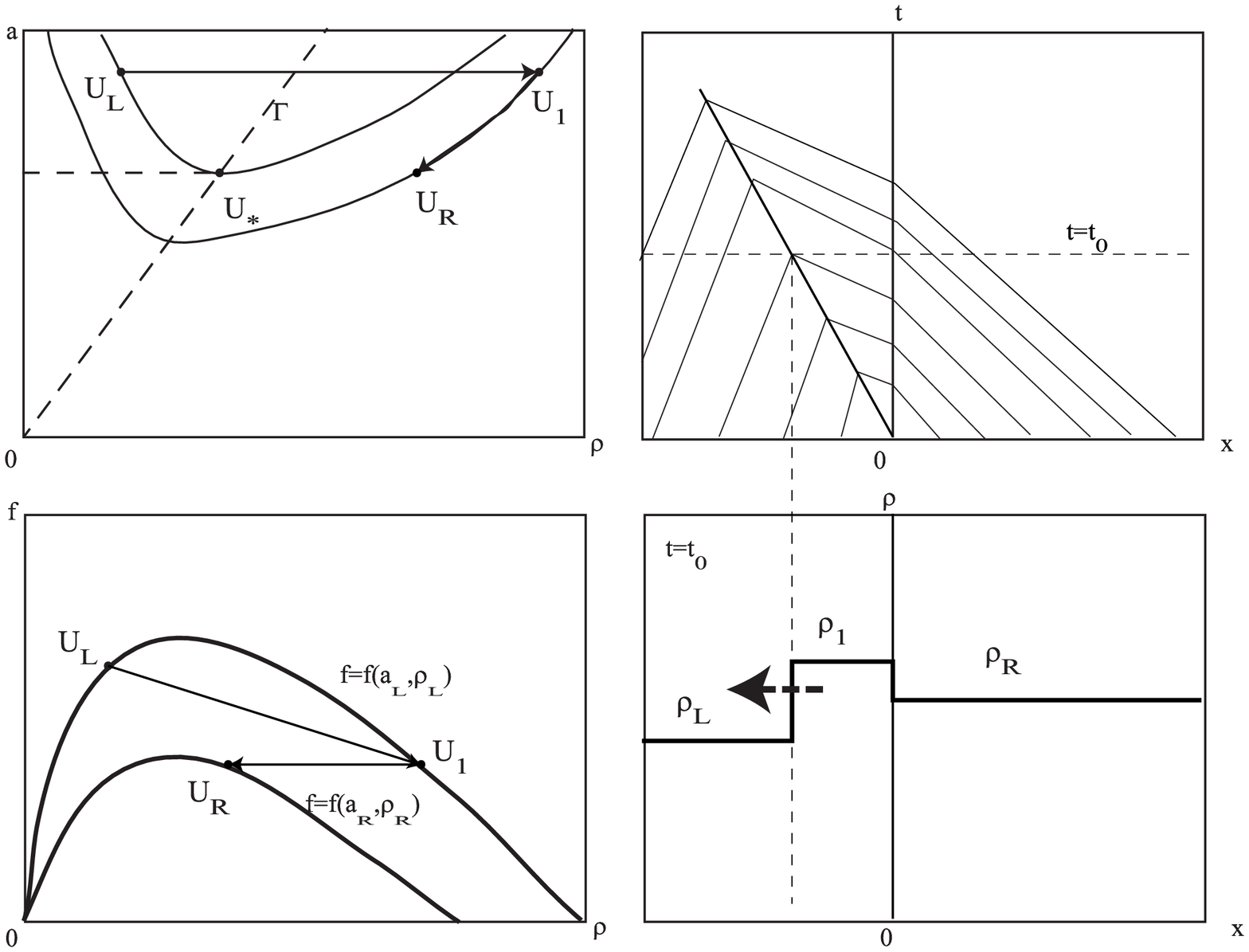}\ec\caption{One example for wave solutions of type 3 for \refe{system} with initial conditions \refe{inh.ini}}\label{F_case3}\efg
\bfg\bc\includegraphics[width=9cm] {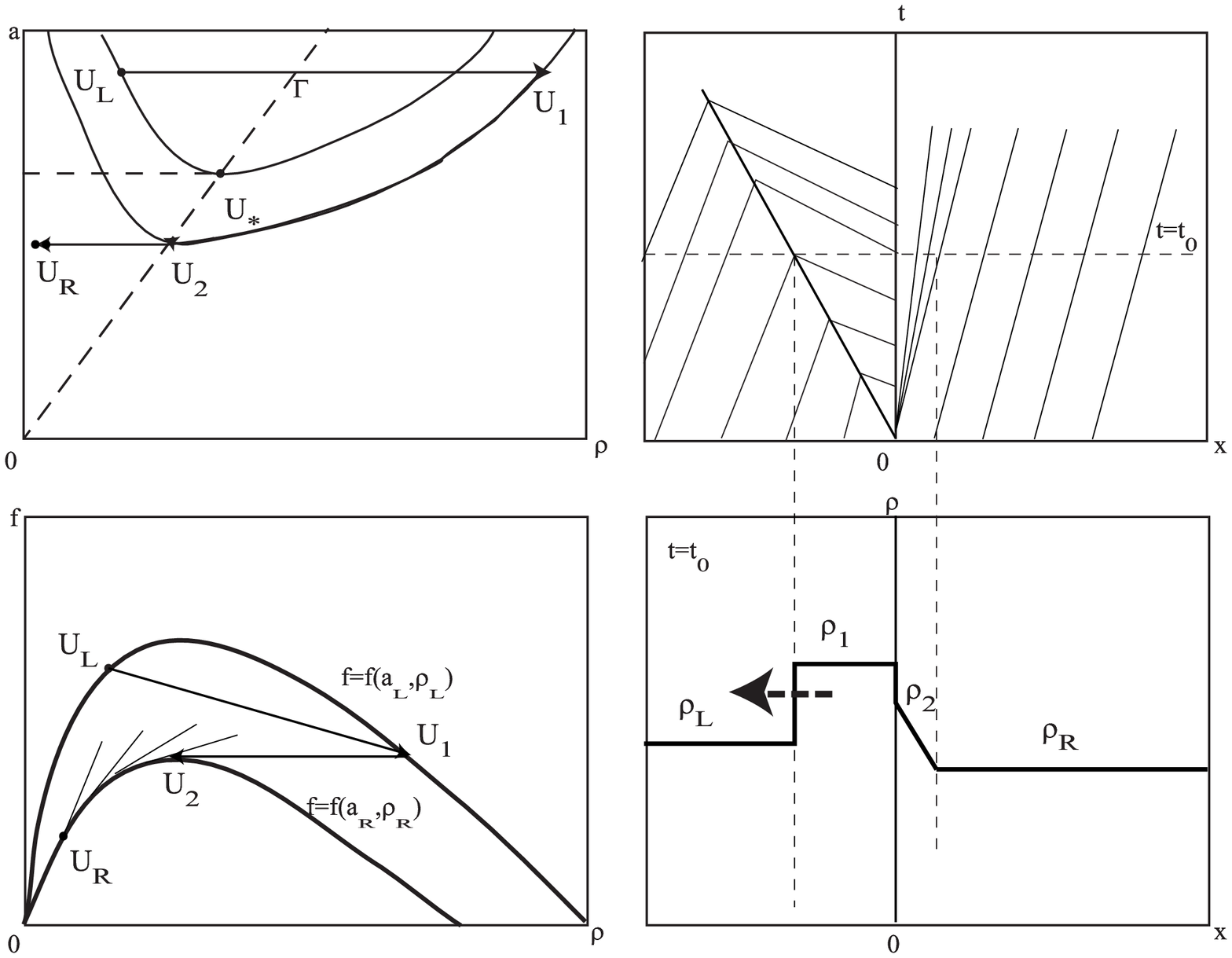}\ec\caption{One example for wave solutions of type 4 for \refe{system} with initial conditions \refe{inh.ini}}\label{F_case4}\efg

When $U_L=(a_L,\r_L)$ is undercritical; i.e., $\r_L/a_L<\alpha$, where $\alpha$ is defined in \refe{trancur}, we denote the special critical point on standing wave passing $U_L$ as $U_{\ast}$. Thus, as shown in \reff{F_riemann1}, the $U$-space is partitioned into three regions by $DU_{\ast}$, $OU_{\ast}$ and $U_{\ast}C$, where $DU_{\ast}=\{(a,\r)|a=a_{\ast},\r<\r_{\ast}\}$, $OU_{\ast}=\Gamma\cap\{0\leq \r\leq \r_{\ast}\}$ and $U_{\ast}C=\{(a,\r)|f(a,\r)=f(U_L),\r>\r_{\ast}\}$.  Related to different positions of the right state $U_R$ in the $U$-space, the Riemann problem for \refe{system} with initial conditions \refe{inh.ini} has the following six types of wave solutions. After discussion for each type of solutions we  provide formula for calculating the associated boundary flux $f^{\ast}_0$.

\bi
\item [Type 1] When $U_R$ is in region $ABU_LU_{\ast}DA$ shown in \reff{F_riemann1}; i.e.,
\bqn
f(U_R)<f(U_{\ast})=f(U_L),\quad \r_R/a_R<\alpha \m{ and } a_R\geq a_{\ast},
\eqn
wave solutions to the Riemann problem are of type 1.  These solutions consist of two basic waves with an intermediate state $U_1=(a_R,\r_1|_{f(a_R,\r_1)=f(U_{\ast})=f(U_L)})$. Of these two waves, the left $(U_L, U_1)$ is a standing wave, and the right $(U_1,U_R)$ is a rarefaction wave with characteristic velocity $\l_1(a,\r)>0$.

From \reff{F_riemann1}, we can see that the Riemann problem may admit this type of solutions when $a_L>a_R$ or $a_L\leq a_R$; i.e., when the road merges or diverges at $x=0$. Here we present an example of this type of solutions in \reff{F_case1}, where the roadway merges at $x=0$; i.e., $a_L>a_R$. For the case  when the roadway diverges at $x=0$, we can find similar solutions.

From \reff{F_case1},  we obtain the boundary flux $f^{\ast}_0=f(U_L)=f(U_{\ast})$ for wave solutions of type 1.

\item [Type 2]  When $U_R$ is in region $BU_LU_{\ast}CB$ shown in \reff{F_riemann1}; i.e.,
\bqn
f(U_R)\geq f(U_{\ast})=f(U_L), \label{type2}
\eqn
wave solutions to the Riemann problem are of type 2.  These solutions consist of two basic waves with an intermediate state $U_1=(a_R,\r_1|_{f(a_R,\r_1)=f(U_{\ast})=f(U_L)})$. Of these two waves, the left $(U_L, U_1)$ is a standing wave, and the right $(U_1,U_R)$ is a shock wave with positive speed $\sigma=\frac{f(U_R)-f(U_{\ast})}{\r_R-\r_1}>0$.

From \reff{F_riemann1}, we can see that the Riemann problem  may admit this type of solutions when the downstream traffic condition $U_R$ is undercritical or overcritical, or the roadway merges or diverges at $x=0$. Here we present an example of this type of solutions in \reff{F_case2}, where the downstream traffic condition overcritical and the roadway merges at $x=0$. For other situations as long as \refe{type2} is satisfied, we can find similar solutions.

From \reff{F_case2},  we obtain the boundary flux $f^{\ast}_0=f(U_L)=f(U_{\ast})$ for wave solutions of type 2. Here we have the same formula as that for wave solutions of type 1.

\item [Type 3] When $U_R$ is in region $OU_{\ast}CO$ shown in \reff{F_riemann1}; i.e.,
\bqn
f(U_R)<f(U_{\ast})=f(U_L),\qquad \r_R/a_R \geq \alpha,
\eqn
wave solutions to the Riemann problem are of type 3. These solutions consist of two basic waves with an intermediate state $U_1=(a_L,\r_1|_{f(a_L,\r_1)=f(U_R)})$.  Of these two waves, the left one $(U_L, U_1)$ is a shock wave with negative speed $\sigma=\frac{f(U_1)-f(U_L)}{\r_1-\r_L}<0$, and the right one $(U_1,U_R)$ is a standing wave.

From \reff{F_riemann1}, we can see that the Riemann problem may admit this type of solutions when the roadway is merges or diverges at $x=0$. Here we present an example of this type of solutions in \reff{F_case3}, where the roadway merges at $x=0$. For the case when the roadway diverges at $x=0$, we can find similar solutions.

From \reff{F_case3},  we obtain the boundary flux $f^{\ast}_0=f(U_R)$ for wave solutions of type 3.

\item [Type 4] When $U_R$ is in region $OU_{\ast}DO$ shown in \reff{F_riemann1}; i.e.,
\bqn
f(U_R)<f(U_{\ast})=f(U_L),\quad \r_R/a_R<\r_{\ast}/a_{\ast} \m{ and }a_R<a_{\ast},
\eqn
wave solutions to the Riemann problem  are of type 4. These solutions consist of three basic waves with two intermediate states: $U_1=(a_L,\r_1|_{f(a_L,\r_1)=f(U_2)})$ and $U_2=(a_R,\r_2|_{\r_2/a_R=\alpha})$. Of these three waves, the left one $(U_L, U_1)$ is a shock wave with negative speed $\sigma=\frac{f(U_1)-f(U_L)}{\r_1-\r_L}<0$, the middle one $(U_1,U_2)$ is a standing wave with zero speed, and the right one $(U_2,U_R)$ is a rarefaction wave with characteristic velocity $\l_1(a,\r)>0$.

From \reff{F_riemann1}, we can see that this type of solutions are admitted only when the roadway merges at $x=0$. Here we present an example of this type of  solutions in \reff{F_case4}.

From \reff{F_case4},  we obtain the boundary flux $f^{\ast}_0=f(U_2)$ for wave solutions of type 4.

\ei

\bfg\bc\includegraphics[width=9cm] {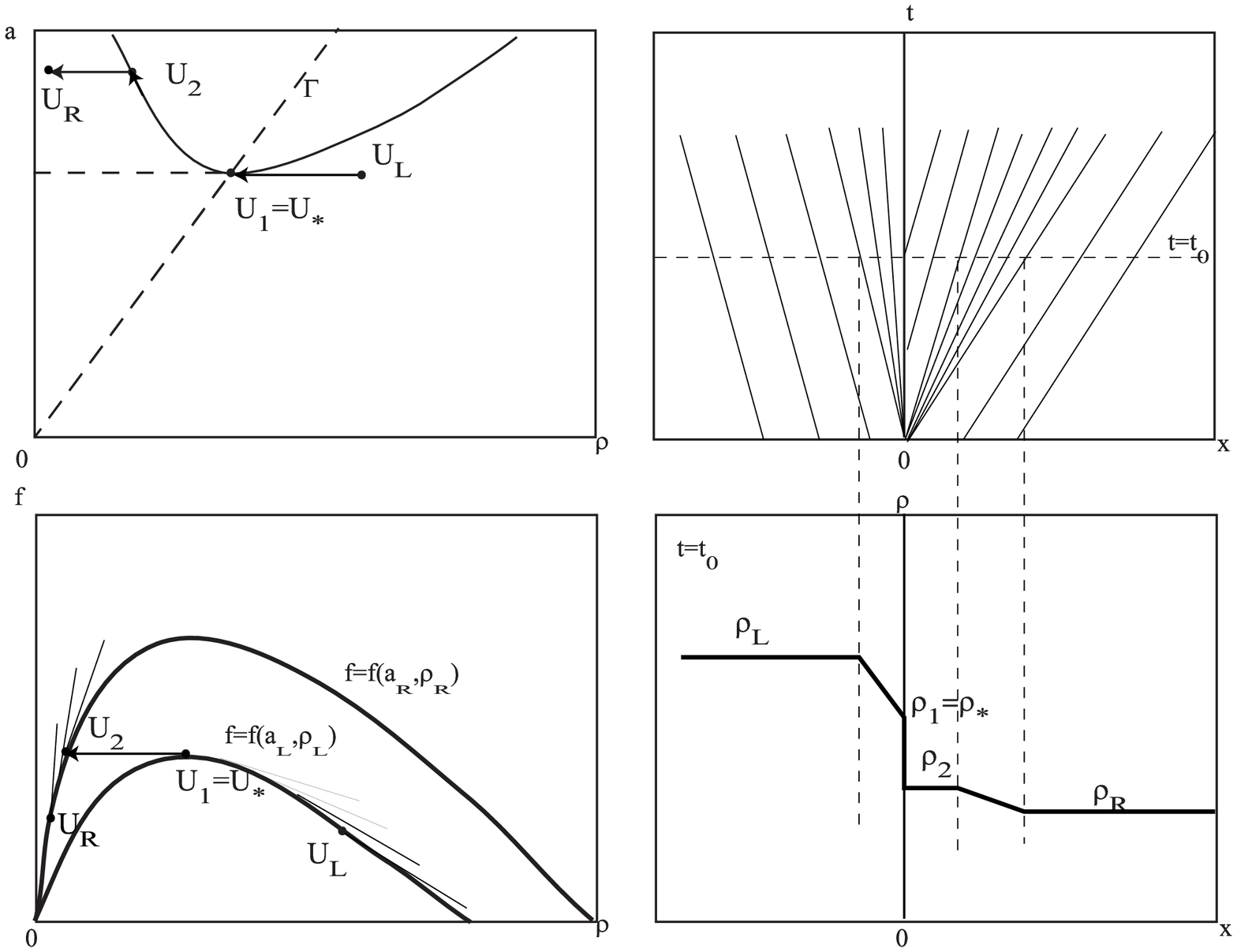}\ec\caption{One example for wave solutions of type 5 for \refe{system} with initial conditions \refe{inh.ini}}\label{F_case5}\efg
\bfg\bc\includegraphics[width=9cm] {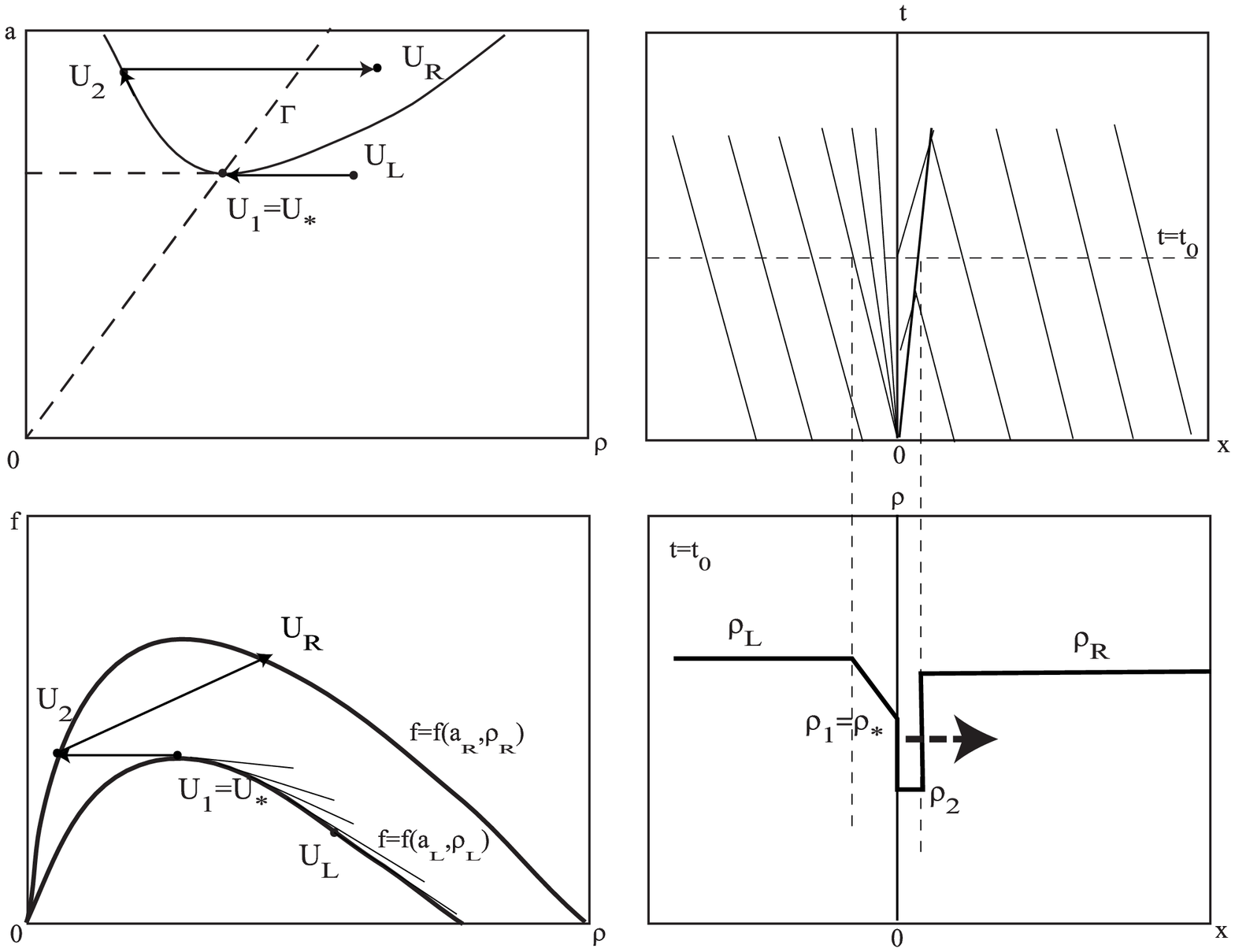}\ec\caption{One example for wave solutions of type 6 for \refe{system} with initial conditions \refe{inh.ini}}\label{F_case6}\efg
\bfg\bc\includegraphics[width=9cm] {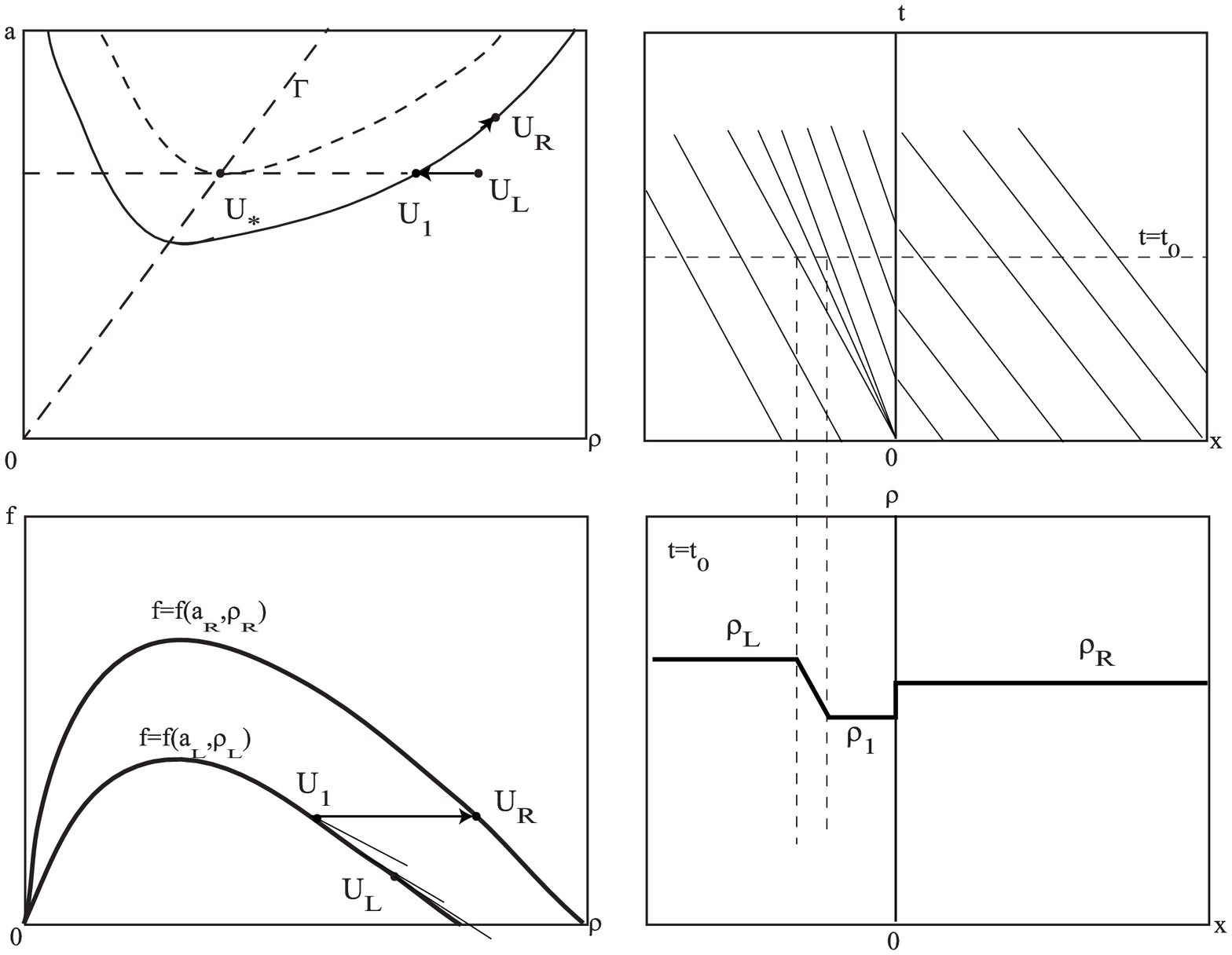}\ec\caption{One example for wave solutions of type 7 for \refe{system} with initial conditions \refe{inh.ini}}\label{F_case7}\efg
\bfg\bc\includegraphics[width=9cm] {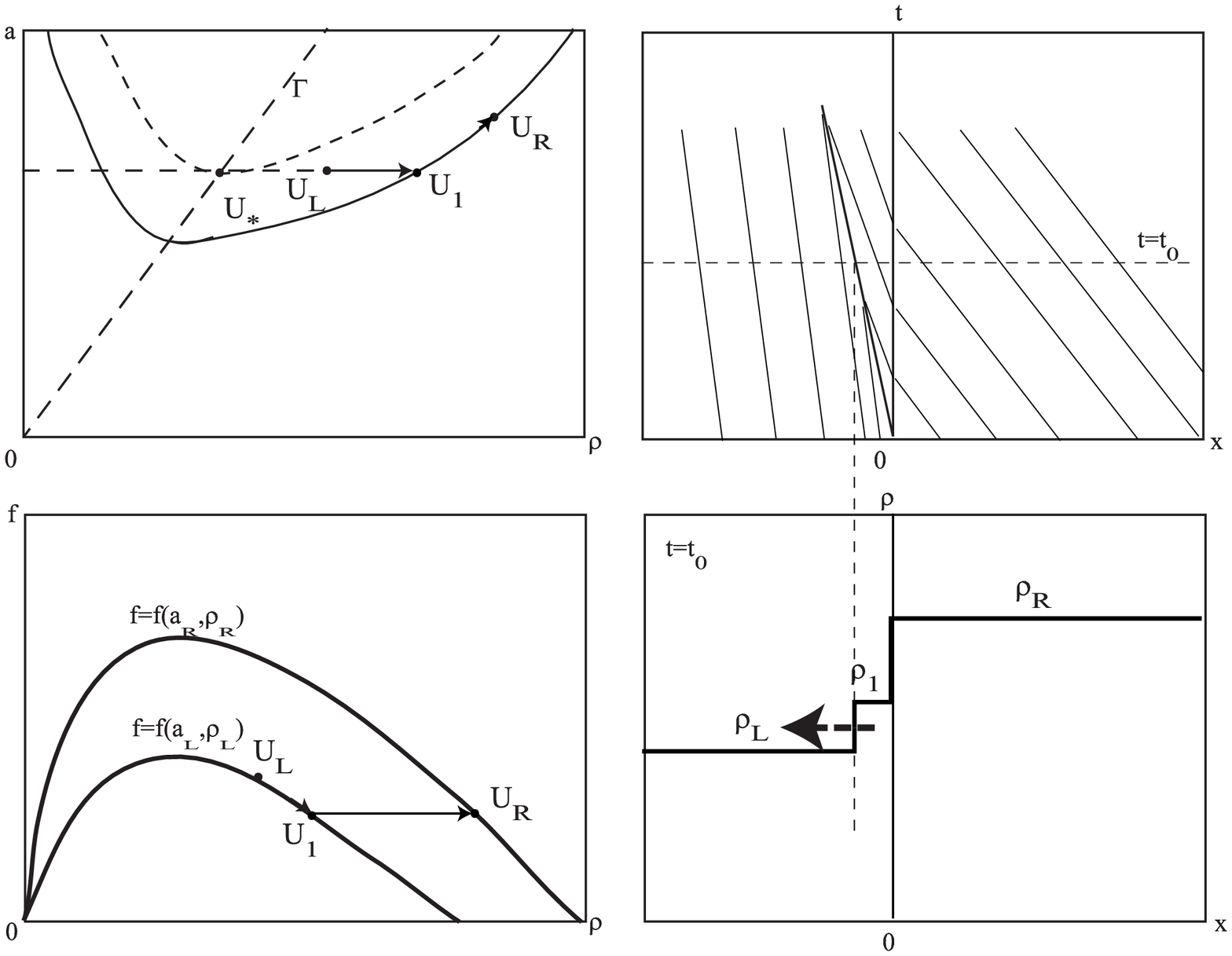}\ec\caption{One example for wave solutions of type 8 for \refe{system} with initial conditions \refe{inh.ini}}\label{F_case8}\efg
\bfg\bc\includegraphics[width=9cm] {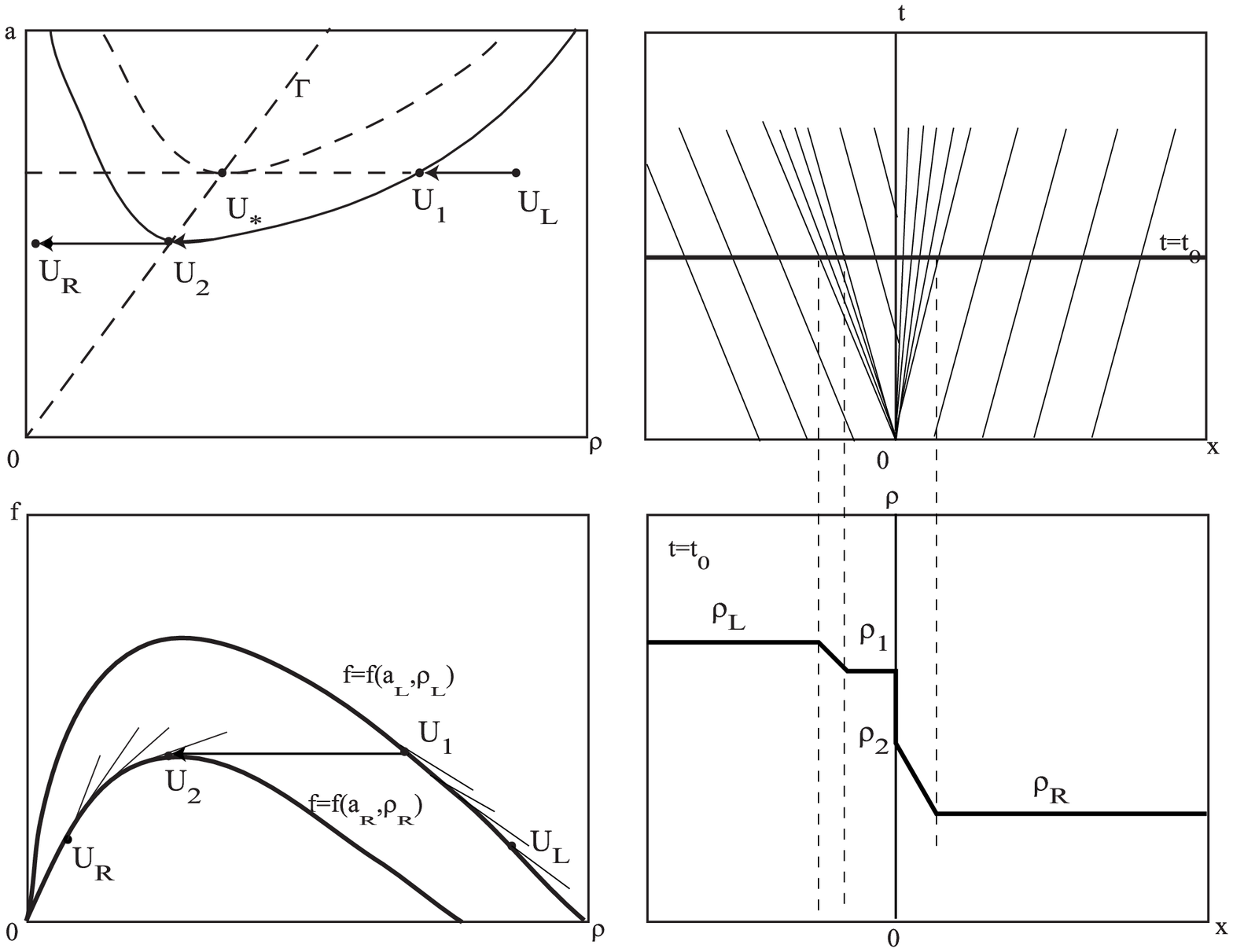}\ec\caption{One example for wave solutions of type 9 for \refe{system} with initial conditions \refe{inh.ini}}\label{F_case9}\efg
\bfg\bc\includegraphics[width=9cm] {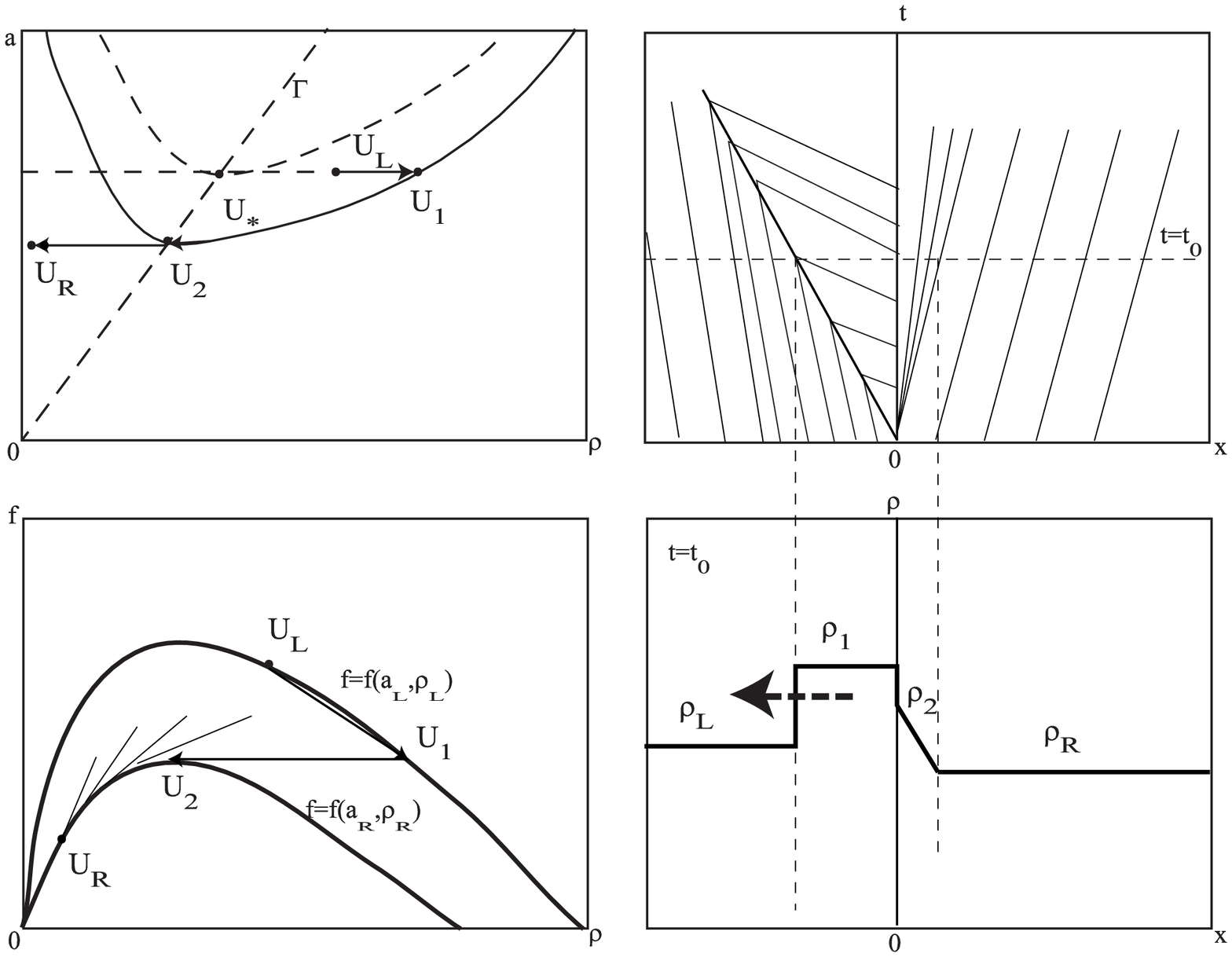}\ec\caption{One example for wave solutions of type 10 for \refe{system} with initial conditions \refe{inh.ini}}\label{F_case10}\efg

When $U_L=(a_L,\r_L)$ is overcritical; i.e., $\r_L/a_L>\alpha$, where $\alpha$ is defined in \refe{trancur}, we denote the special critical point on 1-wave curve passing $U_L$ as $U_{\ast}$; i.e., $U_{\ast}=(a_L,\r_{\ast}|_{\r_{\ast}/a_L=\alpha})$. Thus, as shown in \reff{F_riemann2}, the $U$-space is partitioned into three regions by  three curves $DU_{\ast}=\{a=a_{\ast}=a_L,0\leq \r\leq \r_{\ast}\}$, $OU_{\ast}=\{0\leq a\leq a_{\ast},\r=a\alpha\}$ and $U_{\ast}C=\{a\geq a_{\ast},f(a,\r)=f(U_{\ast})$. Related to different positions of the right state $U_R$ in the $U$-space, the Riemann problem for \refe{system} with initial conditions \refe{inh.ini} has the following six types of wave solutions. After discussion for each type of solutions we  provide formula for calculating the associated boundary flux $f^{\ast}_0$.

\bi
\item [Type 5] When $U_R$ resides in region $ABU_{\ast}DA$ shown in \reff{F_riemann2}; i.e.,
\bqn
f(U_R)<f(U_{\ast}), \quad \r_R/a_R<\alpha \m{ and } a_R \geq a_{\ast}=a_L,
\eqn
wave solutions to the Riemann problem are of type 5. These solutions consist of three basic waves with two intermediate states: $U_1=U_{\ast}$ and $U_2=(a_R,\r_2|_{f(U_2)=f(U_{\ast})})$. Of these three waves, the left one $(U_L,U_1)$ is a rarefaction wave with negative characteristic wave velocity $\l_1(a,\r)$, the middle one $(U_1,U_2)$ is a standing wave and the right one $(U_2,U_R)$ is a rarefaction wave with positive characteristic velocity $\l_1(a,\r)$.

From \reff{F_riemann2}, we can see that this type of solutions are admitted only when the roadway diverges at $x=0$; i.e., $a_R>a_L$. Here we present an example of this type of solutions in \reff{F_case5}.

From \reff{F_case5}, we obtain the boundary flux $f^{\ast}_0=f(U_2)$ for wave solutions of type 5.

\item [Type 6] When $U_R$ resides in region $BU_{\ast}CB$As shown in \reff{F_riemann2}; i.e.,
\bqn
f(U_R)\geq f(U_{\ast}),
\eqn
solutions to the Riemann problem are of type 6. These solutions consist of three basic waves with two intermediate states:  $U_1=U_{\ast}$ and $U_2=(a_R,\r_2|_{f(U_2)=f(U_{\ast})})$. Of these three waves, the left one $(U_L,U_1)$ is a rarefaction wave with negative characteristic velocity $\l_1(a,\r)$, the middle one $(U_1,U_2)$ is a standing wave and the right one $(U_2,U_R)$ is a shock wave with positive speed $\sigma=\frac{f(U_R)-f(U_2)}{\r_R-\r_2}$.

From \reff{F_riemann2}, we can see that this type of solutions may be admitted when the downstream traffic condition is undercritical or overcritical; However, they are admitted only when the roadway diverges at $x=0$. Here we present an example of this type of solutions in \reff{F_case6}, where the downstream traffic condition is overcritical. For the case when the downstream traffic condition is undercritical, we can find similar solutions.

From \reff{F_case6}, we obtain the boundary flux $f^{\ast}_0=f(U_2)$ for this type of wave solutions. Here we have the same formula as that for wave solutions of type 5.

\item [Type 7] When $U_R$ resides in region $CU_{\ast}FU_LEC$ shown in \reff{F_riemann2}; i.e.,
\bqn
f(U_L)\leq f(U_R)<f(U_{\ast}) \m{ and }\r_R/a_R\geq\alpha,
\eqn
wave solutions to the Riemann problem are of type 7. These solutions consist of two basic waves with an intermediate state $U_1=(a_L,\r_1|_{f(U_1)=f(U_R)})$. Of these two waves, the left one $(U_L,U_1)$ is a rarefaction with negative characteristic velocity $\l_1(a,\r)$, and the right one $(U_1,U_R)$ is a standing wave.

From \reff{F_riemann2}, we can see that the Riemann problem may admit this type of solutions when the roadway merges or diverges at $x=0$. Here we present an example of this type of solutions in \reff{F_case7}, where the roadway diverges at $x=0$; i.e., $a_R>a_L$. For the case when the roadway merges, we can find similar solutions.

From \reff{F_case7}, we obtain the boundary flux $f^{\ast}_0=f(U_R)$ for wave solutions of type 7.

\item [Type 8] When $U_R$ locates in region $FU_LEOF$ shown in \reff{F_riemann2}; i.e.,
\bqn
f(U_R)<f(U_L)<f(U_{\ast}) \m{ and } \r_R/a_R \geq \alpha,
\eqn
wave solutions to the Riemann problem are of type 8. These solutions consist of two basic waves with an intermediate state $U_1=(a_L,\r_1|_{f(U_1)=f(U_R)})$. Of these two waves, the left on $(U_L,U_1)$ is a shock with negative speed $\sigma=\frac{f(U_L)-f(U_1)}{\r_L-\r_1}$, and the right one $(U_1,U_R)$ is a standing wave.

From \reff{F_riemann2}, we can see that the Riemann problem may admit this type of solutions when the roadway merges or diverges at $x=0$. Here we present an example of this type of solutions in \reff{F_case8}, where the roadway diverges at $x=0$; i.e., $a_R>a_L$. For the case when the roadway merges, we can find similar solutions.

From \reff{F_case8}, we obtain the boundary flux $f^{\ast}_0=f(U_R)$ for wave solutions of type 8. Here we have the same formula as that for wave solutions of type 7.

\item [Type 9] When $U_R$ resides in region $DU_{\ast}FGD$ shown in \reff{F_riemann2}; i.e.,
\bqn
f(U_L)\leq f(U_R)<f(U_{\ast}),\quad \r_R/a_R<\alpha \m{ and } a_R<a_{\ast}=a_L,
\eqn
wave solutions to the Riemann problem are of type 9. These solutions consist of three basic waves with two intermediate states: $U_1=(a_L,\r_1|_{f(U_1)=f(U_2)})$ and $U_2=(a_R,\r_2|_{\r_2/a_R=\alpha})$. Of these three waves, the left on $(U_L,U_1)$ is a rarefaction with negative characteristic velocity $\l_1(a,\r)$, the middle one $(U_1,U_2)$ is a standing wave, and the right one $(U_2,U_R)$ is a rarefaction with positive speed $\l_1(a,\r)$.

From \reff{F_riemann2}, we can see that this type of solutions are admitted only when the roadway merges at $x=0$; i.e., $a_R<a_L$. Here we present an example of this type of solutions in \reff{F_case9}.

From \reff{F_case9}, we obtain the boundary flux $f^{\ast}_0=f(U_2)$ for wave solutions of type 9.

\item [Type 10] When $U_R$ resides in region $GFOG$ shown in \reff{F_riemann2}; i.e.,
\bqn
f(U_R)<f(U_L)<f(U_{\ast}),\quad \r_R/a_R<\alpha \m{ and }a_R<a_{\ast}=a_L,
\eqn
wave solutions to the Riemann problem are of type 10. These solutions consist of three basic waves with two intermediate states: $U_1=(a_L,\r_1|_{f(U_1)=f(U_2)})$ and $U_2=(a_R,\r_2|_{\r_2/a_R=\alpha})$. Of these three waves, the left one $(U_L,U_1)$ is a shock with negative speed, the middle one $(U_1,U_2)$ is a standing wave, and the right one $(U_2,U_R)$ is a rarefaction wave with positive characteristic velocity $\l_1(a,\r)$.

From \reff{F_riemann2}, we can see that this type of solutions are admitted only when the roadway merges at $x=0$; i.e., $a_R<a_L$. Here we present an example of this type of solutions in \reff{F_case10}.

From \reff{F_case10}, we obtain the boundary flux $f^{\ast}_0=f(U_2)$ for wave solutions of type 10. Here we have the same formula as that for wave solutions of type 9.
\ei

\subsection{Summary}
In the above subsection, we have discussed 10 types of wave solutions. For each type of solutions, the boundary flux $f^{\ast}_0$ is equal to one of the following four quantities: the upstream flow rate $f(U_L)$, the downstream flow rate $f(U_R)$, the capacity of the upstream roadway $f^{max}_L$ and the capacity of the downstream roadway $f^{max}_R$.  For wave solutions of type 1 and 2, the boundary flux is equal to the upstream traffic flow rate; i.e., $f^{\ast}_0=f(U_L)$. For wave solutions of type 3, 7 and 8, the boundary flux is equal to the downstream traffic flow rate; i.e., $f^{\ast}_0=f(U_R)$. For wave solutions of type 4, 9 and 10, the boundary flux is equal to the capacity of the downstream roadway; i.e., $f^{\ast}_0=f^{max}_R$. For wave solutions of type 5 and 6, the boundary flux is equal to the capacity of the upstream roadway; i.e., $f^{\ast}_0=f^{max}_L$. In Table \ref{Table0}, the boundary fluxes are listed for the 10 types of wave solutions to the Riemann problem, as well as the conditions when the Riemann problem admit those solutions.

Note that when $a_L=a_R$; i.e., when \refe{inh_1st} becomes a homogeneous LWR model, wave solutions and the solutions of the boundary fluxes provided here are the same as those for the homogeneous LWR model.

\btb
\begin{tabular}{c||l|l||c}
No.&left state $U_L$ & right state $U_R$ & boundary flux $f^{\ast}_0$ \\\hline
1&undercritical &$f(U_R)<f(U_L)$, $a_R>a_{\ast}$, $\r_R/a_R<\alpha$ & $f(U_L)$ \\\hline
2&undercritical &$f(U_R)>f(U_L)$ & $f(U_L)$ \\\hline
3&undercritical & $f(U_R)<f(U_L)$, $\r_R/a_R>\alpha$ &$f(U_R)$ \\\hline
4&undercritical & $f(U_R)<f(U_L)$, $\r_R/a_R<\alpha$, $a_R<a_{\ast}$ & $f^{max}_R$ \\\hline
5&overcritical & $f(U_R)<f^{max}_L$, $a_R>a_L$, $\r_R/a_R<\alpha$ & $f^{max}_L$ \\\hline
6&overcritical & $f(U_R)>f^{max}_L$ & $f^{max}_L$ \\\hline
7&overcritical & $f(U_L)<f(U_R)<f^{max}_L$, $\r_R/a_R>\alpha$ & $f(U_R)$ \\\hline
8&overcritical & $f(U_R)<f(U_L)$, $\r_R/a_R>\alpha$ & $f(U_R)$ \\\hline
9&overcritical & $f(U_L)<f(U_R)<f^{max}_L$, $\r_R/a_R<\alpha$, $a_R<a_L$ & $f^{max}_R$ \\\hline
10&overcritical & $f(U_R)<f(U_L)$, $\r_R/a_R<\alpha$, $a_R<a_L$ & $f^{max}_R$
\end{tabular}
\caption {Solutions of the Boundary Fluxes}\label{Table0}
\etb

\btb
\begin{tabular}{l|c||l|c}
Conditions& Solutions by Lebacque & Types &Our solutions\\\hline
$a_L\leq a_R$, $U_L$ uc, $U_R$ uc& $f(U_L)$& 1  & $f(U_L)$ \\\hline
$a_L\leq a_R$, $U_L$ uc, $U_R$ oc& $\min \{f(U_L),f(U_R)\}$ & 2 or 3 & $f(U_L)$ or $f(U_R)$\\\hline
$a_L\leq a_R$, $U_L$ oc, $U_R$ uc & $f^{max}_L$ & 5 or 6 & $f^{max}_L$ \\\hline
$a_L\leq a_R$, $U_L$ oc, $U_R$ oc & $\min \{f^{max}_L, f(U_R)\}$ &6, 7 or 8&$f^{max}_L$ or $f(U_R)$\\\hline
$a_L\geq a_R$, $U_L$ uc, $U_R$ uc & $\min \{f^{max}_R, f(U_L)$ & 1 or 4 &$f(U_L)$, $f^{max}_R$ \\\hline
$a_L\geq a_R$, $U_L$ uc, $U_R$ oc & $\min \{f(U_L),f(U_R)\}$ & 2 or 3 & $f(U_L)$ or $f(U_R)$  \\\hline
$a_L\geq a_R$, $U_L$ oc, $U_R$ uc & $f^{max}_R$ & 9 or 10 & $f^{max}_R$\\\hline
$a_L\geq a_R$, $U_L$ oc, $U_R$ oc & $f(U_R)$ & 7 or 8 & $f(U_R)$
\end{tabular}
\caption{Comparison with Lebacque's results}\label{table1}
\etb

The Riemann problem for an inhomogeneous LWR model was also studied by Lebacque (1995\nocite{lebacque1995}). He discussed the Riemann problem for \refe{system} with general inhomogeneity, and presented empirical solutions for it.
He categorized the solutions according to two criteria. The first criterion is to consider the relationship between the upstream capacity and the downstream capacity. For the roadway with variable number of lanes, this criterion is equivalent to considering the relationship between the number of lanes of the upstream and downstream roadway, since the roadway with greater number of lanes has larger capacity. The second criterion is to consider whether the upstream and downstream traffic conditions are undercritical or overcritical. With these criteria, he discussed 8 types of waves solutions to the Riemann problem and obtained the formula for the boundary flux related to each type of solutions. The conditions for those types of wave solutions as well as the formulae related to those types of solutions are listed in Table \ref{table1}, in which oc and uc stand for overcritical and undercritical respectively. Under each of those conditions, the Riemann problem may admit different types of solutions discussed the above subsection \ref{bf}. The types of solutions and our related formulae for the boundary flux are also presented in Table \ref{table1}. From this table, we can see that our results are consistent with those provided by Lebacque, although the Riemann problem is solved through different approaches.

The consistency of our results with existing results can also be shown by introducing a simple formula for the boundary flux.
If we define the upstream demand as
\bqn
f^{\ast}_L&=&\cas{{ll} f(U_L) & \r_L/a_L<\alpha \\f^{max}_L & \r_L/a_L\geq \alpha}
\eqn
and define the downstream supply as
\bqn
f^{\ast}_R&=&\cas{{ll} f^{max}_R & \r_R/a_R<\alpha \\f(U_R) &\r_R/a_R\geq \alpha}
\eqn
then the boundary flux can be simply computed as
\bqn
f^{\ast}_0&=&\min\{f^{\ast}_L, f^{\ast}_R\}. \label{ds}
\eqn
Note that $f^{\ast}_L=f(U_{\ast})$.
Formula \refe{ds} was also provided by Daganzo (1994\nocite{Daganzo94}, 1995\nocite{Daganzo95}) and Lebacque (1995\nocite{lebacque1995}).

\section{Numerical solution method}\label{go}
Since the inhomogeneous LWR model can be written in a conservation form \refe{system}, Godunov's method is efficient for its numerical solutions.
In this section, we describe the Godunov method for solving \refe{system} with general initial and boundary conditions.

In a Godunov's method for \refe{system}, the roadway is partitioned into $N$ zones and a duration of time is discretized into $M$ time steps. In a zone $i$, we approximate the continuous equation \refe{inh_1st} with a finite difference equation
\bqn
\frac {\r_i^{m+1}-\r_i^m} {\dt}+\frac {f^{\ast}_{i-1/2}-f^{\ast}_{i+1/2}}{\dx}&=&0, \label{dis.1}
\eqn
where $\r_i^m$ denotes the average of $\r$ in zone $i$ at time step $m$, similarly $\r_i^{m+1}$ is the average at time step $m+1$; $f^{\ast}_{i+1/2}$ denotes the flux through the upstream boundary of zone $i$, and similarly $f^{\ast}_{i+1/2}$ denotes the downstream boundary flux of zone $i$. In \refe{dis.1}, the boundary flux $f^{\ast}_{i-1/2}$ is related to solutions to a Riemann problem for \refe{system} with the following initial conditions:
\bqn
U(x=x_{i-1/2}, t=t_{m}) &=&\cas{{ll} U^m_{i-1} & x <x_{i-1/2}\\U^m_i & x>x_{i-1/2}}.
\eqn
The wave solutions to the Riemann problem and the associated formula for the boundary flux have been discussed in section \ref{ri}. Then according to \refe{dis.1} we can find $\r$ at time step $m+1$, given traffic conditions at time step $m$ and values of $\r$ at the road boundaries.

\section{Conclusions}
In this chapter, we find that the inhomogeneous LWR model \refe{system} exhibits nonlinear resonance. For this system, we discuss the Riemann problem and find 10 types of wave solutions and the formula for the boundary flux related to each type of wave solutions. Then we obtain a general formula for the boundary flux with the definition of upstream demand and downstream supply. We also conclude that the results obtained here are consistent with those in existing literature. Since the inhomogeneous LWR model can be written in a conservation form \refe{system}, Godunov's method is efficient for its numerical solutions. Godunov's method for the inhomogeneous LWR model is briefly described in section \ref{go}.

In this chapter, we consider the inhomogeneity as variable number of lanes. Similarly we can solve the inhomogeneous LWR model with other inhomogeneities, by changing the functions $v(a,x)$ and $f(a,x)$ corresponding to the new inhomogeneities. For all kinds of inhomogeneity, equation \refe{ds} is always the general formula for computing the boundary flux.

In this chapter, we consider an inhomogeneous LWR model, which is a first-order traffic flow model. However, we want to point out that the method provided here  can be extended to solve  inhomogeneous higher-order traffic flow models such as Zhang's model discussed in Chapter 3 and the PW model discussed in Chapter 4.

\newpage
\pagestyle{myheadings} 
\markright{  \rm \normalsize CHAPTER 6. \hspace{0.5cm}
 A First-Order Multi-Commodity Model and Its Numerical Simulations}
\large
\chapter{A First-Order Multi-Commodity Model and Its Numerical Simulations} 
\section{Introduction}
All the models discussed in previous chapters are for link flows, in which all vehicles have the same contribution to the flow dynamics.  When one considers a network flow, he/she has to deal with vehicles of different origins, destinations, or other attributes. These attributes of vehicles play a role in determining their choice of routes, and therefore affect the flow dynamics. Hence, to model a network flow, one has to take these attributes into account and disaggregate traffic flow into different components. Such a network flow model considering disaggregated traffic flow is called a multi-commodity model.

In literature, there have been several models for multi-commodity flow, including the model by Vaughan, Hurdle and Hauer (\nocite{vaughan84}1984), Jayakrishnan's model (\nocite{jayakrishnan91}1991) and Daganzo's cell transmission model (\nocite{Daganzo94}1994; \nocite{Daganzo95}1995).

Vaughan, Hurdle and Hauer (\nocite{vaughan84}1984) suggested a continuous model consisting of a `local equation' and a `history equation'. In their model, the basic element of traffic flow is a vehicle, and for each vehicle a label for its trajectory is introduced. In this model, the local equation ensures traffic conservation, and the history equation computes the trajectory of each vehicle. It has been shown that this model is consistent with the traditional LWR model at the aggregate level.  

Jayakrishnan (\nocite{jayakrishnan91}1991) introduced another multi-commodity model. This is a discrete model, in which each link is partitioned into a number of zones. In each zone, vehicles close to each other and with the same origin, destination or other common commodity-specific characteristics are considered as a ``macroparticle". To determine the position of a macroparticle at next time step,  Jayakrishnan considered its travel speed and the length of the zone where it stays. This model doesn't always preserve traffic conservation; i.e., this model may not be consistent with the LWR model, since the macroparticles are considered separately. 

Daganzo (\nocite{Daganzo94}1994; \nocite{Daganzo95}1995) introduced a multi-commodity model based on his {\em cell transmission model}. In this discrete multi-commodity model, traffic flow in a zone is also disaggregated into macroparticles. In every zone, the macroparticles are ordered according to the waiting time. In this order, the macroparticles are moved, and those macroparticles can proceed to a downstream zone if their waiting times are greater than a threshold minimum waiting time. In Daganzo's model, the minimum waiting time is set to the waiting time of a macroparticle if the total number of vehicles which are in the same zone and in front of this macroparticle is equal to the number of vehicles moving into the downstream zones. In this model, the number of vehicles moving into the downstream zones are computed as the boundary flux times the length of a time step based on solutions to the Riemann problem for the LWR model, which has been discussed in Chapter 5.

Of these three models, Jayakrishnan's and Daganzo's are simpler in computation since a macroparticle instead of a vehicle is considered as the basic element of traffic flow. 

Both the model by Vaughan et al. and Daganzo's model consider the dynamics of a network flow at two levels: aggregate level and disaggregated level. At aggregate level, vehicles in a network flow is considered to be identical particles. At this level, both models are consistent with the LWR model, and therefore they satisfy traffic conservation. In the model by Vaughan et al., the continuous LWR model is used, and while Daganzo's model uses the discrete form of the LWR model.

In these three models the vehicles are ordered. In the model by Vaughan et al. vehicles are ordered according to the trajectory labels, Jayakrishnan ordered macroparticles according to location and Daganzo ordered the macroparticles according to time. They all assume that the vehicles observe the First-In-First-Out (FIFO) discipline when they travel on the roadway. In particular, the model by Vaughan et al. assumes the vehicle trajectories don't cross each other at any time $t$ and any location $x$, Jayakrishnan's model assumes that vehicles always keep the order in location, and Daganzo's model assumes that vehicles always keep the order in time. Of all these models, the one by Vaughan et al. uses the full information of order in time and location, while the other two use part of the order information of either time or location. 

In this chapter we introduce a new multi-commodity model. Our model is based on the concept of ``macroparticle", and has also a two-level structure. In our model, we use the discrete form of the LWR model for traffic flow at the aggregate level. In our new model, we will interpret the FIFO discipline in a new way so that we can simplify the dynamics at disaggregated level.

This chapter is organized as follows. In section 2, we discuss the new multi-commodity model and define the network structure, data structure and program structure. In section 3, we present some numerical simulations. In section 4, we discuss the future plans on this topic. 

\section{A two-level multi-commodity model}
The dynamics of a network flow is considered at two levels: aggregate level and disaggregated level.

For traffic flow at the aggregate level, we use the discrete form of an inhomogeneous LWR model. As discussed in Chapter 5, the inhomogeneous LWR model can be written as
\bqn
\r_t+f(a,\r)_x&=&0,\label{lwrm}
\eqn
in which $a(x)$ is an inhomogeneity factor. By partitioning each link into $N$ zones, and discretizing the time interval into $M$ time steps, we obtain the Godunov-type finite difference equation for \refe{lwrm}:
\bqn
\frac{\r_i^{j+1}-\r_i^j}{\dt}+\frac{f^{j\ast}_{i-1/2}-f^{j\ast}_{i+1/2}}{\dx}&=&0, \label{fi}
\eqn
where $\dx$ is the length of zone $i$, $\dt$ is the time from time step $j$ to time step $j+1$, and the choice of $\frac {\dt}{\dx}$ is governed by the CFL condition. In equation \refe{fi}, $\r_i^j$ is the average of traffic density $\r$ in zone $i$ at time step $j$, similarly $\r_i^{j+1}$ is the average of $\r$ at time step $j+1$; $f^{j\ast}_{i-1/2}$ is the flux through the upstream boundary of zone $i$ from time step $j$ to time step $j+1$, and similarly $f^{j\ast}_{i+1/2}$ is the downstream boundary flux. Given traffic conditions at time step $j$, we can  calculate the traffic density in zone $i$ at time step $j+1$ as 
\bqn
\r_i^{j+1}&=&\r_i^j+\dxt(f^{j\ast}_{i-1/2}-f^{j\ast}_{i+1/2}), \label{fi2}
\eqn
where computation of the boundary fluxes $f^{j\ast}_{i-1/2}$ and $f^{j\ast}_{i+1/2}$ has been discussed in Chapter 5. 

Defining $N_i^j=\r_i^j \dx$ as the number of vehicles in zone $i$ at time step $j$, $N_i^{j+1}=\r_i^{j+1}\dx$ as the number of vehicles at time step $j+1$, $F_{i-1/2}^j=\dt f(\r^{j\ast}_{i-1/2})$ as the number of vehicles flowing into zone $i$ from time step $j$ to $j+1$, and $F_{i+1/2}^j$ as the number of vehicles flowing out of zone $i$, equation \refe{fi2} can be written as:
\bqn
N_i^{j+1}&=&N_i^j+F_{i-1/2}^j-F_{i+1/2}^j, \label{diffe}
\eqn
which is a conservation form. 

For all the link flows in a network, equation \refe{diffe} is an efficient model at aggregate level, which captures the dynamics and preserves the conservation of traffic flow. However, how to model the dynamics for the flows at merges and diverges is still an open question. 

There have been two empirical treatments of flows at merges by Daganzo and Lebacque. For a merge with $K$ upstream zones, Daganzo (1995a\nocite{Daganzo95}) solve the flux through the merge using traffic supply of the downstream zone, and the sum of traffic demands of these $K$ upstream zones as traffic demand. This flux is distributed to all the upstream zones according to some pre-defined distribution fractions. This treatment has an underlying assumption that the upstream flows have the same contributions to the merge. This assumption limits the application of this treatment for a merge where the upstream zones consist of mainline highway and on-ramps, since  on-ramps have lower priority than mainline highway. Another approach suggested by Lebacque (1995\nocite{lebacque1995}) first splits traffic supply of the downstream zone into $K$ parts according to some pre-defined fractions, and forms a pair of demand and supply for each of the $K$ upstream zones. Then for each pair of demand and supply, a flux through the merge can be solved. The sum of these $K$ fluxes, considered by Lebacque, is then the total flux through the merge. In our model, we will use Lebacque's method to treat the merges. 

For traffic flow at diverges, the approaches suggested by Daganzo and Lebacque can still be used. In Daganzo's method, the traffic supply for a diverge with $K$ downstream zones is the sum of $K$ traffic supplies. Then the flux through the diverge is computed with this traffic supply and the traffic demand of the upstream zone and distributed to the downstream zones according to pre-defined fractions. In Lebacque's method, for each of $K$ downstream zones, the traffic demand is a fraction of the upstream traffic demand. With this traffic demand and its traffic supply, the flux is computed. The sum of the $K$ fluxes is the total flux through the diverge. In our model, we use a treatment which is adjusted from Lebacque's. Here, we split the upstream demand based on the route choices instead of some pre-defined fractions; i.e., the demand for each downstream zone is equal to the number of vehicles choosing that zone. This treatment may better model the route choices at diverges. 

At disaggregated level, we partition vehicles in a zone into macroparticles. In a macroparticle, vehicles are close to each other and have the same disaggregated information such as destination. The macroparticles in a zone are ordered in location as a queue, with those at downstream as the head of the queue and those at upstream as the tail of the queue. Since the CFL condition guarantees that no macroparticle can cross a zone in one time step, our question is whether a macroparticle can be moved into a downstream zone and what position it is in that zone.

To determine the movement of macroparticles, we introduce a new concept -- `` boundary connector". A boundary connector stands for a boundary between linked zones in a link of a traffic network. It is not a physical entity, doesn't have length and cannot store vehicles. With this concept, zones are not considered to be linked to each other directly, but are linked to boundary connectors. Therefore, in a traffic network, zones reflect traffic conditions, and boundary connectors determine the network structure. 

For each boundary connector, we determine the movement of macroparticles at every time step as follows. From each upstream zone connected to the boundary connector, we pick out macroparticles starting from the head of the queue of macroparticles and store them in the boundary connector temporarily, until the total number of vehicles moved is equal to the flux from that zone, which is determined by the computation at aggregate level. When there are more than one zones connected to the boundary connector, the macroparticles moved into the connector are ordered according to the times they enter. However, in reality, we do not know the exact order of their arrival times and a random process may be used to order them. Then every macroparticle is moved into the downstream zone leading to its destination. If a macroparticle has more than one downstream zones to choose, corresponding fractions may be used. After a macroparticle is moved into a zone, it is attached to the tail of the queue in that zone. In our model, the FIFO principle is also applied. However, this principle is interpreted in a different way than in Vaughan et al., Daganzo or Jayakrishnan. Here, FIFO means those macroparticles in front (location) in a zone are ahead (time) in a boundary connector, and similarly those ahead in a boundary connector are in front in a zone. 

The model suggested here has been implemented for a simple network. In the remaining part of this chapter, we design the network, data and program structures for it and carry out some numerical simulations. 

\section{Network, data and program structures for a specific traffic network}
In this section a simple traffic network, shown in \reff{F_network}, is considered. This network consists of a mainline highway, one on-ramp and one off-ramp. Therefore, there are two origins and two destinations in this network. For modeling purpose, the mainline highway is partitioned into 20 zones, labeled from $1$ to $20$, and 21 boundary connectors are used to connect these 20 zones, highway origin, highway destination and two ramps. Of these 21 boundary connectors, boundary connector $1$ is the upstream highway boundary and boundary connector $21$ is the downstream highway boundary. Besides 20 highway zones, we use artificial zone $0$ to denote highway origin 1, zone $21$ to denote highway destination 1, zone $22$ to denote on-ramp, and zone $23$ to denote off-ramp.

\bfg
\bc\includegraphics{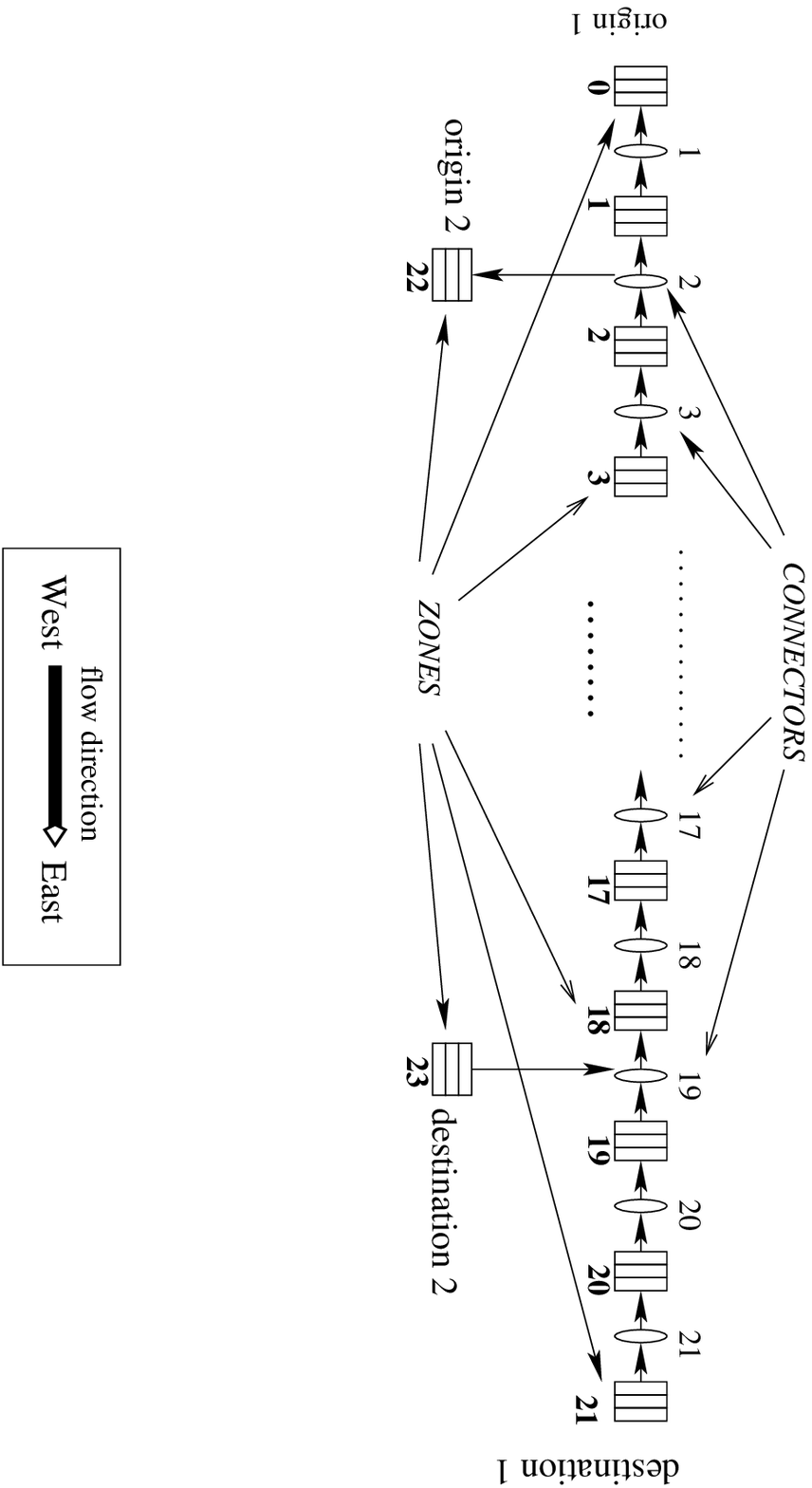}\ec
\caption{An example of a highway network}\label{F_network}
\efg

\bfg
\bc\includegraphics[height=16cm]{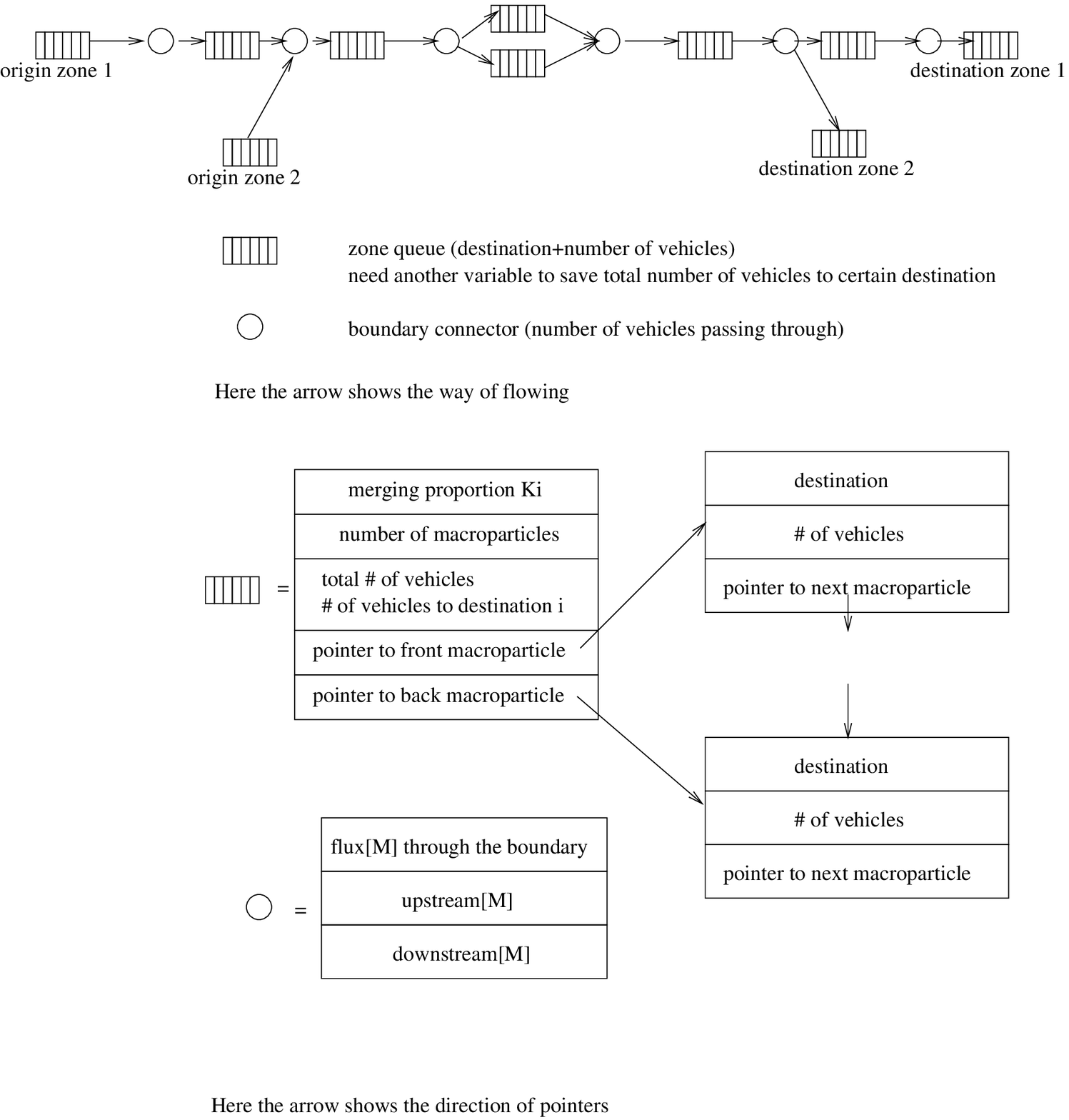}\ec
\caption{Data structure of traffic network}\label{F_data}
\efg

The data used to reflect network traffic flow are shown in \reff{F_data}. At aggregate level the data include:
\ben
\item Traffic measurements of zone $i$ at time step $j$ -- $\rho_i^j$, $v_i^j$. The number of vehicles $N_i^j=\rho_i^j \Delta x_i$, where $\Delta x_i$ is the length of zone $i$. These measurements are given initially and are updated at each time step. The traffic conditions of the origin zones 1 and 22, the destination zones 21 and 23 may be given in certain type of boundary conditions.
\item Fluxes through boundary $i+1/2$ between time step $j$ and time step $j+1$. We solve the Riemann problem to get $\rho^{j\ast}_{i+1/2}$, $v^{j\ast}_{i+1/2}$. The number of vehicles through the boundary $F^j_{i+1/2}=\rho^{j\ast}_{i+1/2}v^{j\ast}_{i+1/2}\Delta t^j$, where $\Delta t^j$ is length of the time interval from $j$ to $j+1$. The fluxes at boundary $1$ and boundary $21$ may be given with certain type of boundary conditions.
\een

At disaggregated level the data include
\ben
\item Array of queues.  Each zone is a dynamically allocated queue. Each node of a queue stands for a macroparticle. The parameters for a macroparticle include  its destination, the number of vehicles in it and a pointer to its upstream macroparticle. The parameters of a zone include  the number of macroparticles in this zone, the number of vehicles to each destination, the total number of vehicles in the zone and some other related information of the zone.
\item Array of boundary connectors. Each boundary connector stores the upstream  and downstream zones that are connected to it. For different types of boundary connectors, e.g., merge or diverge, the treatments are different, therefore the type of a boundary connector is a parameter. Another parameter is the number of vehicles passing through a boundary connector at a time step.  Under some conditions, these fluxes through some specific boundaries may be set, for example, at boundary connector 2, the flux from zone 22 may be set to a value if this on-ramp is metered at a rate.
\een

\bfg\bc\includegraphics[height=8cm]{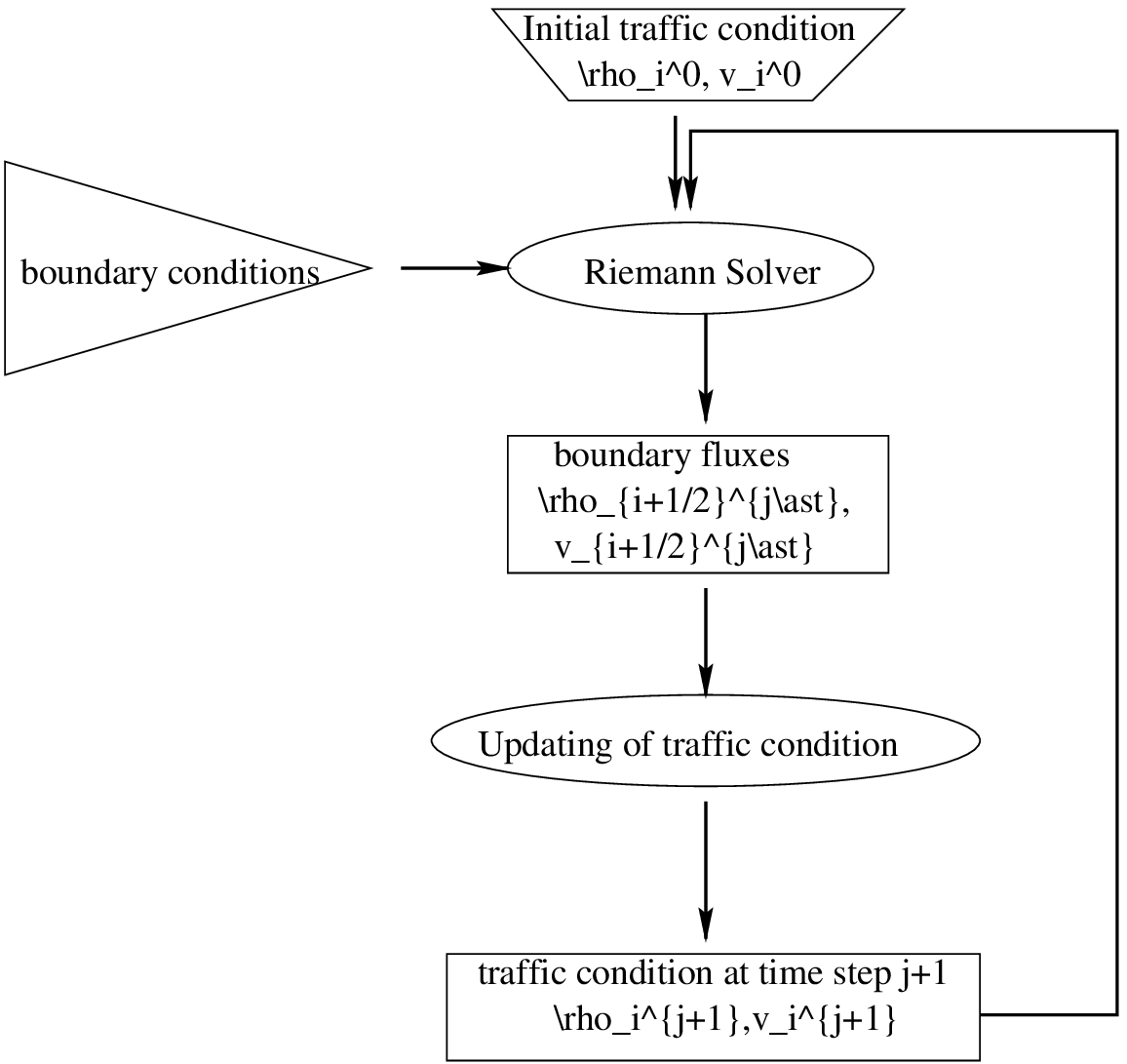}\ec\caption{The program flow chart for aggregate flows} \label{firstlevel}\efg

The program flow chart for computing aggregate flows is shown in \reff{firstlevel}. This program consists of the following operators:
\ben
\item The Riemann solver. The Riemann solver is used to calculate the fluxes through the boundaries. The solver can be of first- or second-order.
\item Updating of traffic conditions. This is done according to the finite difference equations. For the LWR model, traffic conditions are updated according to \refe{diffe}.
\een

\bfg\bc\includegraphics[height=8cm] {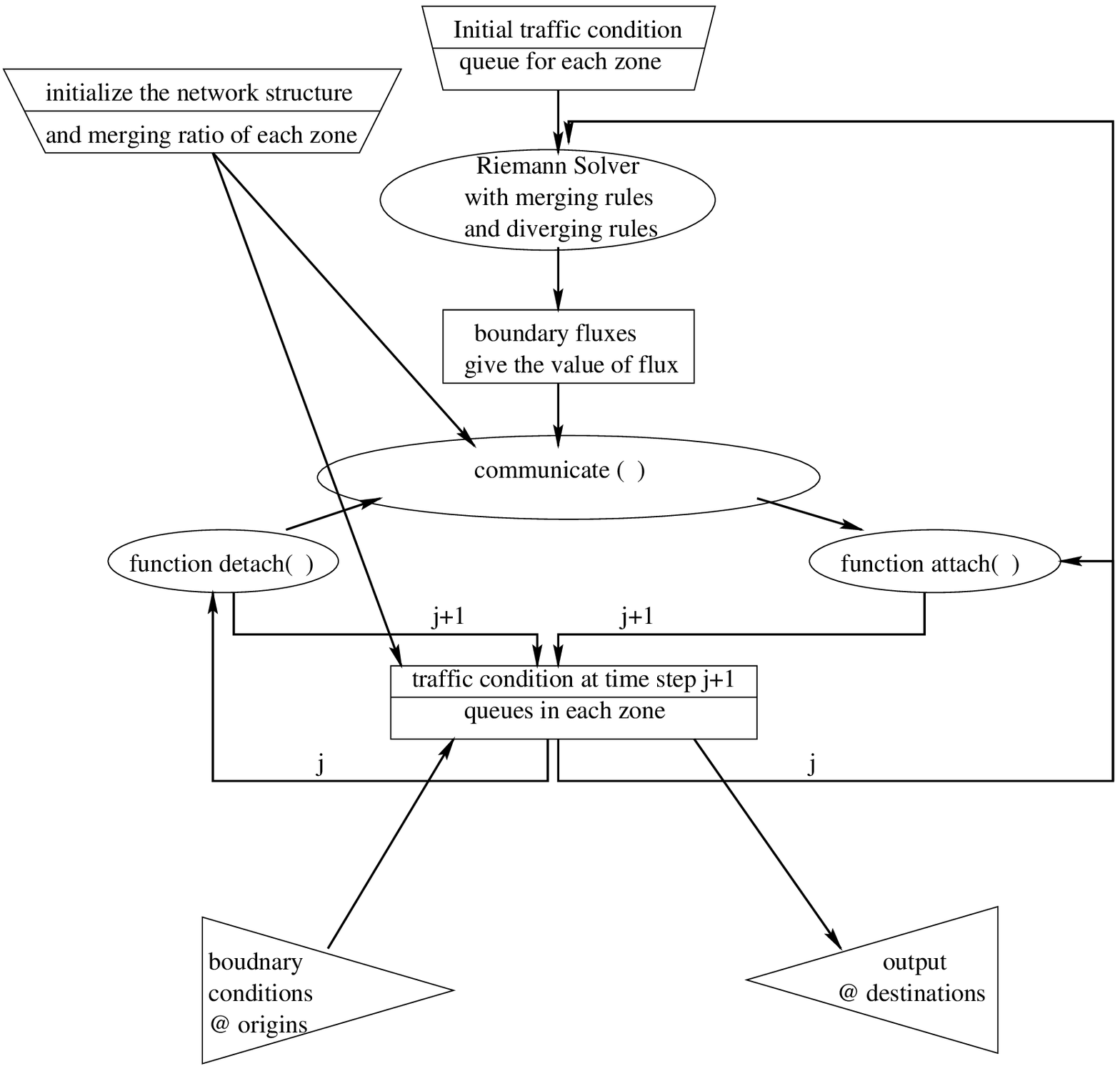}\ec\caption{The program flow chart for two-level flows}\label{secondlevel}\efg

When disaggregated traffic flow is considered, the program flow chart is shown in \reff{secondlevel}. This program consists of the following operators:
\ben
\item Initialization of traffic condition, which provides the initial traffic measurements to all the zones and macroparticles. One simple way to start is to assume that the network has no traffic in the beginning.
\item Initialization of network structure, which provides values for the parameters of each boundary connector.
\item The Riemann solver. The solver is used to compute the aggregate traffic flux. It gives the number of vehicles through the boundary connector during a time interval.
\item Updating of traffic conditions around a boundary connector. At each time step, the traffic conditions of those zones connected to a boundary connector are updated with the flux through the boundary. We retrieve macroparticles from upstream zones, and distribute them to downstream zones according to the destination information carried by each macroparticle. This process stops when the number of vehicles passing the boundary is equal to the flux calculated by the Riemann solver. The retrieving from upstream zones is done by the function \texttt{detach}; and the distributing to downstream zones is done by the function \texttt{attach}. The function \texttt{communicate} governs these two operators according to the O-D information carried by macroparticles. This function uses certain merging rules and route choice rules which are used at aggregate level. Some ad-hoc treatment on these rules was discussed in previous section.
\item Providing demands at origin zones. The demand at each origin zone is given as a queue of well-ordered macroparticles, although in reality this information is usually unknown.
\item Computing output at destination zones. The traffic conditions of each destination zone are the outputs of interest.
\een

\section{Numerical Simulations}
\bfg\bc\includegraphics[height=10cm]{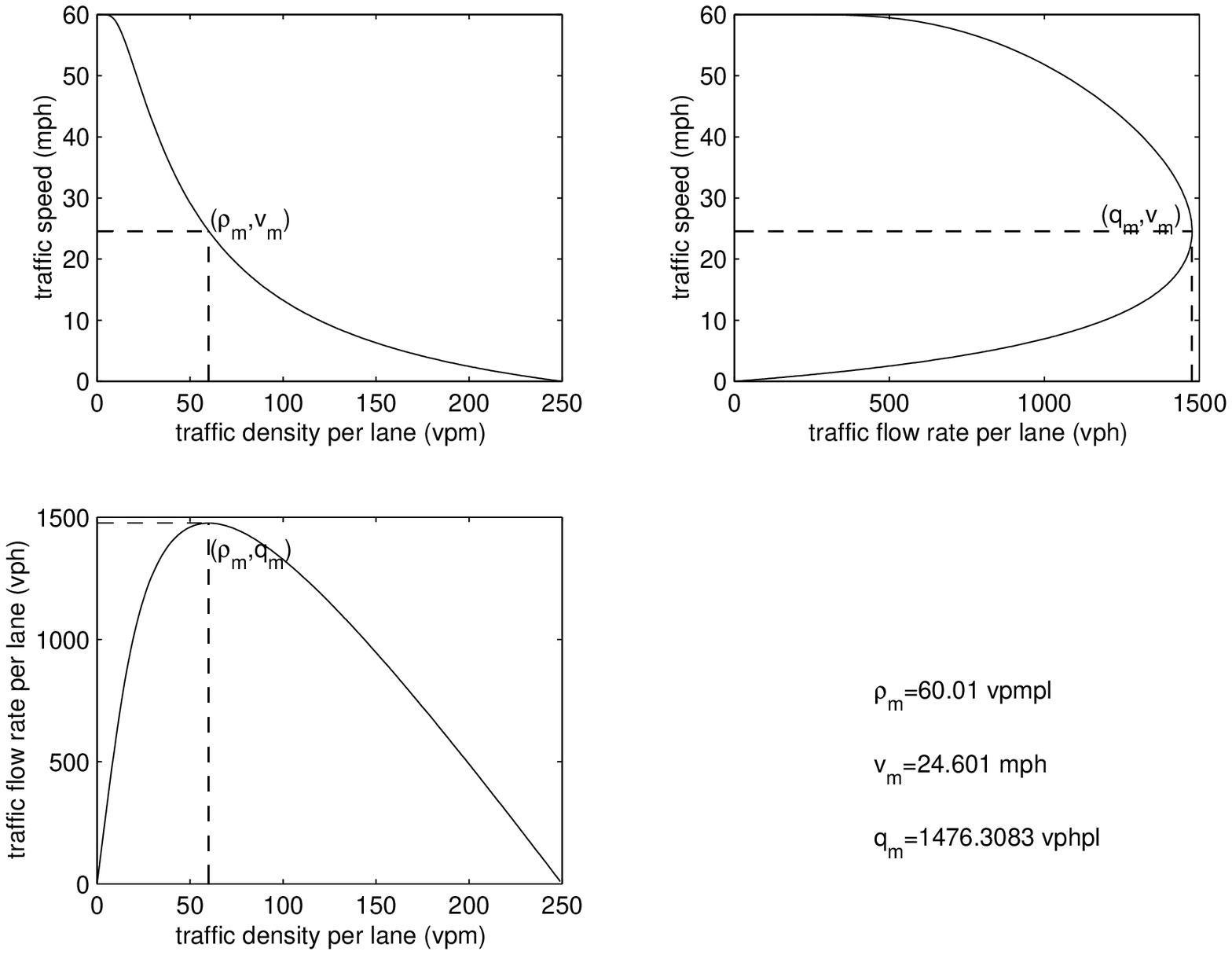}\ec\caption{Newell's Fundamental Diagram}\label{newell2}\efg
In this section, we carry out one numerical simulation for the network shown in \reff{F_network}, in which zone 2 and zone 18 have 4 lanes, zone 22 and zone 23 have 1 lane; all the other zones have 3 lanes.

The LWR model is used for the aggregate flow with Newell's fundamental diagram $f(\r)= \r v_f[1-\exp(\frac {|c_j|}{v_f}(1-\rho_j/\rho)) ]$, shown in \reff{newell2}. The free flow speed $v_f=60$ mph, the wave speed of jam density $c_j=-10$ mph and the jam density $\rho_j=250$ vpm. The characteristic speed of the LWR model $\lambda_{\ast}=v_f*(1-e^{|c_j|/v_f(1-\rho_j/\rho)})-c_j\rho_j/\rho e^{|c_j|/v_f(1-\rho_j/\rho)}$. The length of each zone $\Delta x=0.6$ mile; the length of a time step $\Delta t=30$ sec. The CFL number is $60 \cdot 5/600 /0.6=0.8444$. The maximum number of vehicles in a zone is 250$\times$0.6=150.  The maximum number of vehicles through a boundary in $\dt$ is about $12$ vehicles, when $\rho=60$vpm, $f=1476$vph.

We assume origin zone 0 is jammed with a repeating sequence of platoons with 25 vehicles to destination 1 followed by 10 vehicles to destination 2. Origin zone 22 is assumed to be jammed with a repeating sequence of platoons with 5 vehicles to destination 1 followed by 2 vehicles to destination 2. We assume destination zone 21 has the same number of vehicles as zone 20, and destination zone 23 has the number of vehicles to destination 2 in zone 18.

Using the merging rule introduced by Lebacque (1995\nocite{lebacque1995}), two Riemann problems are solved at boundary connector 2 and the boundary flux is calculated. At boundary connector 19, two fluxes to destination zone 23 and zone 19 are obtained by solving two Riemann problems. The sum of them is the total flux through the boundary.

\bfg\bc\includegraphics[height=8cm]{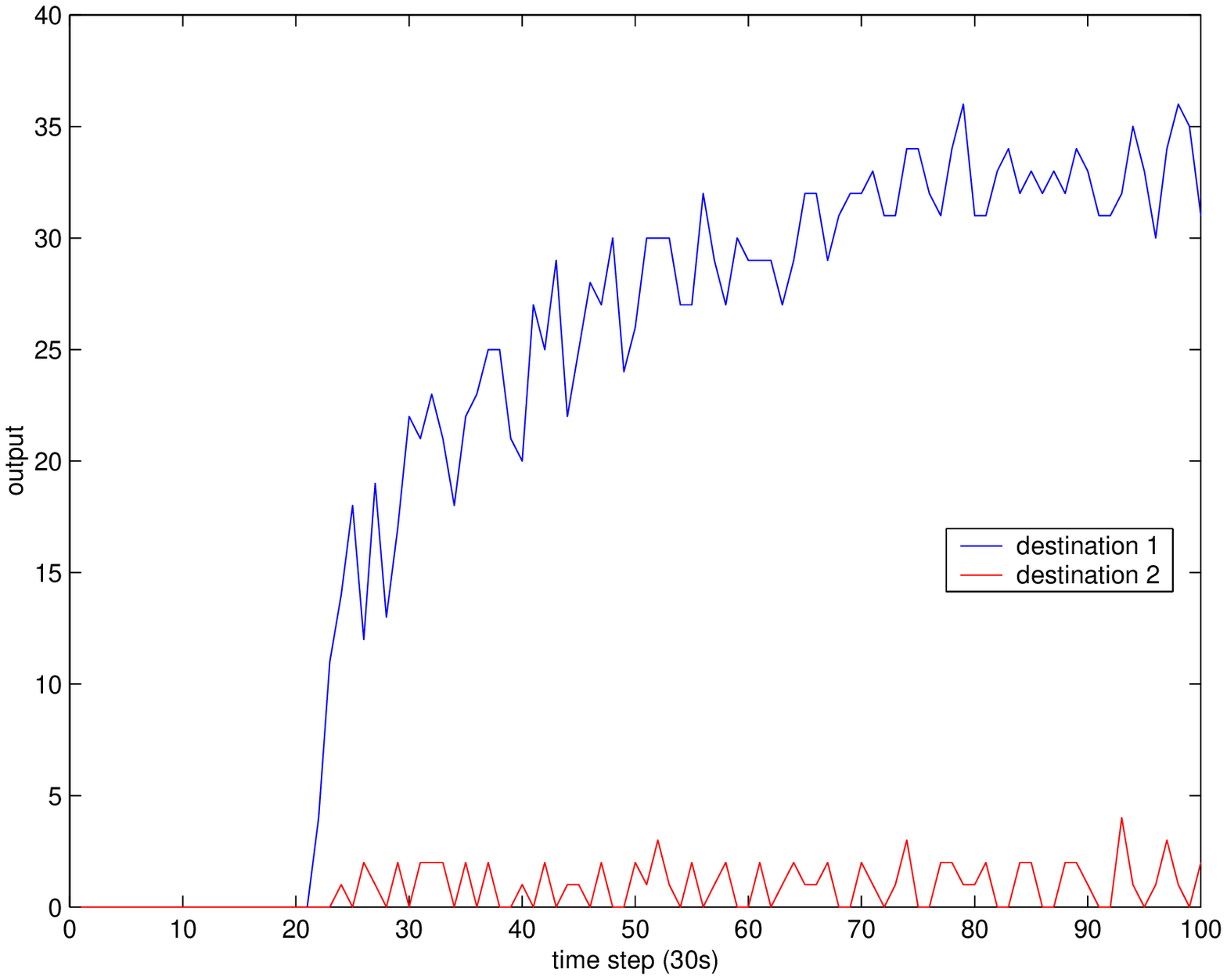}\ec\caption{Output at two destinations}\label{oupt5}\efg
\bfg\bc\includegraphics[height=8cm]{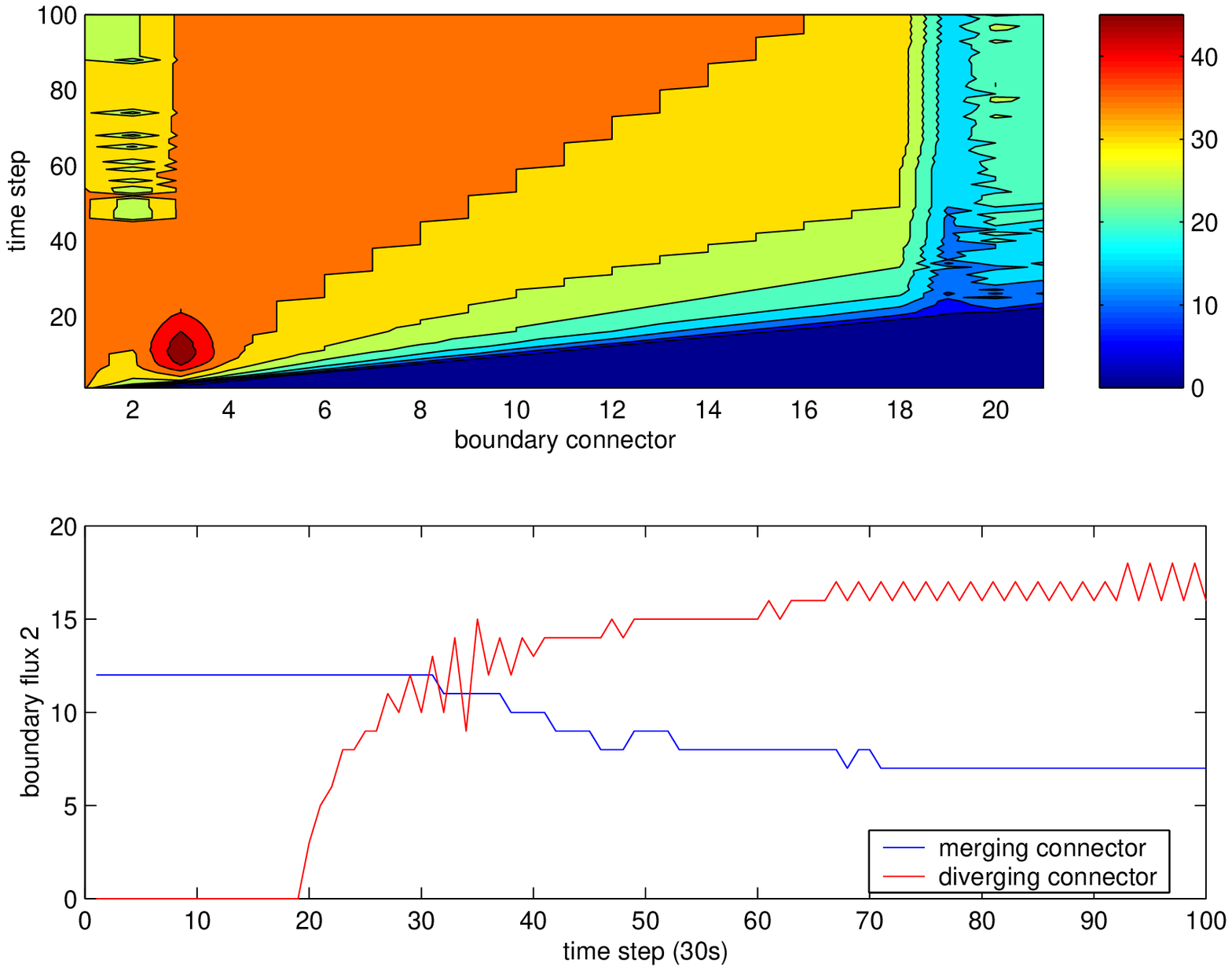}\ec\caption{The number of cars passing the boundary connectors}\label{conn5}\efg
\bfg\bc\includegraphics[height=8cm]{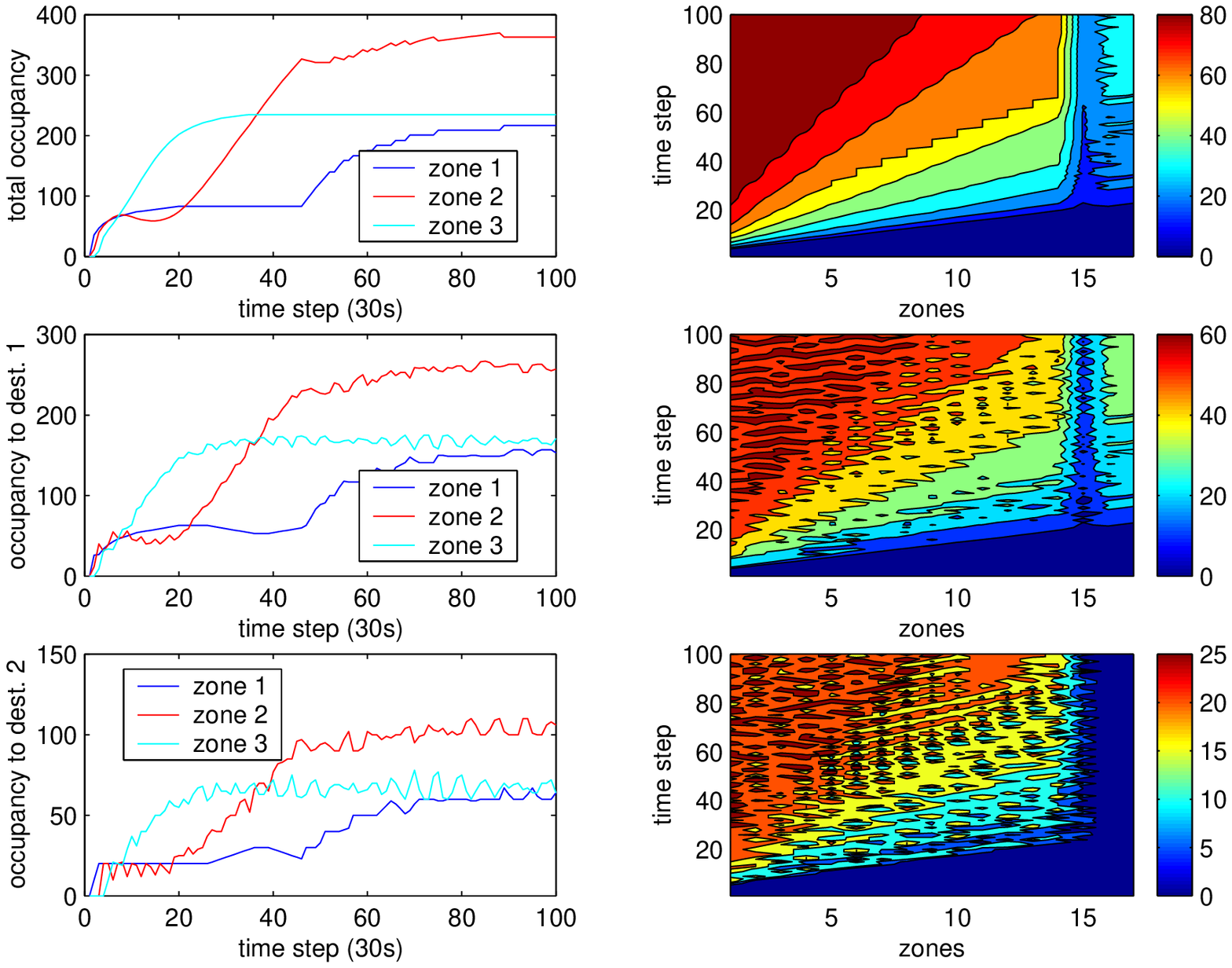}\ec\caption{Above: the total number of vehicles; Middle: the number of vehicles to destination 1; Below: the number of vehicles to destination 2}\label{zone5}\efg
\bfg\bc\includegraphics[height=8cm]{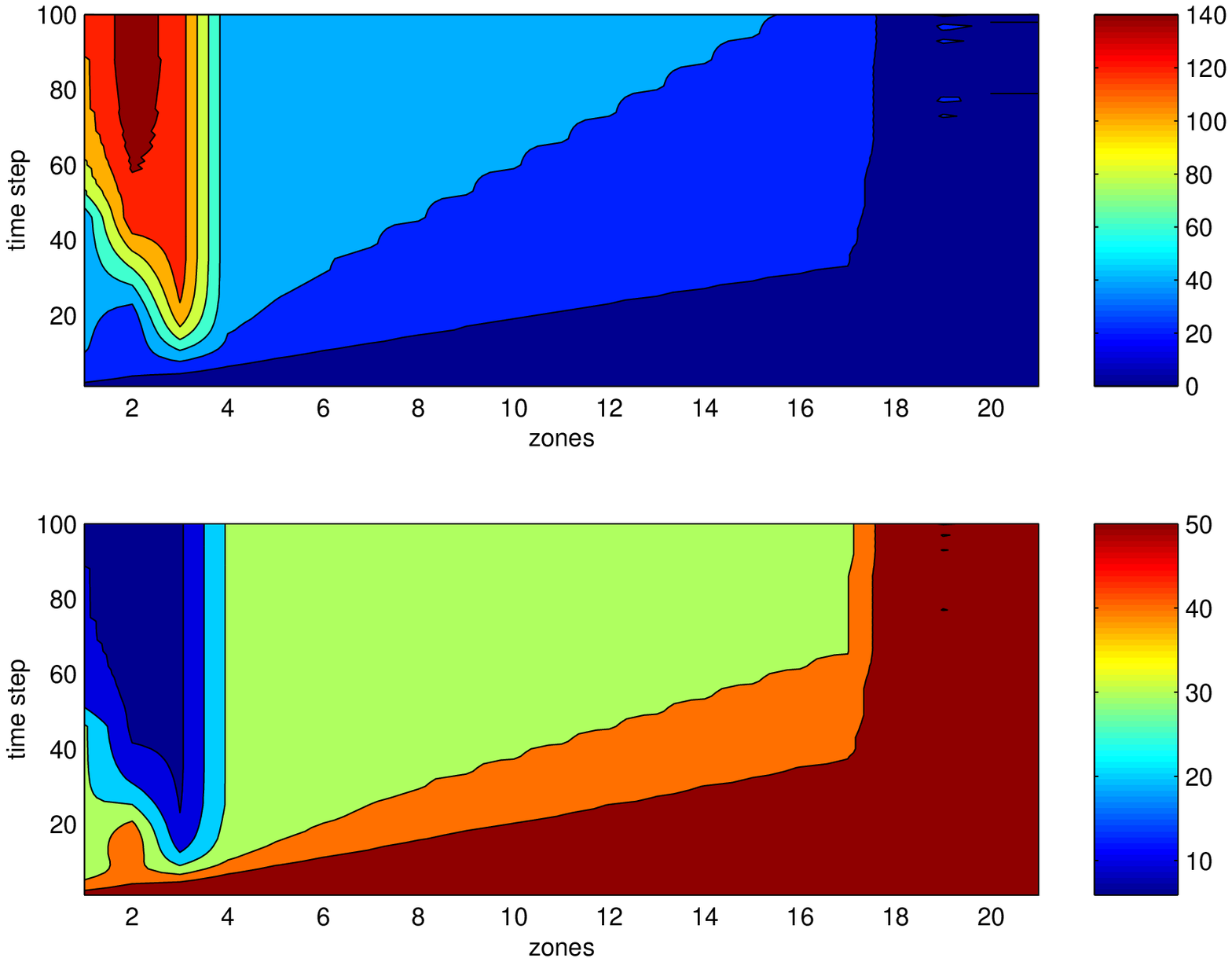}\ec\caption{Above: The densities in zone 1 to zone 21; Below: The velocities in zone 1 to 21}\label{rhov5}\efg

We get the following numerical results:
\ben
\item The output of two destinations are given in \reff{oupt5}.

\item The fluxes through all the boundary connectors are given in \reff{conn5}

\item The number of vehicles of all the zones are given in \reff{zone5}

\item The densities and velocities in zones 1 to 21 are given in \reff{rhov5}

\een

According to the results, a rarefaction wave still forms. However, when traffic flow moves to boundary connector 19, all the vehicles to destination 2 move out of the mainline of highway. (Refer to the below part of \reff{zone5}) Zone 2 becomes jammed after about $60\times 30$ seconds, since traffic has to merge from  4 lanes in zone 2 to 3 lanes in zone 3. These are consistent to experiences. However we do not have field observation data to justify our model at this time.

\section{Conclusions}
In this chapter, a new multi-commodity model is developed. In this model, by interpreting the FIFO principle in a new way and introducing the concept of boundary connectors, we develop a more efficient approach for disaggregated flows. Therefore, this new multi-commodity model is promising for modeling complex network flows.

For traffic flows at aggregate level, higher-order traffic flow models may be combined in order to capture more dynamics. For example, the PW model or Zhang's model discussed in previous chapters can be combined for link flows.

However, how to model the dynamics at merges or diverges is still an open problem, which is much harder than to model link flows.

Since this model is discrete, and network, data and program structures have been discussed in detail, it is ready for practical test. For example, this multi-commodity model can be used in simulating network traffic conditions in order to verify it validity.

Some possible applications of this network flow model can be seen now. One application is for developing better ramp metering methods. Another application is to determine dynamic assignment in a complex traffic network.

\newpage
\pagestyle{myheadings} 
\markright{  \rm \normalsize CHAPTER 7. \hspace{0.5cm}
 CONCLUSIONS}
\large
\chapter{CONCLUSIONS} 
This chapter provides concluding remarks on the research effort reported in this thesis. The chapter starts with an overall conclusion in section 7.1; In section 7.2 we discuss our contribution to traffic flow models and their numerical solutions. In the last section we discuss the future research directions in some areas related to this topic.

\section {Overall Conclusions}
In this research we studied five traffic flow models: the LWR model, Zhang's model, the PW model, the inhomogeneous LWR model and the multi-commodity model. For these models, we discussed their analytical solutions, developed numerical solution methods and carried out numerical simulations.

All these models preserve conservation of traffic. Since the homogeneous LWR model is the simplest model preserving traffic conservation, it is the basis to understand the wave solutions of these models and develop numerical solution methods for them. The LWR model is a first-order hyperbolic conservation law`. Shock and rarefaction waves are the basic solutions to a Riemann problem for such a conservation law. From the solutions to the Riemann problem, we can easily compute the boundary averages of traffic density and boundary fluxes, which are used in the numerical solution methods -- Godunov's methods. We solve the homogeneous LWR model with a first-order Godunov method and our numerical results are consistent with theoretical predictions.

In chapter 3 we discussed a non-equilibrium traffic flow model -- Zhang's model, which is a second-order hyperbolic system of conservation laws with a source term. Due to the difficulty in solving the Riemann problem for a system with a source term, we studied wave solutions to the Riemann problem for the homogeneous version of this model without considering the source term. The wave solutions are much more complicated than those to the Riemann problem for the LWR model. From these wave solutions, we computed the boundary averages of $\r$ and $v$, and developed a first-order and second-order Godunov methods for Zhang's model. The performance of these methods are examined in Section 3.4. We found that the second-order Godunov method performs better than a first-order method, however, its converge rate is 1, instead of 2 as expected and this is believed to be caused by the effect of the source term.

In Chapter 4 we discussed another non-equilibrium second-order traffic flow model -- the PW model. We discussed wave solutions to the Riemann problem for the homogeneous version of the PW model, as well as wave solutions to the Cauchy problem for the PW model with a source term. For the Cauchy problem, we found that the characteristics of a 1-rarefaction wave are approximated by parabolic curve. However, numerical results didn't show significant improvement in solutions to the PW model by solving Cauchy problems than by solving Riemann problems. For the PW model, we studied a first-order Godunov method, a second-order Godunov method, Pember's method, fractional method and LeVeque's method. All the higher-order methods don't show significantly better convergence rates than the first-order method, due to the effect of the source term. In LeVeque's method, we interpolated the given $\r$ and $v$ with solutions of the standing wave in each cell. However, this method was proved to under-estimate the effect of the source term and is not suitable for the PW model when solutions are far from equilibrium states. Different from Zhang's model, the PW model is unstable in certain density/speed regions. Therefore we examined the stability of the PW model with a first-order Godunov method in Section 4.4.1. The results obtained are consistent with theoretical predictions. 

In Chapter 5 we studied the inhomogeneous LWR model. By introducing a profile of inhomogeneities the inhomogeneous LWR model can be written as a $2\times 2$ non-strictly hyperbolic system of conservation laws. Then in Section 5.3 we rigorously solved the Riemann problem for the inhomogeneous LWR model for a roadway with variable number of lanes.  We found that the results are consistent with those empirical results found in literature. However our method is easier to be extended to higher-order inhomogeneous models. 

The multi-commodity model discussed in Chapter 6 is different from other models since it has a two-level structure and uses a discrete form of the LWR model at aggregate level. In Chapter 6, we made clear of the two-level structure of all multi-commodity models. By introducing boundary connectors, we suggested a more efficient method for computing the dynamics at disaggregated level. We also designed the network, data and program structures for a specific network in Section 6.3, and in Section 6.4, we presented some numerical simulation results which are consistent with theoretical expectations.  

Through discussions on different types of models and extensive numerical simulations, this thesis builds a solid base for  modeling traffic flow macroscopically. It helps to understand the wave properties of the continuum models. The numerical methods discussed in this work will be applicable to all applications related with those models. 

\section {Research Contributions}
In this subsection, we discuss about our contributions in this thesis.

Godunov methods have been well developed for studying hyperbolic systems of conservation laws. However they have not been extensively applied to the study of traffic flow problems. In this research, we made many efforts to apply Godunov methods and their variants to solve different traffic flow models. Our research shows that, to the continuum traffic flow models, Godunov methods are as useful as to other fluid dynamics.

In traditional Godunov methods, boundary averages are computed from solutions to a Riemann problem for a homogeneous system. In this research (Section 4.2.2), we tried to compute boundary averages from solutions to a Cauchy problem for the system with a source term. Although this new approach in computing boundary averages doesn't appear to improve the numerical accuracy significantly, it may help people to design Godunov-type methods which are more suitable for systems with source terms. 

The inhomogeneous LWR model was considered as a $2\times 2$ nonlinear resonant system in this thesis (in chapter 5), and rigorous solutions to the Riemann problem are studied based on existing theories, which were not known to the area of traffic flow models. This new approach helps us in understanding the wave phenomenon on an inhomogeneous roadway, and also brings us a new choice in solving inhomogeneous traffic flow models. 

In chapter 6, we made the two-level structure of multi-commodity models clearer, and presented a way to deal with traffic flows at different levels more efficiently. The discussion on the network, data and program structures for traffic networks in this research is another contribution.

\section {Future Research}
In this research, we made good progress in understanding some of the well-known continuum models. However, many topics discussed here still require further study. For example, an inhomogeneous model can be introduced for non-equilibrium flow and for multi-commodity flow, whose wave properties could be different than those in the homogeneous models.

Modeling of traffic flow is becoming increasingly important. However, development of new realistic traffic flow models is still a challenging task.  For example, there are few models that accurately capture the dynamics at merges and diverges or the interaction between different lanes on a multi-lane highway, a problem briefly touched in Chapter 6.

Besides developing new models, another major challenge is the validation of traffic flow models that have been proposed, such as the PW model and Zhang's model. It's not simply a matter of computing the results and comparing with observed data, a deeper understanding of the wave phenomena found in traffic flows and how to model them are critical to this endeavor and should be further pursued.

Although Godunov type of methods have been proved to be useful in solving traffic flow model, more efficient Godunov methods are needed to handle the diverse characteristics of these models. Another direction of research in computation may come from the development of parallel algorithms to simulate traffic flow in large networks.


\newpage\begin {thebibliography}{10} 
\bibitem{Bell89} 
{\sc J.B. Bell, P.Colella, and J. Trangenstein}, {\em Higher order Godunov methods for general systems of hyperbolic conservation laws}, J. Comput. Phys., 82(1989), pp. 362-397.
\bibitem{Chen94}
{\sc G.Q. Chen, C.D. Levermore, and T.P. Liu}, {\em Hyperbolic Conservation laws with stiff relaxation terms and entropy}, Communications on Pure and Applied Mathematics, Vol. XLVII, 787-830 (1994).
\bibitem{ColellaPuckett}
{\sc P. Colella and E.G. Puckett}, {\em Modern Numerical Methods for Fluid Flow}, In Draft, egp@math.ucdavis.edu, 2000. 
\bibitem{Daganzo94}
{\sc C.F. Daganzo}, {\em The cell transmission model: a dynamic representation of highway traffic consistent with hydrodynamic theory}, Transpn. Res. B, Vol. 28B, No.4, pp. 269-287, 1994.
\bibitem{Daganzo95}
{\sc C.F. Daganzo}, {\em The cell transmission model, Part II: Network traffic}, Transpn. Res. B, Vol. 29B, No.2, pp. 79-93, 1995. 
\bibitem{Daganzo95a}
{\sc C.F. Daganzo}, {\em A finite difference approximation of the kinematic wave model of traffic flow}, Transpn. Res. B, Vol 29B, No.4, pp. 261-276, 1995.
\bibitem{Daganzo95b}
{\sc C.F. Daganzo}, {\em Requiem for second-order fluid approximations of traffic flow}, Transpn. Res. B, Vol 29B, No.4, pp. 277-286, 1995.
\bibitem{gazis1961}
{\sc D.C. Gazis, R. Herman and R.W. Rothery}, {\em Nonlinear follow-the-leader models of traffic flow}, Operations Research 9(4), 545-567, 1961.
\bibitem{godunov1959}
{\sc S.K. Godunov}, {\em A difference method for numerical calculations of discontinuous solutions of the equations of hydrodynamics}, Mat. Sb., 47(1959), pp. 271-306. (In Russian.)

\bibitem{isaacson_t1992}
{\sc E.I. Isaacson and J.B. Temple}, {\em Nonlinear resonance in systems of conservation laws}, SIAM J. Appl. Math., Vol. 52, No. 5, pp. 1260-1278, October 1992.
\bibitem{jayakrishnan91}
{\sc R. Jayakrishnan}, {\em In-vehicle information systems for network traffic control: a simulation framework to study alternative guidance strategies}, The University of Texas at Austin, 1991.
\bibitem{Jin96}
{\sc S. Jin and C.D. Levermore}, {\em Numerical schemes for hyperbolic conservation laws with stiff relaxation terms}. Journal of Computational Physics, vol.126, (no.2), Academic Press, July 1996. p.449-67. 48.

\bibitem{JWL2000a}
{\sc W.L. Jin and H.M. Zhang}, {\em Evaluations of on-ramp control algorithms}, University of California at Davis, in preparation.

\bibitem{JWL2000b}
{\sc W.L. Jin and H.M. Zhang}, {\em The inhomogeneous LWR model and its numerical solutions}, University of California at Davis, in preparation.

\bibitem{JWL2000c}
{\sc W.L. Jin and H.M. Zhang}, {\em Numerical studies on the PW model}, University of California at Davis, in preparation.
\bibitem{Kerner}
{\sc B.S. Kerner and P. Konh\"a user},{\em Structure and parameters of clusters in traffic flow}, Physical Review E Volume 50, Number 1, pp54-83, July 1994.

\bibitem{lax1972}
{\sc P.D. Lax}, {\em Hyperbolic systems of conservation laws and the mathematical theory of shock waves}, SIAM, Philadelphia, Pennsylvania, 1972.

\bibitem{lebacque1995}
{\sc J.P. Lebacque}, {\em The Godunov scheme and what it means for first order traffic flow models}, \underline{http://cermics.enpc.fr/reports}, Rapports de recherche Research reports, Nov. 1995.

\bibitem{lebacque1996}
{\sc J.P. Lebacque}, {\em The Godunov scheme and what it means for first order traffic flow models}, Transportation and traffic theory, ITTT 1996.

\bibitem{lebacque1998}
{\sc J.P. Lebacque}, {\em Macroscopic traffic flow models: a question of order}, Rapports de recherche Research reports, \underline{http://cermics.enpc.fr/reports}, Sep. 1998.

\bibitem{lebacque1999}
{\sc J.P. Lebacque}, {\em On a simple model of the interaction between bus and traffic flow}, Rapports de recherche Research reports, \underline{http://cermics.enpc.fr/reports}, Jan. 1999.

\bibitem{LeVeque98a}
{\sc R.J. LeVeque and D.S. Bale}, {\em Wave Propagation Methods for Conservation Laws with Source Terms}, Submitted to Proc. 7'th Int'l Conf. on Hyperbolic Problems, Zurich, Feb. 1998.
\bibitem{LeVeque98b}
{\sc R.J. LeVeque}, {\em Balancing Source Terms and Flux Gradients in High-Resolution Godunov Methods:The Quasi-Steady Wave-Propagation Algorithm}, \underline{http://www.amath.washington.edu/people/faculty/leveque/}, May 1998.
\bibitem{lighthill}
{\sc M.J. Lighthill and G.B. Whitham}, {\em On kinematic waves: II. A theory of traffic flow on long crowded roads}, Proc. Royal Society, volume 229 (1178) of A, pages 317-345, 1955.
\bibitem{lin_t_w1995}
{\sc L. Lin, J.B. Temple and J. Wang}, {\em A comparison of convergence rates for Godunov's method and Glimm's method in resonant nonlinear systems of conservation laws}, SIAM J. Numer. Anal., Vol. 32, No. 3, pp. 824-840, June 1995.

\bibitem{Liu79}
{\sc T.P. Liu}, {\em Quasi-linear Hyperbolic Systems}, Commun. Math. Phys. 68, pp 141-172, 1979.
\bibitem{Liu87}
{\sc T.P. Liu}, {\em Hyperbolic conservation laws with relaxation}, Comm. Math. Phys., 108, pp. 153-175, 1987.
\bibitem{newell93}
{\sc G.F. Newell}, {\em A simplified theory of kinematic waves in highway traffic, part I: general theory}, Transpn. Res. B, Vol 27B, No.4, pp. 281-287, 1993.
\bibitem{payne}
{\sc H.J. Payne}, {\em Models of freeway traffic and control}, in Beeky,G.A., editor, Mathematical Models of Public Systems, Vol. 1 of Simulation Councils Proc. Ser., pp 51-60, 1971.
\bibitem{Pember93a}
{\sc R.B. Pember}, {\em Numerical Methods for Hyperbolic Conservation Laws with Stiff Relaxation I. Spurious Solutions}, Siam J. Applied Mathematics, Vol 53,No. 5, pp 1293-1330, 1993.
\bibitem{Pember93b}
{\sc R.B. Pember}, {\em Numerical Methods for Hyperbolic Conservation Laws with Stiff Relaxation II. Higher Order Godunov Methods}, Siam J. Scientific Computing,Vol 14, No. 4, pp 824-859, 1993.
\bibitem{richards}
{\sc P.I. Richards}, {\em Shock waves on the highway}, Operations Research, 4:42-51, 1956.

\bibitem{schochet1988}
{\sc S. Schochet}, {\em The instant response limit in Whitham's nonlinear traffic model: uniform well-posedness and global existence}, Asymptotic Analysis 1: 263-282, 1988.
\bibitem{smoller1983}
{\sc J. Smoller}, {\em Shock waves and reaction-diffusion equations}, New York, Springer-Verlag, 1983.

\bibitem{vaughan84}
{\sc R. Vaughan, V. F. Hurdle and E. Hauer}, {\em A traffic flow 
model with time dependent O-D patterns}, Ninth International 
Symposium on Transportation and Traffic Theory, VNU Science Press, pp. 155-178, 1984.
\bibitem{Whitham59}
{\sc G.B. Whitham}, {\em Some comments on wave propagation and shock wave structure with application to magnetohydronamics}, Comm. Pure Appl. Math.,XII,pp. 113-158, 1959.
\bibitem{Whitham}
{\sc G.B. Whitham}, {\em Linear and Nonlinear Waves}, John Wiley and Sons, New York, 1974.
\bibitem{zhang96rm}
{\sc H.M. Zhang, S.G. Ritchie and W.W. Recker}, {\em Some General
 Results on the Optimal Ramp Control Problem}, Transpn Res. -C, Vol. 4, No.2, pp. 51-69, 1996.
\bibitem {LocANN}
{\sc H.M. Zhang and S.G. Ritchie}, {\em Freeway ramp
metering using artificial neural networks}, Transportation Research C, Vol. 5, No. 5, pp. 273-286, 1997.
\bibitem{Zhang1998}
{\sc H.M. Zhang}, {\em A theory of nonequilibrium traffic flow}, Transportation Research B, Vol. 32, No. 7, pp. 485-498, 1998.

\bibitem{Zhang1999a}
{\sc H.M. Zhang}, {\em An analysis of the stability and wave properties of a new continuum theory}, Transportation Research B 33(6): 387-398, 1999.

\bibitem{zhang99rm}
{\sc H.M. Zhang and W.W. Recker}, {\em On optimal freeway ramp control policies for congested traffic corridors}, Transportation Research Part B 33 (1999) 417-436.
\bibitem{Zhang2000a}
{\sc H.M. Zhang}, {\em Structural properties of solutions arising from a non-equilibrium traffic flow theory}, Transportation Research B 34, 583-603, 2000.
\bibitem{Zhang2000b}
{\sc H.M. Zhang}, {\em A finite difference model of nonequilibrium traffic flow}, to appead in Transportation Research B, 2000.

\end {thebibliography}
\end {document}